\documentclass[11pt]{amsart}

\usepackage{amsmath}
\usepackage{amssymb}
\usepackage{graphicx}

\usepackage{url}
\usepackage{dsfont}
\usepackage{fancyhdr}
\usepackage{graphicx}
\usepackage{subfigure}
\usepackage{amssymb}
\usepackage{amsmath}
\usepackage{amsfonts}
\usepackage{adjustbox}
\usepackage{mathrsfs}
\usepackage{mathtools}
\usepackage{bm}
\usepackage{tabularx}
\usepackage{color}
\usepackage{setspace}
\usepackage{exscale}
\usepackage{relsize}
\usepackage{bm}

\newcommand\eref[1]{(\ref{#1})}

\newcommand*\xbar[1]{%
  \hbox{%
    \vbox{%
      \hrule height 0.5pt 
      \kern0.4ex
      \hbox{%
        \kern-0.05em
        \ensuremath{#1}%
        \kern-0.05em
      }%
    }%
  }%
}

\newcommand{\mF}{\bm{F}}

\newcommand{\mG}{\bm{G}}

\newcommand{\mU}{\bm{U}}

\newcommand{\dx}{\Delta x}
\newcommand{\dy}{\Delta y}

\newcommand{\hf}{{\frac{1}{2}}}

\newcommand{\jph}{{j+\frac{1}{2}}}
\newcommand{\jmh}{{j-\frac{1}{2}}}
\newcommand{\iph}{{i+\frac{1}{2}}}
\newcommand{\imh}{{i-\frac{1}{2}}}
\newcommand{\kph}{{k+\frac{1}{2}}}
\newcommand{\kmh}{{k-\frac{1}{2}}}

\newenvironment{DA}{{\flushleft \bf Declarations:}}{}
\def\softd{{\leavevmode\setbox1=\hbox{d}%
		\hbox to 1.05\wd1{d\kern-0.4ex{\char039}\hss}}}

\renewcommand{\S}{Section~}

\newtheorem{theorem}{Theorem}[section]

\theoremstyle{definition}

\newtheorem{remark}[theorem]{Remark}
\newtheorem{defn}[theorem]{Definition}
\def\softd{{\leavevmode\setbox1=\hbox{d}%
		\hbox to 1.05\wd1{d\kern-0.4ex{\char039}\hss}}}

\numberwithin{equation}{section}

\title[Numerical Study of Dissipative Weak Solutions]{Numerical Study of Dissipative Weak Solutions for the Euler Equations of Gas Dynamics}

\author[S. Chu]{Shaoshuai Chu}
\address[S. Chu]{Department of Mathematics, RWTH Aachen University, 52056 Aachen, Germany}
\email{{\tt chu@igpm.rwth-aachen.de}}

\author[M. Herty]{Michael Herty}
\address[M. Herty]{Department of Mathematics, RWTH Aachen University, 52056 Aachen, Germany; Department of Mathematics and Applied
Mathematics, University of Pretoria, Hatfield, 0028, South Africa}
\email{\tt herty@igpm.rwth-aachen.de}

\author[A. Kurganov]{Alexander Kurganov}
\address[A. Kurganov]{Department
of Mathematics and Shenzhen International Center for Mathematics, Southern University of Science and Technology, Shenzhen, 518055, China}
\email{\tt alexander@sustech.edu.cn}

\author[M.  Luk\'{a}\v{c}ov\'{a}-Medvi\v{d}ov\'{a}]{M\'{a}ria Luk\'{a}\v{c}ov\'{a}-Medvi\v{d}ov\'{a}}
\address[M.  Luk\'{a}\v{c}ov\'{a}-Medvi{\softd}ov\'{a}]{Institute of Mathematics, Johannes Gutenberg-University Mainz, 55128 Mainz, Germany}
\email{\tt lukacova@uni-mainz.de}

\author[C. Yu]{Changsheng Yu}
\address[C. Yu]{Institute of Mathematics, Johannes Gutenberg-University Mainz, 55128 Mainz, Germany}
\email{\tt chayu@uni-mainz.de}



\begin{document}

\maketitle
{\centerline{\emph{Dedicated to Professor Eitan Tadmor}}}

\begin{abstract}
We study dissipative weak solutions of the Euler equations of gas dynamics using the first-, second-, third-, fifth-, seventh-, and ninth-order local characteristic decomposition-based central-upwind, low-dissipation central-upwind, and viscous finite volume methods, whose higher-order extensions are obtained via the framework of the alternative weighted essentially non-oscillatory (A-WENO) schemes. These methods are applied to several benchmark problems, including several two-dimensional Riemann problems and a Kelvin-Helmholtz instability test. The numerical results demonstrate that for methods converging only weakly in space and time, the limiting solutions are generalized dissipative weak solutions, approximated in the sense of ${\mathcal K}$-convergence and dependent on the numerical scheme. For all of the studied methods, we compute the associated Young measures and compare the dissipative weak solutions using entropy production and energy defect criteria. 
\end{abstract}


\section{Introduction}
We consider the two-dimensional (2-D) Euler equations of gas dynamics, which read as
\begin{equation}
\begin{aligned}
&\rho_t+(\rho u)_x+(\rho v)_y=0,\\
&(\rho u)_t+(\rho u^2 +p)_x+(\rho uv)_y=0,\\
&(\rho v)_t+(\rho uv)_x+(\rho v^2+p)_y=0,\\
&E_t+\left[u(E+p)\right]_x+\left[v(E+p)\right]_y=0.
\end{aligned}
\label{1.1}
\end{equation}
Here, $x$ and $y$ are spatial variables, $t$ the time, $\rho$ is the density, $u$ and $v$ are the velocities in the $x$- and $y$-directions,
respectively, $E$ is the total energy, $p$ is the pressure. The system \eref{1.1} is closed using the following polytropic equation of
state:
\begin{equation}
\begin{aligned}
& e(\rho,\theta)=c_v\theta,\quad  c_v=\frac{1}{\gamma-1},\\
& p(\rho,\theta)=(\gamma-1)\rho e=
\rho\theta=(\gamma-1)\Big[E-\frac{1}{2}\rho(u^2+v^2)\Big],
\end{aligned}
\label{1.2}
\end{equation}
where $e$ is the specific internal energy, $\theta$ is the absolute temperature, $\bm m=(\rho u,\rho v)^\top$ is the momentum vector, and
$\gamma>1$ is a constant representing the adiabatic coefficient. The second law of thermodynamics is expressed by the entropy inequality
for the total entropy $S$:
\begin{equation}
S_t+(Su)_x+(Sv)_y\ge0,\quad S:=c_v\rho\ln\left(\frac{p}{\rho^\gamma}\right).
\label{1.3}
\end{equation}
We consider a bounded domain $\Omega\subset\mathbb R^2$, with either periodic or no-flux boundary conditions, and the initial conditions
\begin{equation}\label{1.4}
\begin{aligned}
&\rho(x,y,0)=\rho_0(x,y),~~\bm m(x,y,0)=\bm m_0(x,y),\\
&S(x,y,0)=S_0(x,y),~~E_0=\rho_0e(\rho_0,S_0)+\frac{|\bm m_0|^2}{2\rho_0}.
\end{aligned}
\end{equation}

Note that the  energy $E$,  expressed as a function of the density $\rho,$ momentum $\bm m$, and entropy $S,$
\begin{equation}
\label{ee}
E = \left\{ \begin{array}{l} \frac{|\bm m|^2}{2\rho}  + c_v \rho^{\gamma} \exp\left(\frac{S}{c_v \rho}\right) \mbox{ if }\ \rho > 0,\ \bm m \in \Bbb R^2, S \in \Bbb R,\\ \\
0 \ \mbox{ if }\ \rho =0, \bm m = 0, S\leq 0,\\ \\ 
\infty \ \mbox{ if }\ \rho < 0 \mbox{ or } \rho = 0, \bm m \ne 0, S > 0,
\end{array} \right.
\end{equation}
is a convex lower semi-continuous function on $ \Bbb R^4.$

It is well-known that solutions of \eref{1.1}--\eref{1.2} may develop complex wave structures, such as shock waves, rarefactions, and
contact discontinuities, even when the initial data are infinitely smooth.
Consequently, distributional weak solutions of \eref{1.1}--\eref{1.2} are typically considered. However, it was proved in
\cite{de2010admissibility} that the multidimensional compressible Euler equations admit infinitely many weak entropy solutions satisfying
\eref{1.3} for a wide range of initial data; see also \cite{buckmaster2017onsager,chiodaroli2015global,feireisl2020density} and references
therein. 

As it was demonstrated in many simulations reported in, e.g.,
\cite{CCHKL_22,CCK23_Adaptive,feireisl2021computing,fjordholm2016computation,KX_22}, no strong convergence is observed when numerical
methods are applied to shear flows appearing, for instance, in Kelvin-Helmholtz (KH) instability problems. Consequently, a natural question
arises: {\em What is the weak limit of those numerical approximations?} To answer this question, one has to introduce a generalized solution
concept, which is amenable to numerical methods. As advocated in \cite{diperna1985measure}, measure-valued solutions are suitable to
represent low-regularity or rough solutions. Measure-valued solutions for multidimensional hyperbolic conservation laws were also studied in
\cite{fjordholm2016computation} and references therein. In this paper, we will concentrate on consistent approximations of dissipative weak
(DW) solutions \cite{Lukacova_book} of the Euler equations of gas dynamics.

According to the results reported in \cite{feireisl2021computing,fjordholm2016computation}, the limit of consistent approximations is not a
weak solution in the sense of distributions, if the convergence is only weak. However, a weak convergence of numerical solutions can lead to
the strong convergence to a DW solution by averaging the numerical solutions over different mesh resolutions. In
\cite{Lukacova_book,lukavcova2023limit}, the Ces\` aro averages over different meshes were considered for low-order numerical methods and a
strong convergence of the Ces\`aro averages to a DW solution was observed; see also \cite{feireisl2021computing} as well as \cite{CHIIKL},
where a fifth-order alternative weighted essentially non-oscillatory (A-WENO) scheme was used to compute the Ces\`aro averages.

In this paper, we study the convergence of several modern finite volume (FV) and finite difference (FD) methods of different orders of
spatial accuracy. Although a rigorous consistency analysis is not available for high-order methods, the numerical study presented in this
paper indicates that the numerical approximations of different orders converge weakly. This demonstrates that the consistency errors decay
when the mesh is refined. The sequence of numerical solutions behaves in the following way. When a weak solution of the problem exists, all
different numerical methods converge to this weak solution. However, when the solutions exhibit oscillations and develop turbulent
structures, different numerical methods may yield different limiting DW solutions. In light of this observation, we discuss possible
selection criteria for choosing an admissible and physically relevant solution among those oscillatory numerical solutions.

The paper is organized as follows. In \S\ref{sec2}, we give a brief overview of consistent approximations and DW solutions. In \S\ref{sec3},
we present a description of FV and FD methods used in the conducted numerical experiments. \S\ref{sec4} contains simulations for several
problems, including several configurations of the 2-D Riemann problems and the KH instability. We analyze the density profiles of the
obtained numerical solutions, along with the averages over solutions of different orders and their time averages. For the KH problem, we
also compute the underlying Young measure in specific subdomains. Additionally, we assess entropy production and energy defects of various
methods to evaluate them with respect to the selection criteria. Finally, \S\ref{sec5} is devoted to a brief conclusion.

\section{Consistent Approximations and Dissipative Weak Solutions}\label{sec2}
To analyze a given numerical method and establish its convergence, it is essential to investigate the local consistency error. For this
purpose, the notion of a consistent approximation seems to be particularly useful; see, e.g., \cite{Lukacova_book}.

\begin{defn}[\textbf{Consistent approximation}]\label{D3}
A sequence $\left\{\rho_\ell,\bm m_\ell,S_\ell\right\}_{\ell=1}^\infty$ is a {\em consistent approximation} of \eref{1.1}--\eref{1.4} in
$\Omega\times(0,T)$ if

\medskip
\noindent
$\bullet$ There is a sequence $\left\{\rho_{0,\ell},\bm m_{0,\ell},S_{0,\ell}\right\}_{\ell=1}^\infty$, that weakly approximates the
initial data \eref{1.4}, that is,
\begin{equation*}
\begin{aligned}
&\rho_{0,\ell}\to\rho_0~\mbox{weakly in}~L^1(\Omega),\quad\bm m_{0,\ell}\to\bm m_0~\mbox{weakly in}~L^1(\Omega;\mathbb R^2),\\
& S_{0,\ell}\to S_0~\mbox{weakly in}~L^1(\Omega),\\
&\int\limits_\Omega\left(\frac{|\bm m_{0,\ell}|^2}{2\rho_{0,\ell}}+\rho_{0,\ell}e(\rho_{0,\ell},S_{0,\ell})\right){\rm d}x{\rm d}y\to
\int\limits_\Omega\left(\frac{|\bm m_0|^2}{2\rho_0}+\rho_0e(\rho_0,S_0)\right){\rm d}x{\rm d}y,
\end{aligned}
\end{equation*}
and satisfies the {\bf energy inequality}
\begin{equation*}
\resizebox{\linewidth}{!}{$
\int\limits_\Omega\left(\frac{|\bm m_\ell|^2}{2\rho_\ell}+\rho_\ell e(\rho_\ell,S_\ell)\right)(\cdot,t)\,{\rm d}x{\rm d}y\le
\int\limits_\Omega\left(\frac{|\bm m_{0,\ell}|^2}{2\rho_{0,\ell}}+\rho_{0,\ell}e(\rho_{0,\ell},S_{0,\ell})\right){\rm d}x{\rm d}y+e^1_\ell
$}
\end{equation*}
for a.a. $t\in[0,T]$ with $e^1_\ell\to0$ as $\ell\to\infty$;

\medskip
\noindent
$\bullet$ The {\bf equation of continuity}
\begin{equation*}
\int\limits_0^T\int\limits_\Omega\left(\rho_\ell\varphi_t+\bm m_\ell\!\cdot\!\nabla\varphi\right){\rm d}x{\rm d}y{\rm d}t=
-\int\limits_\Omega\rho_{0,\ell}\varphi(\cdot,0)\,{\rm d}x{\rm d}y+ e^2_\ell[\varphi]
\end{equation*}
holds for any $\varphi\in C^1_c(\Omega\times[0,T))$ with $e^2_\ell[\varphi]\to0$ as $\ell\to\infty$ for any
$\varphi\in C^2_c(\Omega\times[0,T))$;

\medskip
\noindent
$\bullet$ The {\bf momentum equation}
\begin{equation*}
\begin{aligned}
&\int\limits_0^T\int\limits_\Omega\left(\bm m_\ell\cdot\bm\varphi_t+\mathds 1_{\rho_\ell>0}\frac{\bm m_\ell\otimes\bm m_\ell}{\rho_\ell}:
\bm\nabla\bm\varphi+\mathds 1_{\rho_\ell>0}p(\rho_\ell,S_\ell)\bm{\nabla}\!\cdot\!\bm\varphi\right){\rm d}x{\rm d}y{\rm d}t\\
&\quad=-\int\limits_\Omega\bm m_{0,\ell}\bm\varphi(\cdot,0)\,{\rm d}x{\rm d}y+e^3_\ell[\bm\varphi]
\end{aligned}
\end{equation*}
holds for any $\bm\varphi\in C^1_c(\Omega\times[0,T);\mathbb R^2)$ with $e^3_\ell[\bm\varphi]\to0$ as $\ell\to\infty$ for any
$\bm\varphi\in C^2_c(\Omega\times[0,T);\mathbb R^2)$;

\medskip
\noindent
$\bullet$ The {\bf entropy inequality}
\begin{equation*} 
\int\limits_0^T\int\limits_\Omega\left(S_\ell\varphi_t+\mathds 1_{\rho_\ell>0}\left(S_\ell\frac{\bm m_\ell}{\rho_\ell}\right)\cdot
\nabla\varphi\right){\rm d}x{\rm d}y{\rm d}t\le-\int\limits_\Omega S_{0,\ell}\varphi(\cdot,0)+e^4_\ell[\varphi]
\end{equation*}
holds for any $\varphi\in C^1_c(\Omega\times[0,T))$, $\varphi\ge0$ with $e^4_\ell[\varphi]\to0$ as $\ell\to\infty$ for any
$\varphi\in C^2_c(\Omega\times[0,T))$;

\medskip
\noindent
$\bullet$ The {\bf minimum entropy principle} is satisfied, namely, there exists $\underline s\in\mathbb R$ such that
$S_\ell\ge\rho_\ell\underline s$ a.e. in $\Omega\times(0,T)$.
\end{defn}

A consistent approximation can be obtained through various strategies, such as vanishing viscosity methods or structure-preserving numerical
schemes. Different numerical methods on a sequence of refined meshes denoted as $\{{\mathcal I}_{h_n}\}_{n=1}^\infty$, where $h_n>0$ represents
a mesh parameter, were considered in \cite{Lukacova_book}. This approach yields consistent approximations denoted by
$$
\left\{\rho_\ell,\bm m_\ell,S_\ell\right\}_{\ell=1}^\infty=\left\{\rho_{h_n},\bm m_{h_n},S_{h_n}\right\}_{n=1}^\infty,~
h_n\to0~\mbox{and}~\ell=\ell(h_n)\to\infty~\mbox{as}~n\to\infty.
$$

The following result was established in \cite{feireisl2021computing}.
\begin{theorem}\label{T1}
Let $\{\rho_\ell,\bm m_\ell,S_\ell\}_{\ell=1}^\infty$ be a family of consistent approximations of \eref{1.1}--\eref{1.4} in the sense of
Definition \ref{D3}, such that
$$
\begin{aligned}
&\rho_\ell\to\rho&&\mbox{weakly-(*) in}~L^\infty(0,T;L^\gamma(\Omega)),\\
&\bm m_\ell\to\bm m&&\mbox{weakly-(*) in}~L^\infty(0,T;L^{\frac{2\gamma}{\gamma+1}}(\Omega;\mathbb R^2)),\\
&S_\ell\to S&&\mbox{weakly-(*) in}~L^\infty(0,T;L^\gamma(\Omega)) \mbox{ for } \ell\to\infty. 
\end{aligned}
$$
Suppose 
the limit $(\rho,\bm m,S)$ is a weak solution of \eref{1.1}--\eref{1.4}. Then,  there exists a
subsequence $\{\rho_{\ell_k},\bm m_{\ell_k},S_{\ell_k}\}$ such that
$$
\begin{aligned}
&\rho_{\ell_k}\to\rho&&\mbox{in}~L^q(0,T;L^\gamma(\Omega)),\\
&\bm m_{\ell_k}\to\bm m&&\mbox{in}~L^q(0,T;L^{\frac{2\gamma}{\gamma+1}}(\Omega; \mathbb R^2)),\\
&S_{\ell_k}\to S&&\mbox{in}~L^q(0,T;L^\gamma(\Omega)) \qquad \mbox{ for } \ell_k\to\infty, \  1\le q<\infty.
\end{aligned}
$$
\end{theorem}

Next, we provide a definition of DW solutions and state two convergence results proved in
\cite{feireisl2021computing,feireisl2020finite,Lukacova_book}.
\begin{defn}[{\bf DW solution}]\label{D2}
Let the initial data satisfy
$$
\resizebox{\linewidth}{!}{$
\rho_0\in L^\gamma(\Omega),\quad\bm m_0\in L^\frac{2\gamma}{\gamma+1}(\Omega; \mathbb R^2),\quad S_0 \in L^\gamma(\Omega),~~\mbox{and}~~
\int\limits_{\Omega}E(\rho_0,\bm m_0,S_0)\,{\rm d}x{\rm d}y<\infty.
$}
$$
We say that $(\rho,\bm m,S)$ is a {\em DW solution} of \eref{1.1}--\eref{1.4} in
$\Omega\times[0,T)$, $0<T<\infty$, if it satisfies the following properties.

\medskip
\noindent
$\bullet$ {\bf Regularity}
$$
\begin{aligned}
&\rho\in C_{\rm weak}([0,T];L^\gamma(\Omega)),\quad\bm m\in C_{\rm weak}([0,T];L^{\frac{2\gamma}{\gamma+1}}(\Omega; \mathbb R^2)),
\\
&S\in L^\infty(0,T;L^\gamma(\Omega))\cap BV_{\rm weak}([0,T];L^\gamma(\Omega)),\\
&\int\limits_{\Omega}E(\rho,\bm m,S)(\cdot,t)\,{\rm d}x{\rm d}y\le\int\limits_{\Omega}E(\rho_0,\bm m_0,S_0)\,{\rm d}x{\rm d}y~~
\mbox{for any}~t\in[0,T).
\end{aligned}
$$

\medskip
\noindent
$\bullet$ The {\bf equation of continuity}
$$
\int\limits_0^T\int\limits_{\Omega}\left(\rho\varphi_t+\bm m\!\cdot\!\nabla\varphi\right){\rm d}x{\rm d}y{\rm d}t=
-\int\limits_{\Omega}\rho_0\varphi(\cdot,0)\,{\rm d}x{\rm d}y
$$
holds for any $\varphi\in C^1_c(\Omega\times[0,T))$.

\medskip
\noindent
$\bullet$ The {\bf momentum equation}
$$
\begin{aligned}
\int\limits_0^T&\int\limits_{\Omega}\left(\bm m\cdot\bm\varphi_t+\mathds 1_{\rho>0}\frac{\bm m\otimes\bm m}{\rho}: \bm\nabla\bm\varphi+
p(\rho,S)\bm{\nabla}\!\cdot\!\bm\varphi\right){\rm d}x{\rm d}y{\rm d}t\\
&=-\int\limits_0^T\int\limits_{\Omega}\bm\nabla\bm\varphi: {\rm d}\mathfrak R(t)\,{\rm d}t-
\int\limits_{\Omega}\bm m_0\cdot\bm\varphi(\cdot,0)\,{\rm d}x{\rm d}y
\end{aligned}
$$
holds for any $\bm\varphi\in C^1_c(\Omega\times[0,T); \mathbb R^2)$, where
$\mathfrak R\in L^\infty(0,T;{\mathcal M}^+(\Omega;\mathbb R^{2\times2}_{\rm sym}))$ is the Reynolds defect stress.

\medskip
\noindent
$\bullet$ The {\bf entropy inequality}
$$
\resizebox{\linewidth}{!}{$
\int\limits_\Omega\Big[S(\cdot,\tau_2+)\varphi(\cdot,\tau_2+)-S(\cdot,\tau_1-)\varphi(\cdot,\tau_1-)\Big]{\rm d}x{\rm d}y\ge
\int\limits_{\tau_1}^{\tau_2}\int\limits_\Omega\left(S\varphi_t+\left<{\mathcal V}_{x,y,t};\mathds 1_{\widehat\rho\,>0}
\big(\widehat S\frac{\widehat{\bm m}}{\widehat \rho}\big)\right>\cdot\nabla\varphi\right){\rm d}x{\rm d}y{\rm d}t,
$}
$$
with $S(\cdot,0-)=S_0$ holds for any $0\le\tau_1\le\tau_2<T$, any $\varphi\in C^1_c(\Omega\times[0,T))$, $\varphi\ge0$, where
$\mathcal{V} \equiv \{{\mathcal V}_{x,y,t}\}_{(x,y,t)\in\Omega\times(0,T)}$ is a parametrized probability (Young) measure:
$$
\begin{aligned}
&{\mathcal V}_{x,y,t} \in L^\infty(\Omega\times(0,T); {\mathcal P}(\mathbb R^4)),\quad
\mathbb R^4=\left\{\widehat\rho\in\mathbb R,\widehat{\bm m}\in\mathbb R^2,\widehat S\in\mathbb R\right\};\\
&\left<{\mathcal V}_{x,y,t};\widehat\rho\,\right>=\rho,\quad\left<{\mathcal V}_{x,y,t};\widehat{\bm m}\right>=\bm m,\quad\big<{\mathcal V}_{x,y,t};\widehat S\,\big>=S.
\end{aligned}
$$

\medskip
\noindent
$\bullet$ The {\bf energy inequality and compatibility of the energy and Reynolds stress defects.}
$$
\int\limits_\Omega E\big(\rho(\cdot,t),\bm m(\cdot,t),S(\cdot,t)\big)\,{\rm d}x{\rm d}y+\int\limits_{\overline{\Omega}}{\rm d}\mathfrak E(t)\le
\int\limits_\Omega E(\rho_0,\bm m_0,S_0)\,{\rm d}x{\rm d}y,
$$
where $\mathfrak E\in L^\infty(0,T;{\mathcal M}^+(\Omega))$ is the energy defect satisfying
$$
2\min\big\{1,\gamma-1\big\}\,\mathfrak E\le{\rm trace}[\mathfrak R]\le2\max\big\{1,\gamma-1\big\}\,\mathfrak E.
$$
\end{defn}
\begin{remark}\label{R1}
A DW solution $(\rho,\bm m,S)$ can be interpreted as the expected value of a Young measure ${\mathcal V}_{x,y,t}$, which is generated by a
consistent approximation of \eref{1.1}--\eref{1.4}; see \cite[Chapter 5]{Lukacova_book}. Note that a dissipative-strong uniqueness principle
for DW solutions holds. Consequently, as long as the strong solution exists, all DW solutions coincide with the strong solution.
\end{remark}
\begin{theorem}[\textbf{Weak versus strong convergence}]\label{T3}
Let the initial data $\{\rho_{0,\ell},\bm m_{0,\ell},E_{0,\ell}\}_{\ell=1}^\infty$ satisfy
$$
\rho_{0,\ell}>0,\quad E_{0,\ell}-\frac{|\bm m_{0,\ell}|^2}{2\rho_{0,\ell}}>0,\quad\forall\ell.
$$
Let $\{\rho_\ell,\bm m_\ell,S_\ell\}_{\ell=1}^\infty$ be a consistent approximation of \eref{1.1}--\eref{1.4} in the sense of Definition
\ref{D3} and let it be uniformly bounded, that is, there exists a positive constant $C$ such that
$$
	\|(\rho_\ell,\bm m_\ell,S_\ell)\|_{L^\infty(\Omega\times(0,T))}\le C,\quad\forall\ell.
$$
Then, a subsequence of the consistent approximation $\{\rho_{\ell_k},\bm m_{\ell_k},S_{\ell_k}\}$ generates a DW solution
$(\rho,\bm m,S)$ in the sense of Definition \ref{D2}, that is,
$$
\begin{aligned}
&(\rho_{\ell_k},\bm m_{\ell_k},S_{\ell_k})\to(\rho,\bm m,S)~\,\mbox{weakly-(*) in}~L^\infty(\Omega\times(0,T); \mathbb R^4),\\
&E(\rho_{\ell_k},\bm m_{\ell_k},S_{\ell_k})\to\big<{\mathcal V}_{x,y,t}; E(\widehat\rho,\widehat{\bm m},\widehat S)\big>~\,\mbox{weakly-(*) in}~
L^\infty(\Omega\times(0,T)),
\end{aligned}
$$
as ${\ell_k}\to\infty$. Moreover, if the Euler system admits a strong solution  that is Lipschitz continuous, that is,~
$(\rho,\bm m,E) \in W^{1,\infty}(\Omega \times (0,T))$, then
\begin{equation*}
\rho_\ell\to\rho,\quad\bm m_\ell\to\bm m,\quad E_\ell\to E\quad\mbox{in}~L^q(\,\Omega\times(0,T))~\mbox{as}~\ell\to\infty,~
\forall q\in[1,\infty).
\end{equation*}
\end{theorem}
\begin{theorem}[\textbf{${\mathcal K}$-convergence}]\label{T4}
Let the assumptions of Theorem \ref{T3} hold and $\{\rho_\ell,\bm m_\ell,S_\ell\}_{\ell=1}^\infty$ be a consistent and bounded approximation
of \eref{1.1}--\eref{1.4}. Then, there exists a subsequence $\{\rho_{\ell_k},\bm m_{\ell_k},S_{\ell_k}\}$, which converges strongly to a DW
solution $(\rho,\bm m,S)$ in the following sense:
\begin{equation*}
\resizebox{\linewidth}{!}{$
\begin{aligned}
&\big(\widetilde\rho_n,\widetilde{\bm m}_n,\widetilde S_n\big):=\frac{1}{n}\sum_{k=1}^n(\rho_{\ell_k},\bm m_{\ell_k},S_{\ell_k})\to 
(\rho,\bm m,S)~~\mbox{in}~L^q(\Omega\times(0,T);\mathbb R^4),~\forall q\in[1,\infty),\\
&\widetilde E_n:=\frac{1}{n}\sum_{k=1}^nE(\rho_{\ell_k},\bm m_{\ell_k},S_{\ell_k})\to
\big<{\mathcal V}_{x,y,t};E(\widehat\rho,\widehat{\bm m},\widehat S)\big>~~\mbox{in}~L^q(\Omega\times(0,T)),~\forall q\in[1,\infty),
\end{aligned}
$}
\end{equation*}	
as $n\to\infty$.
\end{theorem}

\section{Numerical Schemes}\label{sec3}
In this section, we provide a brief overview of the 2-D semi-discrete schemes used to conduct the numerical experiments in this paper.

We introduce uniform rectangular cells $C_{j,k}:=\big[x_\jmh,x_\jph\big]\times\big[y_\kmh,y_\kph\big]$ centered at $(x_j,y_k)$ with
$x_j=\big(x_\jmh+x_\jph\big)/2$ and $y_k=\big(y_\kph+y_\kmh\big)/2$. The spatial mesh sizes are $x_\jph-x_\jmh\equiv\dx$ and
$y_\kph-y_\kmh\equiv\dy$, respectively. Below, we will describe low-order FV methods (\S\ref{sec31}) and high-order FD A-WENO schemes
(\S\ref{sec32}).

In the FV methods, the computed discrete quantities are the cell averages, which are denoted by 
$\,\xbar{\bm U}_{\!j,k}(t):\approx\frac{1}{\dx\dy}\int_{C_{j,k}}\mU(x,y,t)\,{\rm d}x{\rm d}y$ and assumed to be available at a certain time
level $t$. Prior to evolving them in time, we will perform a conservative piecewise constant (for first-order methods) or piecewise linear
(for second-order methods) reconstruction
\begin{equation}
\widetilde{\mU}(x,y)=\,\xbar{\bm U}_{\!j,k}+(\mU_x)_{j,k}(x-x_j)+(\mU_y)_{j,k}(y-y_k),\quad(x,y)\in C_{j,k},
\label{3.1}
\end{equation}
whose one-sided point values at the midpoints of the cell interfaces are
\begin{equation}
\begin{aligned}
&\mU_{\jmh,k}^+=\,\xbar\mU_{\!j,k}-\frac{\dx}{2}(\mU_x)_{j,k},&&\mU_{\jph,k}^-=\,\xbar\mU_{\!j,k}+\frac{\dx}{2}(\mU_x)_{j,k},\\
&\mU_{j,\kmh}^+=\,\xbar\mU_{\!j,k}-\frac{\dy}{2}(\mU_y)_{j,k},&&\mU_{j,\kph}^-=\,\xbar\mU_{\!j,k}+\frac{\dy}{2}(\mU_y)_{j,k},
\end{aligned}
\label{3.2}
\end{equation}
where $(\mU_x)_{j,k}=(\mU_y)_{j,k}\equiv0$ for all $j,k$ for first-order methods, while for the second-order methods, the slopes
$(\mU_x)_{j,k}$ and $(\mU_y)_{j,k}$ are to be computed using a nonlinear limiter. In the numerical results reported in \S\ref{sec4}, we
have used a generalized minmod limiter (see, e.g., \cite{Lie03,NT,Sweby84}):
\begin{equation}
\begin{aligned}
&(\mU_x)_{j,k}={\rm minmod}\left(\theta\,\frac{\,\xbar\mU_{\!j+1,k}-\,\xbar\mU_{\!j,k}}{\dx},\,
\frac{\,\xbar\mU_{\!j+1,k}-\,\xbar\mU_{\!j-1,k}}{2\dx},\,\theta\,\frac{\,\xbar\mU_{\!j,k}-\,\xbar\mU_{\!j-1,k}}{\dx}\right),\\
&(\mU_y)_{j,k}={\rm minmod}\left(\theta\,\frac{\,\xbar\mU_{\!j,k+1}-\,\xbar\mU_{\!j,k}}{\dy},\,
\frac{\,\xbar\mU_{\!j,k+1}-\,\xbar\mU_{\!j,k-1}}{2\dy},\,\theta\,\frac{\,\xbar\mU_{\!j,k}-\,\xbar\mU_{\!j,k-1}}{\dy}\right),
\end{aligned}
\label{3.3}
\end{equation}
where the minmod function is defined by
\begin{equation*}
{\rm minmod}(z_1,z_2,\ldots)=\begin{cases}\min(z_1,z_2,\ldots)&\mbox{if}~z_i>0,\forall i,\\
\max(z_1,z_2,\ldots)&\mbox{if}~z_i<0,\forall i,\\0&\mbox{otherwise},\end{cases}
\end{equation*}
and the parameter $\theta$ was set to $1.3$.

Note that most of the indexed quantities in \eref{3.1}--\eref{3.3} as well as in formulae below are time-dependent, but we suppress this
dependence for the sake of brevity.

In the FD A-WENO schemes, the computed discrete quantities are point values $\mU_{j,k}:\approx\mU(x_j,y_k,t)$. Before evolving them in
time, we use one-dimensional WENO interpolants of an appropriate order, implemented in the local characteristic fields (see, e.g.,
\cite{Jiang13,Liu17,Nonomura20,Qiu02,Shu20} and references therein) in the $x$- and $y$-directions to obtain the point values
$\mU_{\jph,k}^\pm$ and $\mU_{j,\kph}^\pm$, respectively.

\subsection{Low-Order FV Schemes}\label{sec31}
In these schemes, the cell averages $\,\xbar\mU_{\!j,k}(t)$ are evolved in time by solving the following system of ODEs:
\begin{equation}
\frac{{\rm d}\xbar\mU_{\!j,k}}{{\rm d}t}=-\frac{\bm{{\mathcal F}}^{\,\rm FV}_{\jph,k}-\bm{{\mathcal F}}^{\,\rm FV}_{\jmh,k}}{\dx}-
\frac{\bm{{\mathcal G}}^{\rm FV}_{j,\kph}-\bm{{\mathcal G}}^{\rm FV}_{j,\kmh}}{\dy},
\label{3.4}
\end{equation}
where $\bm{{\mathcal F}}^{\,\rm FV}_{\jph,k}\big(\bm U_{\jph,k}^-,\bm U_{\jph,k}^+\big)$ and
$\bm{{\mathcal G}}^{\rm FV}_{j,\kph}\big(\bm U_{j,\kph}^-,\bm U_{j,\kph}^+\big)$ are either the local characteristic decomposition based
central-upwind (LCDCU) numerical fluxes from \cite{CCHKL_22}, low-dissipation central-upwind (LDCU) numerical fluxes from \cite{CKX_24}, or
viscous finite volume (VFV) numerical fluxes from \cite{Lukacova_book}. In the rest of this section, we provide a detailed
description of the VFV numerical fluxes, which read as $\bm{{\mathcal F}}^{\,\rm FV}_{\jph,k}=\big({\mathcal F}^{\,{\rm FV},\,\rho}_{\jph,k},
{\mathcal F}^{\,{\rm FV},\,\rho u}_{\jph,k},{\mathcal F}^{\,{\rm FV},\,\rho v}_{\jph,k},{\mathcal F}^{\,{\rm FV},\,E}_{\jph,k}\big)$ and
$\bm{{\mathcal G}}^{\rm FV}_{j,\kph}=\big({\mathcal G}^{{\rm FV},\,\rho}_{j,\kph},{\mathcal G}^{{\rm FV},\,\rho u}_{j,\kph},
{\mathcal G}^{{\rm FV},\,\rho v}_{j,\kph},{\mathcal G}^{{\rm FV},\,E}_{j,\kph}\big)$, where
\begin{equation}\resizebox{\linewidth}{!}{$
\begin{aligned}
{\mathcal F}^{\,{\rm FV},\,\rho}_{\jph,k}&={\mathcal F}^{{\rm upw},\,\rho}_{\jph,k},\quad
{\mathcal F}^{\,{\rm FV},\,\rho u}_{\jph,k}={\mathcal F}^{{\rm upw},\,\rho u}_{\jph,k}+\hf\big(p_{\jph,k}^-+p_{\jph,k}^+\big)-
\big(u_{\jph,k}^+-u_{\jph,k}^-\big),\\
{\mathcal F}^{\,{\rm FV},\,\rho v}_{\jph,k}&={\mathcal F}^{{\rm upw},\,\rho v}_{\jph,k}-\big(v_{\jph,k}^+-v_{\jph,k}^-\big),\\
{\mathcal F}^{\,{\rm FV},\,E}_{\jph,k}&={\mathcal F}^{{\rm upw},\,E}_{\jph,k}+\hf\big(u_{\jph,k}^-p_{\jph,k}^-+u_{\jph,k}^+p_{\jph,k}^+\big)\\
&-\hf\Big(\big(u_{\jph,k}^+\big)^2+\big(v_{\jph,k}^+\big)^2-\big(u_{\jph,k}^-\big)^2-\big(v_{\jph,k}^-\big)^2\Big),
\end{aligned}$}
\label{3.5f}
\end{equation}
and
\begin{equation}
\begin{aligned}
{\mathcal G}^{{\rm FV},\,\rho}_{j,\kph}&={\mathcal G}^{{\rm upw},\,\rho}_{j,\kph},\quad
{\mathcal G}^{{\rm FV},\,\rho u}_{j,\kph}={\mathcal G}^{{\rm upw},\,\rho u}_{j,\kph}-\big(u_{j,\kph}^+-u_{j,\kph}^-\big),\\
{\mathcal G}^{{\rm FV},\,\rho v}_{j,\kph}&={\mathcal G}^{{\rm upw},\,\rho v}_{j,\kph}+\hf\big(p_{j,\kph}^-+p_{j,\kph}^+\big)-
\big(v_{j,\kph}^+-v_{j,\kph}^-\big),\\
{\mathcal G}^{{\rm FV},\,E}_{j,\kph}&={\mathcal G}^{{\rm upw},\,E}_{j,\kph}+\hf\big(v_{j,\kph}^-p_{j,\kph}^-+v_{j,\kph}^+p_{j,\kph}^+\big)\\
&-\hf\Big(\big(u_{j,\kph}^+\big)^2+\big(v_{j,\kph}^+\big)^2-\big(u_{j,\kph}^-\big)^2-\big(v_{j,\kph}^-\big)^2\Big).
\end{aligned}
\label{3.6f}
\end{equation}
In the above formulae,
\begin{equation}\resizebox{0.9\linewidth}{!}{$
\begin{aligned}
\bm{{\mathcal F}}^{\rm upw}_{\jph,k}&=\frac{1}{4}\big(u_{\jph,k}^-+u_{\jph,k}^+-\big|u_{\jph,k}^-+u_{\jph,k}^+\big|-4\big)
\big(\bm U_{\jph,k}^-+\bm U_{\jph,k}^+\big),\\
\bm{{\mathcal G}}^{\rm upw}_{j,\kph}&=\frac{1}{4}\big(v_{j,\kph}^-+v_{j,\kph}^+-\big|v_{j,\kph}^-+v_{j,\kph}^+\big|-4\big)
\big(\bm U_{j,\kph}^-+\bm U_{j,\kph}^+\big),
\end{aligned}
\label{3.7f}$}
\end{equation}
$$
\begin{aligned}
&u_{\jph,k}^\pm=\frac{(\rho u)_{\jph,k}^\pm}{\rho_{\jph,k}^\pm},~~\!v_{\jph,k}^\pm=\frac{(\rho v)_{\jph,k}^\pm}{\rho_{\jph,k}^\pm},\\
&p_{\jph,k}^\pm=(\gamma-1)\!\left[E_{\jph,k}^\pm-
\frac{\big((\rho u)_{\jph,k}^\pm\big)^2+\big((\rho v)_{\jph,k}^\pm\big)^2}{2\rho_{\jph,k}^\pm}\right],\\
&u_{j,\kph}^\pm=\frac{(\rho u)_{j,\kph}^\pm}{\rho_{j,\kph}^\pm},~~\!v_{j,\kph}^\pm=\frac{(\rho v)_{j,\kph}^\pm}{\rho_{j,\kph}^\pm},\\
&\!p_{j,\kph}^\pm=(\gamma-1)\!\left[E_{j,\kph}^\pm-
\frac{\big((\rho u)_{j,\kph}^\pm\big)^2+\big((\rho v)_{j,\kph}^\pm\big)^2}{2\rho_{j,\kph}^\pm}\right].
\end{aligned}
$$
\begin{remark}
We stress that the first-order VFV scheme, which according to \cite{feireisl2020finite} produces consistent approximations, is slightly
different from the scheme \eref{3.4}--\eref{3.7f}. The difference is in the evolution of $\xbar E_{j,k}$, which, in the first-order
scheme, reads as
$$
\begin{aligned}
\frac{{\rm d}\xbar E_{j,k}}{{\rm d}t}=&-\frac{{\mathcal F}^{{\rm upw},\,E}_{\jph,k}-{\mathcal F}^{{\rm upw},\,E}_{\jmh,k}}{\dx}-
\frac{{\mathcal G}^{{\rm upw},\,E}_{j,\kph}-{\mathcal G}^{{\rm upw},\,E}_{j,\kmh}}{\dy}\\
&-u_{j,k}\frac{p_{j+1,k}-p_{j-1,k}}{2\dx}-p_{j,k}\frac{u_{j+1,k}-u_{j-1,k}}{2\dx}\\
&-v_{j,k}\frac{p_{j,k+1}-p_{j,k-1}}{2\dy}-p_{j,k}\frac{v_{j,k+1}-v_{j,k-1}}{2\dy}\\
&+\frac{(u_{j+1,k}^+)^2-2(u_{j,k}^+)^2+(u_{j-1,k}^+)^2+(v_{j+1,k}^+)^2-2(v_{j,k}^+)^2+(v_{j-1,k}^+)^2}{2\dx}\\
&+\frac{(u_{j,k+1}^+)^2-2(u_{j,k}^+)^2+(u_{j,k-1}^+)^2+(v_{j,k+1}^+)^2-2(v_{j,k}^+)^2+(v_{j,k-1}^+)^2}{2\dy}.
\end{aligned}
$$
\end{remark}

\subsection{High-Order FD A-WENO Schemes}\label{sec32}
In this section, we present higher-order extensions of the low-order schemes from \S\ref{sec31} within the FD A-WENO framework introduced in
\cite{Jiang13}; see also \cite{CKX23,Liu17,Gao20}.

In FD A-WENO schemes, computed solutions are realized in terms of the point values $\mU_{j,k}$ (rather than the cell averages
$\,\xbar\mU_{\!j,k}$), which are evolved in time by solving the following system of ODEs:
\begin{equation}
\frac{{\rm d}\mU_{j,k}}{{\rm d}t}=-\frac{{\bm{{\mathcal F}}_{\jph,k}}-{\bm{{\mathcal F}}_{\jmh,k}}}{\dx}-
\frac{{\bm{{\mathcal G}}_{j,\kph}}-{\bm{{\mathcal G}}_{j,\kmh}}}{\dy},
\label{3.5}
\end{equation}
where the numerical fluxes $\bm{{\mathcal F}}_{\jph,k}$ and $\bm{{\mathcal G}}_{j,\kph}$ are defined by
\begin{equation*}\resizebox{\linewidth}{!}{$
\begin{aligned}
\bm{{\mathcal F}}_{\jph,k}&=\bm{{\mathcal F}}^{\,\rm FV}_{\jph,k}-\mu_2(\dx)^2(\mF_{xx})_{\jph,k}+H(r-4)\mu_4(\dx)^4(\mF_{xxxx})_{\jph,k}\\
&-H(r-6)\mu_6(\dx)^6(\mF_{xxxxxx})_{\jph,k}+H(r-8)\mu_8(\dx)^8(\mF_{xxxxxxxx})_{\jph,k},\\
\bm{{\mathcal G}}_{j,\kph}&=\bm{{\mathcal G}}^{\rm FV}_{j,\kph}-\mu_2(\dy)^2(\mG_{yy})_{j,\kph}+H(r-4)\mu_4(\dy)^4(\mG_{yyyy})_{j,\kph}\\
&-H(r-6)\mu_6(\dy)^6(\mG_{yyyyyy})_{j,\kph}+H(r-8)\mu_8(\dy)^8(\mG_{yyyyyyyy})_{j,\kph}.
\end{aligned}$}
\end{equation*}
Here, $\mu_2=1/24$, $\mu_4=7/5760$, $\mu_6=31/967680$, $\mu_8=127/154828800$, $H$ is the Heaviside function, $r=3$, $5$, $7$, or $9$ is the
order of the scheme, $\bm{{\mathcal F}}^{\,\rm FV}_{\jph,k}\big(\bm U_{\jph,k}^-,\bm U_{\jph,k}^+\big)$ and
$\bm{{\mathcal G}}^{\rm FV}_{j,\kph}\big(\bm U_{j,\kph}^-,\bm U_{j,\kph}^+\big)$ are the FV numerical fluxes, which are computed using the
one-sided point values $\bm U_{\jph,k}^\pm$, which are obtained using an $r$th-order WENO interpolation applied to the local characteristic
variables, and $(\mF_{xx})_{\jph,k}$, $(\mF_{xxxx})_{\jph,k}$, $(\mF_{xxxxxx})_{\jph,k}$, $(\mF_{xxxxxxxx})_{\jph,k}$,
$(\mG_{yy})_{j,\kph}$, $(\mG_{yyyy})_{j,\kph}$, $(\mG_{yyyyyy})_{j,\kph}$, and \linebreak $(\mG_{yyyyyyyy})_{j,\kph}$ are the higher-order correction
terms, which can be computed with the help of the FV numerical fluxes $\bm{{\mathcal F}}^{\,\rm FV}_{\jph,k}$ and
$\bm{{\mathcal G}}^{\rm FV}_{j,\kph}$ as it was first done in \cite{CKX23} for the fifth-order A-WENO schemes. The formulae for these terms
depend on the order of the scheme and for different values of $r$ they are given below together with the details on particular WENO
interpolations used for each $r$.\\

\noindent{\bf $r=3$ (third-order schemes)}
\allowdisplaybreaks
\begin{equation*}
\begin{aligned}
&(\mF_{xx})_{\jph,k}=\frac{1}{(\dx)^2}\Big[\bm{{\mathcal F}}^{\,\rm FV}_{\jmh,k}-2\bm{{\mathcal F}}^{\,\rm FV}_{\jph,k}+
\bm{{\mathcal F}}^{\,\rm FV}_{j+\frac{3}{2},k}\Big],\\
&(\mG_{yy})_{j,\kph}=\frac{1}{(\dy)^2}\Big[\bm{{\mathcal G}}^{\rm FV}_{j,\kmh}-
2\bm{{\mathcal G}}^{\rm FV}_{j,\kph}+\bm{{\mathcal G}}^{\rm FV}_{j,k+\frac{3}{2}}\Big].
\end{aligned}
\end{equation*}
The one-sided point values $\mU^\pm_{\jph,k}$ and $\mU^\pm_{j,\kph}$ are computed using the third-order WENO interpolation introduced in
\cite{CH_Third,CHT25}; for details, we refer the reader to \cite[Appendix A]{CHT25}.\\

\noindent\textbf{$r=5$ (fifth-order schemes)} 
\allowdisplaybreaks
$(\mF_{xx})_{\jph,k}$, $(\mF_{xxxx})_{\jph,k}$, $(\mG_{yy})_{j,\kph}$, and $(\mG_{yyyy})_{j,\kph}$
are defined in \cite[(4.4)]{CKX23}, and the one-sided point values $\mU^\pm_{\jph,k}$ and $\mU^\pm_{j,\kph}$ are computed using the
fifth-order WENO-Z interpolation from \cite{Borges08,Castro11,Don13}; for details, we refer the reader to \cite[Appendix A]{CCK23_Adaptive}.\\

\noindent\textbf{$r=7$ (seventh-order schemes)}
\allowdisplaybreaks
\begin{align*}
(\mF_{xx})_{\jph,k}=\frac{1}{180(\dx)^2}\Big[&2\bm{{\mathcal F}}^{\,\rm FV}_{j-\frac{5}{2},k}-27\bm{{\mathcal F}}^{\,\rm FV}_{j-\frac{3}{2},k}+
270\bm{{\mathcal F}}^{\,\rm FV}_{\jmh,k}\\
&-490\bm{{\mathcal F}}^{\,\rm FV}_{\jph,k}+270\bm{{\mathcal F}}^{\,\rm FV}_{j+\frac{3}{2},k}-27\bm{{\mathcal F}}^{\,\rm FV}_{j+\frac{5}{2},k}+
2\bm{{\mathcal F}}^{\,\rm FV}_{j+\frac{7}{2},k}\Big],\\
(\mF_{xxxx})_{\jph,k}=\frac{1}{6(\dx)^4}\Big[&-\bm{{\mathcal F}}^{\,\rm FV}_{j-\frac{5}{2},k}+12\bm{{\mathcal F}}^{\,\rm FV}_{j-\frac{3}{2},k}-
39\bm{{\mathcal F}}^{\,\rm FV}_{\jmh,k}\\
&+56\bm{{\mathcal F}}^{\,\rm FV}_{\jph,k}-39\bm{{\mathcal F}}^{\,\rm FV}_{j+\frac{3}{2},k}+12\bm{{\mathcal F}}^{\,\rm FV}_{j+\frac{5}{2},k}
-\bm{{\mathcal F}}^{\,\rm FV}_{j+\frac{7}{2},k}\Big],\\
(\mF_{xxxxxx})_{\jph,k}=\frac{1}{(\dx)^6}\Big[&\bm{{\mathcal F}}^{\,\rm FV}_{j-\frac{5}{2},k}-6\bm{{\mathcal F}}^{\,\rm FV}_{j-\frac{3}{2},k}+
15\bm{{\mathcal F}}^{\,\rm FV}_{\jmh,k}\\
&-20\bm{{\mathcal F}}^{\,\rm FV}_{\jph,k}+15\bm{{\mathcal F}}^{\,\rm FV}_{j+\frac{3}{2},k}-6\bm{{\mathcal F}}^{\,\rm FV}_{j+\frac{5}{2},k}+
\bm{{\mathcal F}}^{\,\rm FV}_{j+\frac{7}{2},k}\Big],\\
(\mG_{yy})_{j,\kph}=\frac{1}{180(\dy)^2}\Big[&2\bm{{\mathcal G}}^{\rm FV}_{j,k-\frac{5}{2}}-27\bm{{\mathcal G}}^{\rm FV}_{j,k-\frac{3}{2}}+
270\bm{{\mathcal G}}^{\rm FV}_{j,\kmh}\\
&-490\bm{{\mathcal G}}^{\rm FV}_{j,\kph}+270\bm{{\mathcal G}}^{\rm FV}_{j,k+\frac{3}{2}}-27\bm{{\mathcal G}}^{\rm FV}_{j,k+\frac{5}{2}}+
2\bm{{\mathcal G}}^{\rm FV}_{j,k+\frac{7}{2}}\Big],\\
(\mG_{yyyy})_{j,\kph}=\frac{1}{6(\dy)^4}\Big[&-\bm{{\mathcal G}}^{\rm FV}_{j,k-\frac{5}{2}}+12\bm{{\mathcal G}}^{\rm FV}_{j,k-\frac{3}{2}}-
39\bm{{\mathcal G}}^{\rm FV}_{j,\kmh}\\
&+56\bm{{\mathcal G}}^{\rm FV}_{j,\kph}-39\bm{{\mathcal G}}^{\rm FV}_{j,k+\frac{3}{2}}+12\bm{{\mathcal G}}^{\rm FV}_{j,k+\frac{5}{2}}-
\bm{{\mathcal G}}^{\rm FV}_{j,k+\frac{7}{2}}\Big],\\
(\mG_{yyyyyy})_{j,\kph}=\frac{1}{(\dy)^6}\Big[&\bm{{\mathcal G}}^{\rm FV}_{j,k-\frac{5}{2}}-6\bm{{\mathcal G}}^{\rm FV}_{j,k-\frac{3}{2}}+15\bm{{\mathcal G}}^{\rm FV}_{j,\kmh}-20\bm{{\mathcal G}}^{\rm FV}_{j,\kph}\\
&+15\bm{{\mathcal G}}^{\rm FV}_{j,k+\frac{3}{2}}-6\bm{{\mathcal G}}^{\rm FV}_{j,k+\frac{5}{2}}+\bm{{\mathcal G}}^{\rm FV}_{j,k+\frac{7}{2}}\Big].
\end{align*}
The one-sided point values $\mU^\pm_{\jph,k}$ and $\mU^\pm_{j,\kph}$ are computed using the seventh-order WENO-Z interpolation from
\cite{Gao20}.\\

\noindent\textbf{$r=9$ (ninth-order schemes)}
\allowdisplaybreaks
\begin{align*}
(\mF_{xx})_{\jph,k}=\frac{1}{5040(\dx)^2}\Big[&-9\bm{{\mathcal F}}^{\,\rm FV}_{j-\frac{7}{2},k}+128\bm{{\mathcal F}}^{\,\rm FV}_{j-\frac{5}{2},k}-
1008\bm{{\mathcal F}}^{\,\rm FV}_{j-\frac{3}{2},k}\\
&+8064\bm{{\mathcal F}}^{\,\rm FV}_{\jmh,k}-14350\bm{{\mathcal F}}^{\,\rm FV}_{\jph,k}+8064\bm{{\mathcal F}}^{\,\rm FV}_{j+\frac{3}{2},k}\\
&-1008\bm{{\mathcal F}}^{\,\rm FV}_{j+\frac{5}{2},k}+
128\bm{{\mathcal F}}^{\,\rm FV}_{j+\frac{7}{2},k}-9\bm{{\mathcal F}}^{\,\rm FV}_{j+\frac{9}{2},k}\Big],\\
(\mF_{xxxx})_{\jph,k}=\frac{1}{240(\dx)^4}\Big[&7\bm{{\mathcal F}}^{\,\rm FV}_{j-\frac{7}{2},k}-96\bm{{\mathcal F}}^{\,\rm FV}_{j-\frac{5}{2},k}+
676\bm{{\mathcal F}}^{\,\rm FV}_{j-\frac{3}{2},k}-1952\bm{{\mathcal F}}^{\,\rm FV}_{\jmh,k}\\
&+2730\bm{{\mathcal F}}^{\,\rm FV}_{\jph,k}-1952\bm{{\mathcal F}}^{\,\rm FV}_{j+\frac{3}{2},k}+676\bm{{\mathcal F}}^{\,\rm FV}_{j+\frac{5}{2},k}\\
&-96\bm{{\mathcal F}}^{\,\rm FV}_{j+\frac{7}{2},k}+7\bm{{\mathcal F}}^{\,\rm FV}_{j+\frac{9}{2},k}\Big],\\
(\mF_{xxxxxx})_{\jph,k}=\frac{1}{4(\dx)^6}\Big[&-\bm{{\mathcal F}}^{\,\rm FV}_{j-\frac{7}{2},k}+12\bm{{\mathcal F}}^{\,\rm FV}_{j-\frac{5}{2},k}-
52\bm{{\mathcal F}}^{\,\rm FV}_{j-\frac{3}{2},k}+116\bm{{\mathcal F}}^{\,\rm FV}_{\jmh,k}\\
&-150\bm{{\mathcal F}}^{\,\rm FV}_{\jph,k}+116\bm{{\mathcal F}}^{\,\rm FV}_{j+\frac{3}{2},k}-52\bm{{\mathcal F}}^{\,\rm FV}_{j+\frac{5}{2},k}\\
&+12\bm{{\mathcal F}}^{\,\rm FV}_{j+\frac{7}{2},k}-\bm{{\mathcal F}}^{\,\rm FV}_{j+\frac{9}{2},k}\Big],\\
(\mF_{xxxxxxxx})_{\jph,k}=\frac{1}{(\dx)^8}\Big[&\bm{{\mathcal F}}^{\,\rm FV}_{j-\frac{7}{2},k}-8\bm{{\mathcal F}}^{\,\rm FV}_{j-\frac{5}{2},k}+
28\bm{{\mathcal F}}^{\,\rm FV}_{j-\frac{3}{2},k}-56\bm{{\mathcal F}}^{\,\rm FV}_{\jmh,k}\\
&+70\bm{{\mathcal F}}^{\,\rm FV}_{\jph,k}-56\bm{{\mathcal F}}^{\,\rm FV}_{j+\frac{3}{2},k}+28\bm{{\mathcal F}}^{\,\rm FV}_{j+\frac{5}{2},k}-
8\bm{{\mathcal F}}^{\,\rm FV}_{j+\frac{7}{2},k}\\
&+\bm{{\mathcal F}}^{\,\rm FV}_{j+\frac{9}{2},k}\Big],\\
(\mG_{yy})_{j,\kph}=\frac{1}{5040(\dy)^2}\Big[&-9\bm{{\mathcal G}}^{\rm FV}_{j,k-\frac{7}{2}}+128\bm{{\mathcal G}}^{\rm FV}_{j,k-\frac{5}{2}}-
1008\bm{{\mathcal G}}^{\rm FV}_{j,k-\frac{3}{2}}\\
&+8064\bm{{\mathcal G}}^{\rm FV}_{j,\kmh}-14350\bm{{\mathcal G}}^{\rm FV}_{j,\kph}+8064\bm{{\mathcal G}}^{\rm FV}_{j,k+\frac{3}{2}}\\
&-1008\bm{{\mathcal G}}^{\rm FV}_{j,k+\frac{5}{2}}+128\bm{{\mathcal G}}^{\rm FV}_{j,k+\frac{7}{2}}-9\bm{{\mathcal G}}^{\rm FV}_{j,k+\frac{9}{2}}\Big],\\
(\mG_{yyyy})_{j,\kph}=\frac{1}{240(\dy)^4}\Big[&7\bm{{\mathcal G}}^{\rm FV}_{j,k-\frac{7}{2}}-96\bm{{\mathcal G}}^{\rm FV}_{j,k-\frac{5}{2}}+
676\bm{{\mathcal G}}^{\rm FV}_{j,k-\frac{3}{2}}-1952\bm{{\mathcal G}}^{\rm FV}_{j,\kmh}\\
&+2730\bm{{\mathcal G}}^{\rm FV}_{j,\kph}-1952\bm{{\mathcal G}}^{\rm FV}_{j,k+\frac{3}{2}}+676\bm{{\mathcal G}}^{\rm FV}_{j,k+\frac{5}{2}}\\
&-96\bm{{\mathcal G}}^{\rm FV}_{j,k+\frac{7}{2}}+7\bm{{\mathcal G}}^{\rm FV}_{j,k+\frac{9}{2}}\Big],\\
(\mG_{yyyyyy})_{j,\kph}=\frac{1}{4(\dy)^6}\Big[&-\bm{{\mathcal G}}^{\rm FV}_{j,k-\frac{7}{2}}+12\bm{{\mathcal G}}^{\rm FV}_{j,k-\frac{5}{2}}-
52\bm{{\mathcal G}}^{\rm FV}_{j,k-\frac{3}{2}}+116\bm{{\mathcal G}}^{\rm FV}_{j,\kmh}\\
&-150\bm{{\mathcal G}}^{\rm FV}_{j,\kph}+116\bm{{\mathcal G}}^{\rm FV}_{j,k+\frac{3}{2}}-52\bm{{\mathcal G}}^{\rm FV}_{j,k+\frac{5}{2}}\\
&+12\bm{{\mathcal G}}^{\rm FV}_{j,k+\frac{7}{2}}-\bm{{\mathcal G}}^{\rm FV}_{j,k+\frac{9}{2}}\Big],\\
(\mG_{yyyyyyyy})_{j,\kph}=\frac{1}{(\dy)^8}\Big[&\bm{{\mathcal G}}^{\rm FV}_{j,k-\frac{7}{2}}-8\bm{{\mathcal G}}^{\rm FV}_{j,k-\frac{5}{2}}+
28\bm{{\mathcal G}}^{\rm FV}_{j,k-\frac{3}{2}}-56\bm{{\mathcal G}}^{\rm FV}_{j,\kmh}\\
&+70\bm{{\mathcal G}}^{\rm FV}_{j,\kph}-56\bm{{\mathcal G}}^{\rm FV}_{j,k+\frac{3}{2}}+28\bm{{\mathcal G}}^{\rm FV}_{j,k+\frac{5}{2}}-
8\bm{{\mathcal G}}^{\rm FV}_{j,k+\frac{7}{2}}\\
&+\bm{{\mathcal G}}^{\rm FV}_{j,k+\frac{9}{2}}\Big].
\end{align*}
The one-sided point values $\mU^\pm_{\jph,k}$ and $\mU^\pm_{j,\kph}$ are computed using the ninth-order WENO-Z interpolation from
\cite{Gao20}.

\section{Numerical Results}\label{sec4}
In this section, we present the numerical simulations conducted using the schemes described in \S\ref{sec3} for three 2-D Riemann problems
and the KH instability problem. 

In all of the examples, we take a uniform mesh with $\dx=\dy=1/1024$ on the computational domain $[0,1]\times[0,1]$. We numerically
integrate the ODE systems \eref{3.4} and  \eref{3.5} using the three-stage third-order strong stability preserving Runge-Kutta method (see,
e.g., \cite{Gottlieb11,Gottlieb12}) and use the CFL number $0.45$ for the schemes employing the LCDCU and LDCU numerical fluxes and $0.1$,
when the VFV numerical fluxes are utilized.

Below, we will use the following notation for the solutions computed using the same numerical flux, but with different orders of the
resulting scheme: $(\rho_\ell,\bm m_\ell,S_\ell)$, $\ell(r)=1,\dots,6$ are solutions computed by the first-, second-, third-, fifth-,
seventh-, and ninth-order schemes, respectively, that is, $\ell(1)=1$, $\ell(2)=2$, $\ell(3)=3$, $\ell(5)=4$, $\ell(7)=5$, and $\ell(9)=6$.
In accordance with the concept of ${\mathcal K}$-convergence described in Theorem \ref{T4}, we will average  these solutions to obtain
$(\widetilde\rho_n,\widetilde{\bm m}_n,\widetilde S_n)^\top:=\frac{1}{n}\sum_{\ell=1}^n(\rho_\ell,\bm m_\ell,S_\ell)^\top$, where
$n=1,\dots,6$.

In some of the examples below, oscillatory solutions are to be captured. As explained in Remark \ref{R1}, such solutions should be
represented as a Young measure in the phase space. In fact, when the solution is a weak solution, the corresponding Young measure is
expected to be a Dirac $\delta$-function. In order to check whether the obtained solution is a genuine DW solution or a weak solution in the
sense of distributions, we will approximate the probability density function (PDF) of the Young measure $\{{\mathcal V}_{x,y,t}\}$ as follows. 

We take a certain small rectangular subdomain $\widetilde\Omega$ and compute 
$$
\underline\rho:=\min\limits_{(x_j,y_k)\in\widetilde\Omega}\rho_{j,k}\quad\mbox{and}\quad
\xbar\rho:=\max\limits_{(x_j,y_k)\in\widetilde\Omega}\rho_{j,k}.
$$
We then introduce $\Delta\rho:=(\xbar\rho-\underline\rho)/30$ and the following partition of the interval $[\underline\rho,\xbar\rho]$:
$$
\resizebox{\linewidth}{!}{$
[\underline\rho,\xbar\rho]=\bigcup\limits_{i=1}^{30}\big[\rho^{(\imh)},\rho^{(\iph)}\big),~~\mbox{where}~~
\rho^{(\hf)}=\underline\rho,~~\rho^{(\iph)}=\rho^{(\imh)}+\Delta\rho,~~i=2,\dots,30.
$}
$$
The distribution of the dataset formed by the values $\rho_{j,k}$ across all $C_{j,k}\subseteq\widetilde\Omega$ is represented by a function
$\sigma_\ell=\sigma_\ell(\rho)$, where the subscript $\ell=\ell(r)$ indicates that the density values have been computed using the $r$-th
order scheme. More specifically, $\sigma_\ell(\rho^{(i)})$ denotes the number of cells $C_{j,k}$ for which the computed values $\rho_{j,k}$
fall within the interval $\big[\rho^{(\imh)},\rho^{(\iph)}\big)$. Additionally, we apply a normalization such that 
\begin{equation*}
\Delta\rho\sum\limits_{i=1}^{30}\sigma_\ell\big(\rho^{(i)}\big)=1. 
\end{equation*}
In order to obtain a strong convergence, we average the sequence $\{\sigma_\ell\}$ and introduce
$\widetilde\sigma_n:=\frac{1}{n}\sum\limits_{\ell=1}^n\sigma_\ell$, $n=1,\dots,6$.
 
\subsection{2-D Riemann Problems}\label{sec41}
We begin by considering three 2-D Riemann problems from \cite{Kurganov02} (see also \cite{Schulz93,Schulz93a,Zheng01}): Configurations 2--4.
The details on the initial data, subject to the free boundary conditions, are specified below.\\

\textbf{Configuration 2.} In the first example, the initial data are
\begin{equation*}
(\rho,u,v,p)(x,y,0)=\begin{cases}
(1,0,0,1),&x>0.5,~y>0.5,\\
(0.5197,-0.7259,0,0.4),&x<0.5,~y>0.5,\\
(1,-0.7259,-0.7259,1),&x<0.5,~y<0.5,\\
(0.5197,0,-0.7259,0.4),&x>0.5,~y<0.5,
\end{cases}
\end{equation*}
and far from the center of the computational domain, the solution consists of four rarefaction waves.

\begin{figure}[ht!]
\centerline{\includegraphics[trim=1.3cm 0.5cm 1.6cm 0.2cm, clip, width=4.cm]{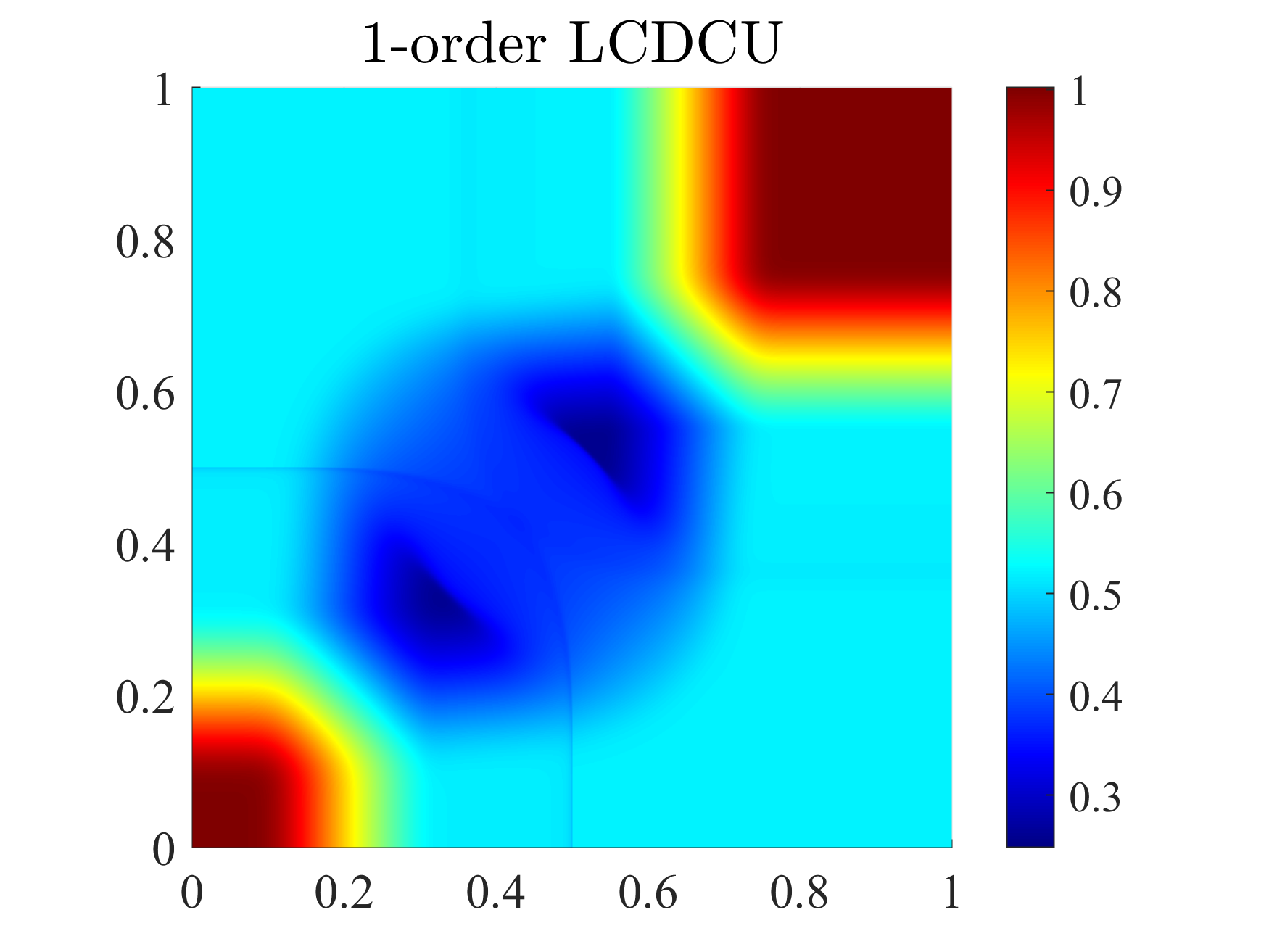}\hspace*{0.2cm}
            \includegraphics[trim=1.3cm 0.5cm 1.6cm 0.2cm, clip, width=4.cm]{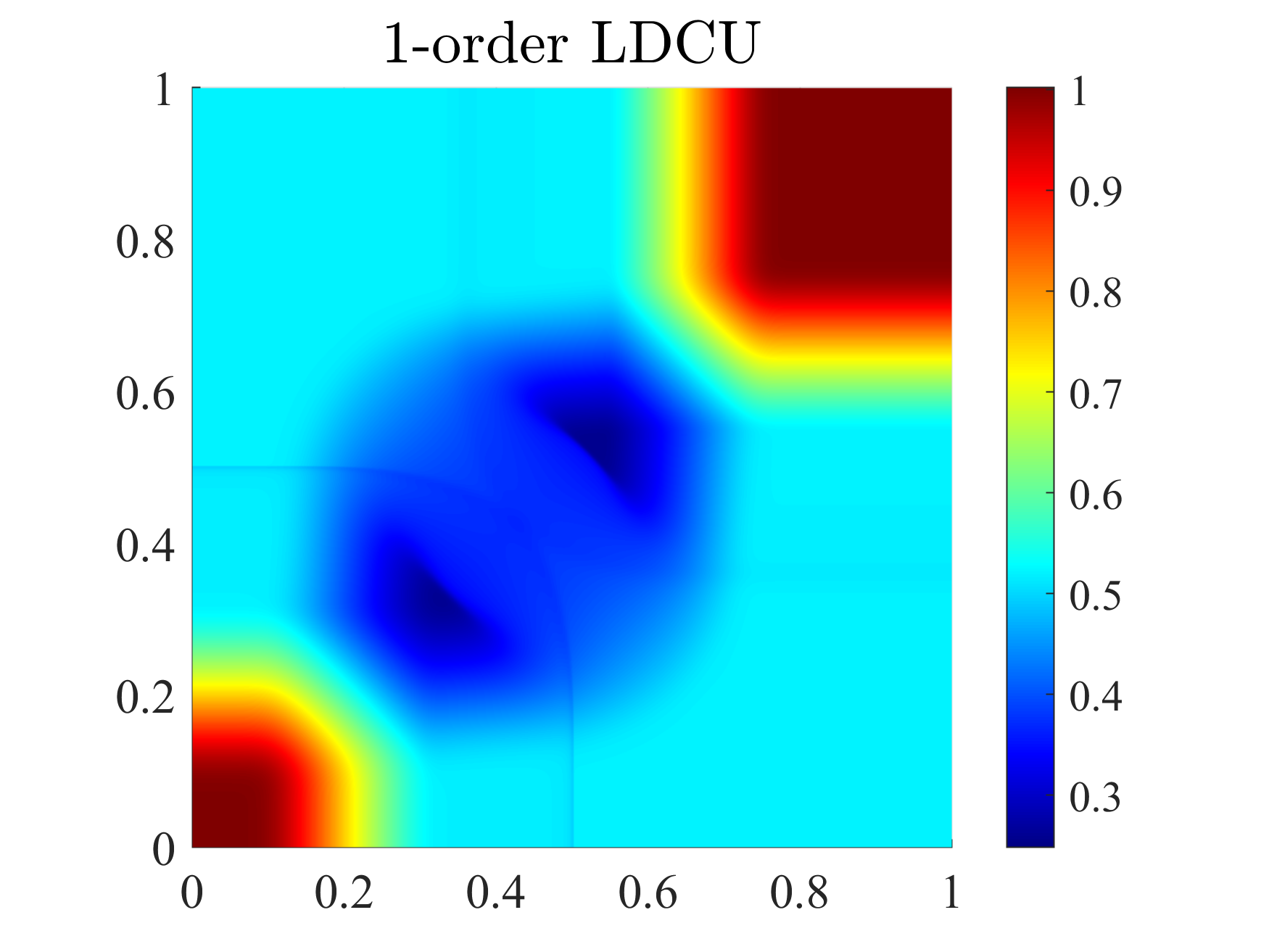}\hspace*{0.0cm}
            \includegraphics[trim=1.3cm 0.5cm 1.6cm 0.2cm, clip, width=4.2cm]{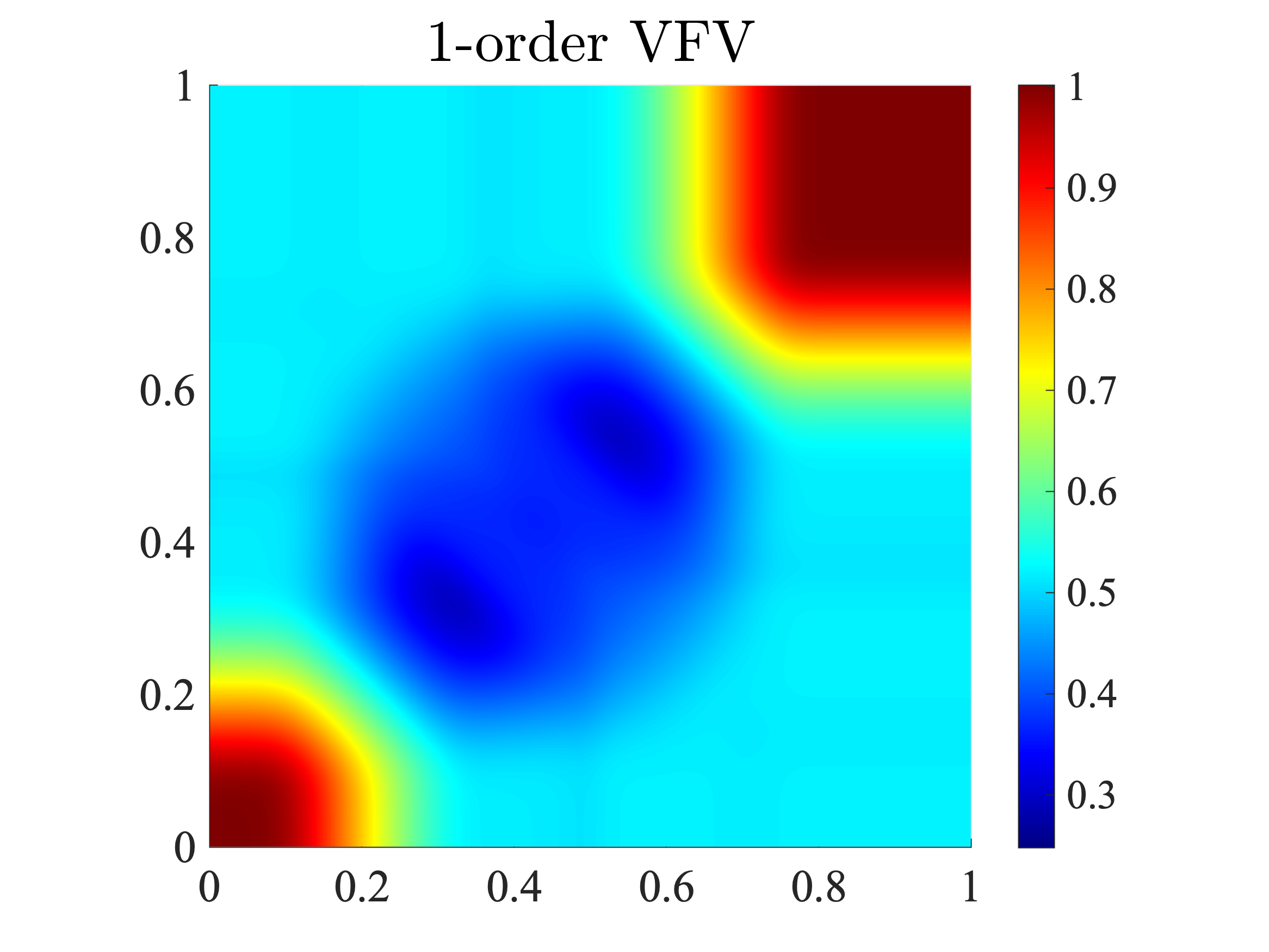}}
\vskip8pt
\centerline{\includegraphics[trim=1.3cm 0.5cm 1.6cm 0.2cm, clip, width=4.cm]{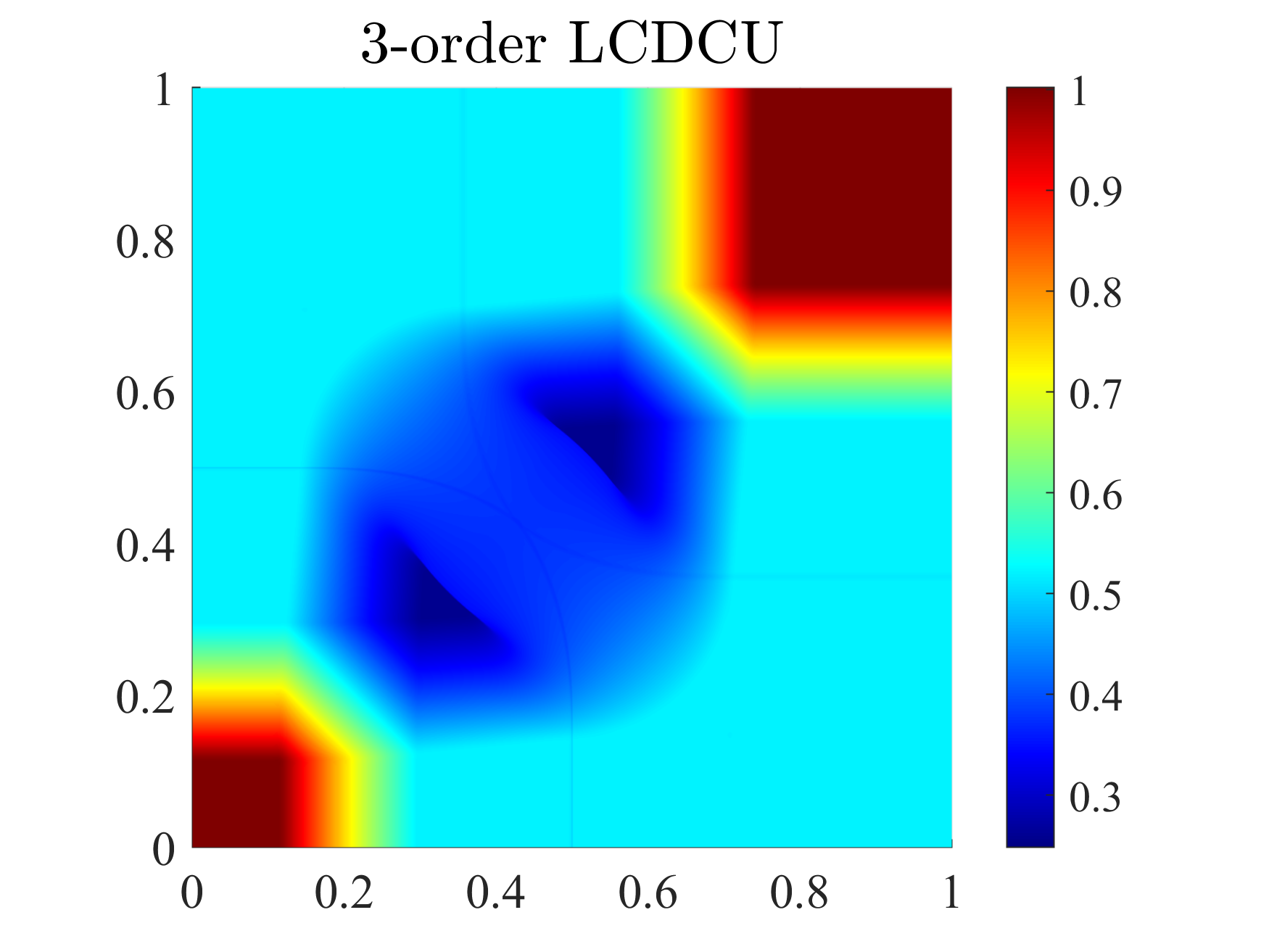}\hspace*{0.2cm}
	        \includegraphics[trim=1.3cm 0.5cm 1.6cm 0.2cm, clip, width=4.cm]{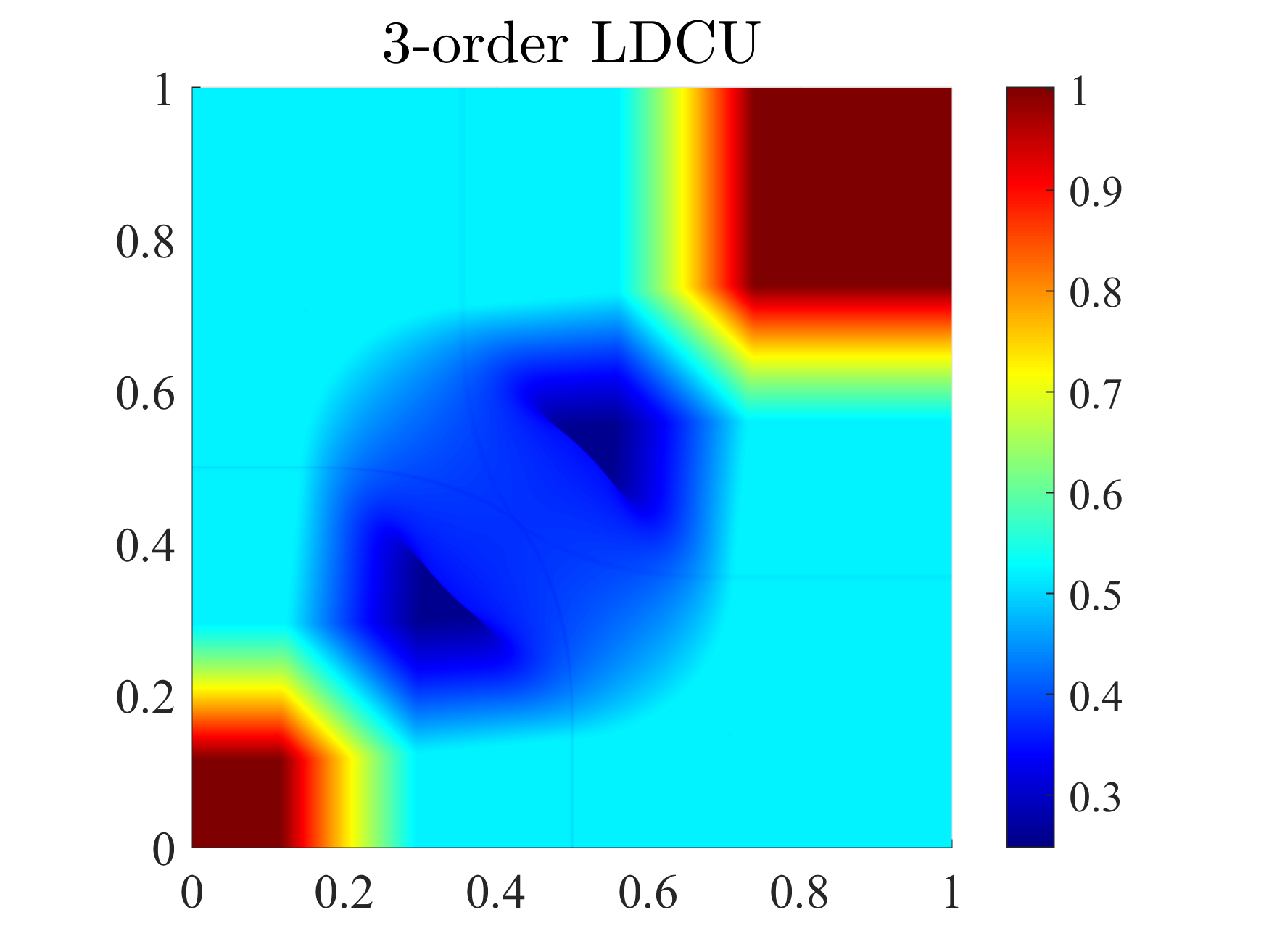}\hspace*{0.0cm}
	        \includegraphics[trim=1.3cm 0.5cm 1.6cm 0.2cm, clip, width=4.2cm]{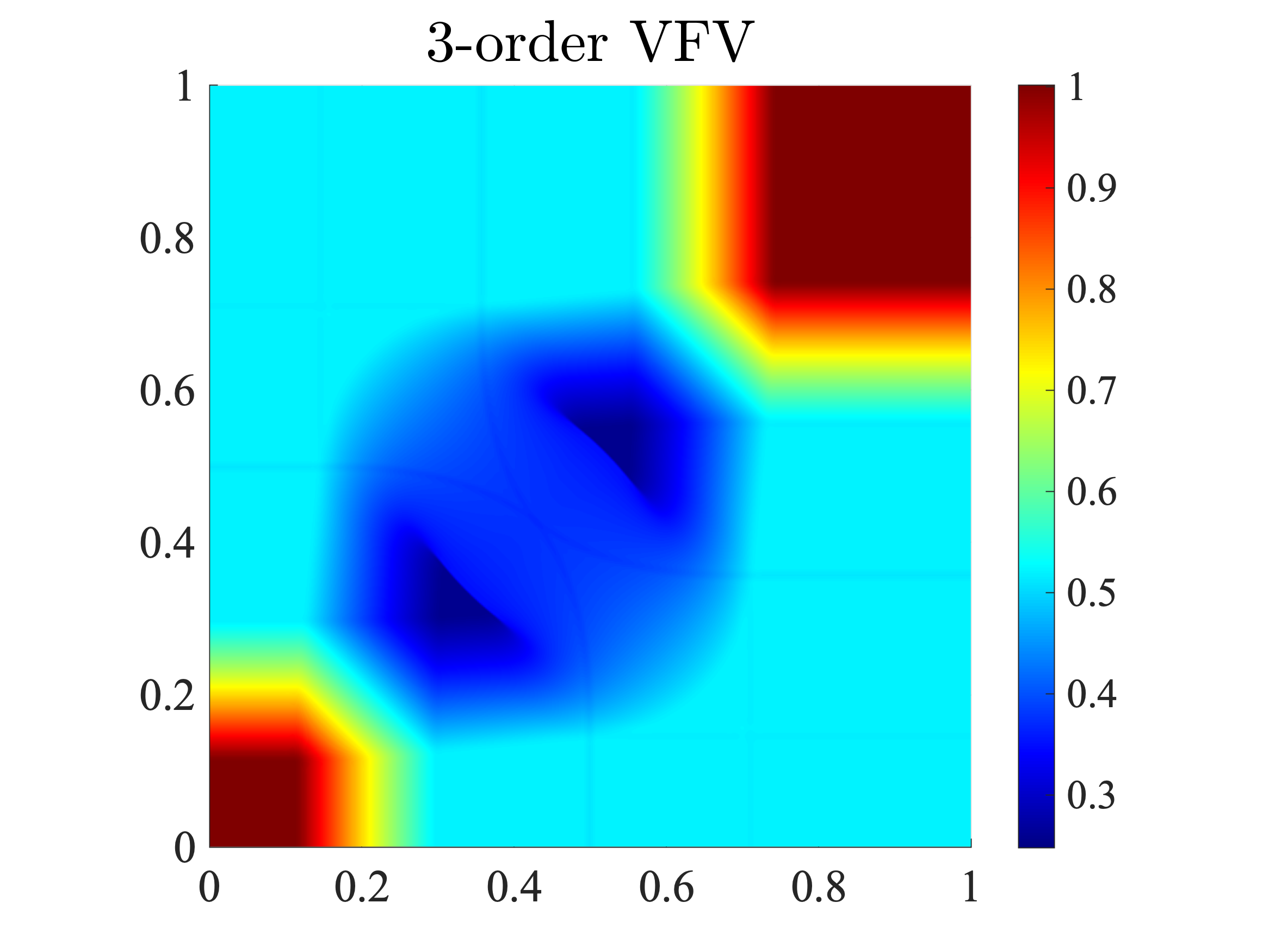}}
\vskip8pt
\centerline{\includegraphics[trim=1.3cm 0.5cm 1.6cm 0.2cm, clip, width=4.cm]{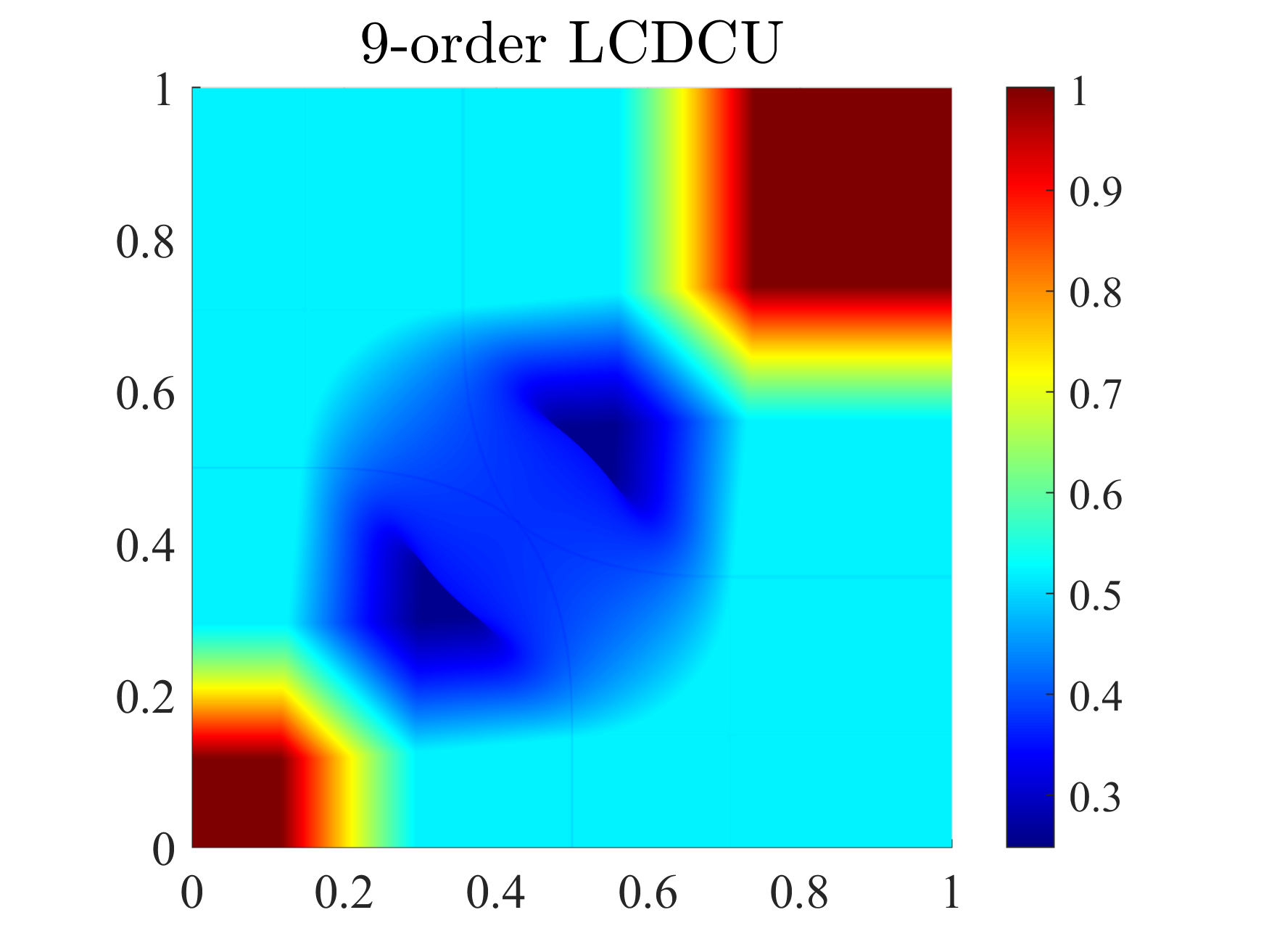}\hspace*{0.2cm}
	        \includegraphics[trim=1.3cm 0.5cm 1.6cm 0.2cm, clip, width=4.cm]{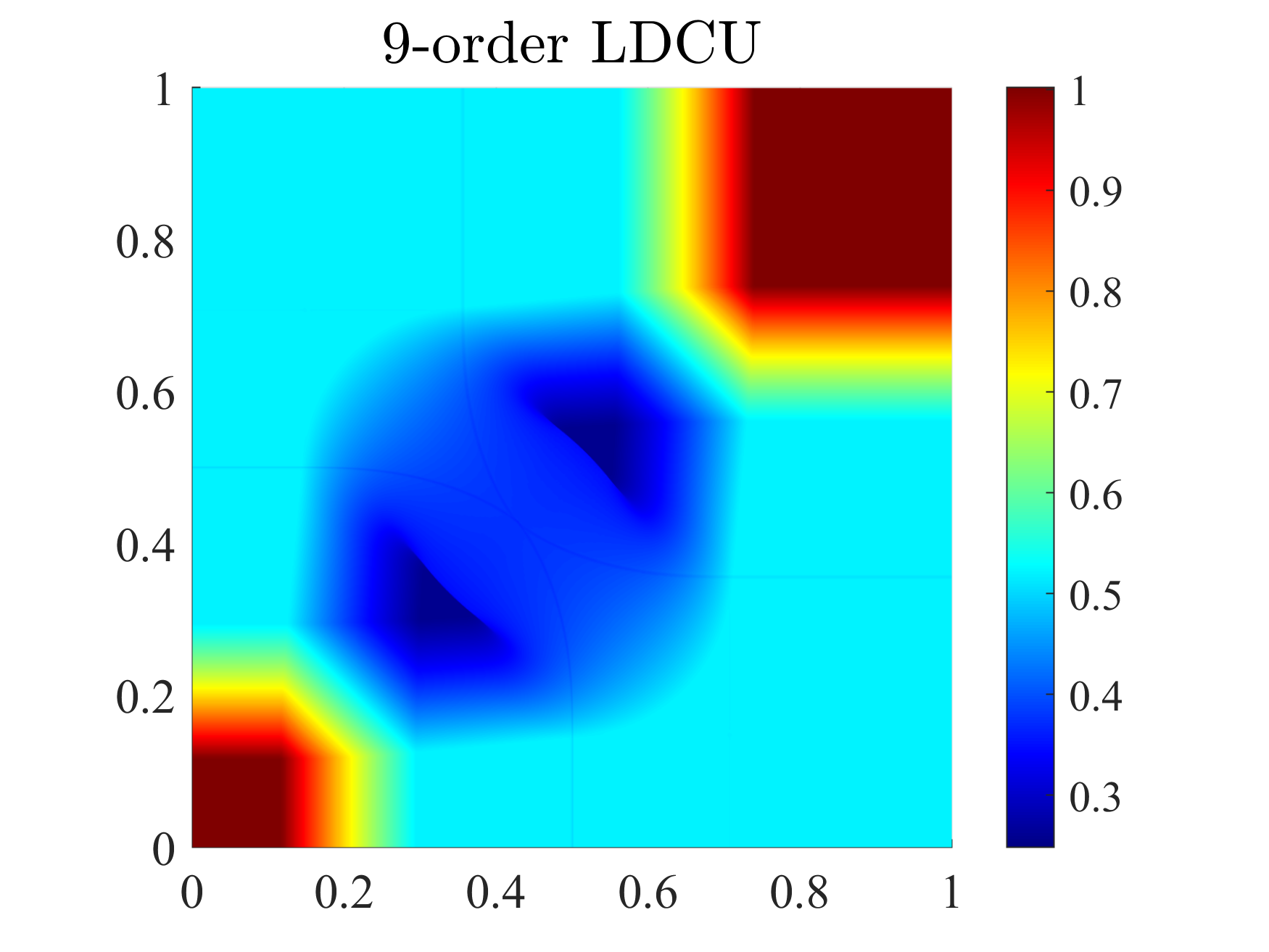}\hspace*{0.0cm}
	        \includegraphics[trim=1.3cm 0.5cm 1.6cm 0.2cm, clip, width=4.2cm]{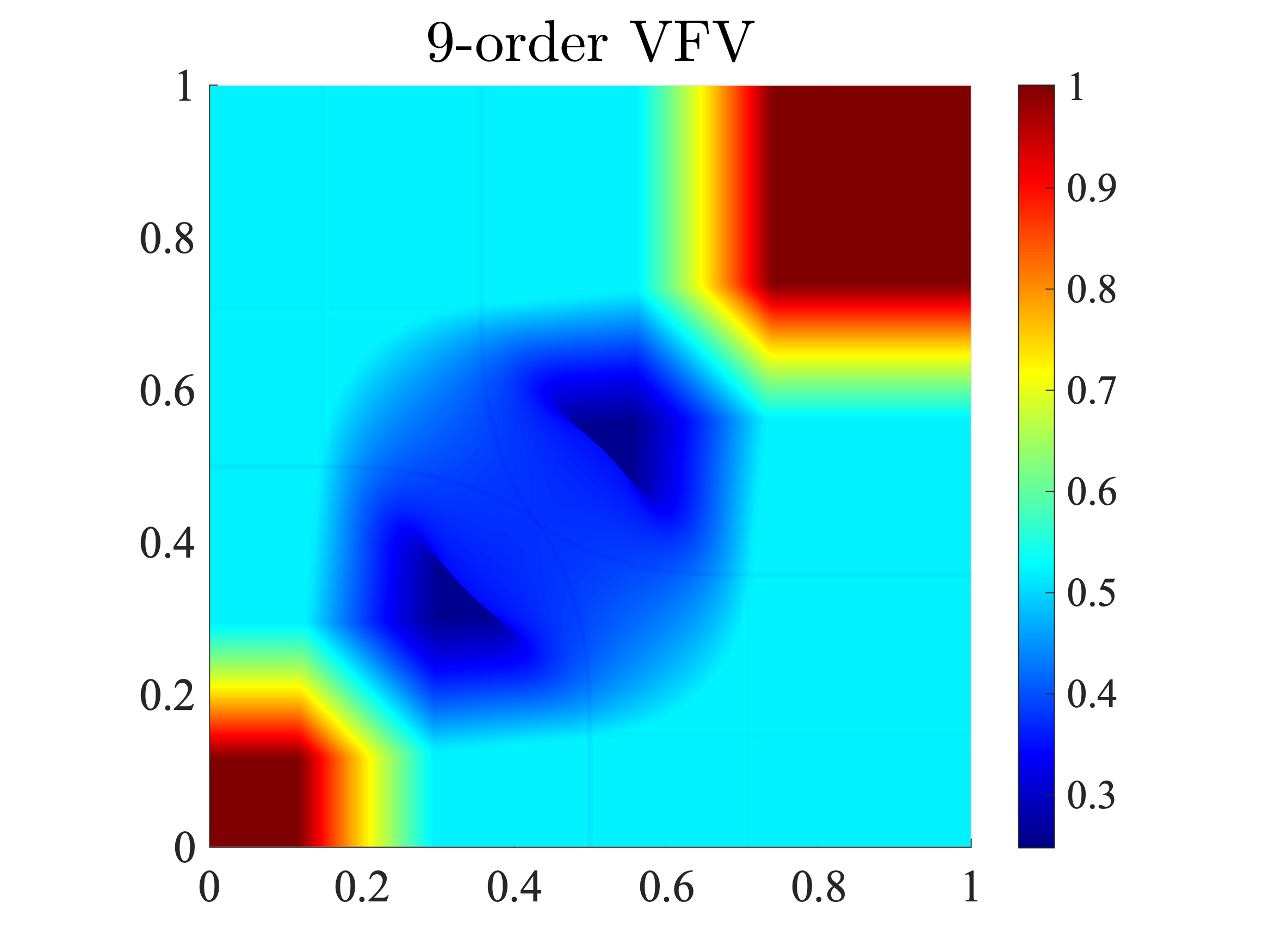}}
\caption{\sf Configuration 2: Density computed by the first- (top row), third- (middle row), and ninth-order (bottom row) LCDCU (left
column), LDCU (middle column), and VFV (right column) schemes.\label{fig7}}
\end{figure}

We compute the numerical solutions until the final time $T=0.2$ using the studied schemes and show the results, obtained by the first-,
third-, and ninth-order schemes in  Figure \ref{fig7}. As one can see, the use of all three numerical fluxes leads to similar solutions,
especially when the ninth-order results are compared. The absence of oscillations in the obtained numerical solutions indicates that they
converge to the same weak solution.

\medskip
\textbf{Configuration 4.} In this example, the initial data are
\begin{equation*}
(\rho,u,v,p)(x,y,0)=\begin{cases}
(1.1,0,0,1.1),&x>0.5,~y>0.5,\\
(0.5065,0.8939,0,0.35),&x<0.5,~y>0.5,\\
(1.1,0.8939,0.8939,1.1),&x<0.5,~y<0.5,\\
(0.5065,0,0.8939,0.35),&x>0.5,~y<0.5,
\end{cases}
\end{equation*}
and far from the center of the computational domain, the solution consists of four shock waves.

\begin{figure}[ht!]
\centerline{\includegraphics[trim=1.3cm 0.5cm 1.6cm 0.2cm, clip, width=4.cm]{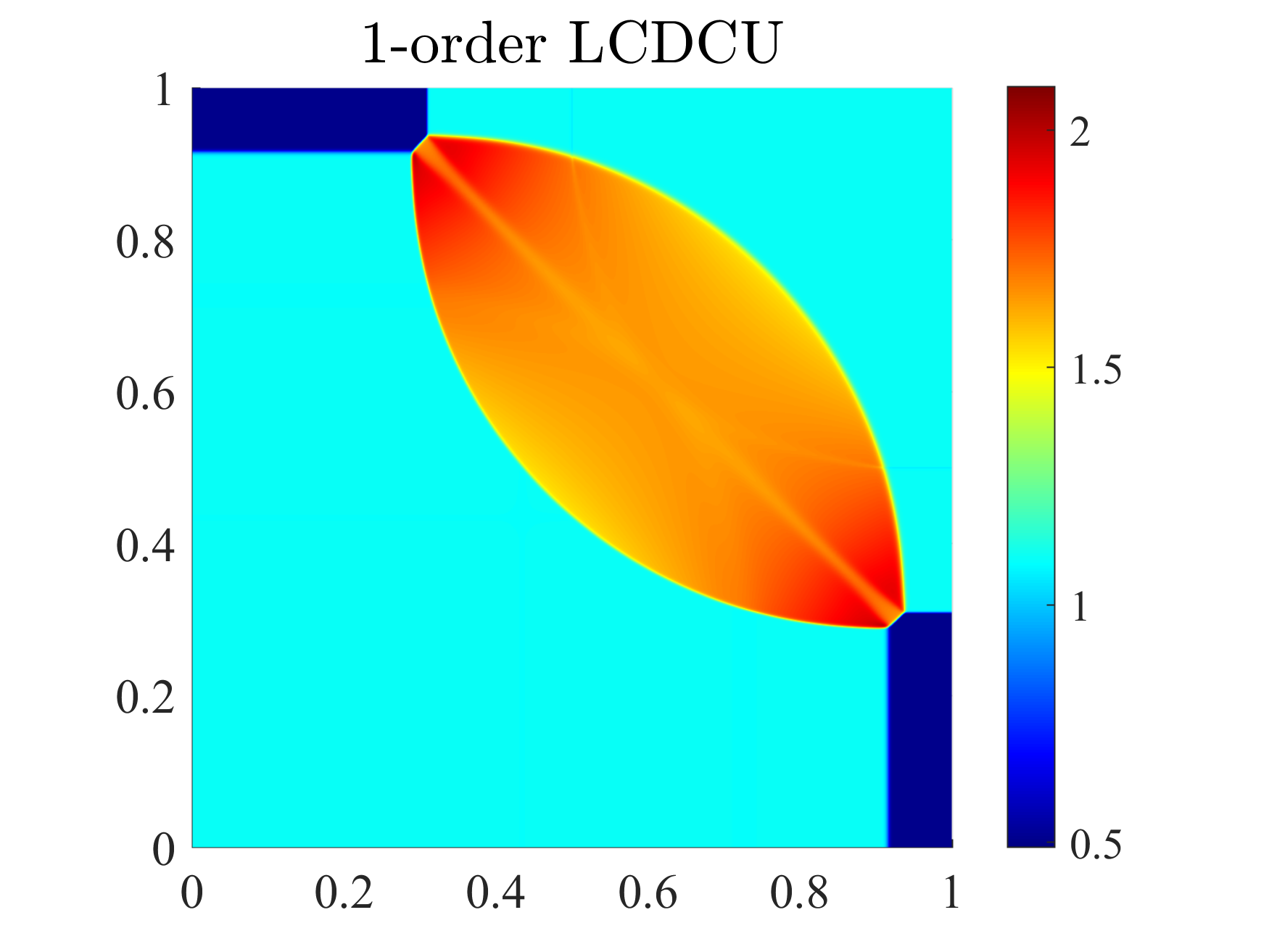}\hspace*{0.2cm}
            \includegraphics[trim=1.3cm 0.5cm 1.6cm 0.2cm, clip, width=4.cm]{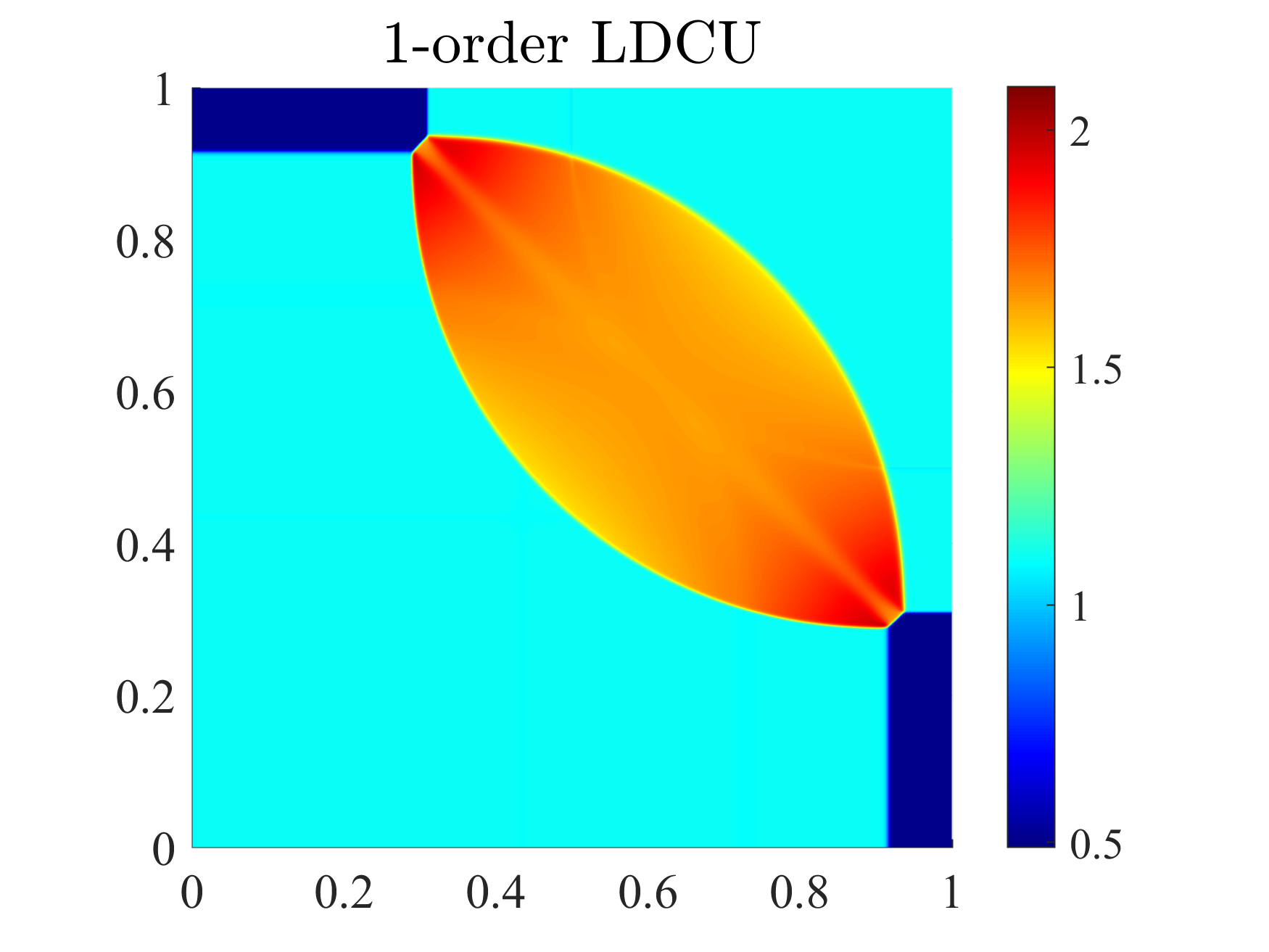}\hspace*{0.0cm}
            \includegraphics[trim=1.3cm 0.5cm 1.6cm 0.2cm, clip, width=4.2cm]{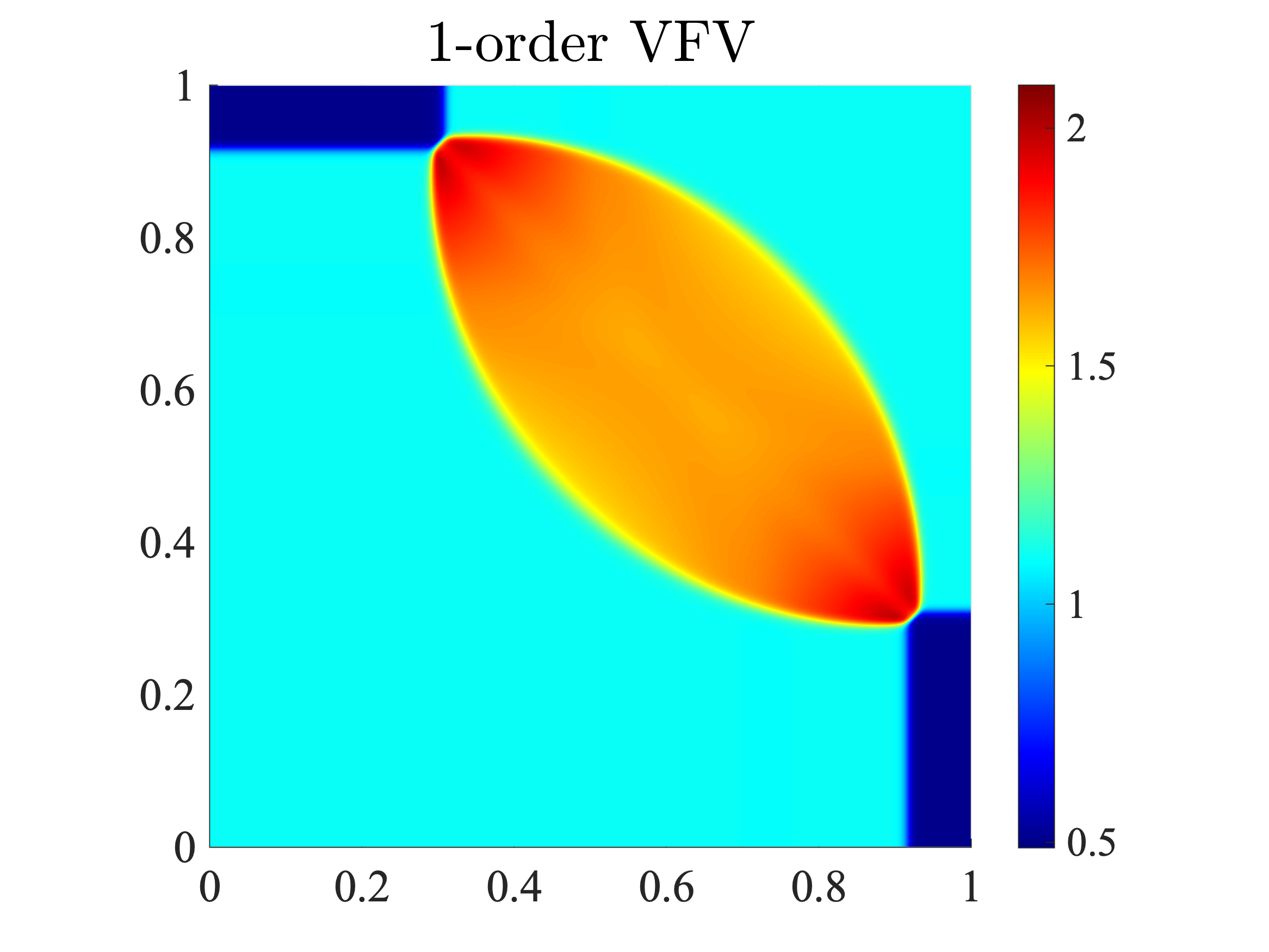}}
\vskip8pt
\centerline{\includegraphics[trim=1.3cm 0.5cm 1.6cm 0.2cm, clip, width=4.cm]{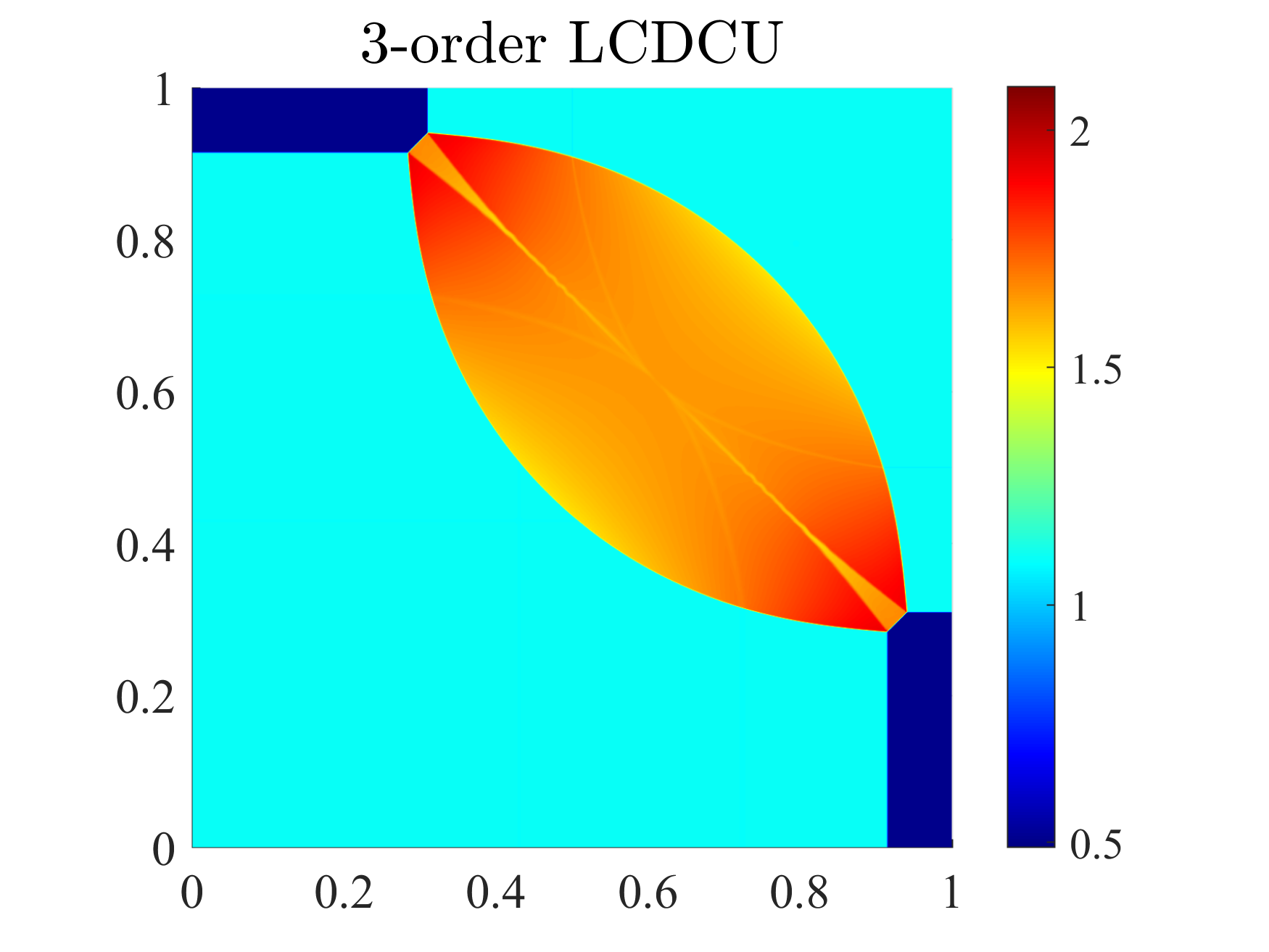}\hspace*{0.2cm}
            \includegraphics[trim=1.3cm 0.5cm 1.6cm 0.2cm, clip, width=4.cm]{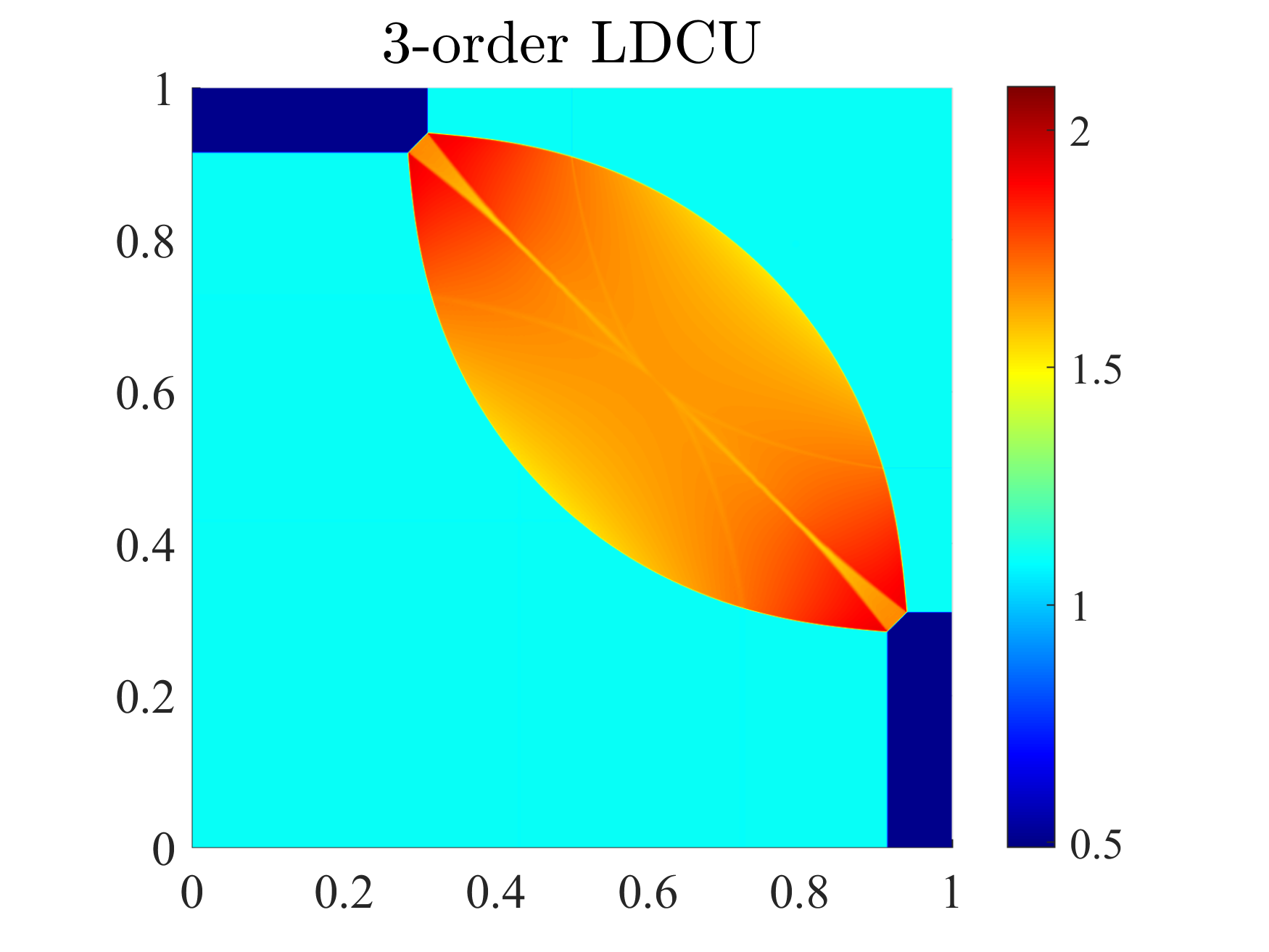}\hspace*{0.0cm}
	        \includegraphics[trim=1.3cm 0.5cm 1.6cm 0.2cm, clip, width=4.2cm]{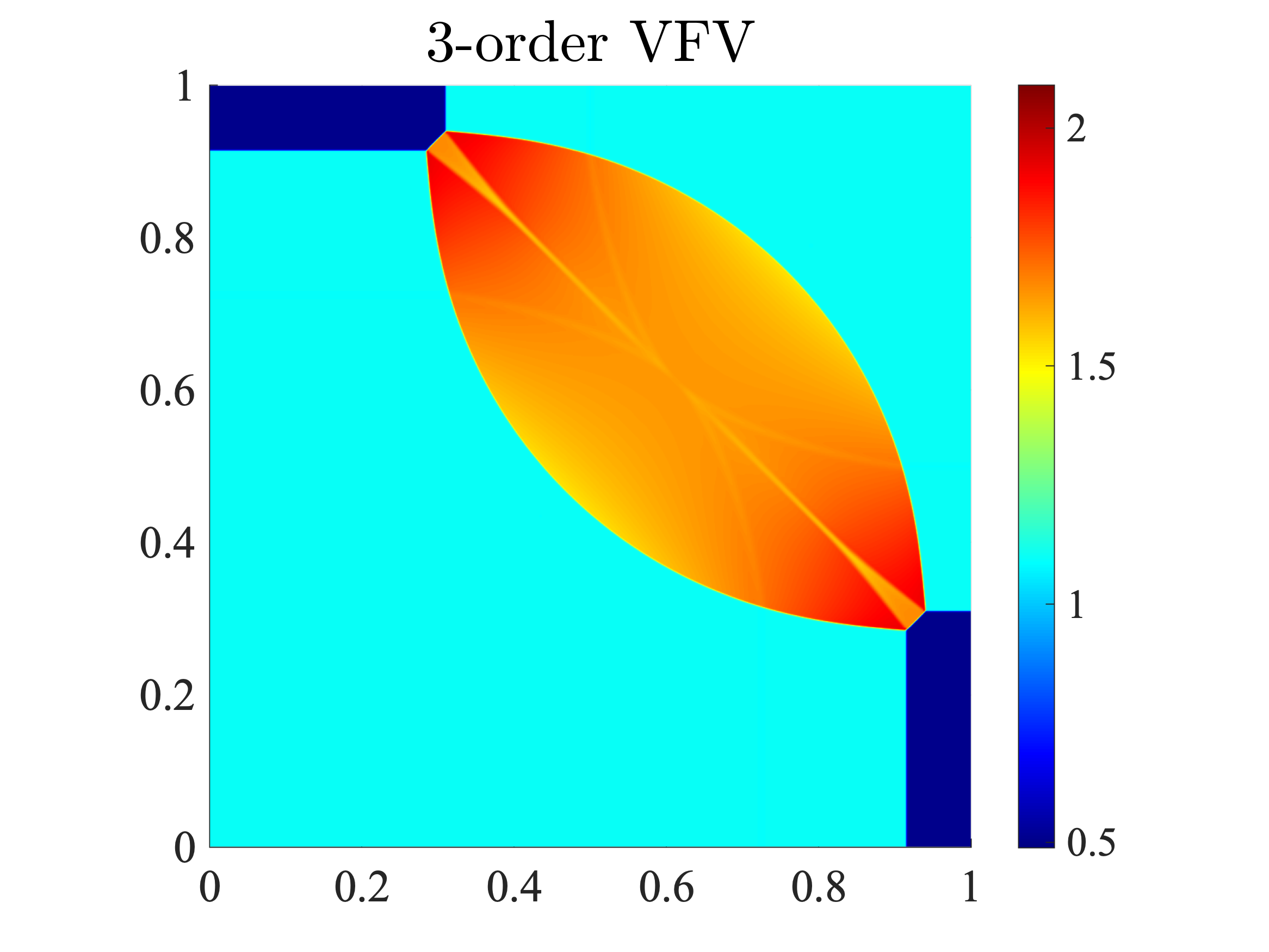}}
\vskip8pt
\centerline{\includegraphics[trim=1.3cm 0.5cm 1.6cm 0.2cm, clip, width=4.cm]{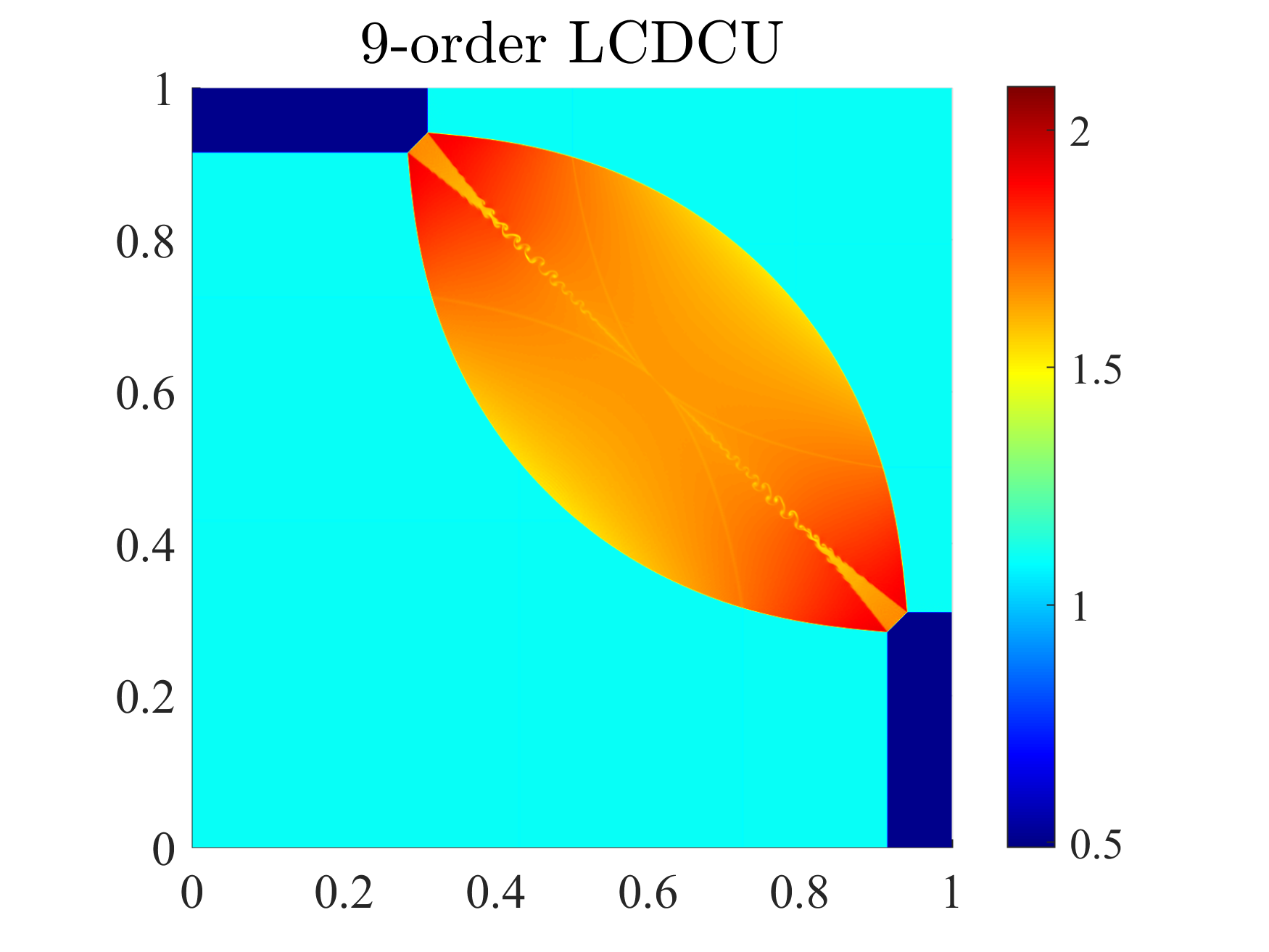}\hspace*{0.2cm}
            \includegraphics[trim=1.3cm 0.5cm 1.6cm 0.2cm, clip, width=4.cm]{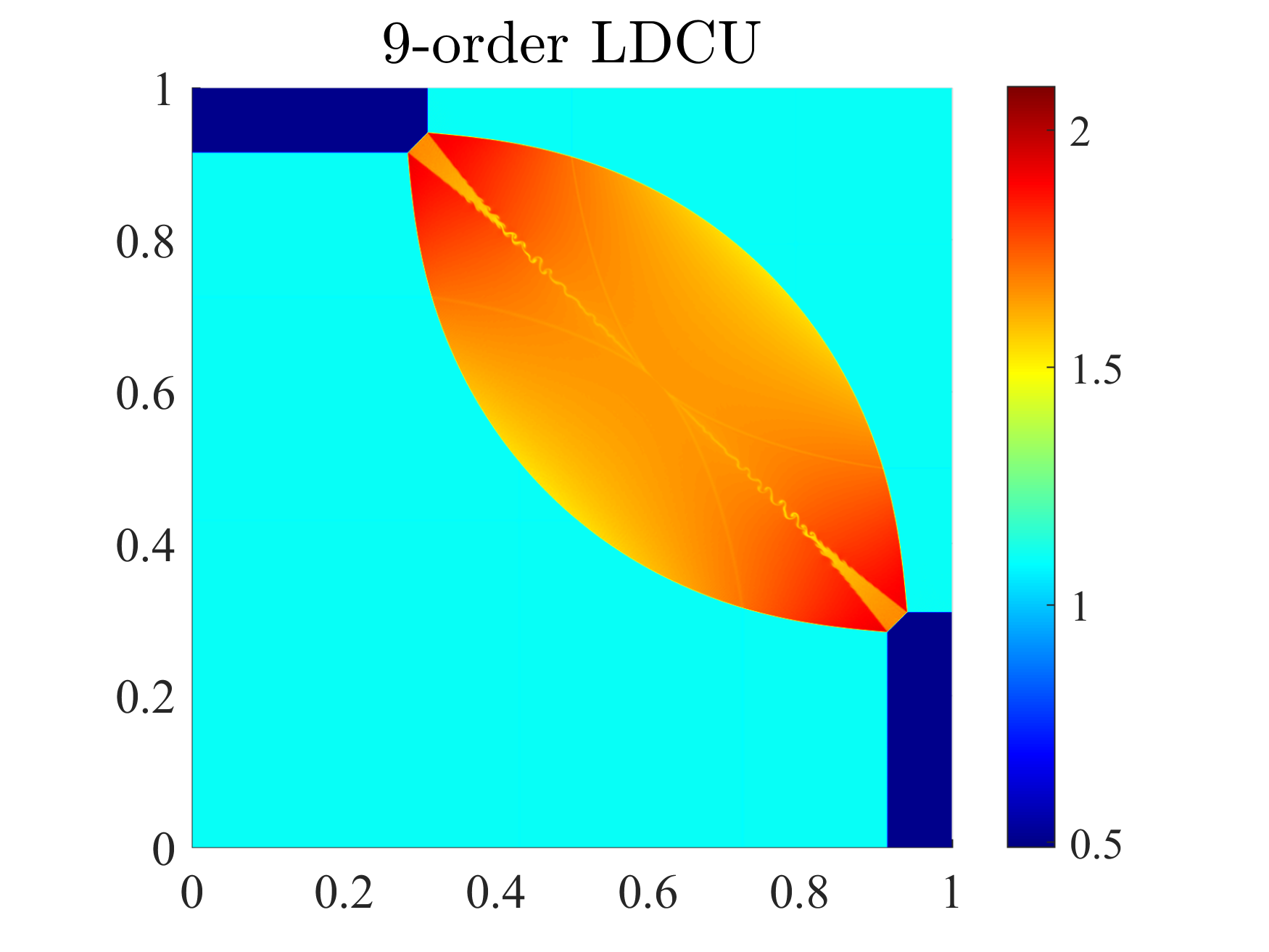}\hspace*{0.0cm}
            \includegraphics[trim=1.3cm 0.5cm 1.6cm 0.2cm, clip, width=4.2cm]{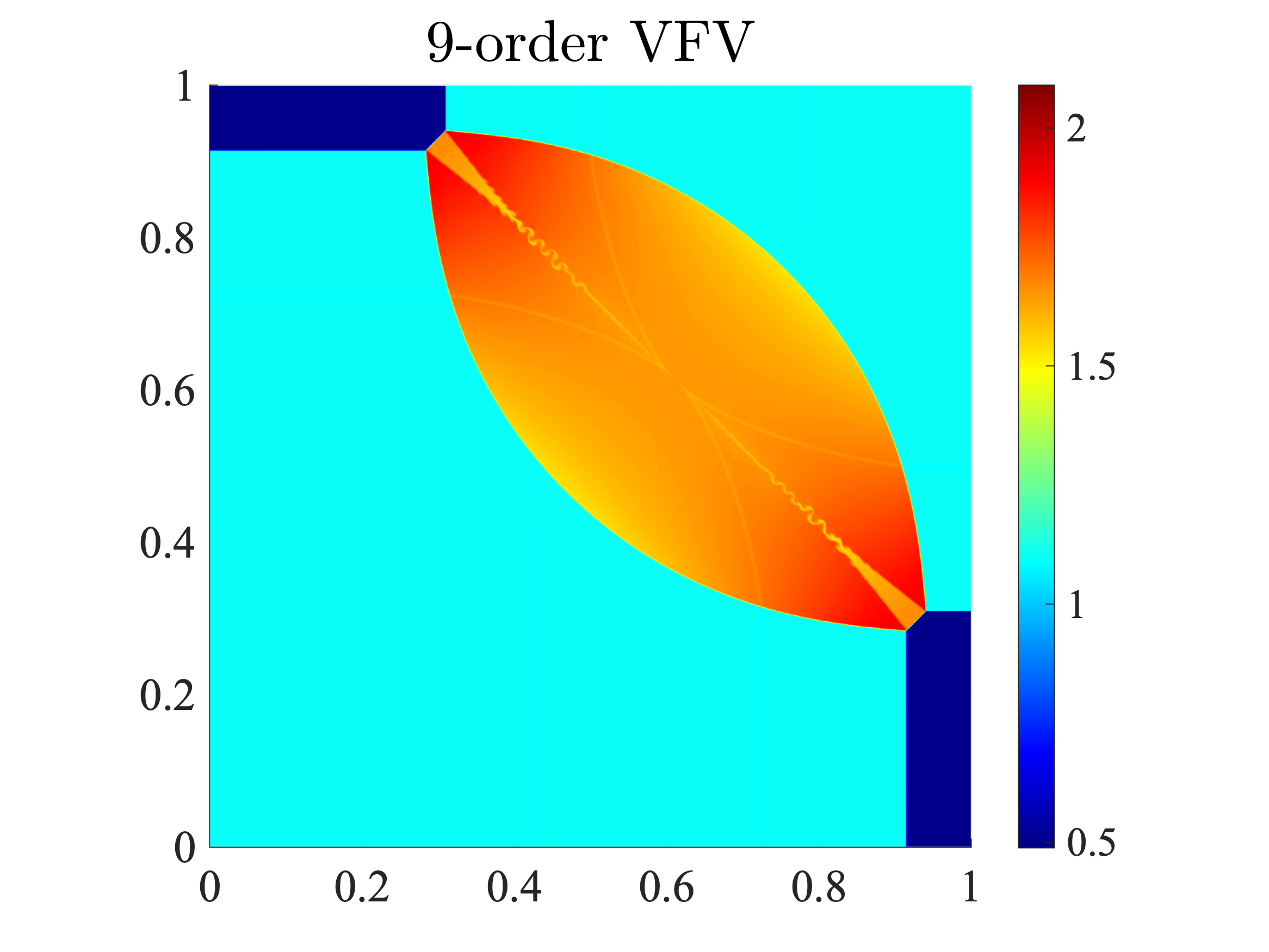}}
\caption{\sf Configuration 4: Density computed by the first- (top row), third- (middle row), and ninth-order (bottom row) LCDCU (left
column), LDCU (middle column), and VFV (right column) schemes.\label{fig11}}
\end{figure}

We compute the numerical solutions until the final time $T=0.25$ using the studied schemes and show the results, obtained by the first-,
third-, and ninth-order schemes in Figure \ref{fig11}. Although the solutions are not smooth, the results exhibit clear similarities,
indicating that the different methods converge to the same weak solution.

\medskip
\textbf{Configuration 3.}
In this example, the initial data are
\begin{equation*}
(\rho,u,v,p)(x,y,0)=\begin{cases}
(1.5,0,0,1.5),&x>0.8,~y>0.8,\\
(0.5323,1.206,0,0.3),&x<0.8,~y>0.8,\\
(0.138,1.206,1.206,0.029),&x<0.8,~y<0.8,\\
(0.5323,0,1.206,0.3),&x>0.8,~y<0.8,
\end{cases}
\end{equation*}
and far from the center of the computational domain, the solution consists of four shock waves.

\begin{figure}[ht!]
\centerline{\includegraphics[trim=1.3cm 0.5cm 1.6cm 0.2cm, clip, width=4.cm]{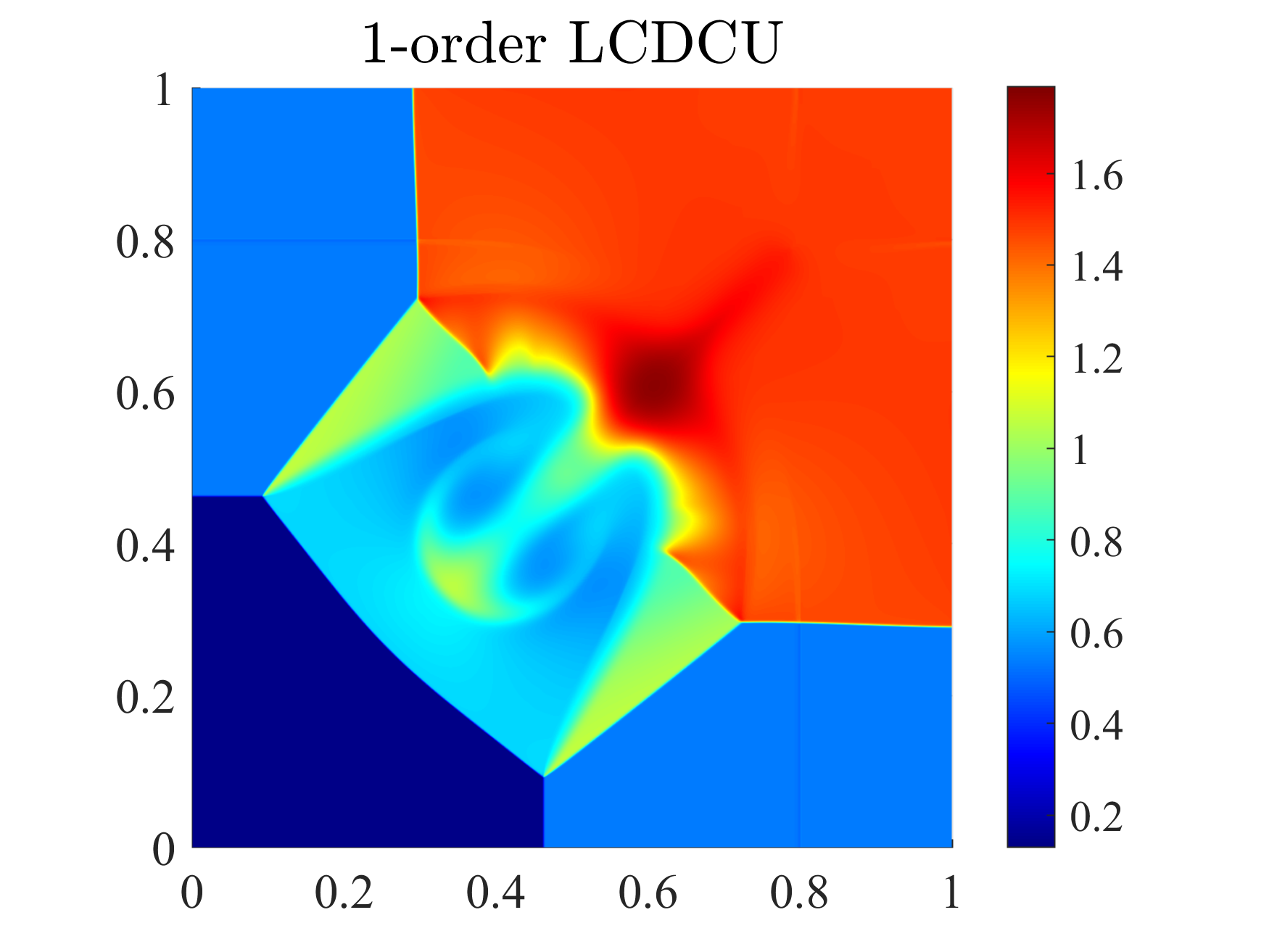}\hspace*{0.2cm}
            \includegraphics[trim=1.3cm 0.5cm 1.6cm 0.2cm, clip, width=4.cm]{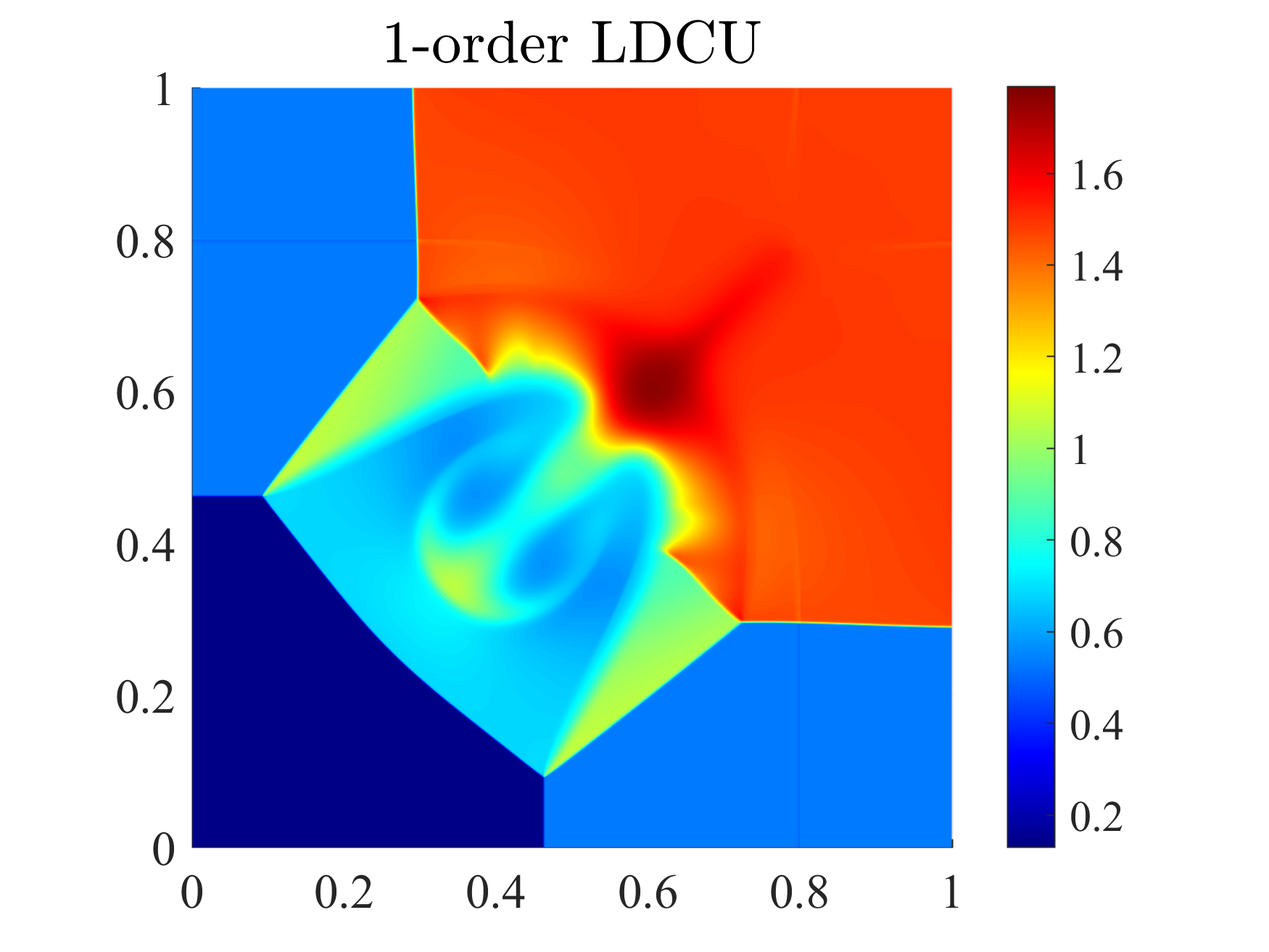}\hspace*{0.0cm}
            \includegraphics[trim=1.3cm 0.5cm 1.6cm 0.2cm, clip, width=4.2cm]{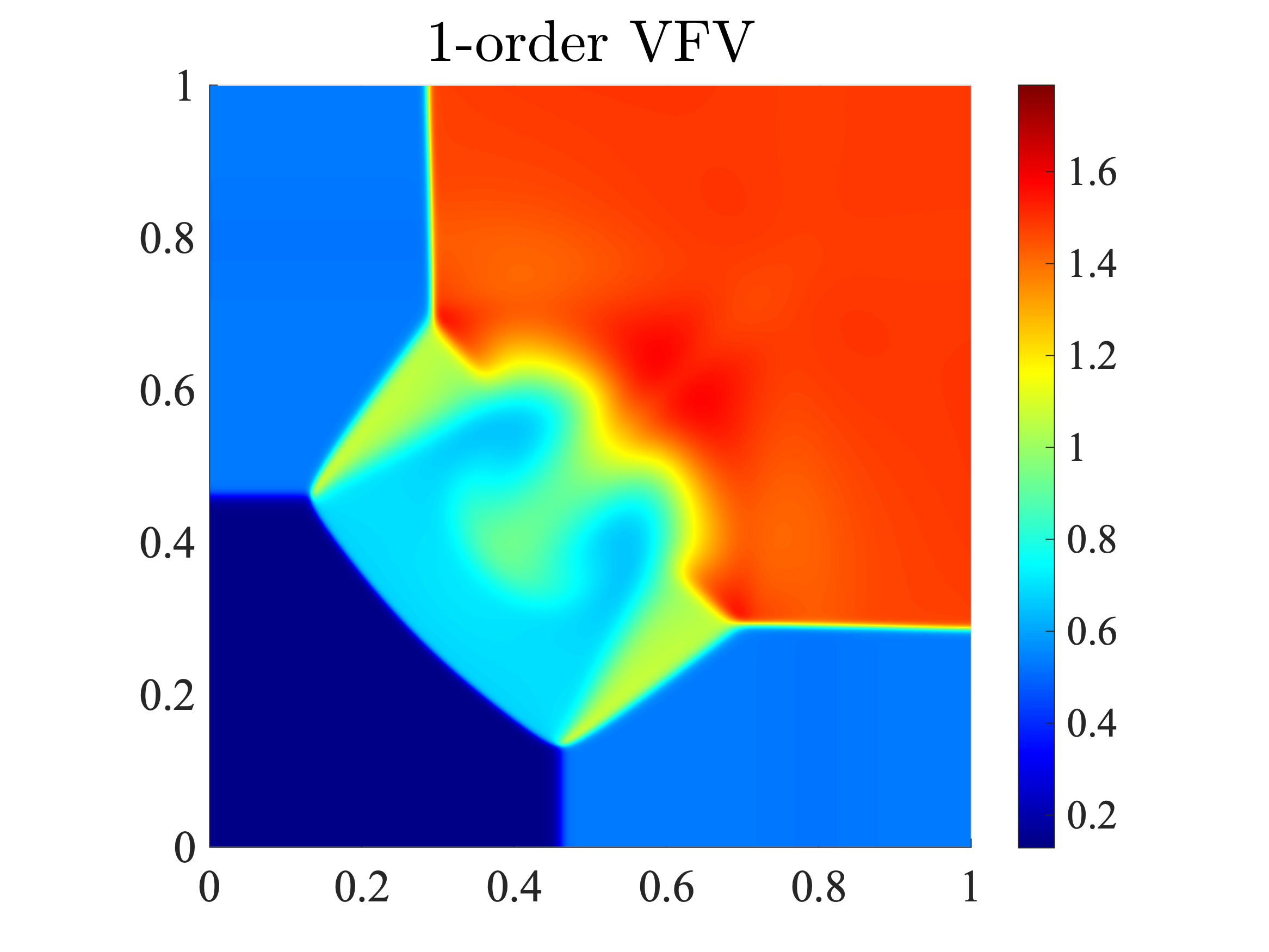}}
\vskip8pt
\centerline{\includegraphics[trim=1.3cm 0.5cm 1.6cm 0.2cm, clip, width=4.cm]{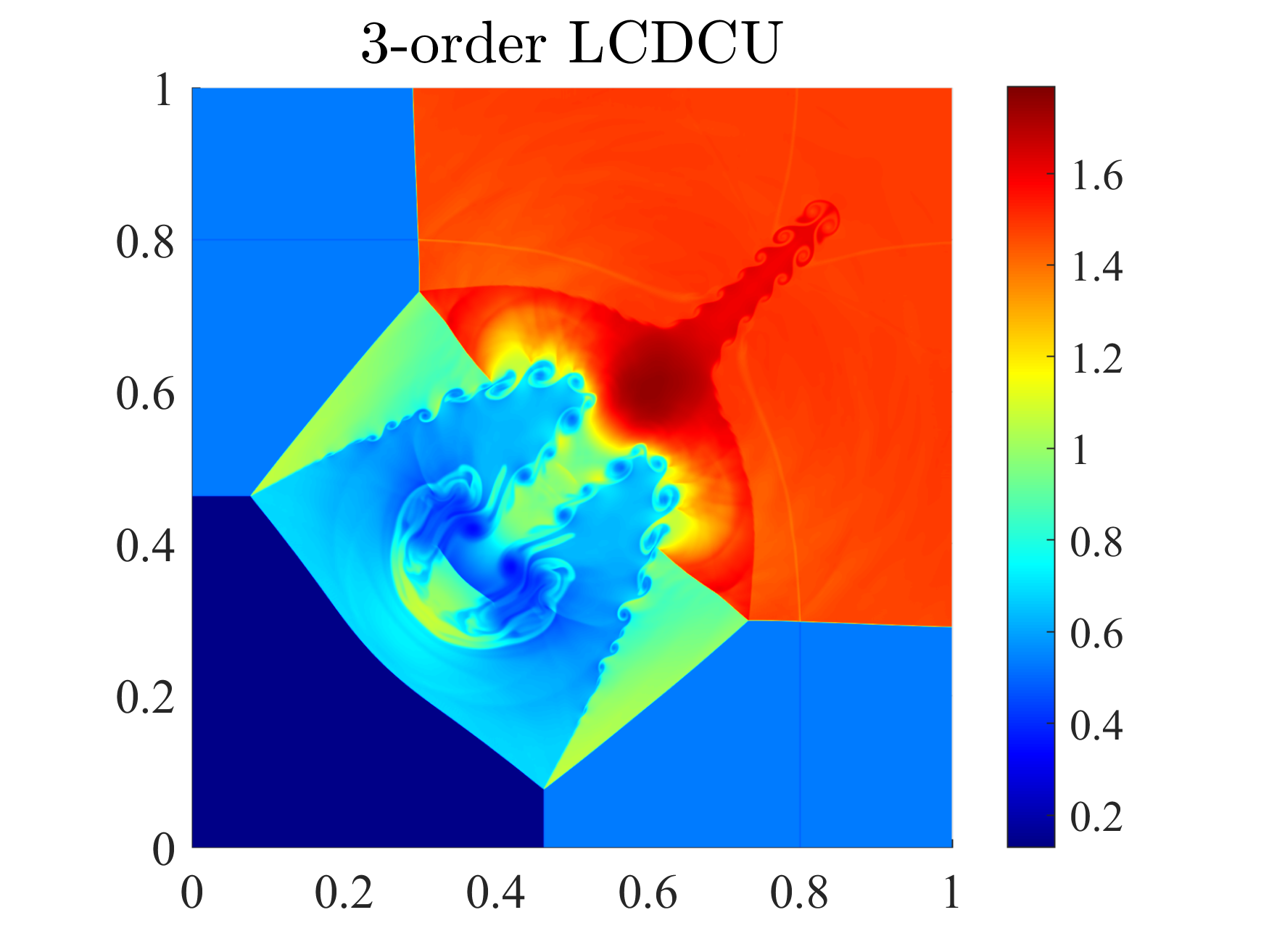}\hspace*{0.2cm}
            \includegraphics[trim=1.3cm 0.5cm 1.6cm 0.2cm, clip, width=4.cm]{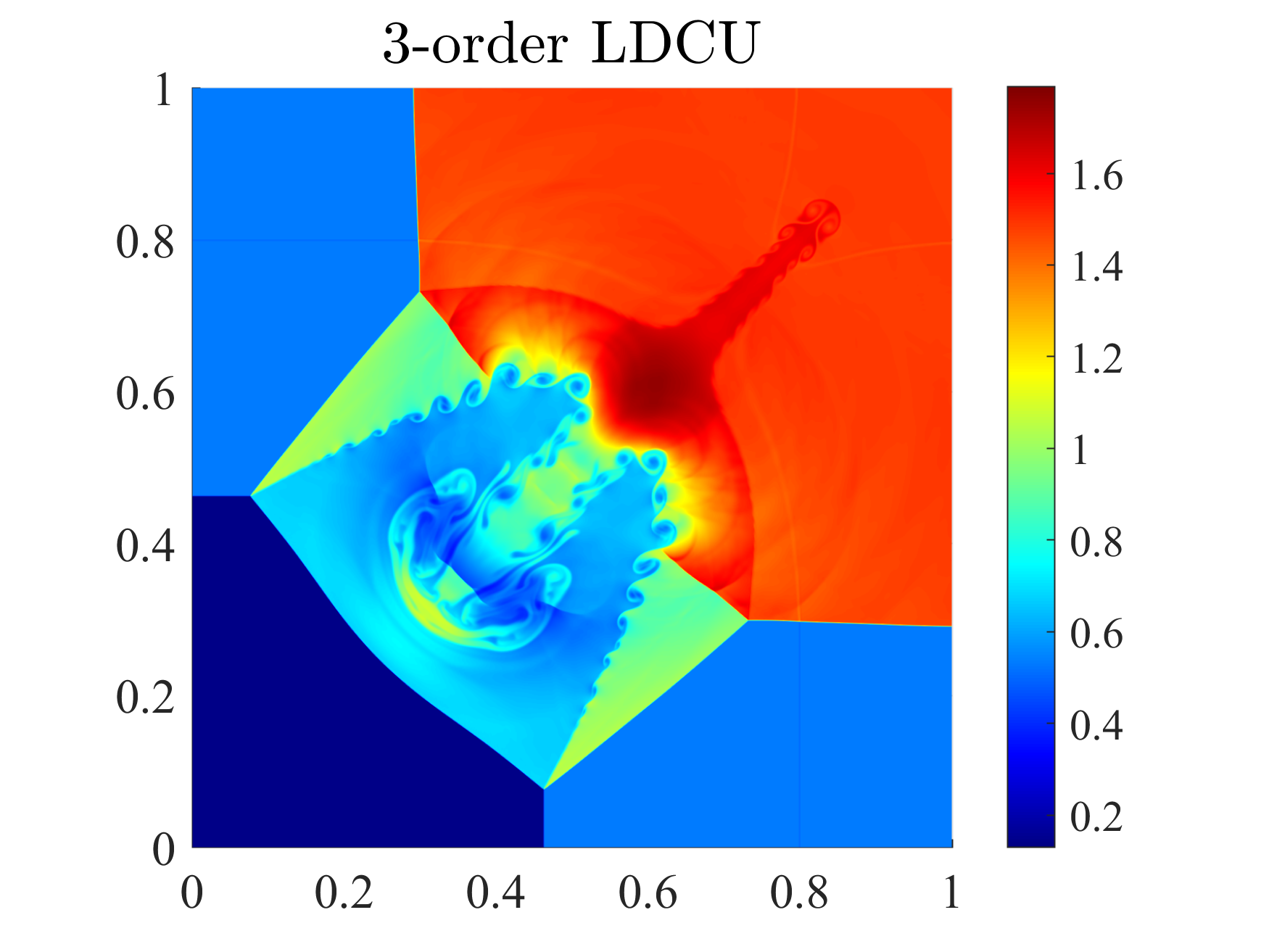}\hspace*{0.0cm}
            \includegraphics[trim=1.3cm 0.5cm 1.6cm 0.2cm, clip, width=4.2cm]{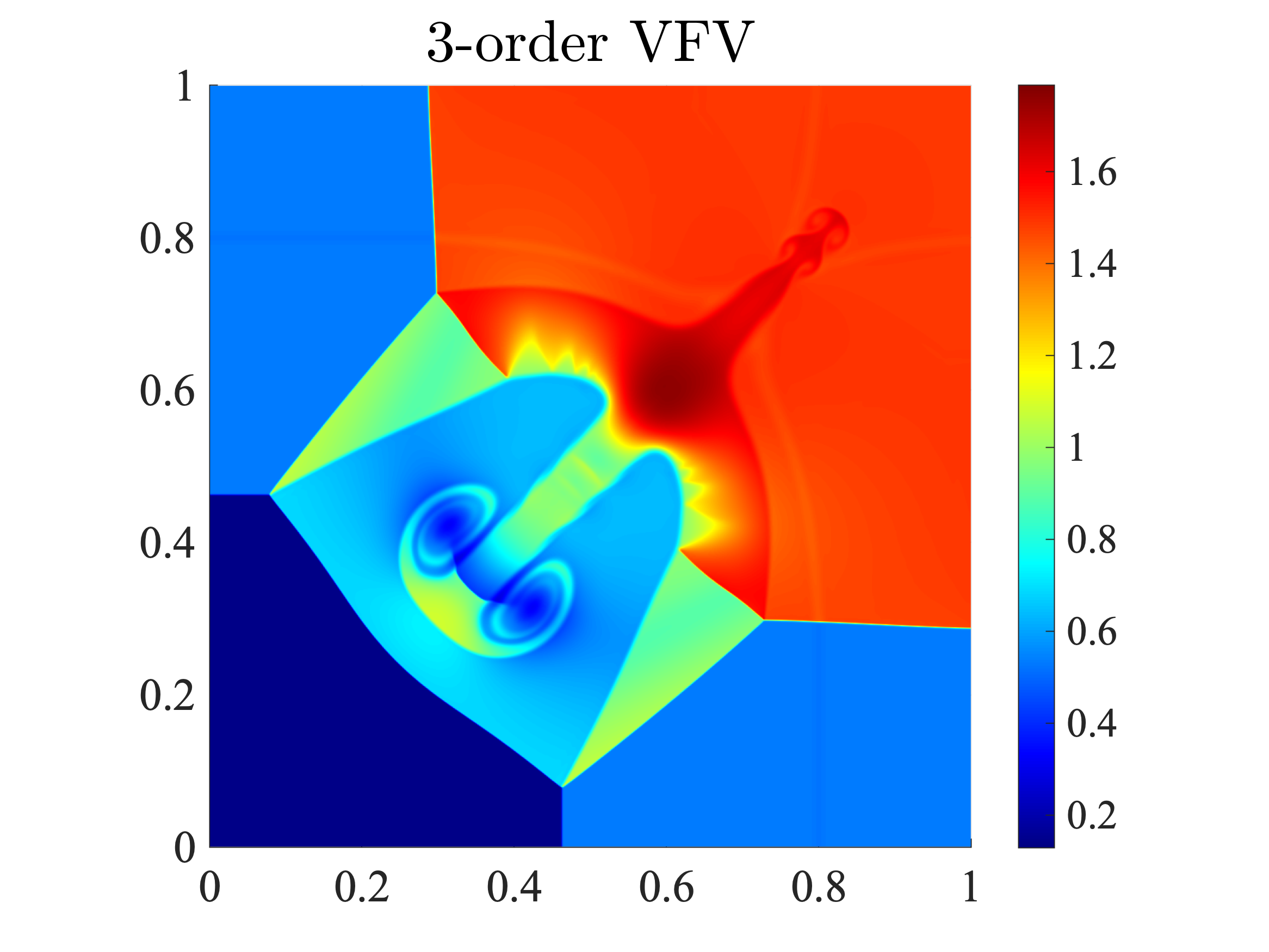}}
\vskip8pt
\centerline{\includegraphics[trim=1.3cm 0.5cm 1.6cm 0.2cm, clip, width=4.cm]{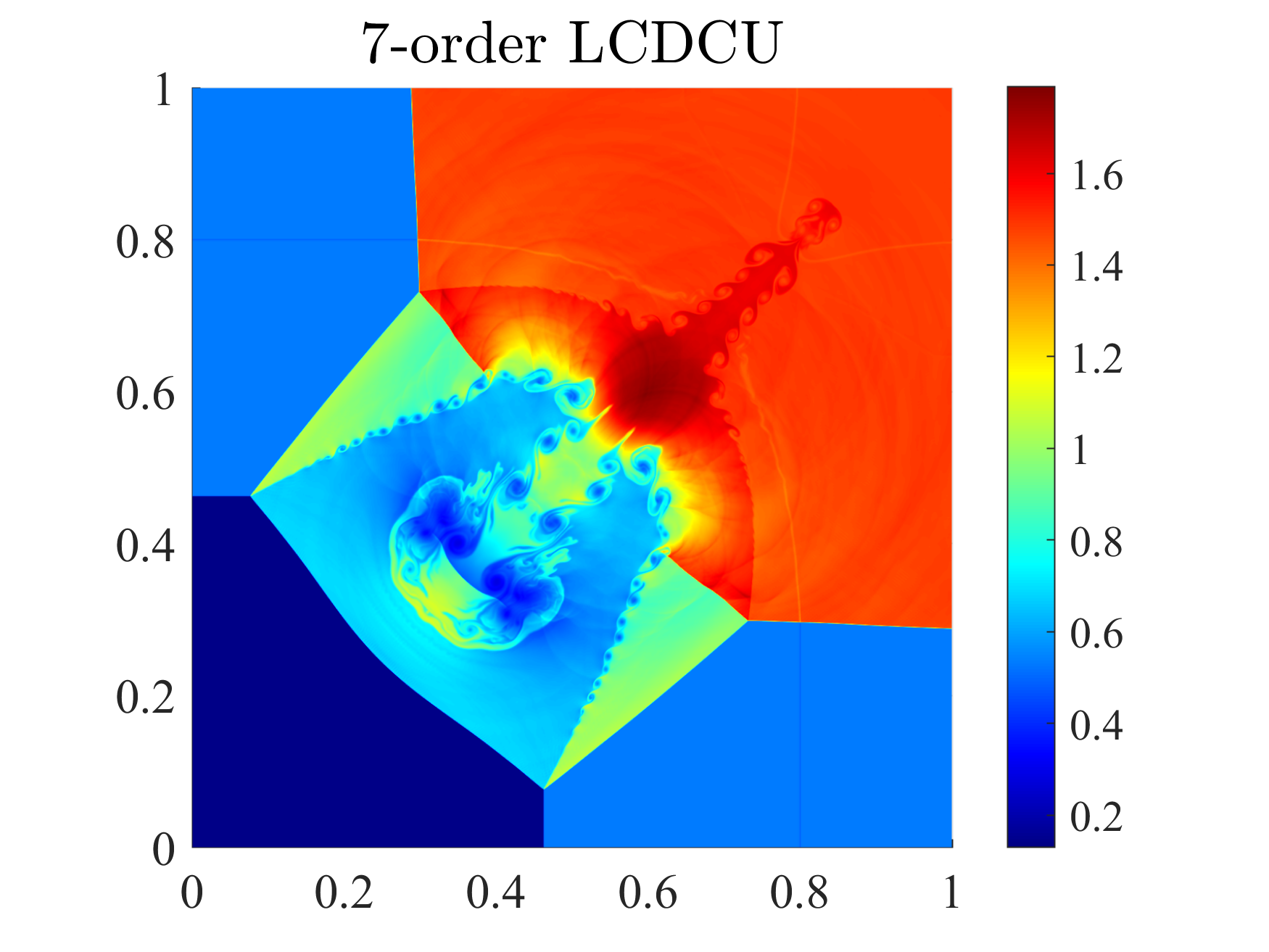}\hspace*{0.2cm}
            \includegraphics[trim=1.3cm 0.5cm 1.6cm 0.2cm, clip, width=4.cm]{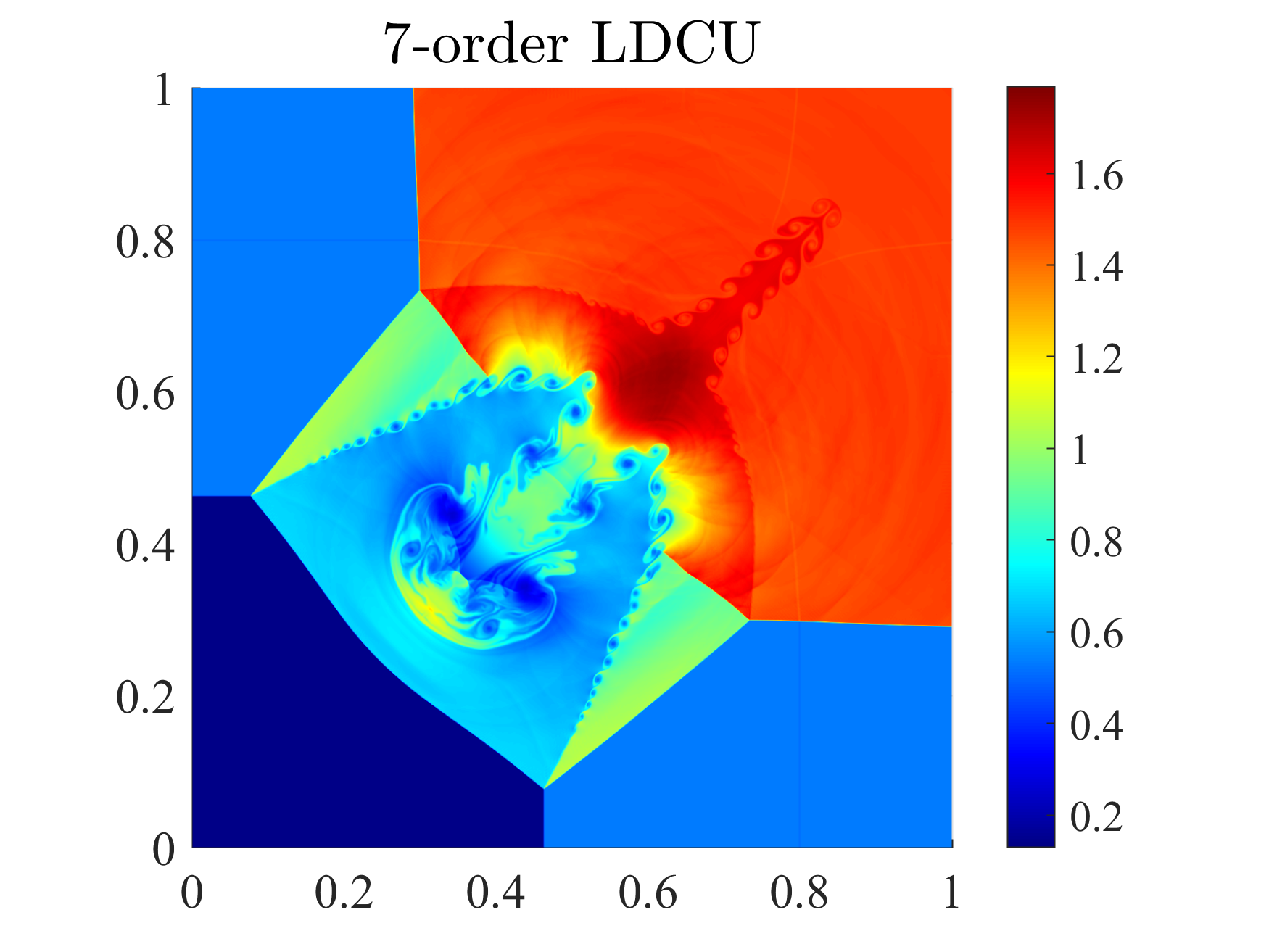}\hspace*{0.0cm}
            \includegraphics[trim=1.3cm 0.5cm 1.6cm 0.2cm, clip, width=4.2cm]{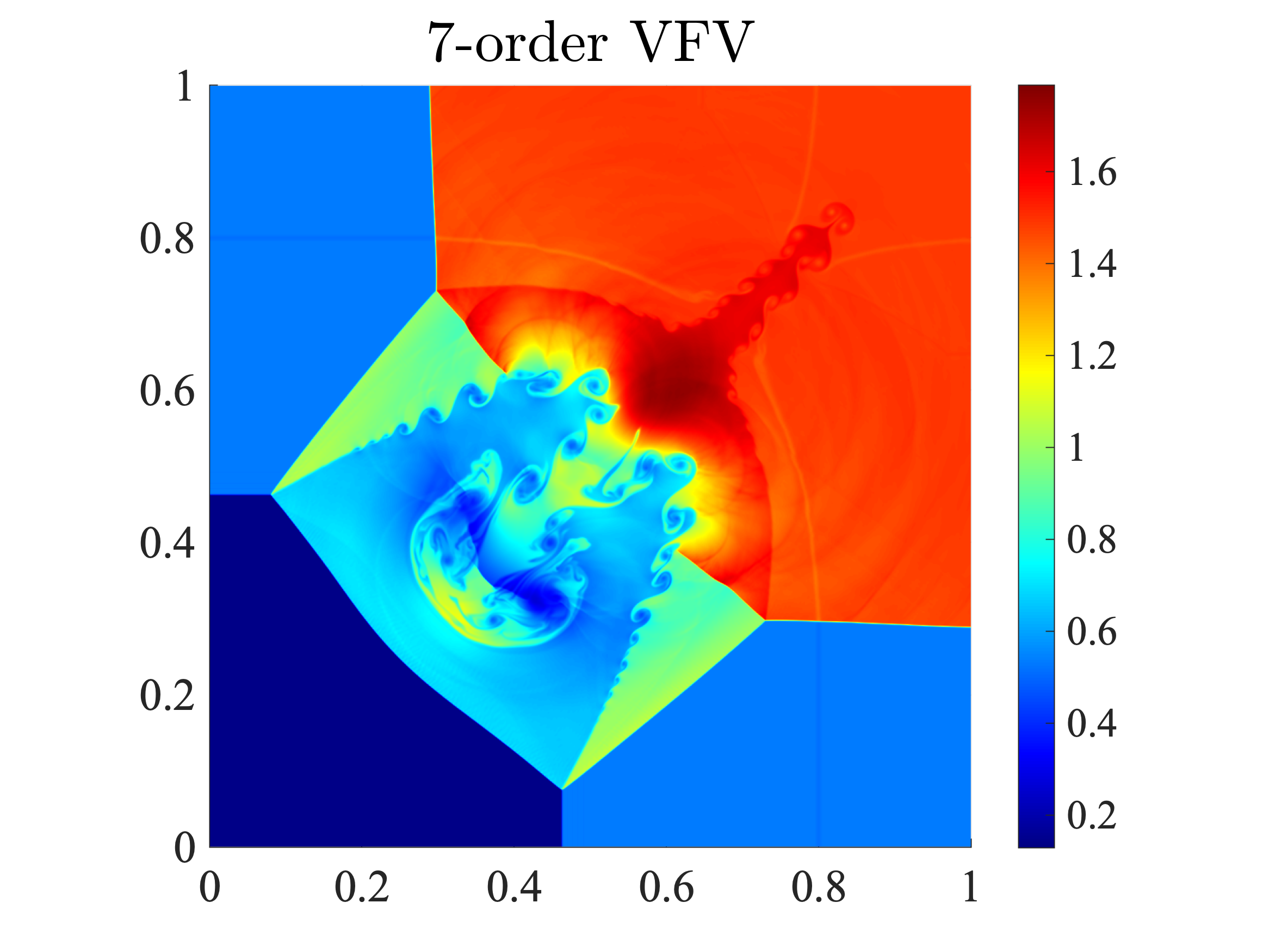}}
\caption{\sf Configuration 3: Density computed by the first- (top row), third- (middle row), and seventh-order (bottom row) LCDCU (left
column), LDCU (middle column), and VFV (right column) schemes.\label{fig22}}
\end{figure}

\begin{figure}[ht!]
\centerline{\includegraphics[trim=1.3cm 0.5cm 1.6cm 0.2cm, clip, width=4.cm]{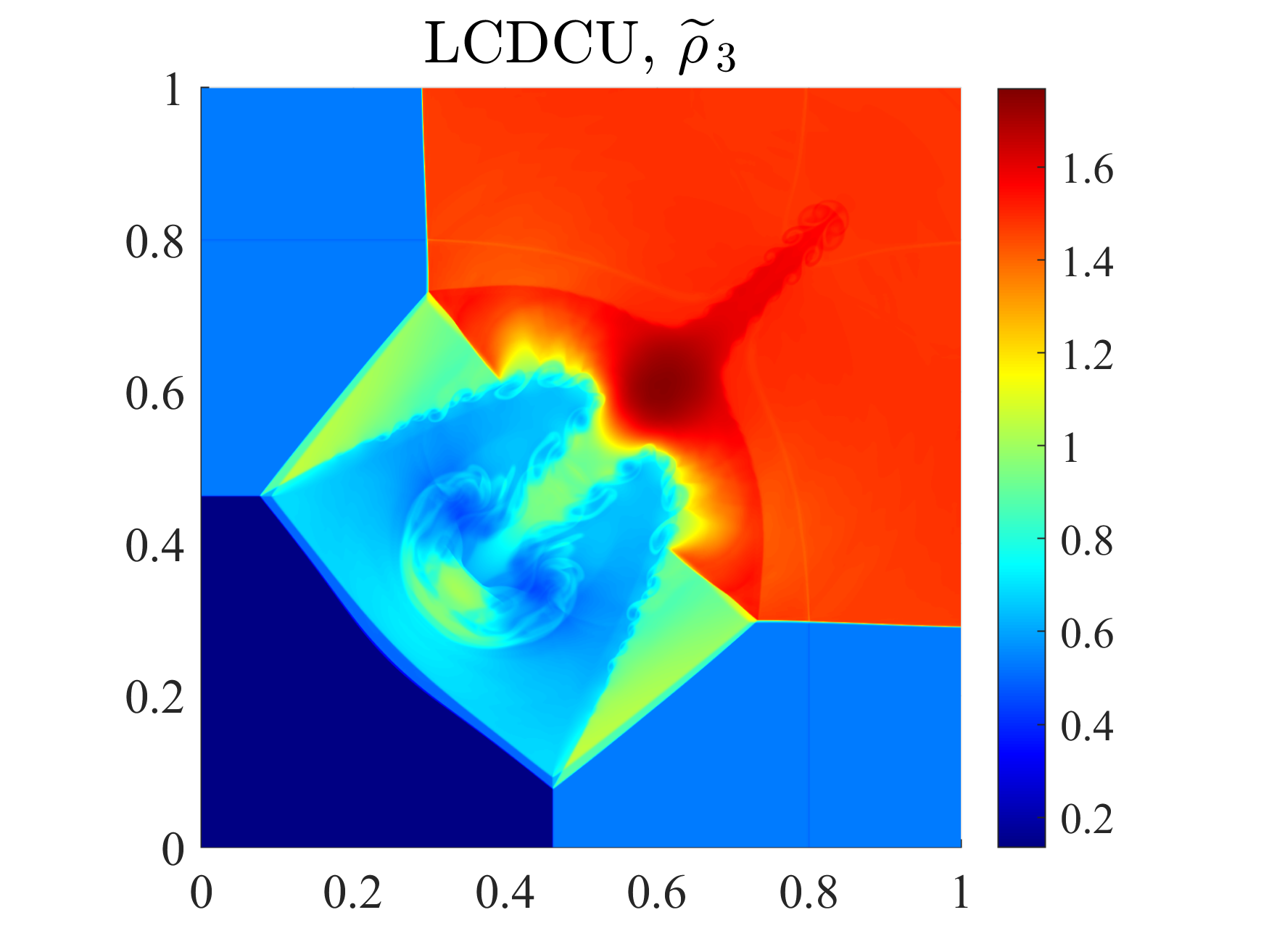}\hspace*{0.2cm}
            \includegraphics[trim=1.3cm 0.5cm 1.6cm 0.2cm, clip, width=4.cm]{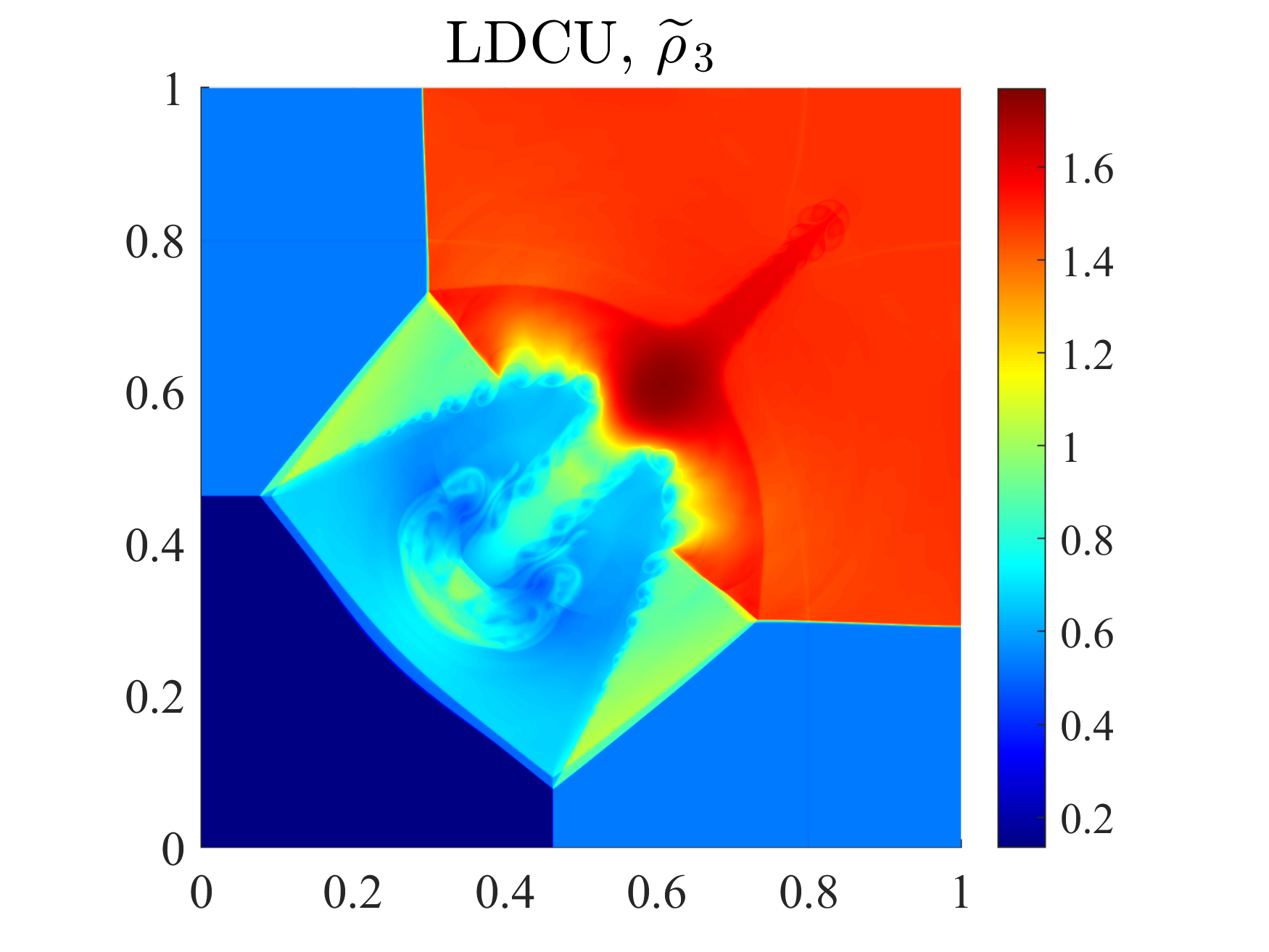}\hspace*{0.0cm}
            \includegraphics[trim=1.3cm 0.5cm 1.6cm 0.2cm, clip, width=4.2cm]{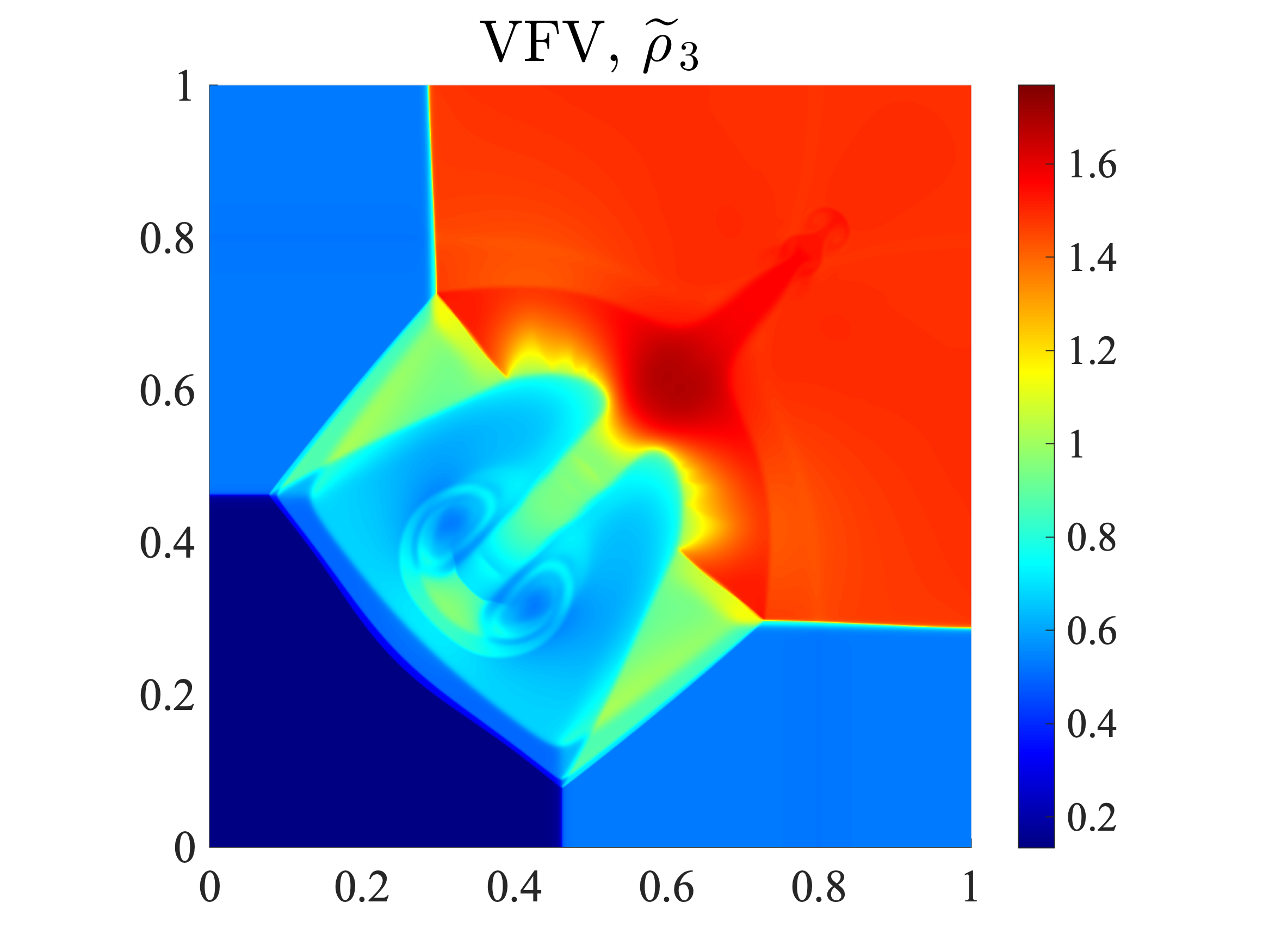}}
\vskip8pt
\centerline{\includegraphics[trim=1.3cm 0.5cm 1.6cm 0.2cm, clip, width=4.cm]{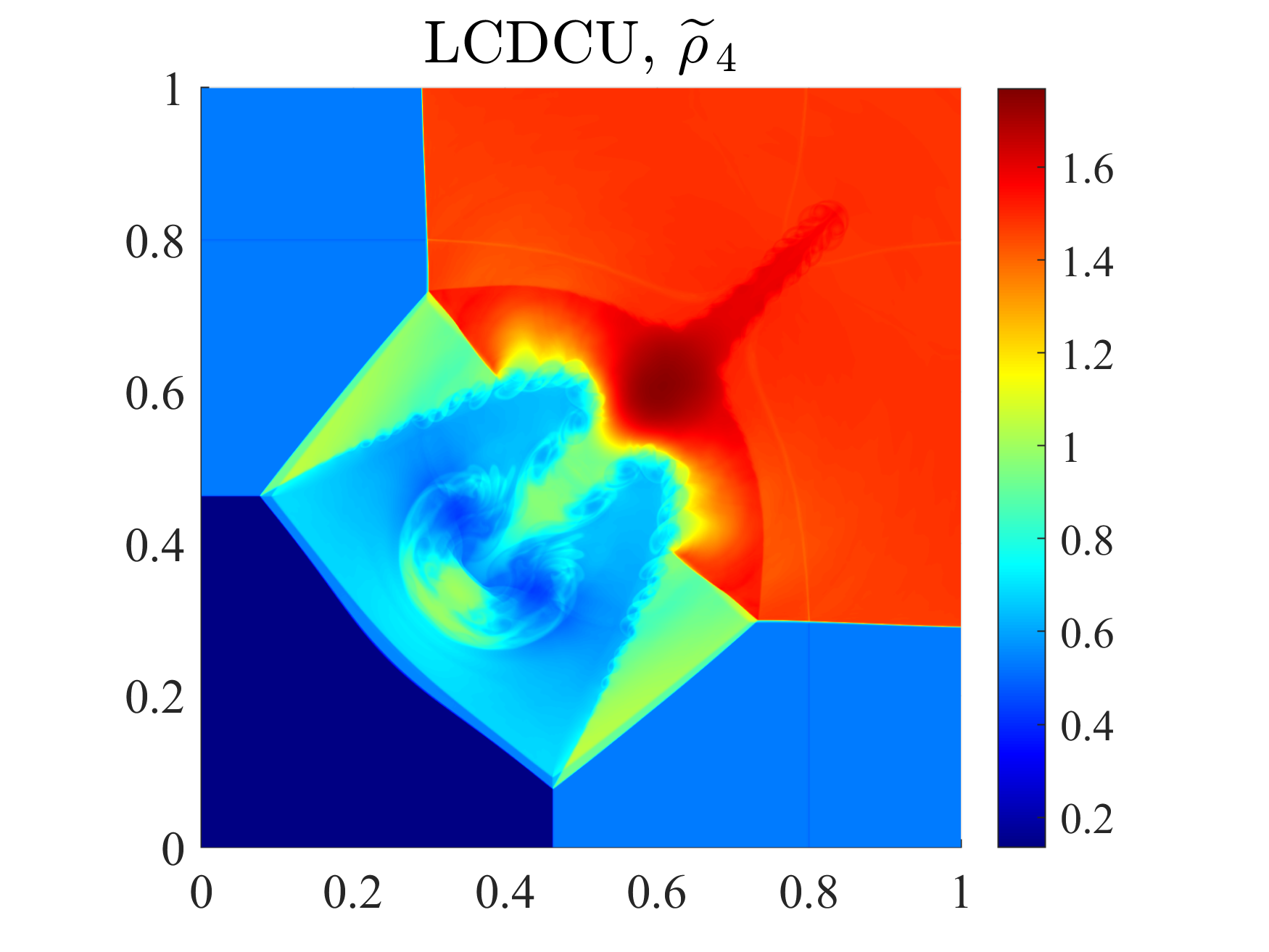}\hspace*{0.2cm}
            \includegraphics[trim=1.3cm 0.5cm 1.6cm 0.2cm, clip, width=4.cm]{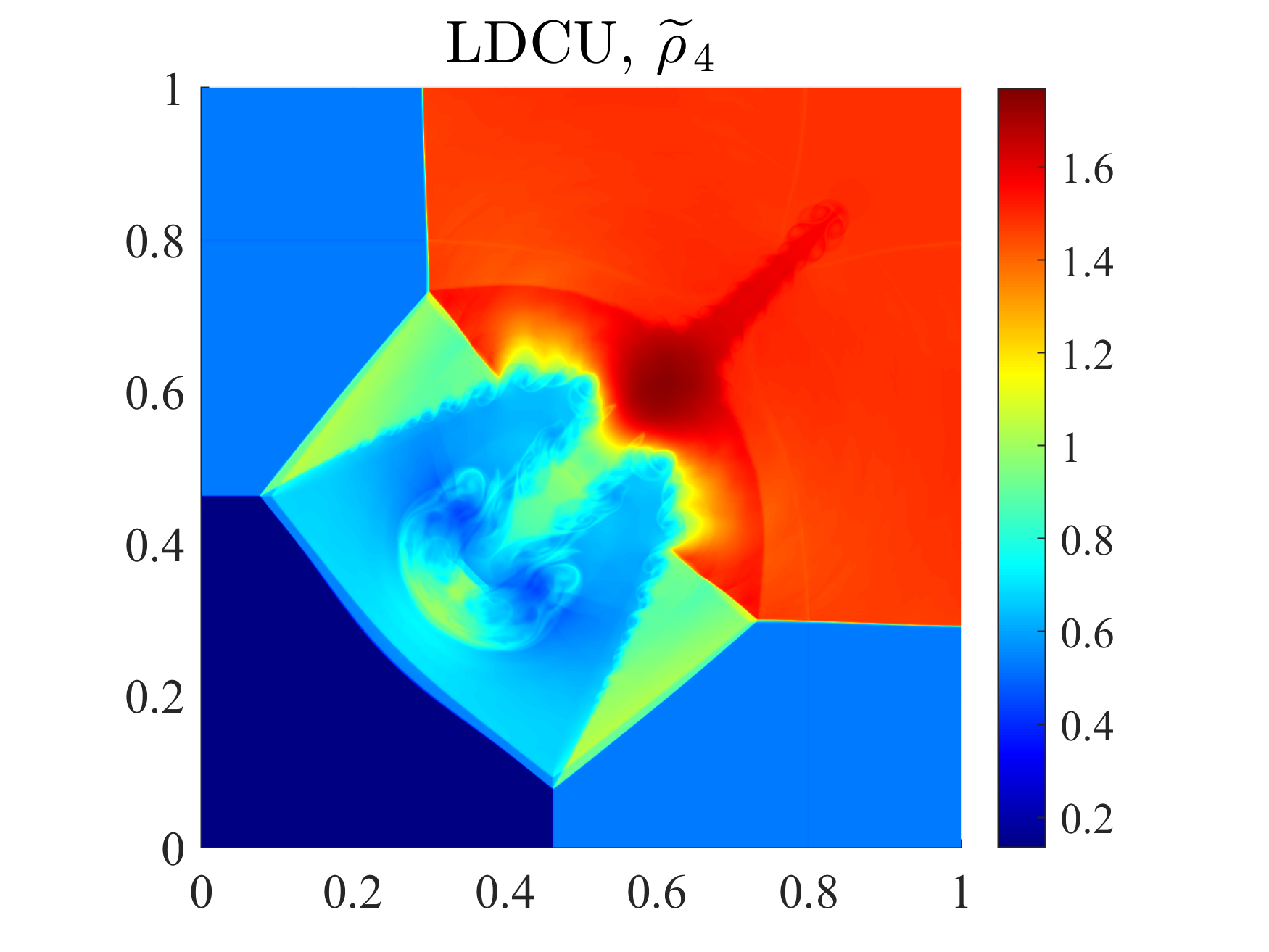}\hspace*{0.0cm}
            \includegraphics[trim=1.3cm 0.5cm 1.6cm 0.2cm, clip, width=4.2cm]{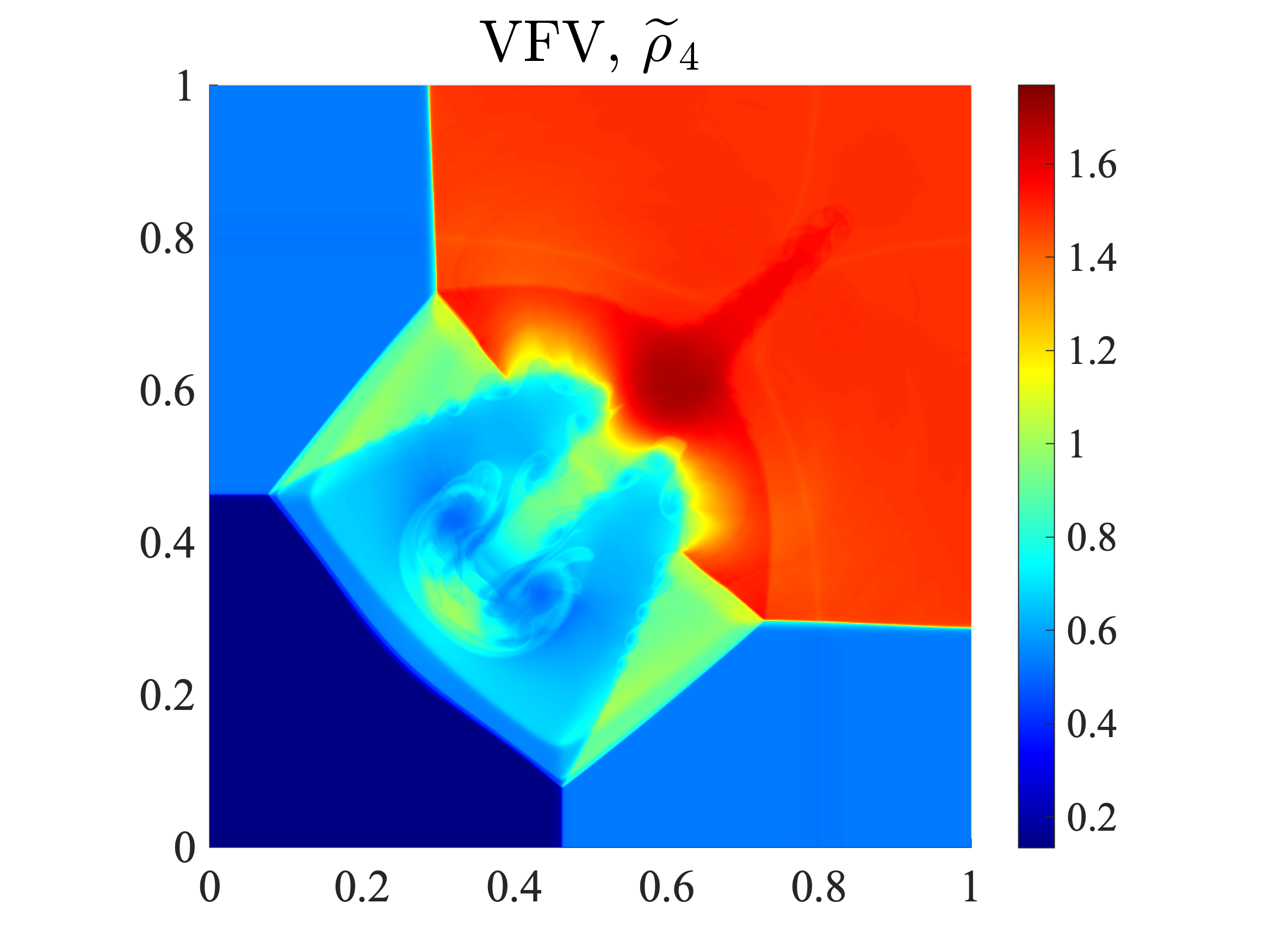}}
\vskip8pt
\centerline{\includegraphics[trim=1.3cm 0.5cm 1.6cm 0.2cm, clip, width=4.cm]{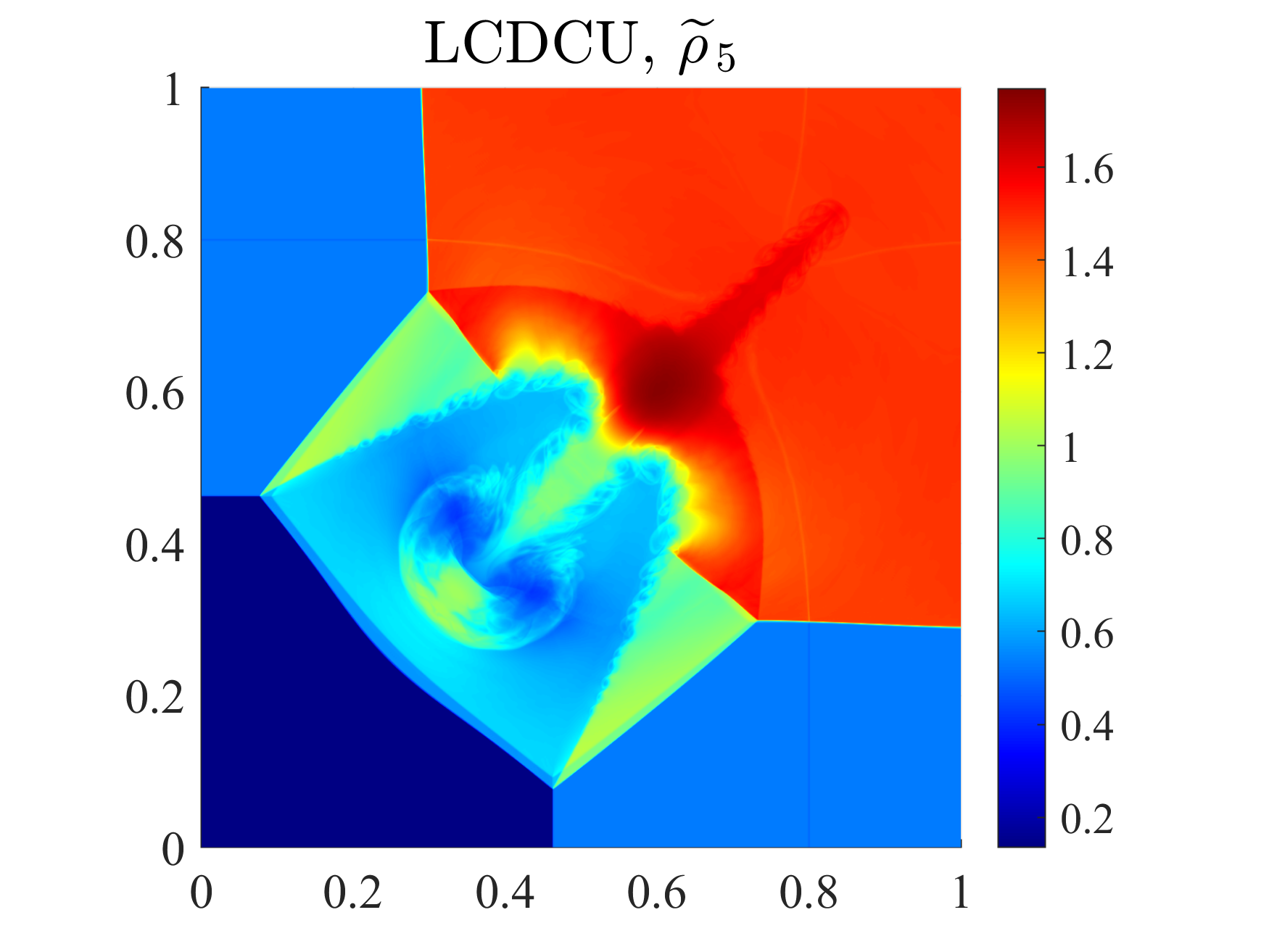}\hspace*{0.2cm}
            \includegraphics[trim=1.3cm 0.5cm 1.6cm 0.2cm, clip, width=4.cm]{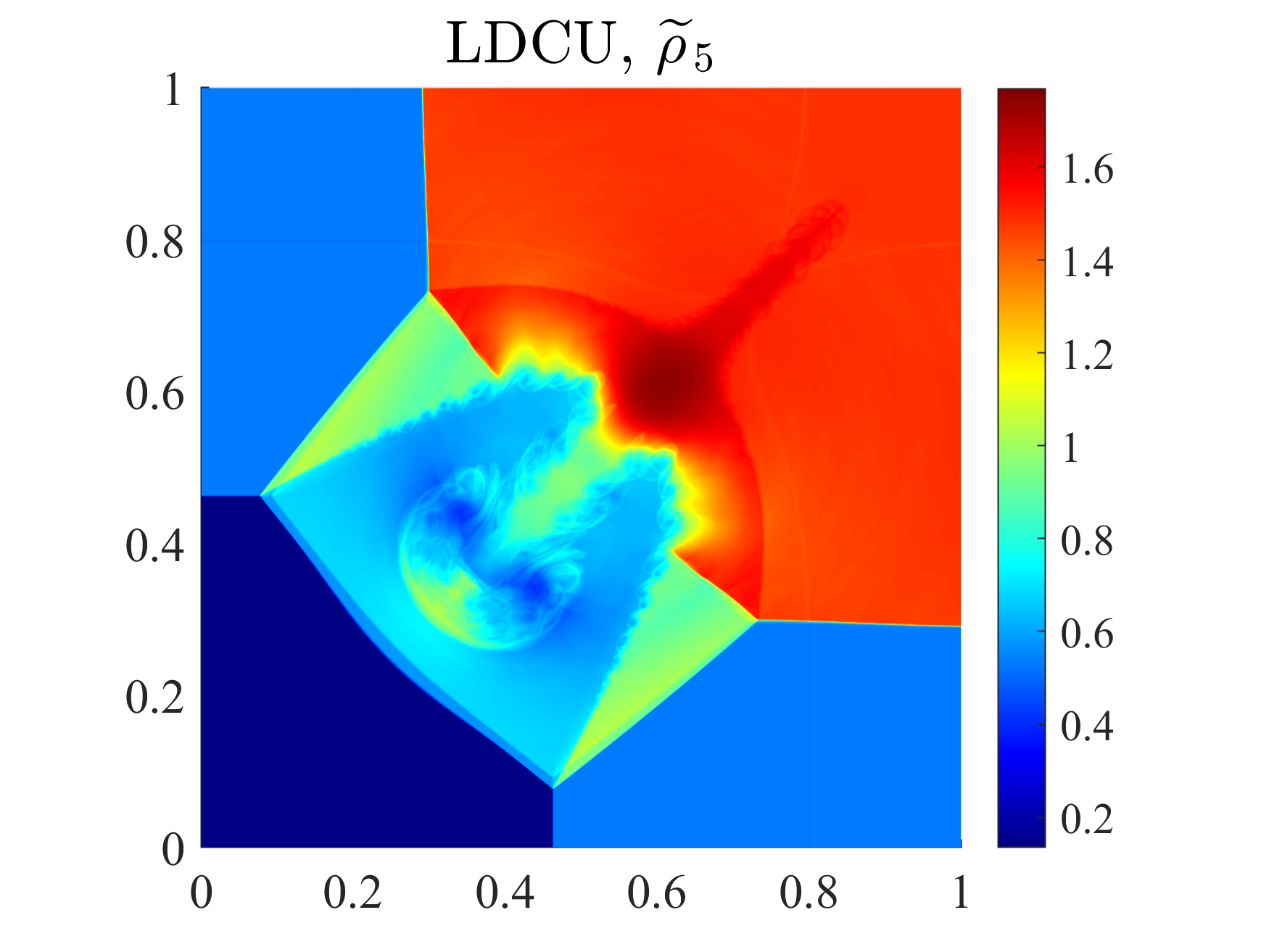}\hspace*{0.0cm}
            \includegraphics[trim=1.3cm 0.5cm 1.6cm 0.2cm, clip, width=4.2cm]{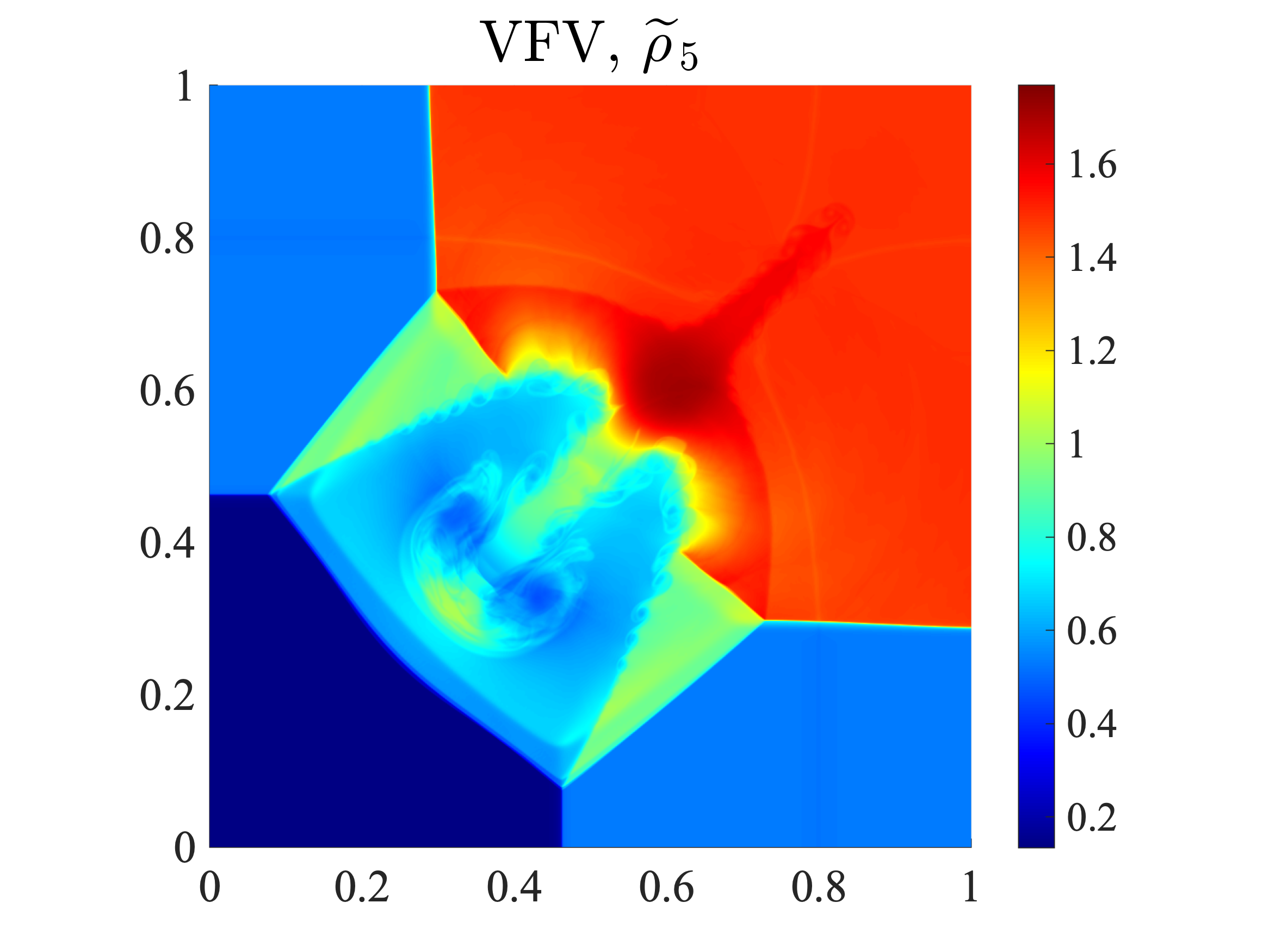}}
\caption{\sf Configuration 3: $\widetilde\rho_3$ (top row), $\widetilde\rho_4$ (middle row), and $\widetilde\rho_5$ (bottom row) computed by
the LCDCU (left column), LDCU (middle column), and VFV (right column) schemes.\label{fig23}}
\end{figure}

We compute the numerical solutions until the final time $T=0.8$ using the studied schemes and show the results, obtained by the first-,
third-, and seventh-order schemes in Figure \ref{fig22}. As one can see, the higher-order schemes can better capture the sideband
instability of the jet in the zones of strong along-jet velocity shear and the instability along the jet's neck. This indicates that the
numerical methods may produce true DW solutions, with convergence being only weak, as indicated by Theorem \ref{T1}. We then investigate the
convergence behavior of the average quantities. In Figure \ref{fig23}, we present the average densities $\widetilde\rho_3$,
$\widetilde\rho_4$, and $\widetilde\rho_5$ for the three studied methods. It can be observed that the average densities exhibit similar
profiles across the different orders of accuracy, indicating a strong convergence of the averaged sequence. Note, however, that the limiting
DW solution for different methods may differ; see the discussion below.

\begin{figure}[ht!]
\centerline{\includegraphics[trim=1.3cm 0.5cm 1.6cm 0.2cm, clip, width=4.cm]{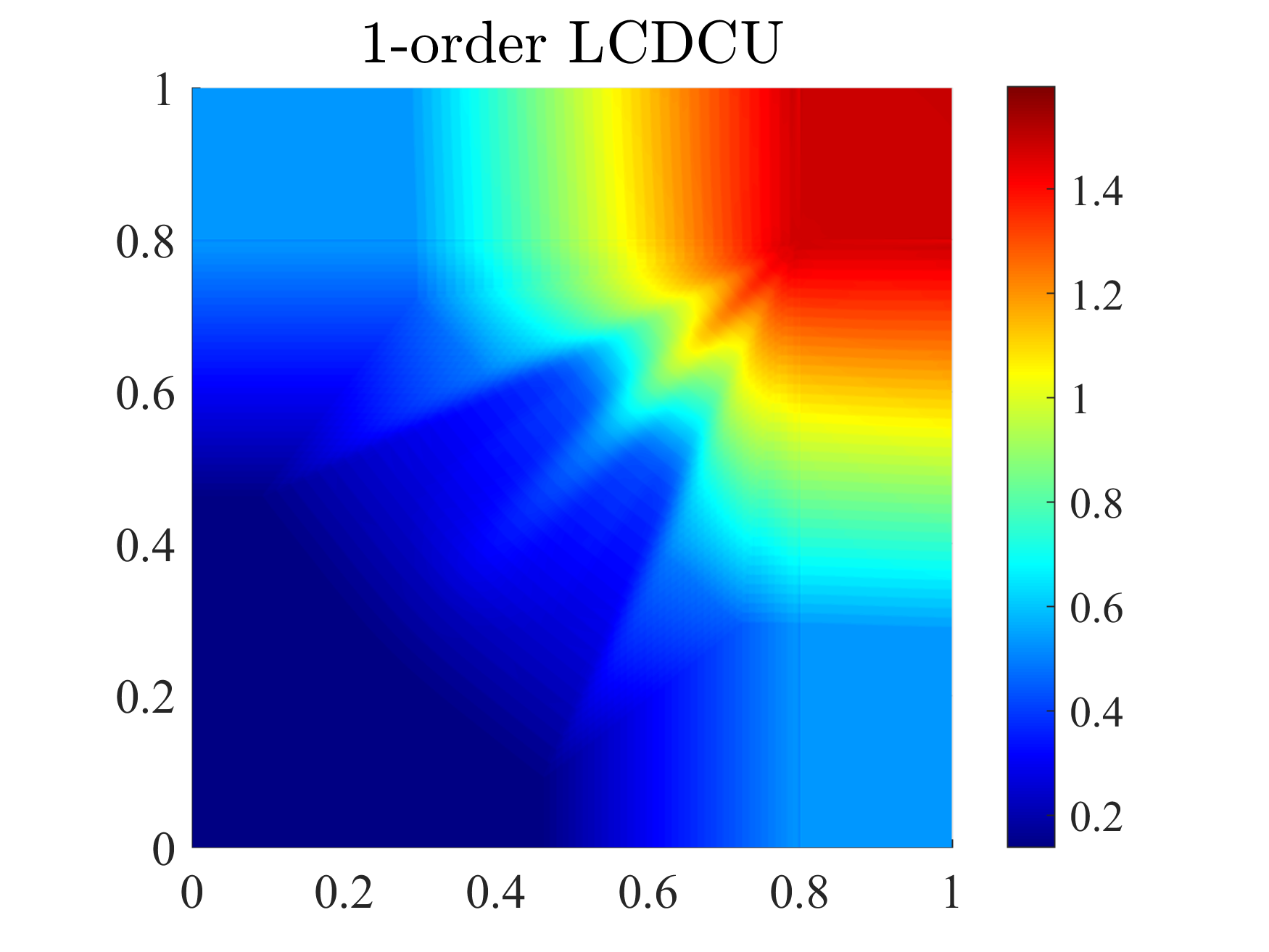}\hspace*{0.2cm}
            \includegraphics[trim=1.3cm 0.5cm 1.6cm 0.2cm, clip, width=4.cm]{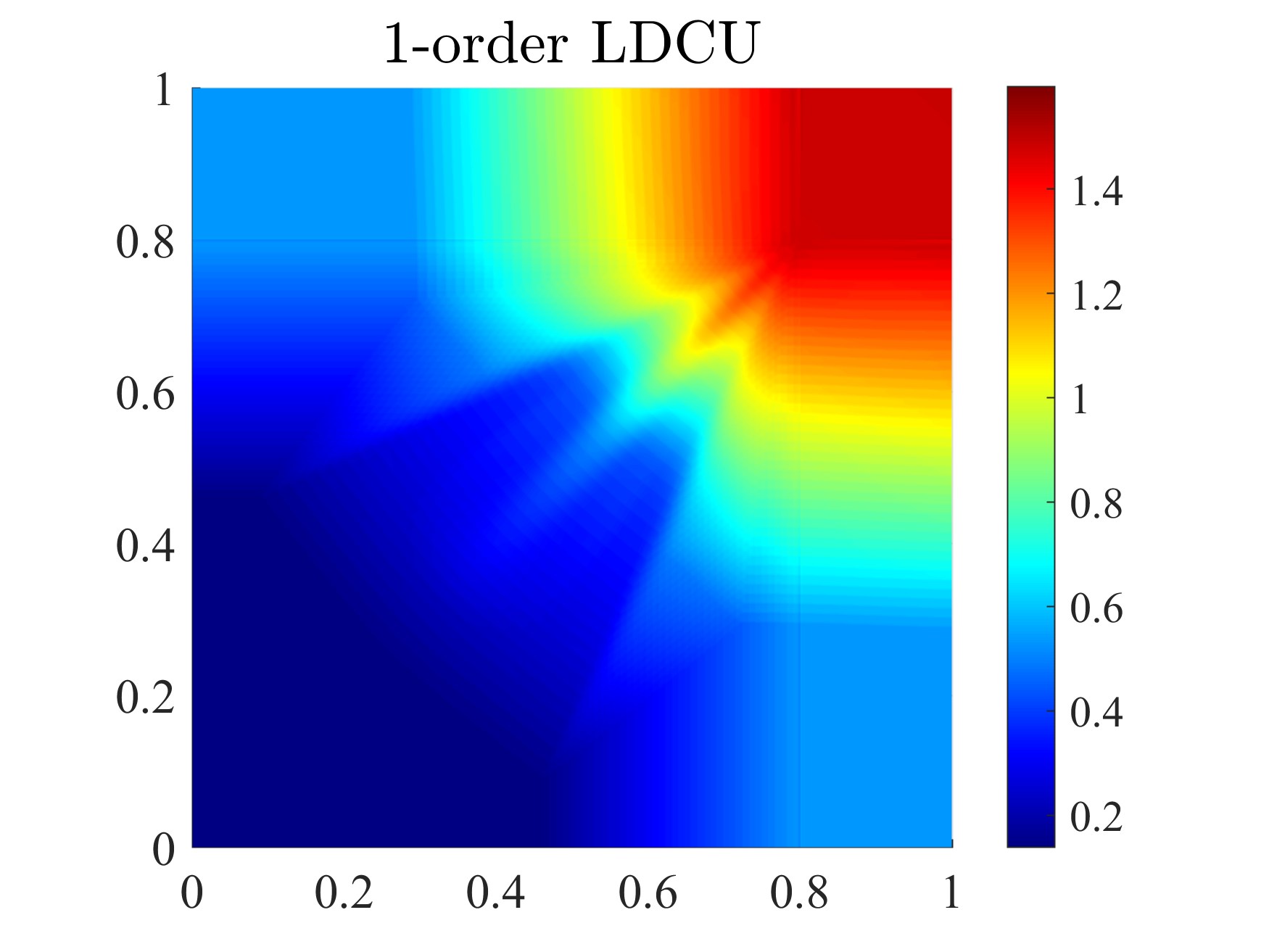}\hspace*{0.0cm}
            \includegraphics[trim=1.3cm 0.5cm 1.6cm 0.2cm, clip, width=4.2cm]{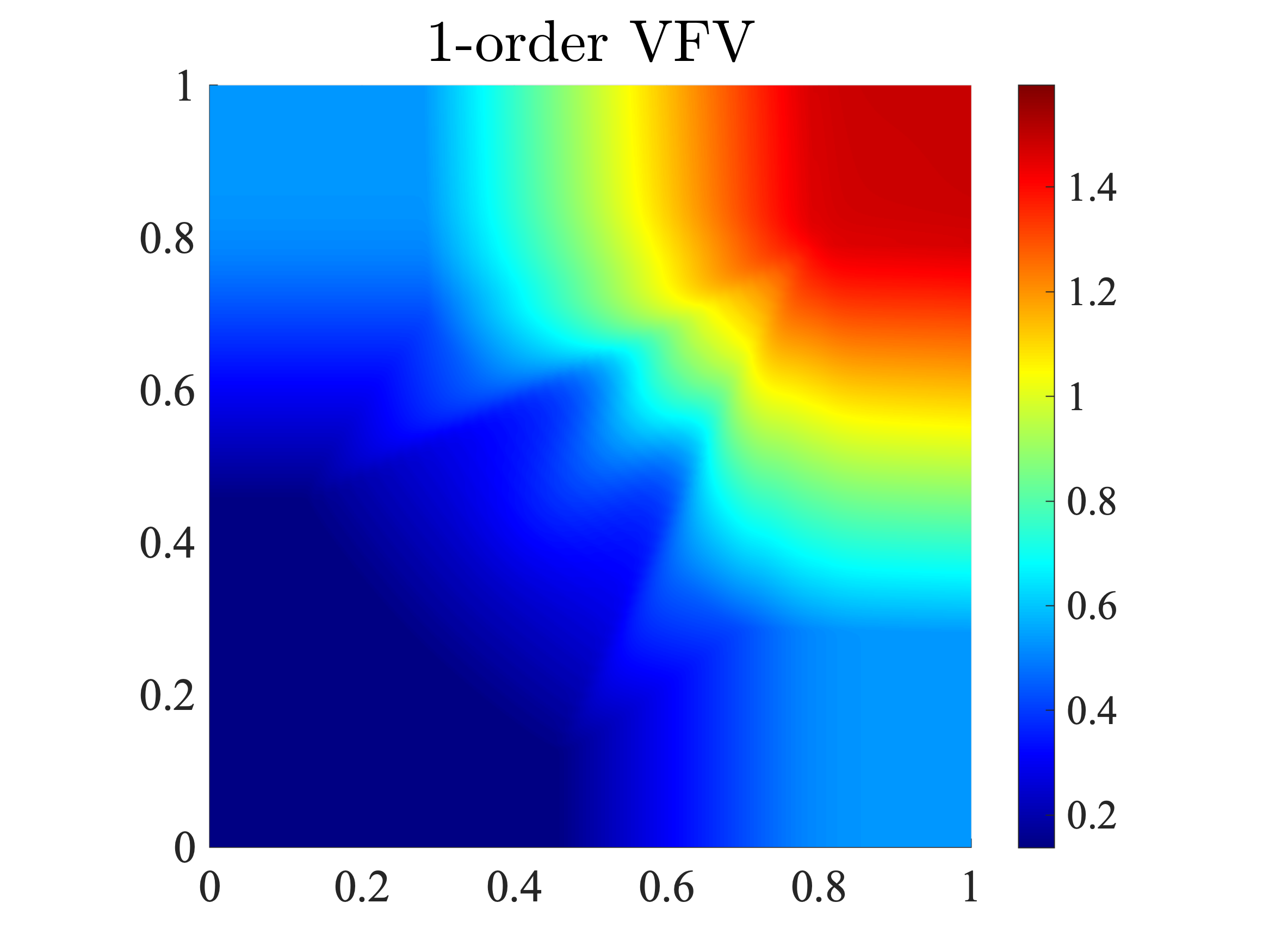}}
\vskip8pt
\centerline{\includegraphics[trim=1.3cm 0.5cm 1.6cm 0.2cm, clip, width=4.cm]{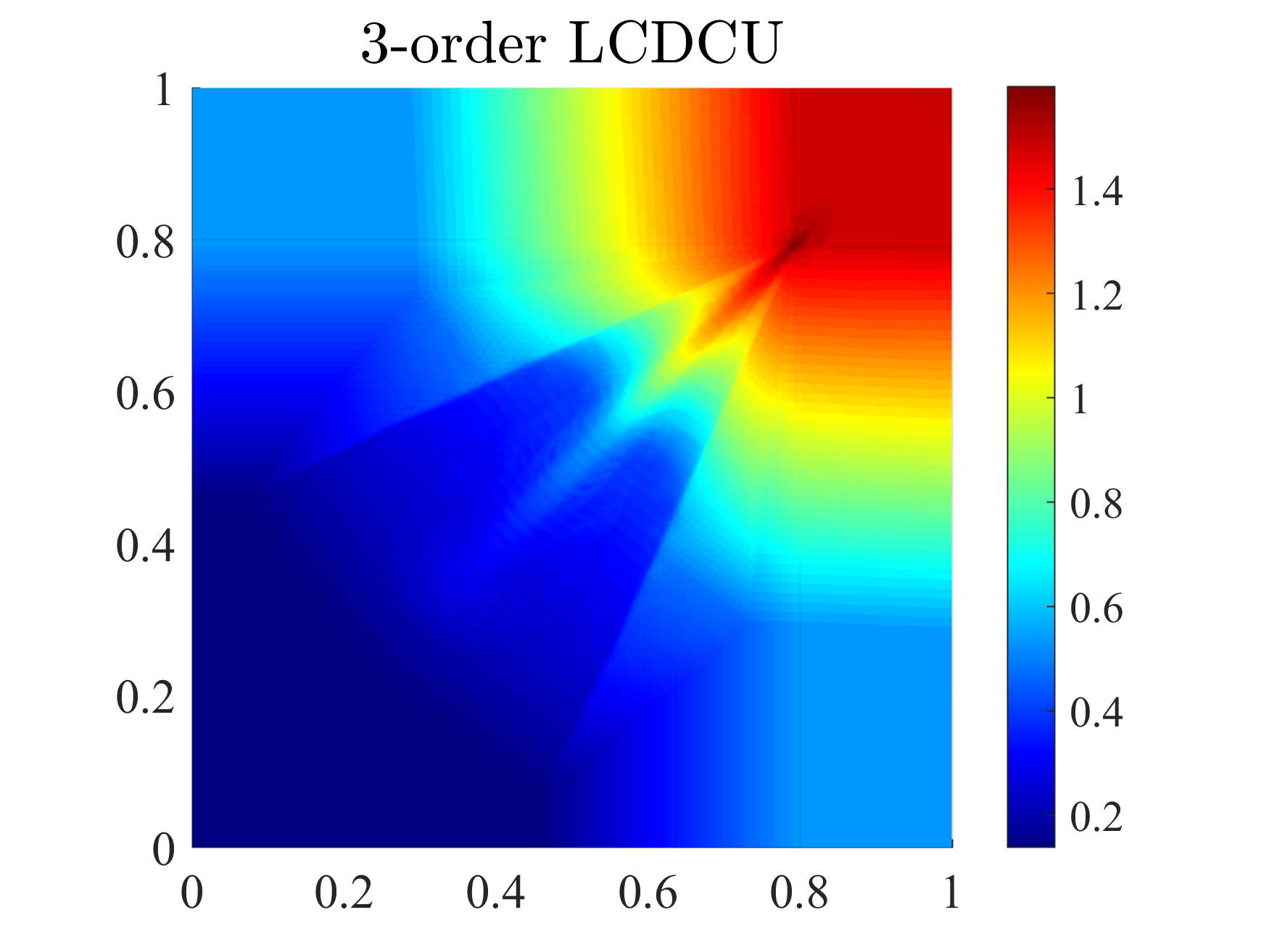}\hspace*{0.2cm}
            \includegraphics[trim=1.3cm 0.5cm 1.6cm 0.2cm, clip, width=4.cm]{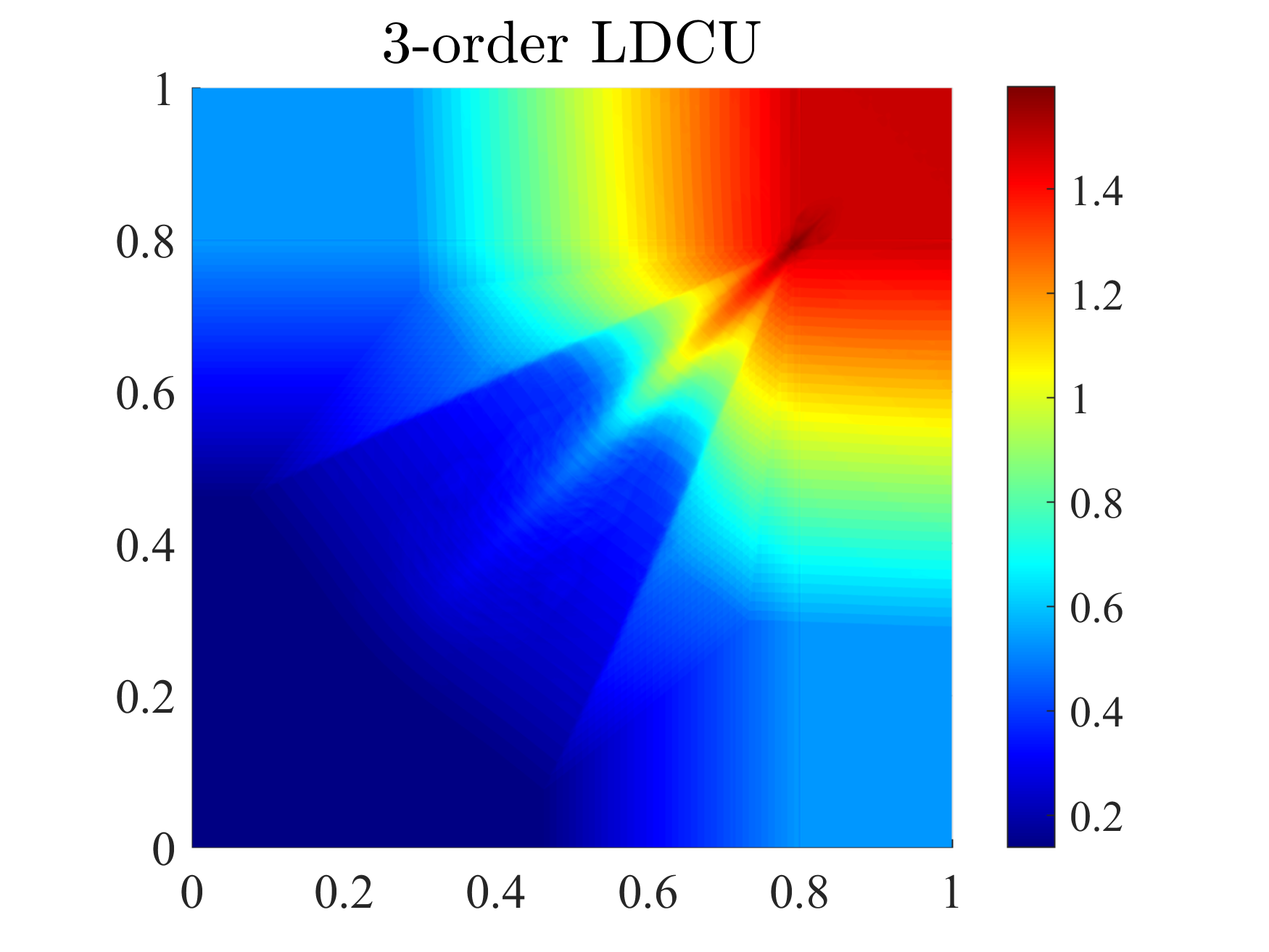}\hspace*{0.0cm}
            \includegraphics[trim=1.3cm 0.5cm 1.6cm 0.2cm, clip, width=4.2cm]{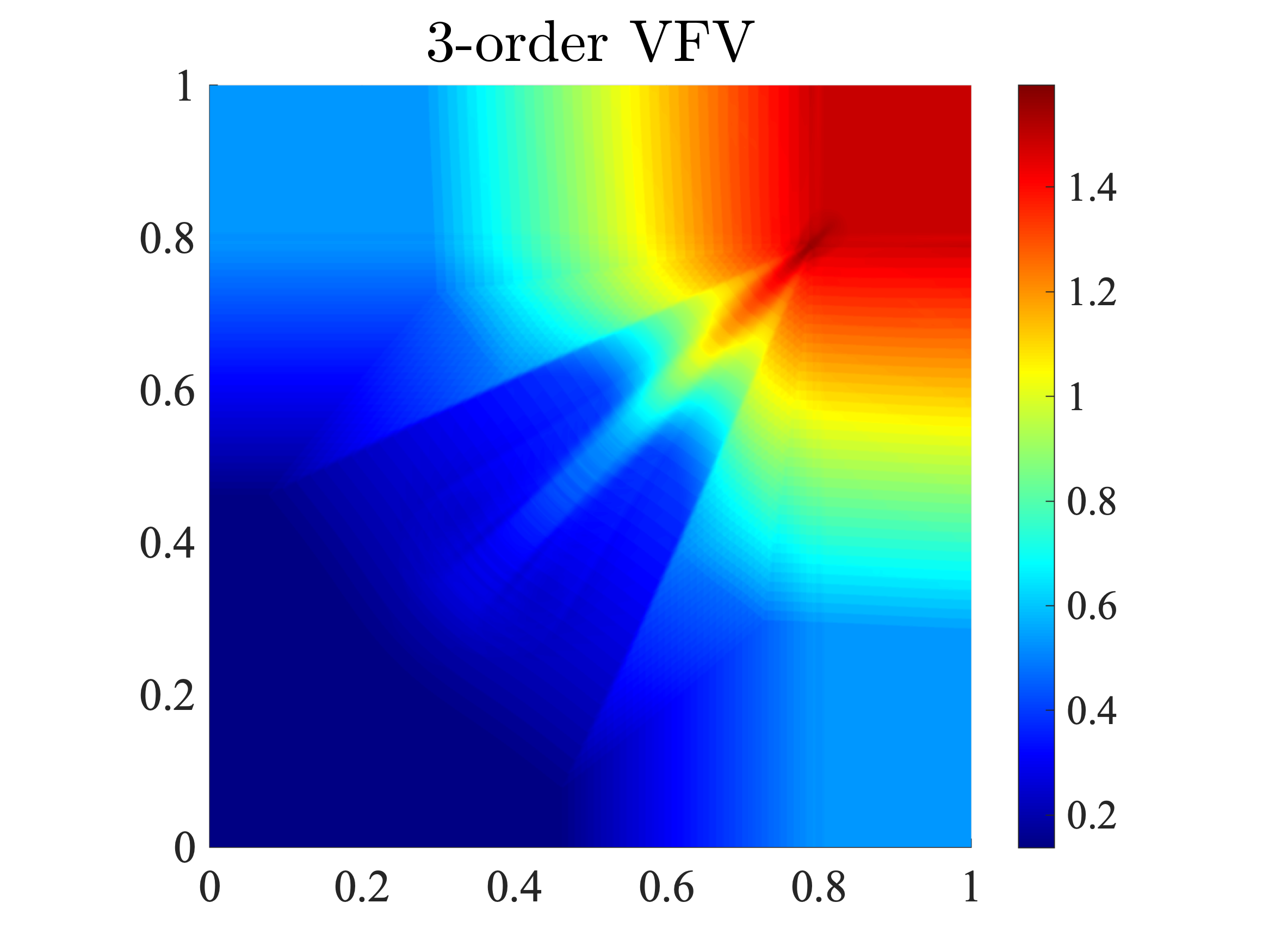}}
\vskip8pt
\centerline{\includegraphics[trim=1.3cm 0.5cm 1.6cm 0.2cm, clip, width=4.cm]{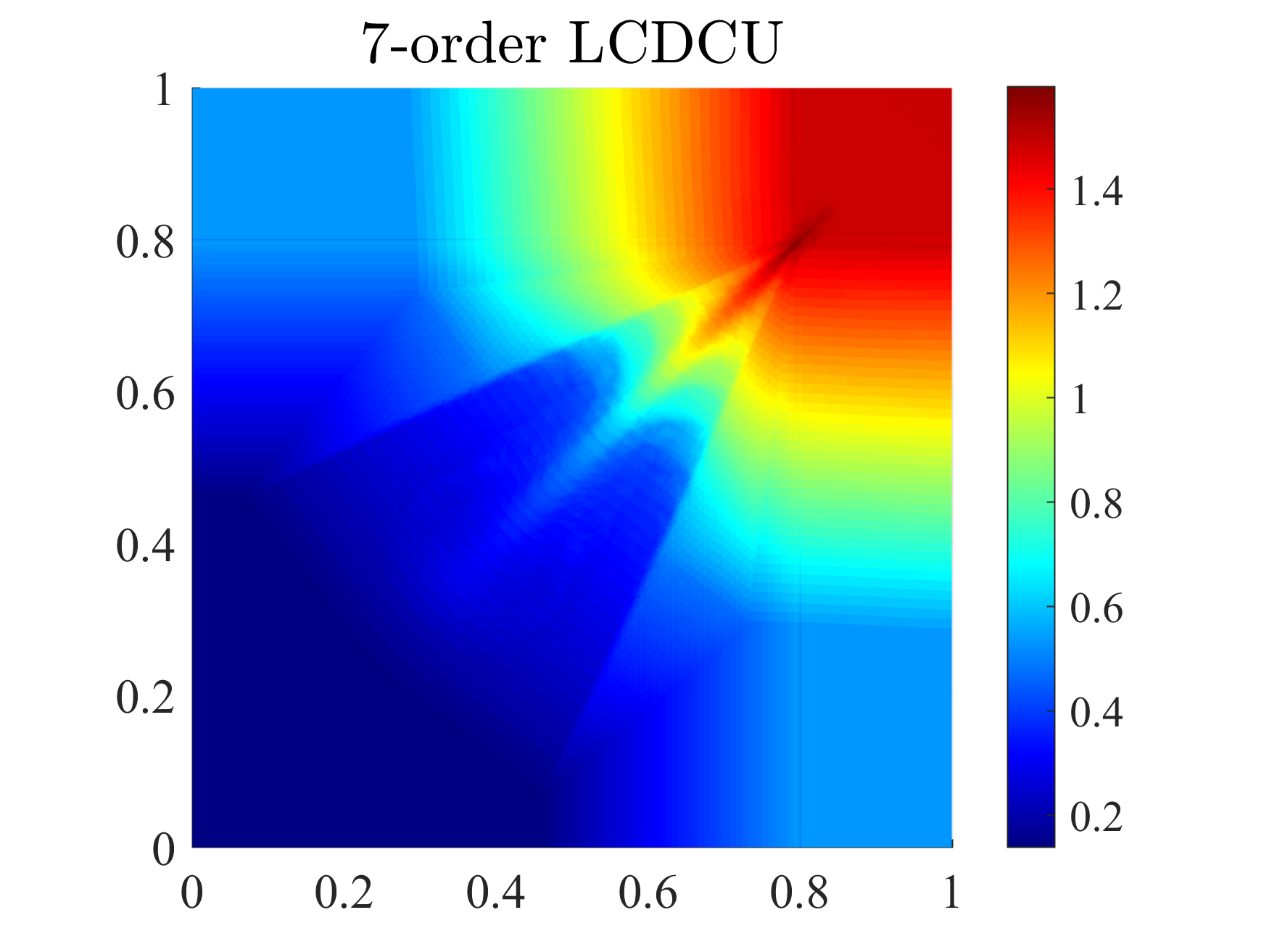}\hspace*{0.2cm}
            \includegraphics[trim=1.3cm 0.5cm 1.6cm 0.2cm, clip, width=4.cm]{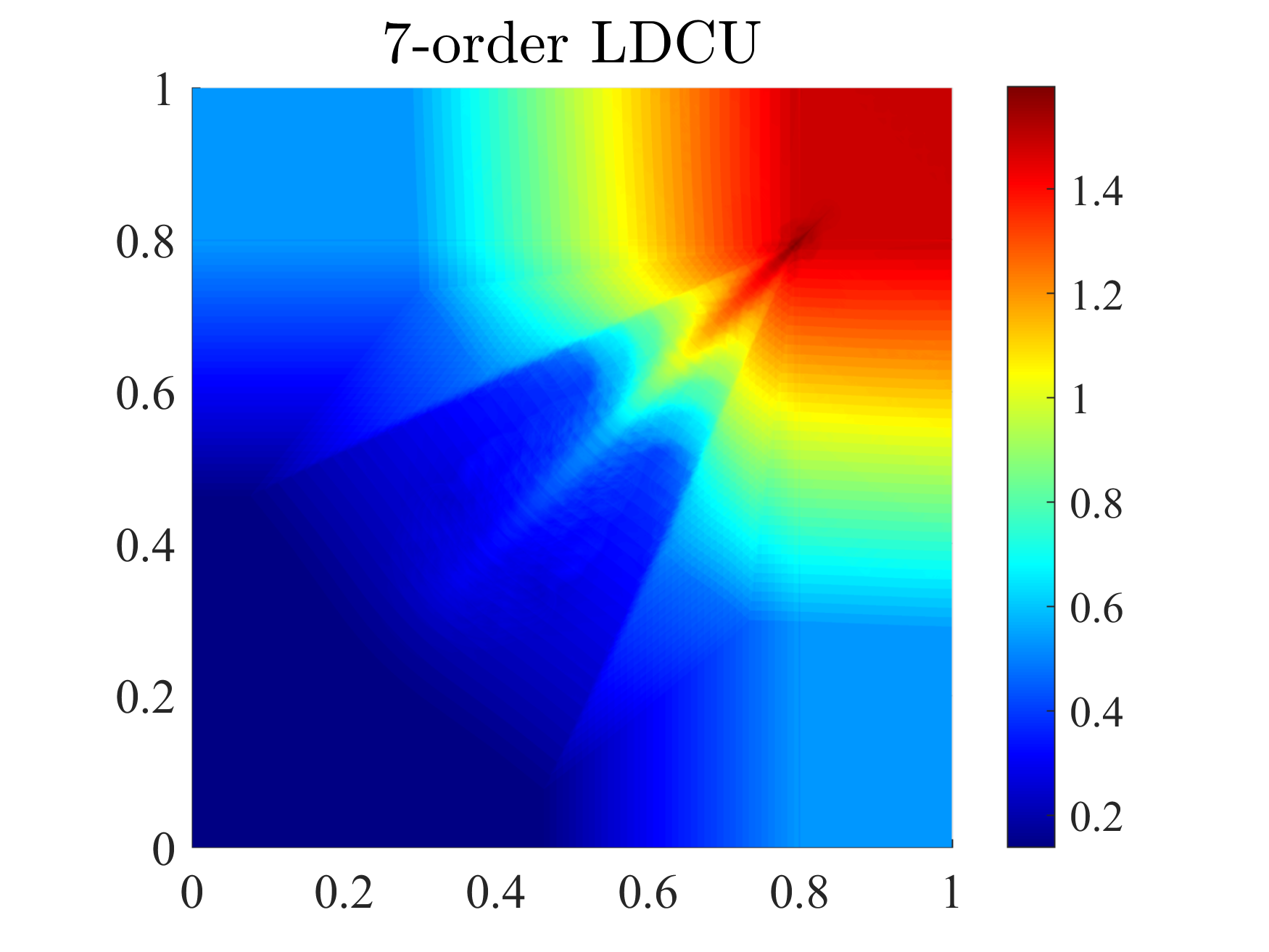}\hspace*{0.0cm}
            \includegraphics[trim=1.3cm 0.5cm 1.6cm 0.2cm, clip, width=4.2cm]{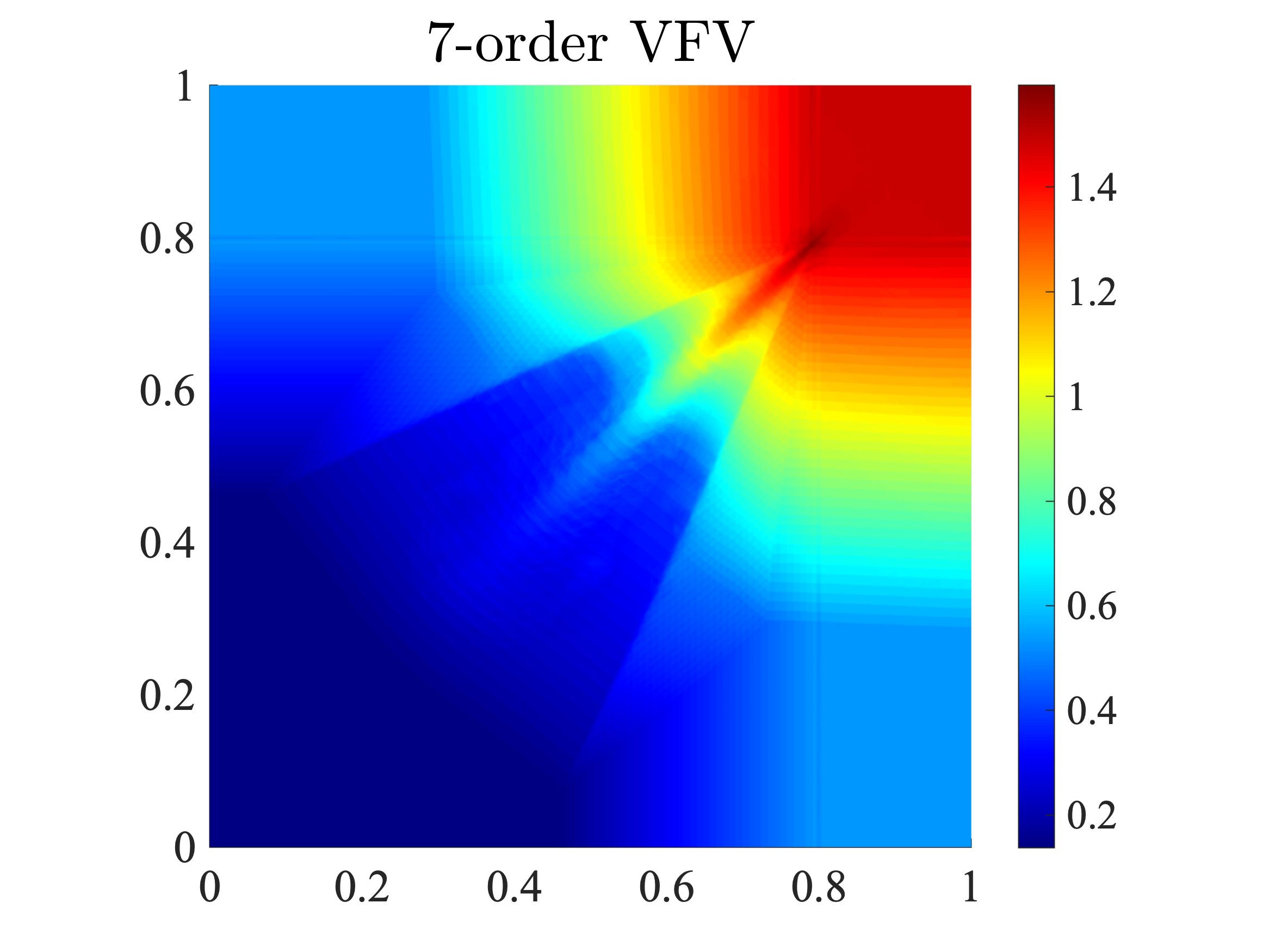}}
\caption{\sf Configuration 3: $\rho^T_1$ (top row), $\rho^T_3$ (middle row), and $\rho^T_5$ (bottom row) computed by the LCDCU (left
column), LDCU (middle column), and VFV (right column) schemes.\label{fig24}}
\end{figure}

\begin{figure}[ht!]
\centerline{\includegraphics[trim=1.3cm 0.5cm 1.6cm 0.2cm, clip, width=4.cm]{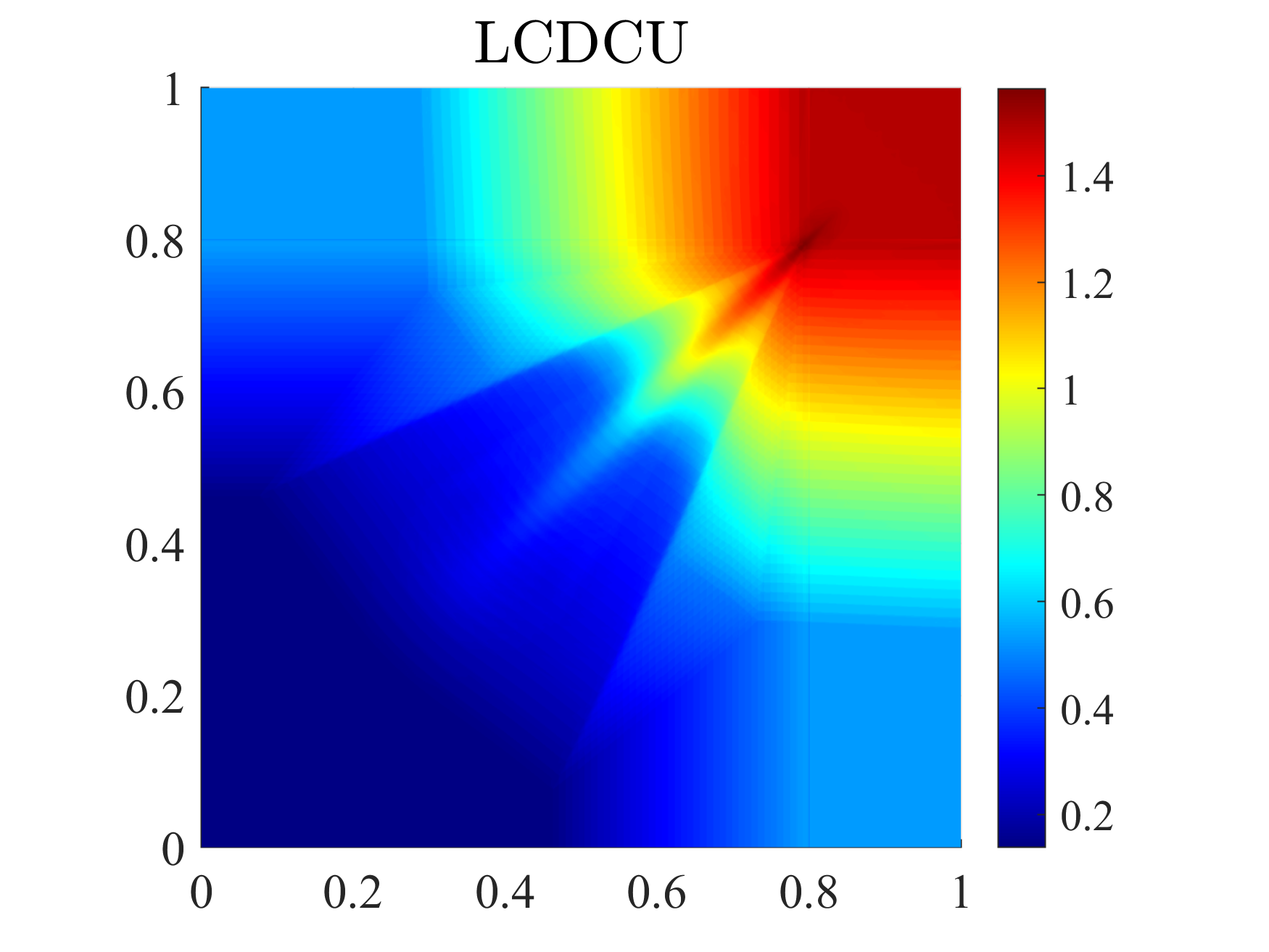}\hspace*{0.2cm}
            \includegraphics[trim=1.3cm 0.5cm 1.6cm 0.2cm, clip, width=4.cm]{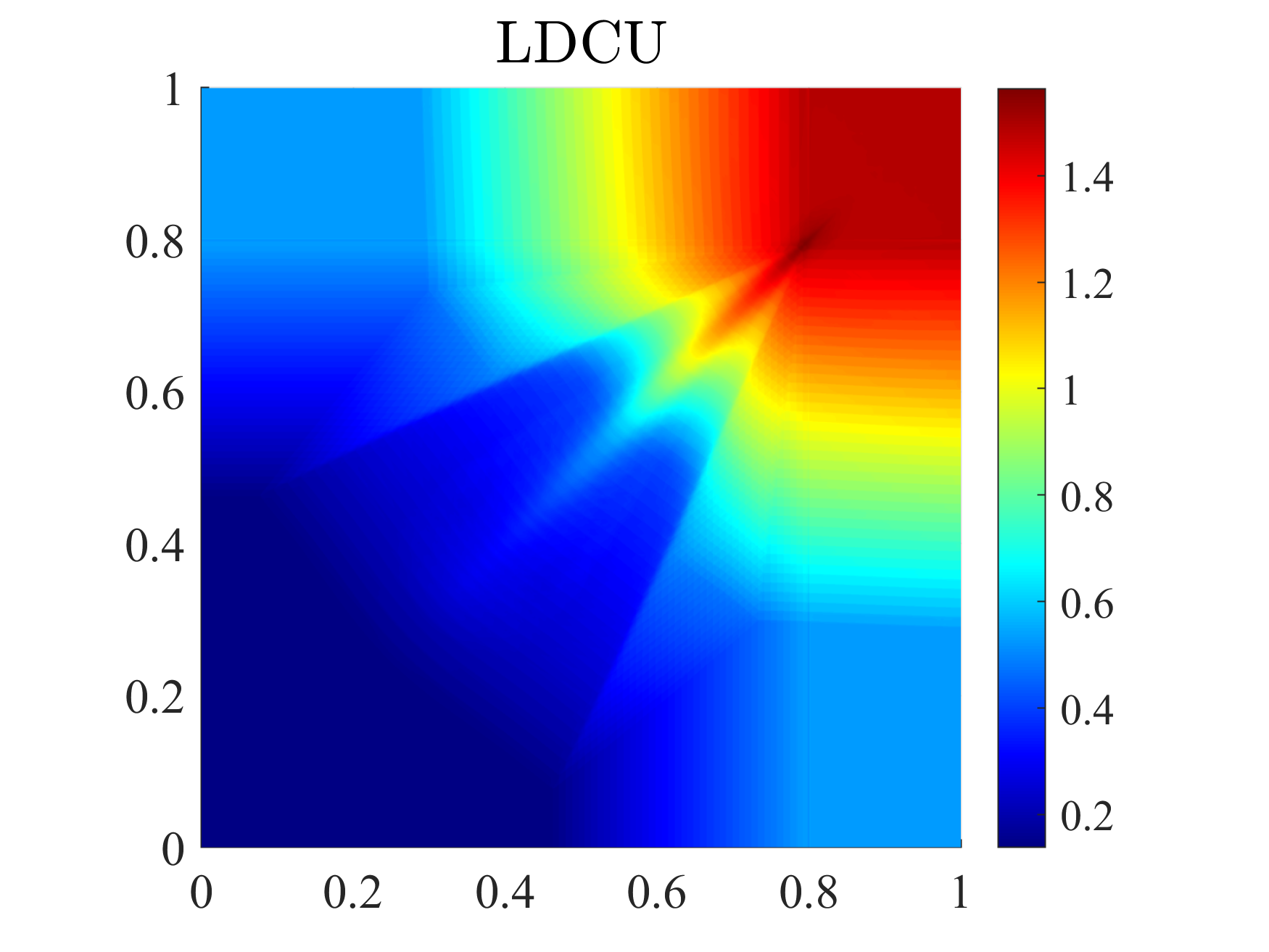}\hspace*{0.0cm}
            \includegraphics[trim=1.3cm 0.5cm 1.6cm 0.2cm, clip, width=4.2cm]{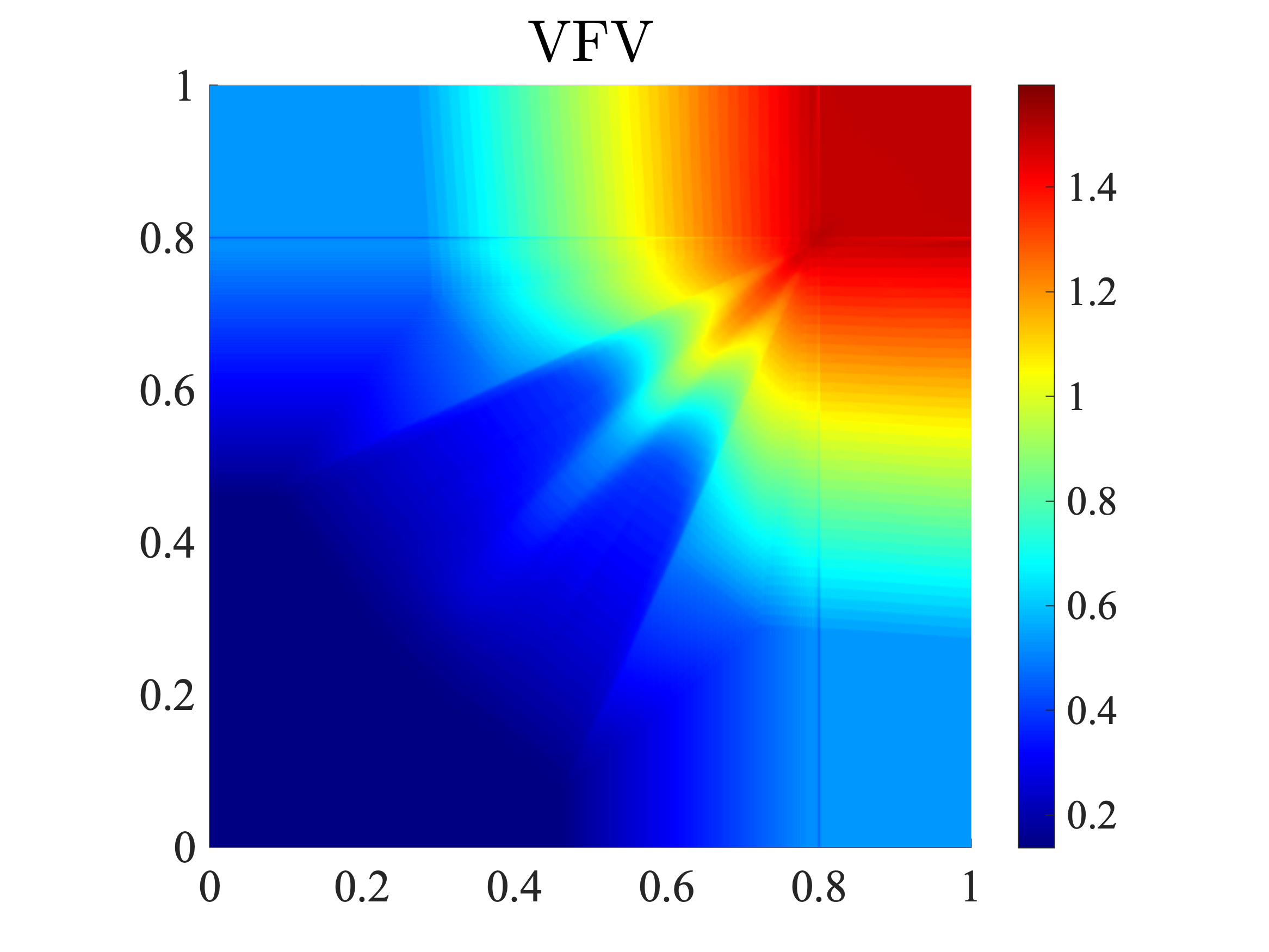}}
\caption{\sf Configuration 3: $\widetilde\rho^{\,T}_5$ computed by the LCDCU (left), LDCU (middle), and VFV (right) schemes.\label{fig25}}
\end{figure}

Next, we approximate the time averages of a DW solution by evaluating the time averages of the computed densities, which are denoted by
\begin{equation*}
\rho^T_\ell:=\frac{1}{T}\int\limits_0^T\rho_\ell(\cdot,t)\,{\rm d}t, 
\end{equation*}
where $T$ is the final time. Figure \ref{fig24} shows $\rho^T_\ell$ for the three studied methods and Figure \ref{fig25} presents the time
averages of $\widetilde{\rho}_5$, which we denote by 
\begin{equation*}
\widetilde\rho^{\,T}_5:=\frac{1}{T}\int\limits_0^T\widetilde\rho_5(\cdot,t)\,{\rm d}t.
\end{equation*}
These figures suggest that there is almost no difference between the time averaged flows computed by different methods. We note that time
averages $\widetilde\rho^{\,T}_5$ are typically used in the analysis of turbulent flows.

\begin{figure}[ht!]
\centerline{\includegraphics[trim=1.2cm 0.5cm 1.4cm 0.1cm, clip, width=3.6cm]{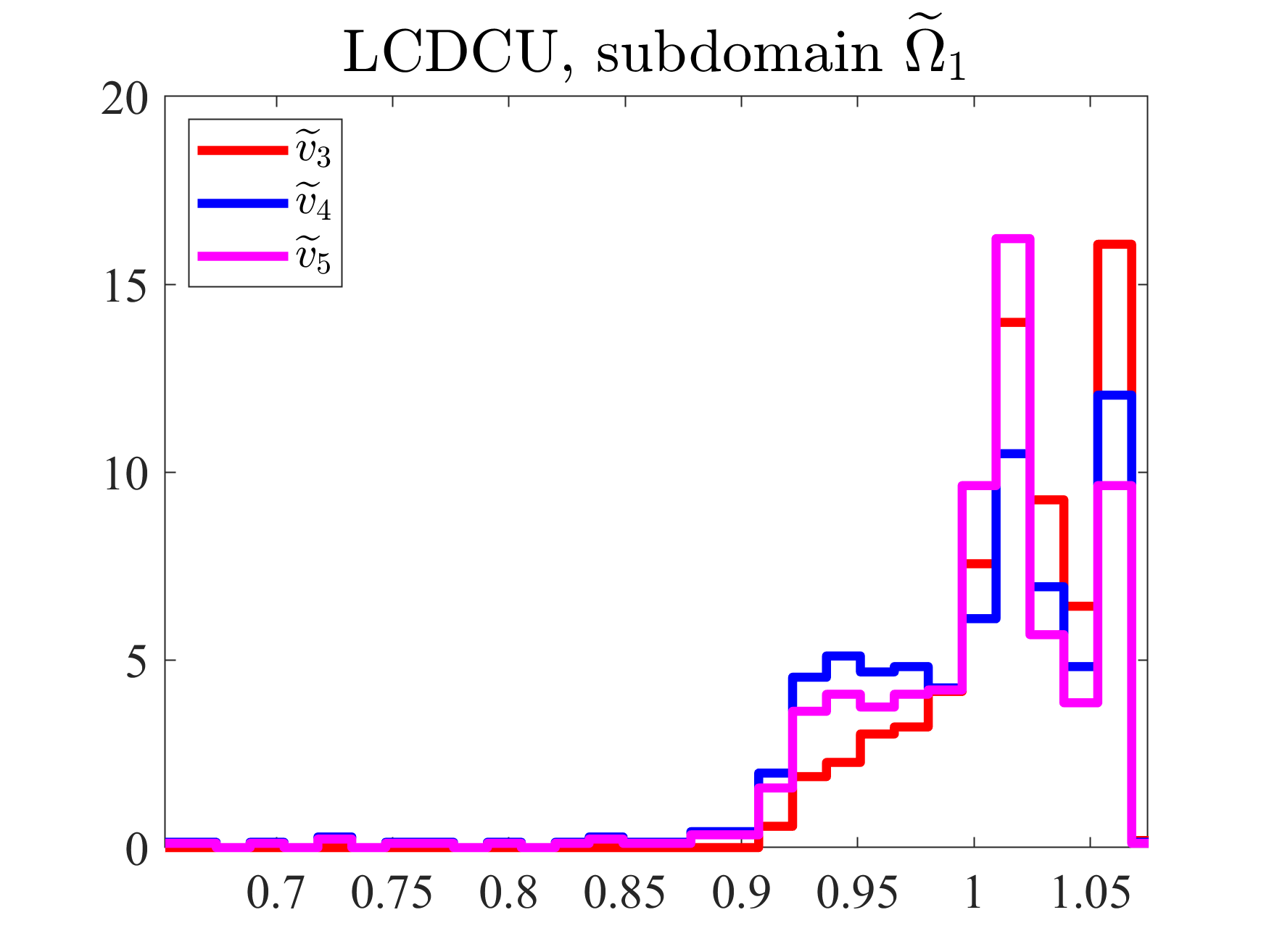}\hspace*{0.6cm}
            \includegraphics[trim=1.2cm 0.5cm 1.4cm 0.1cm, clip, width=3.6cm]{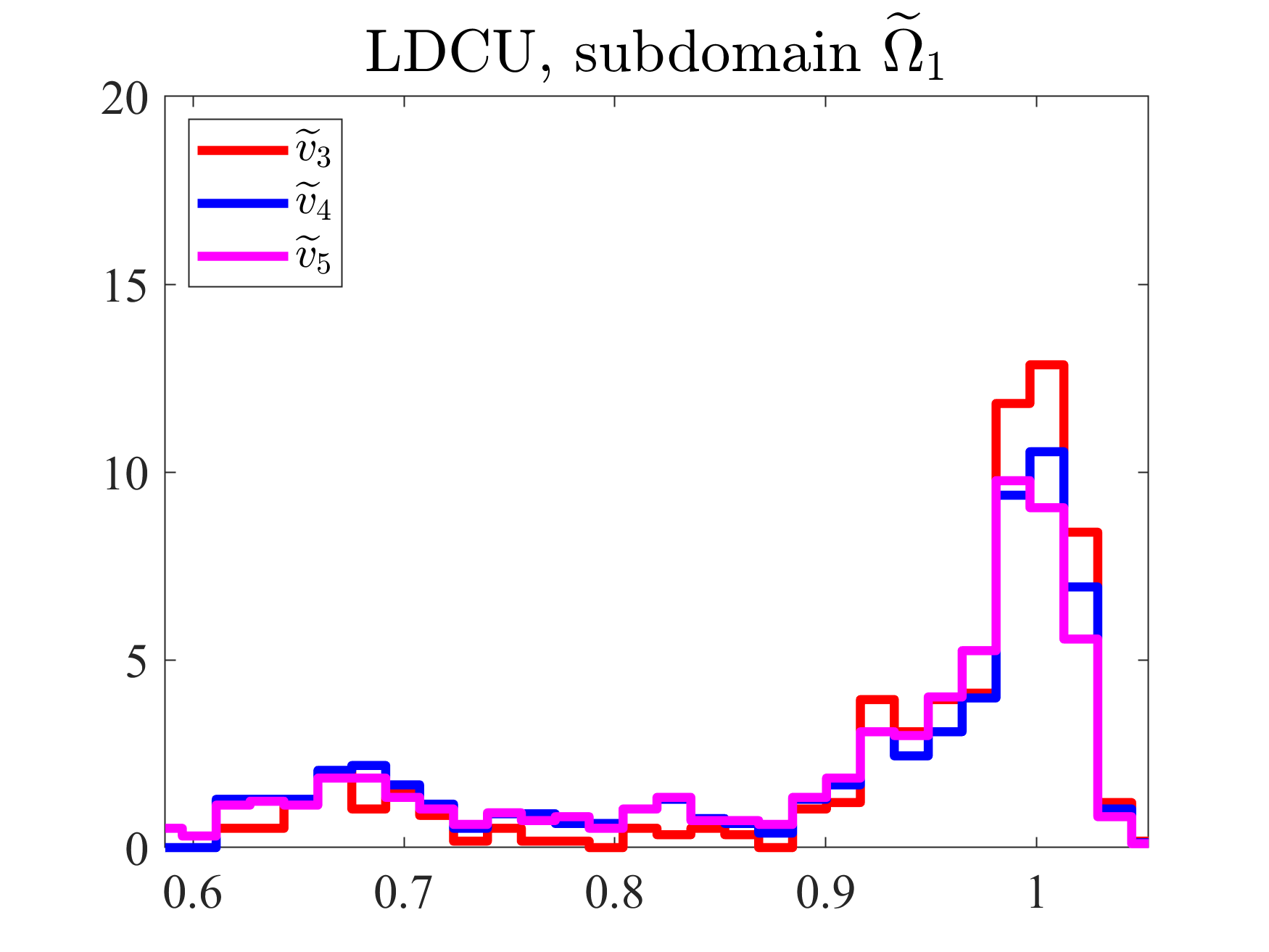}\hspace*{0.4cm}
            \includegraphics[trim=1.2cm 0.5cm 1.4cm 0.1cm, clip, width=3.8cm]{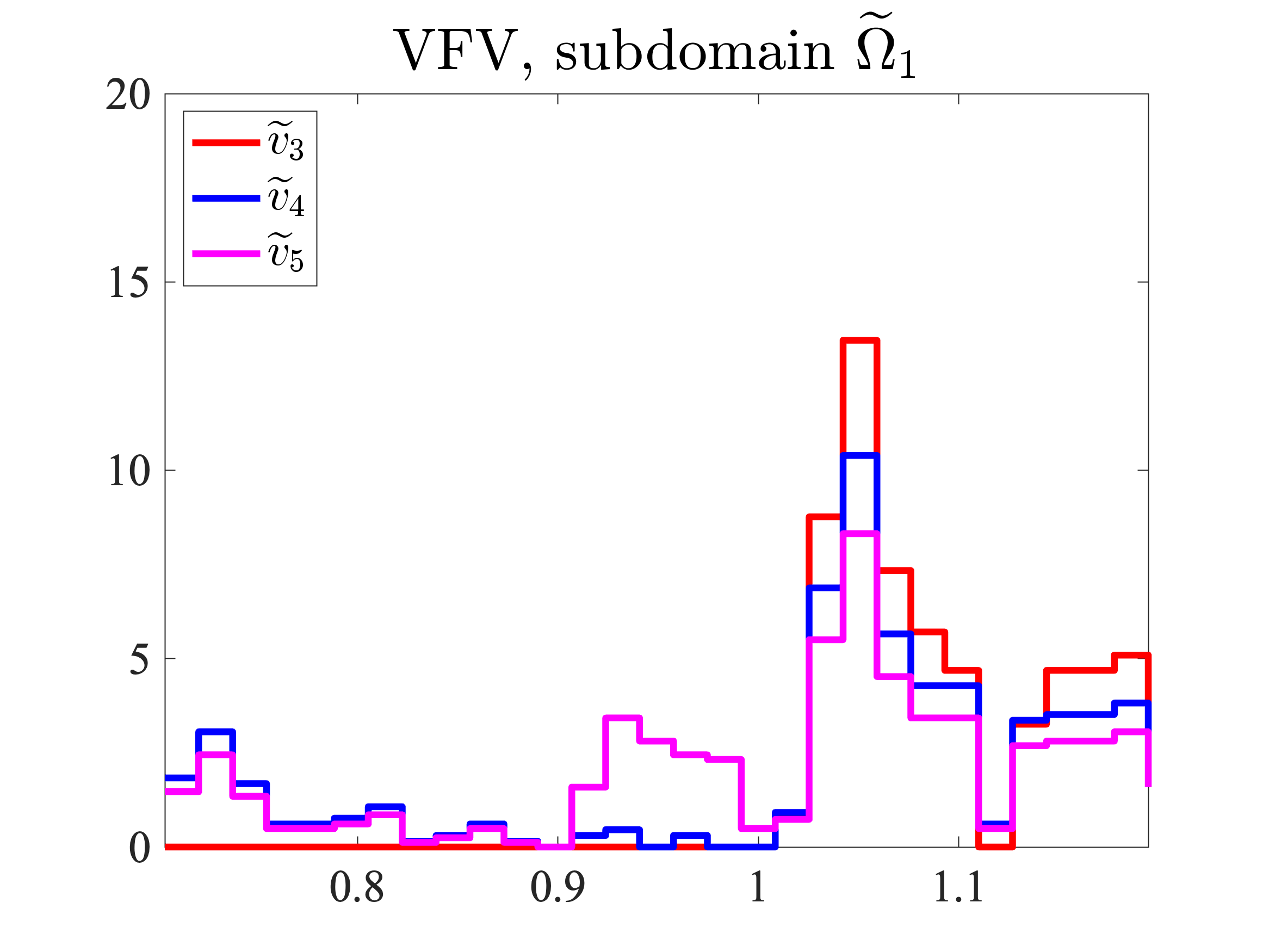}}
\vskip8pt
\centerline{\includegraphics[trim=1.2cm 0.5cm 1.4cm 0.1cm, clip, width=3.6cm]{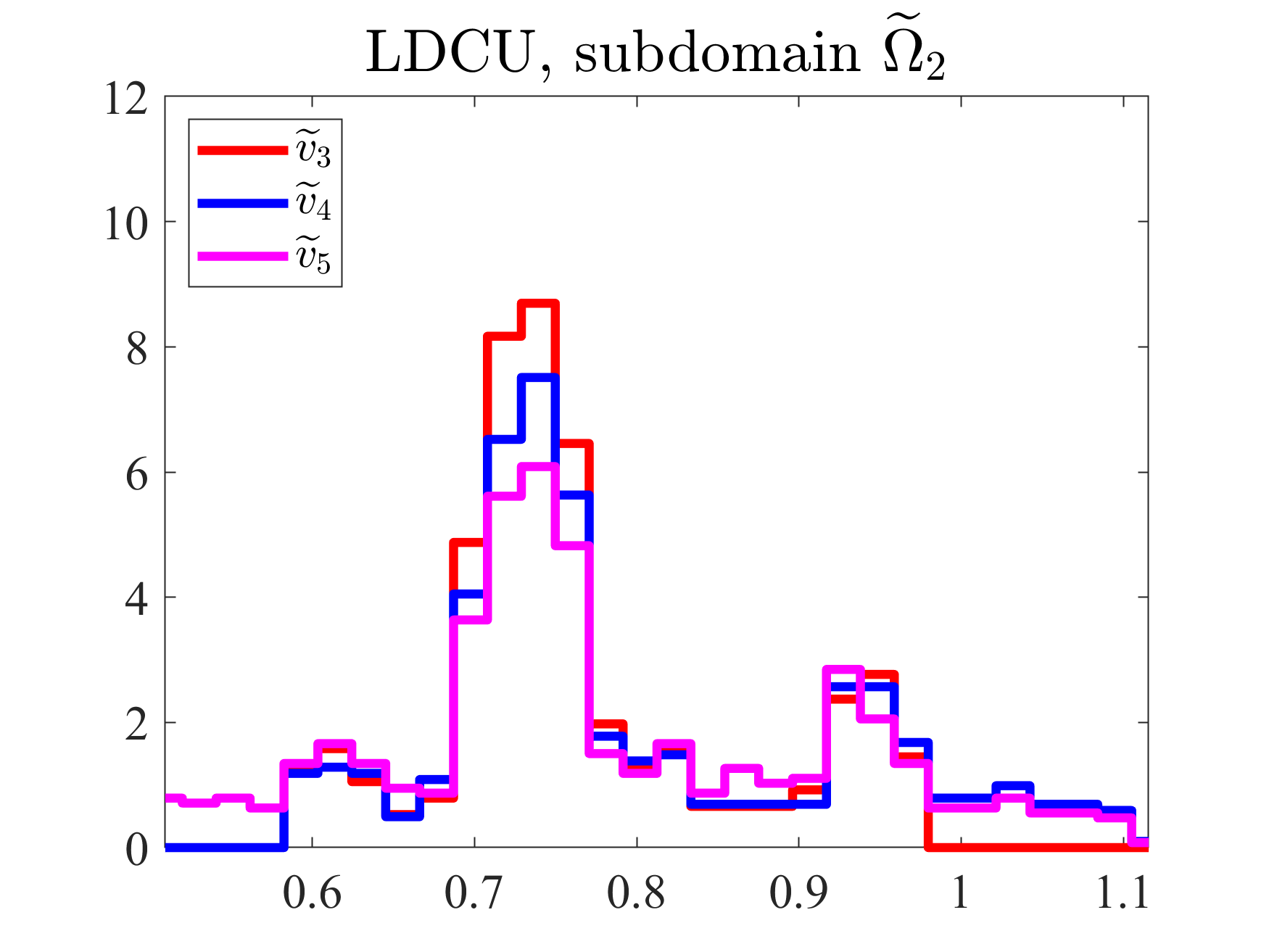}\hspace*{0.6cm}
            \includegraphics[trim=1.2cm 0.5cm 1.4cm 0.1cm, clip, width=3.6cm]{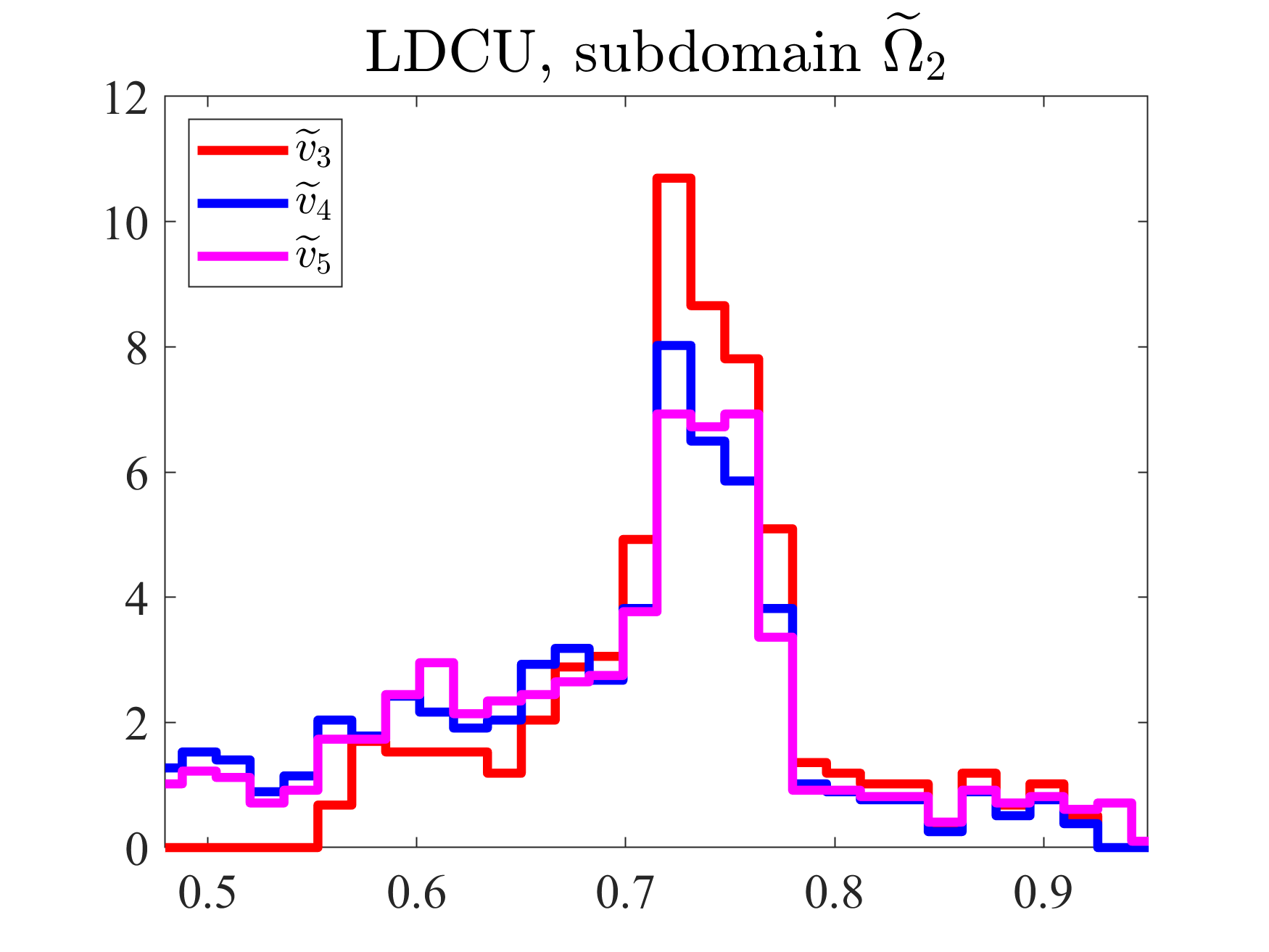}\hspace*{0.4cm}
            \includegraphics[trim=1.2cm 0.5cm 1.4cm 0.1cm, clip, width=3.8cm]{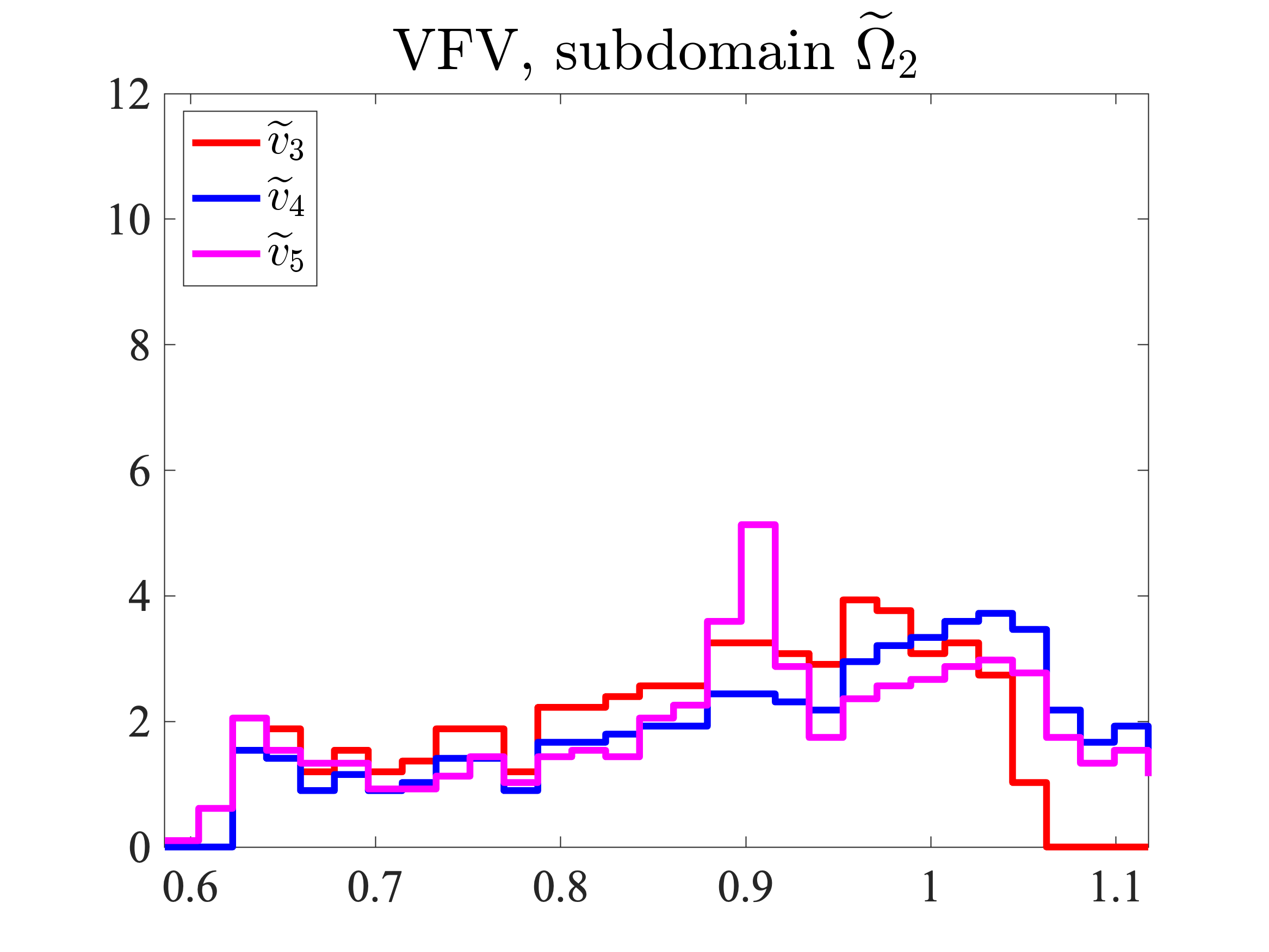}}
\caption{\sf Configuration 3: $\widetilde\sigma_3$, $\widetilde\sigma_4$, and $\widetilde\sigma_5$ computed by the LCDCU (left column), LDCU
(middle column), and VFV (right column) schemes in the subdomains $\widetilde\Omega_1$ (top row) and $\widetilde\Omega_2$ (bottom row).
\label{fig26}}
\end{figure}

Finally, we establish the experimental PDFs of the Young measures for the three studied schemes. In this example, we choose the following
two subdomains:
\begin{equation*}
\widetilde\Omega_1=[0.42,0.43]\times[0.63,0.64]\quad\mbox{and}\quad\widetilde\Omega_2=[0.60,0.61]\times[0.37,0.38],
\end{equation*}
compute $\widetilde\sigma_3$, $\widetilde\sigma_4$, and $\widetilde\sigma_5$, and plot the obtained results in Figure \ref{fig26}. As one
can see, the average sequences from various numerical methods yield different limiting PDFs, particularly for subdomain
$\widetilde\Omega_1$, indicating that the order-averaged solutions converge to different DW solutions.

\subsection{Kelvin-Helmholtz (KH) Instability}
In the last example, we consider the KH instability problem taken from \cite{feireisl2021computing} with the following initial conditions:
$$
(\rho,u,v,p)(x,y,0)=\begin{cases}(2,-0.5,0,2.5),&\mbox{if }0.25+0.01Y_1(x,\omega)<y\\
                                                &\mbox{and }y <0.75+0.01Y_2(x,\omega),\\
(1,0.5,0,2.5),&\mbox{otherwise},\end{cases}
$$
where
$$
Y_j(x,\omega)=\sum_{n=1}^{10}a^n_j(\omega)\cos\big(b^n_j(\omega)+2n\pi x\big),\quad j=1,2,
$$
with $a^n_j=a^n_j(\omega)\in[0,1]$ and $b^n_j=b^n_j(\omega)\in(-\pi,\pi)$ being fixed random numbers, and $a^n_j$ being normalized such that
$\sum_{n=1}^{10}a^n_j=1$. The particular values of $a^n_j$ and $b^n_j$, which have been generated by \texttt{rand} in MATLAB, are given in
Table \ref{table1a}. The boundary conditions in this example are periodic. 
\begin{table}[ht!]
\centering
\begin{adjustbox}{width=\textwidth}
\begin{tabular}{ccccc}
\hline
& $a^n_1$ & $a^n_2$ & $b^n_1$ & $b^n_2$   \\
\hline
1 &6.848086824246653e-08&9.373025805955863e-03&-0.973625473853271& 3.10750325239443\\
2 &4.450348128947341e-03&1.976861219341060e-02& 2.33221742979395 & 1.74829850637860\\
3 &6.156955958786613e-02&1.29928159795144e-01 &-2.57661895600041 &-3.01803367486339\\
4 &1.16481555805349e-01 &2.15403817045303e-01 & 2.43965931651801 &-2.07430785001108\\
5 &1.68204784555961e-01 &1.758678905771403e-02& 1.26278768686501 & 3.10856394704146\\
6 &1.46413246162863e-01 &2.11994809652299e-01 & 1.47373734867445 & 1.59399689987015\\
7 &5.857224323849688e-02&2.134903127162342e-02&-1.25553236484458 &-1.77202310331374\\
8 &1.59868678328304e-01 &1.48283301108284e-01 &-2.82920698380582 & 1.00086476379714\\
9 &1.39003488203723e-01 &1.98161392506709e-01 & 2.56472949845762 & 0.245159802027399\\
10&1.45436027507622e-01 &2.815106156355688e-02&-2.52798558841484 &-0.125265541490568\\
\hline
\end{tabular}
\end{adjustbox}\\[1.ex] 
\caption{\sf The values of $a^n_1$, $a^n_2$, $b^n_1$, and $b^n_2$ for $n=1,\dots,10$ generated by \texttt{rand} in MATLAB.\label{table1a}}
\end{table}

\begin{figure}[ht!]
\centerline{\includegraphics[trim=1.3cm 0.5cm 1.6cm 0.2cm, clip, width=4.cm]{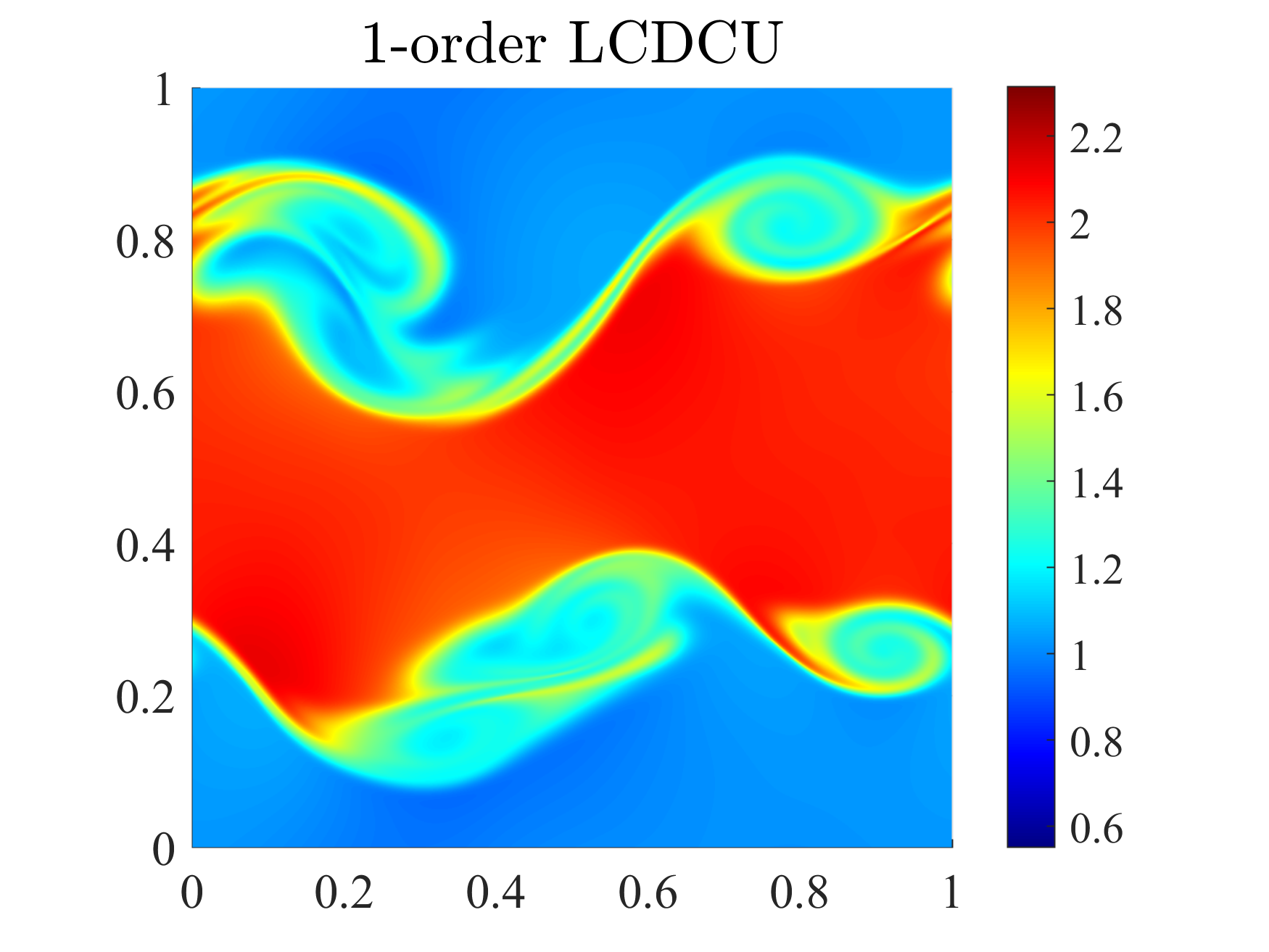}\hspace*{0.2cm}
            \includegraphics[trim=1.3cm 0.5cm 1.6cm 0.2cm, clip, width=4.cm]{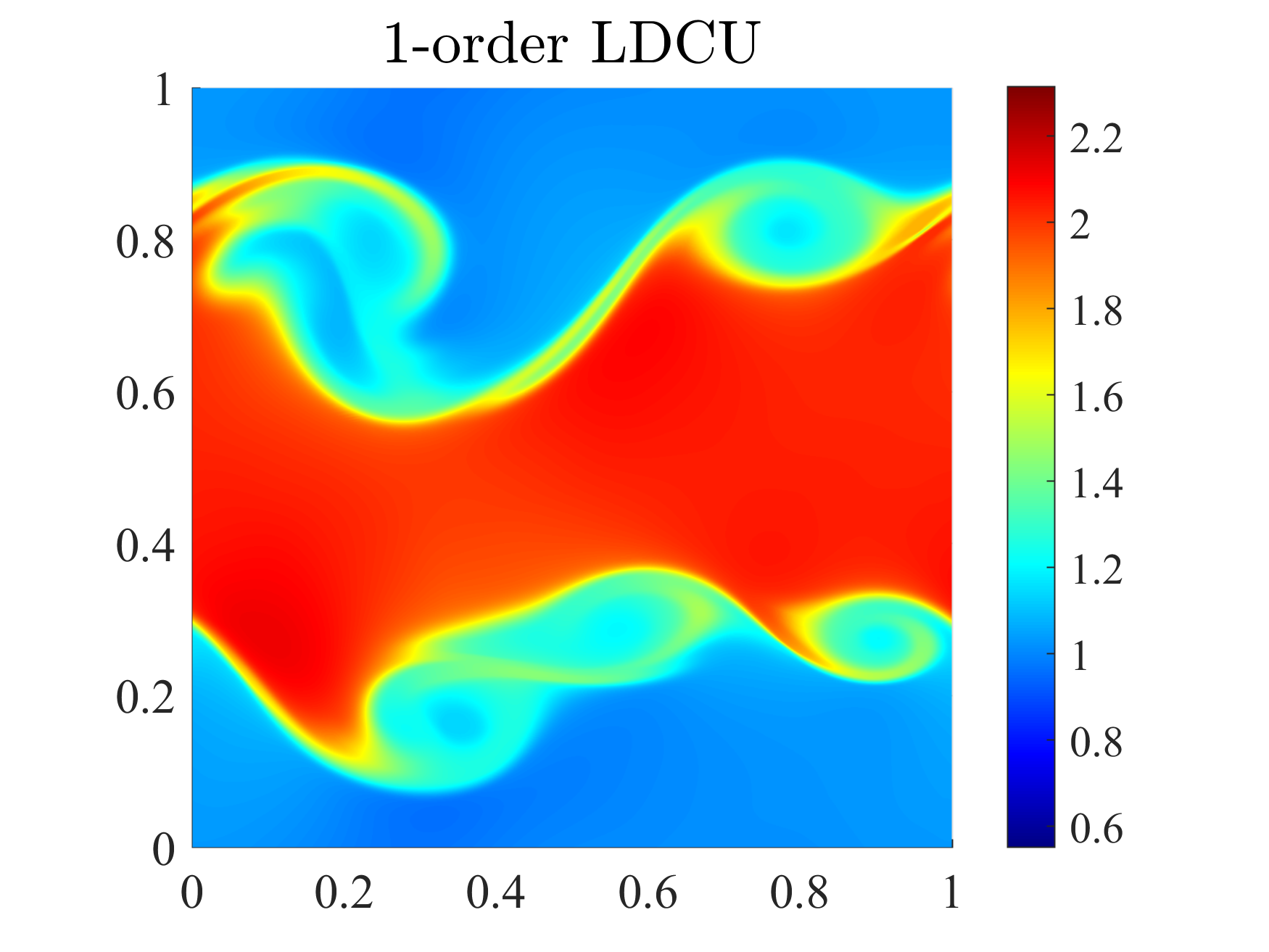}\hspace*{0.0cm}
            \includegraphics[trim=1.3cm 0.5cm 1.6cm 0.2cm, clip, width=4.2cm]{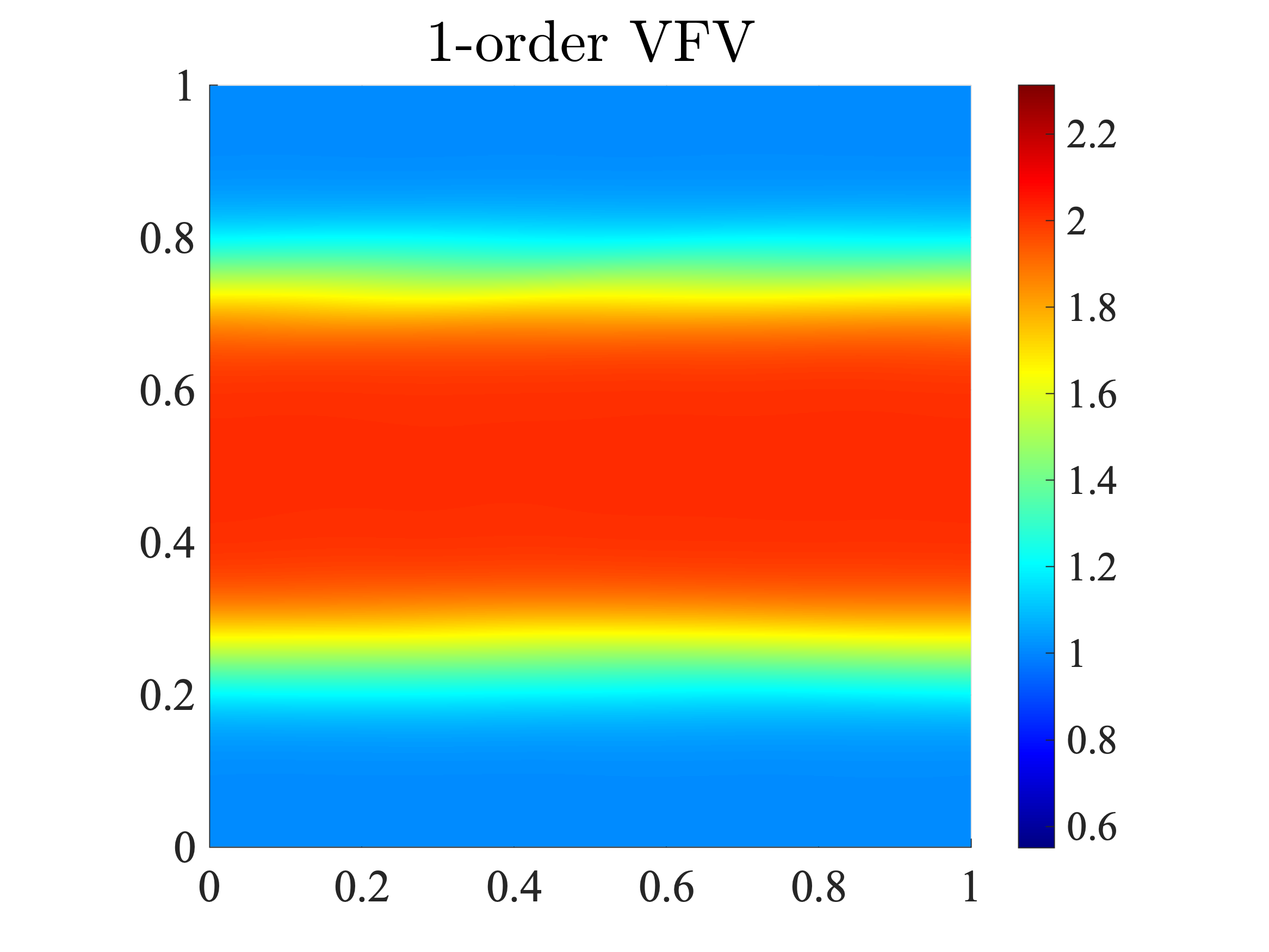}}
\vskip8pt
\centerline{\includegraphics[trim=1.3cm 0.5cm 1.6cm 0.2cm, clip, width=4.cm]{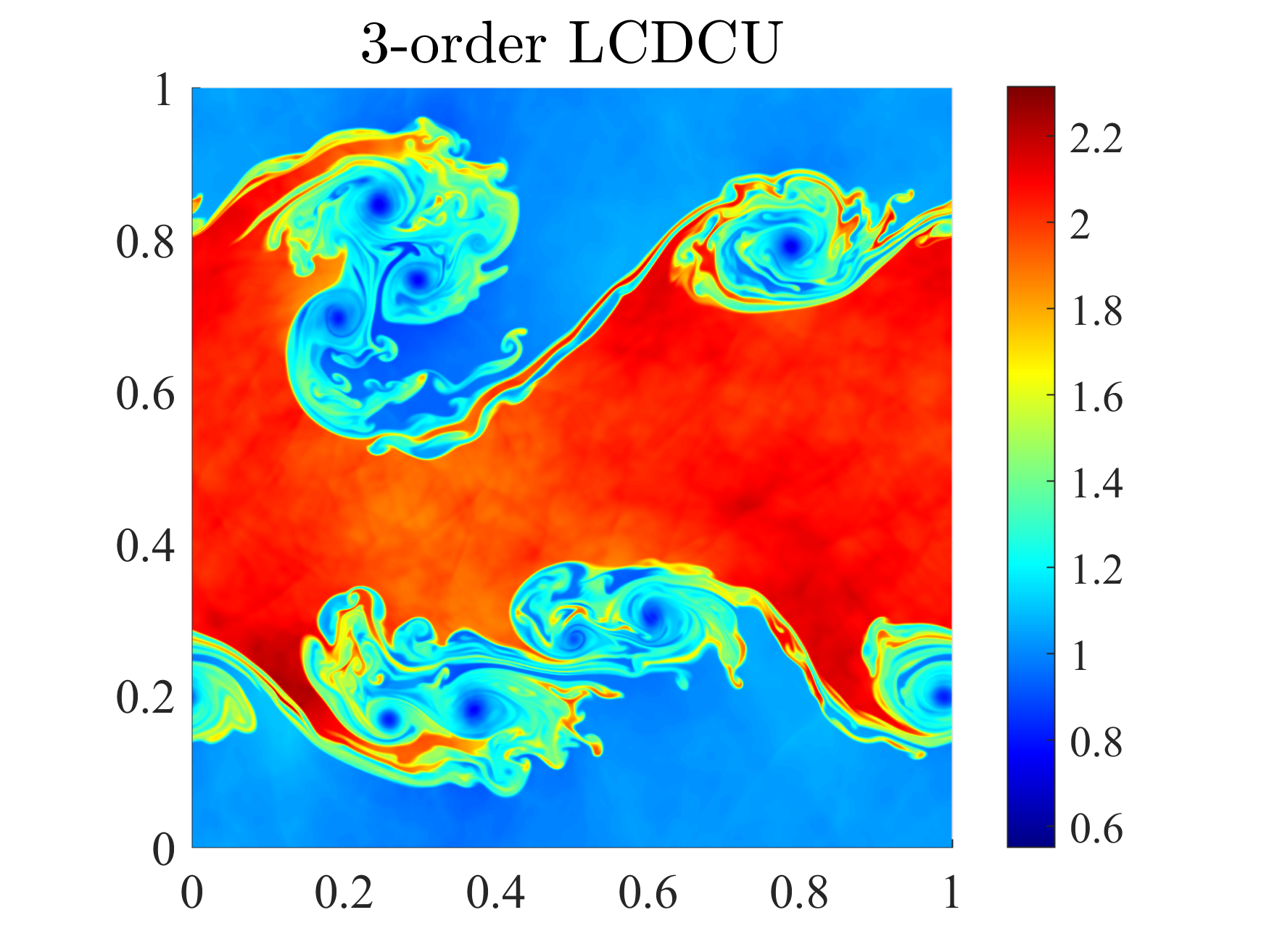}\hspace*{0.2cm}
            \includegraphics[trim=1.3cm 0.5cm 1.6cm 0.2cm, clip, width=4.cm]{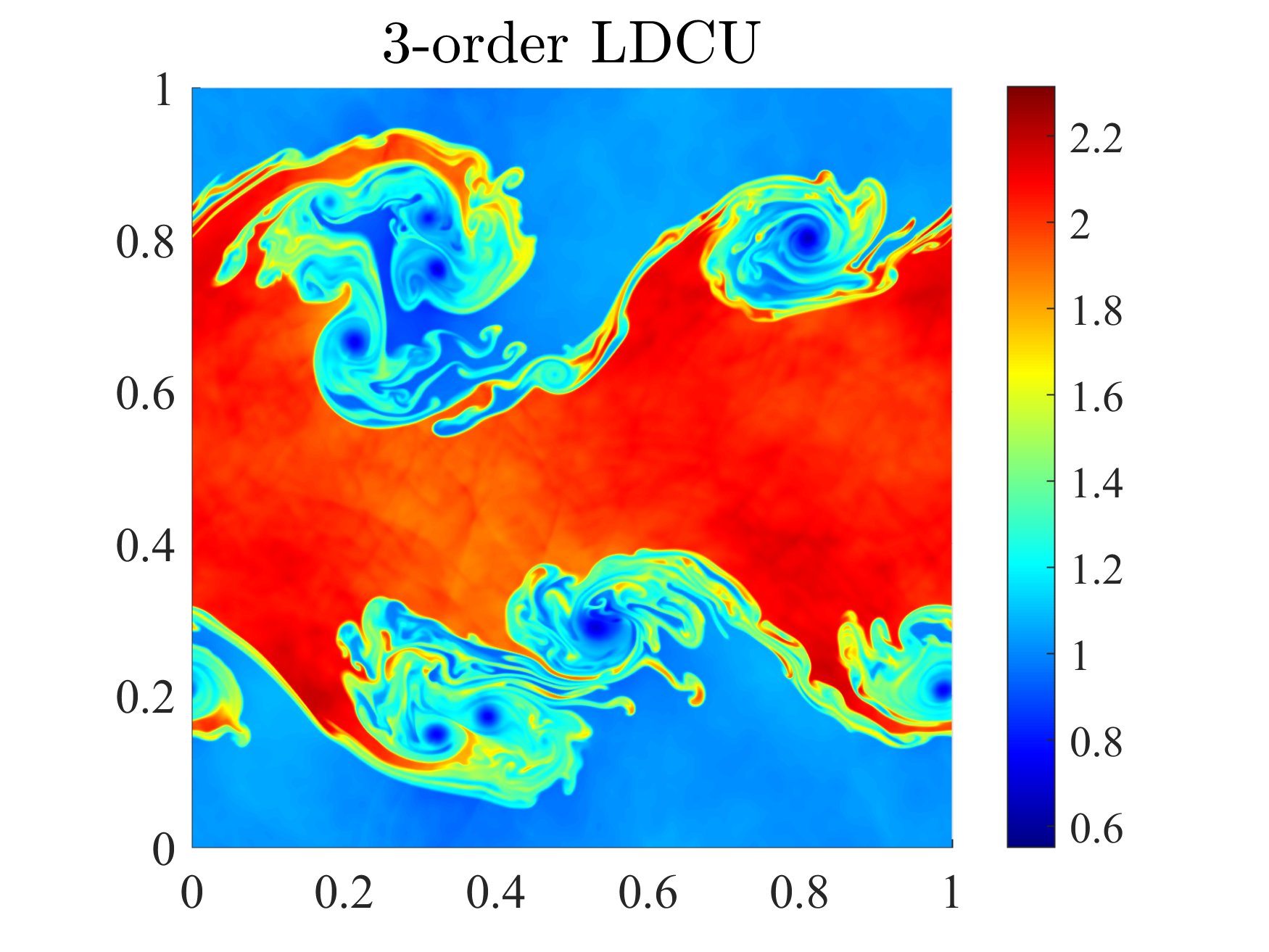}\hspace*{0.0cm}
            \includegraphics[trim=1.3cm 0.5cm 1.6cm 0.2cm, clip, width=4.2cm]{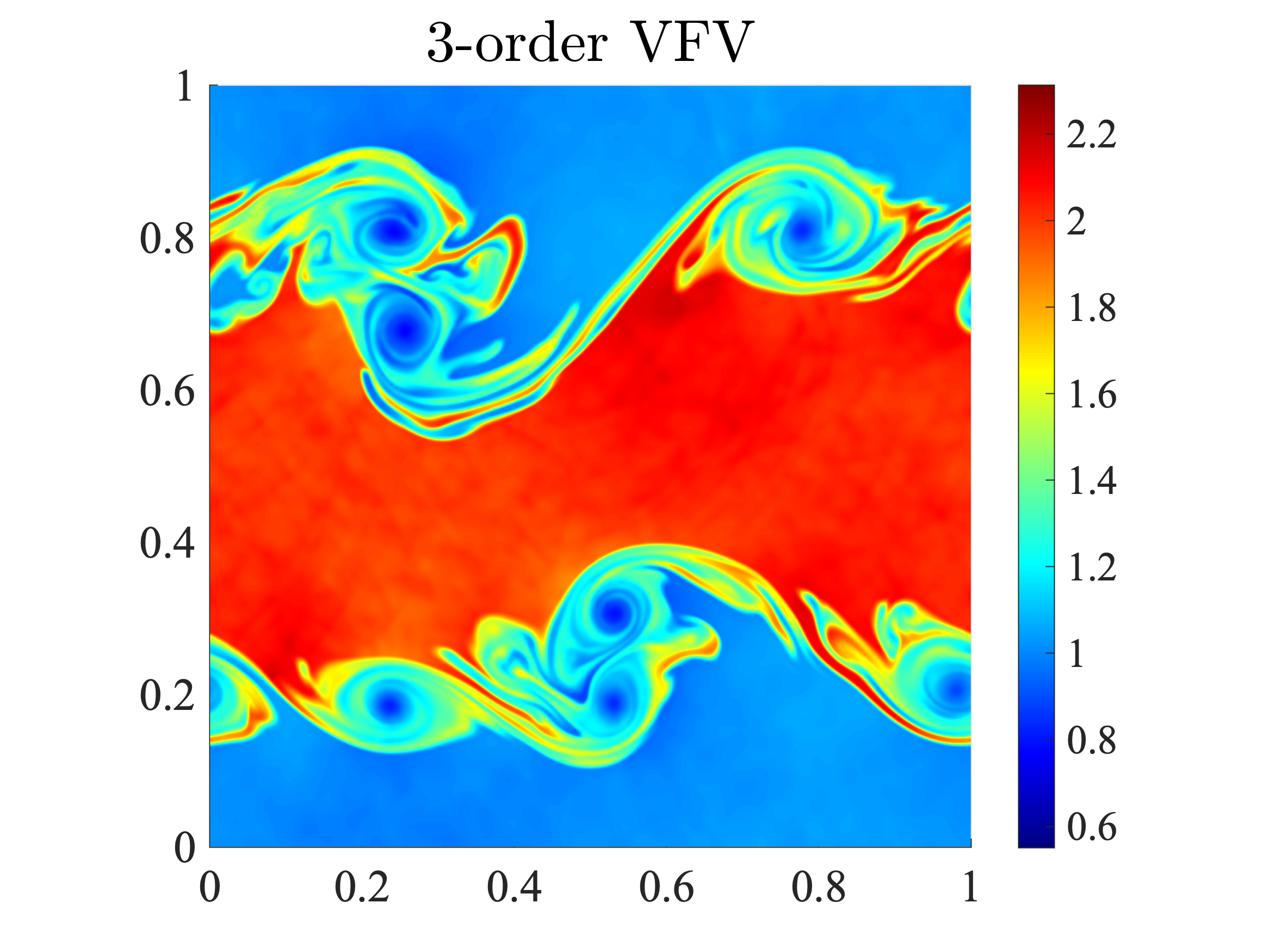}}
\vskip8pt
\centerline{\includegraphics[trim=1.3cm 0.5cm 1.6cm 0.2cm, clip, width=4.cm]{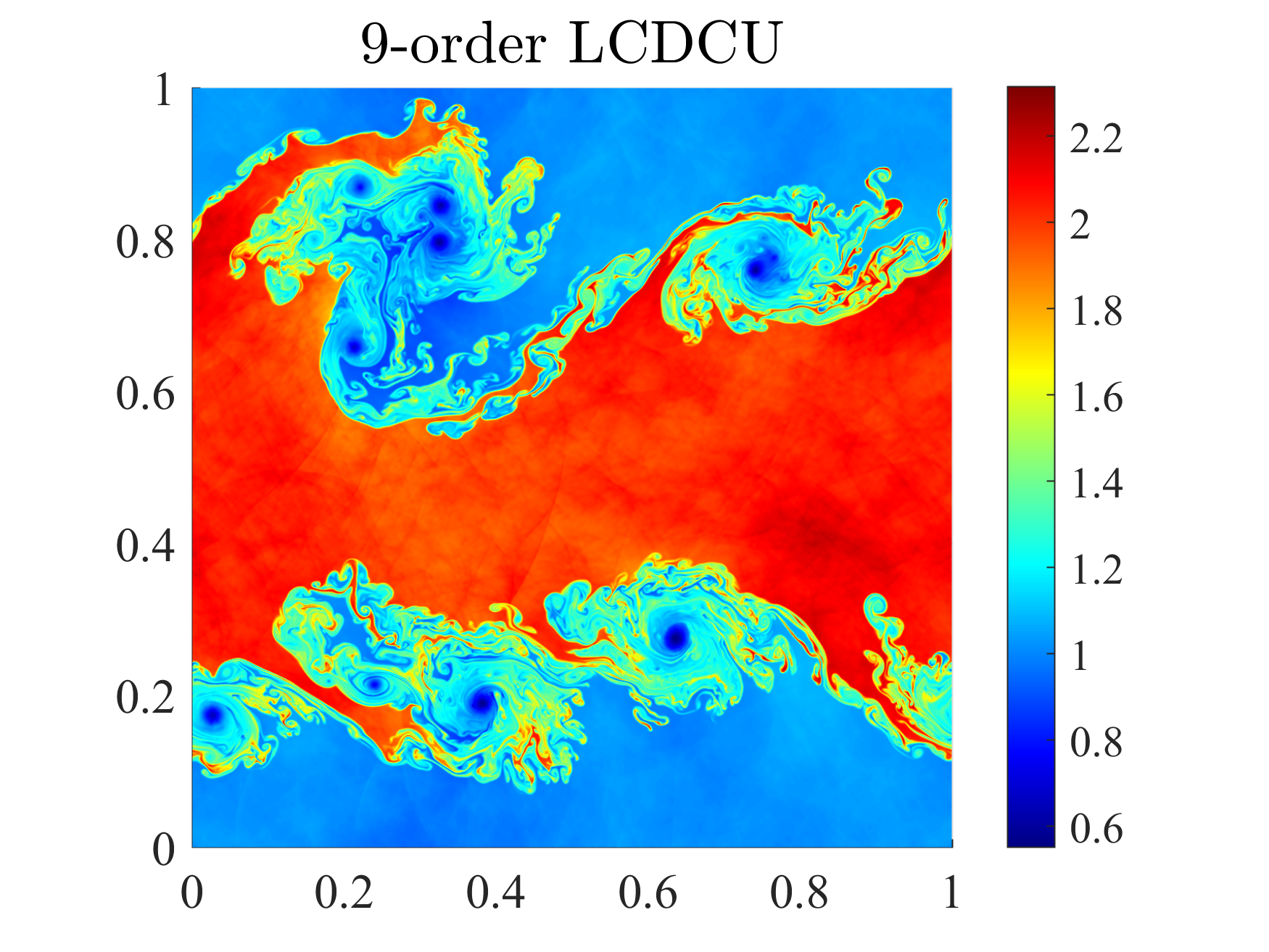}\hspace*{0.2cm}
            \includegraphics[trim=1.3cm 0.5cm 1.6cm 0.2cm, clip, width=4.cm]{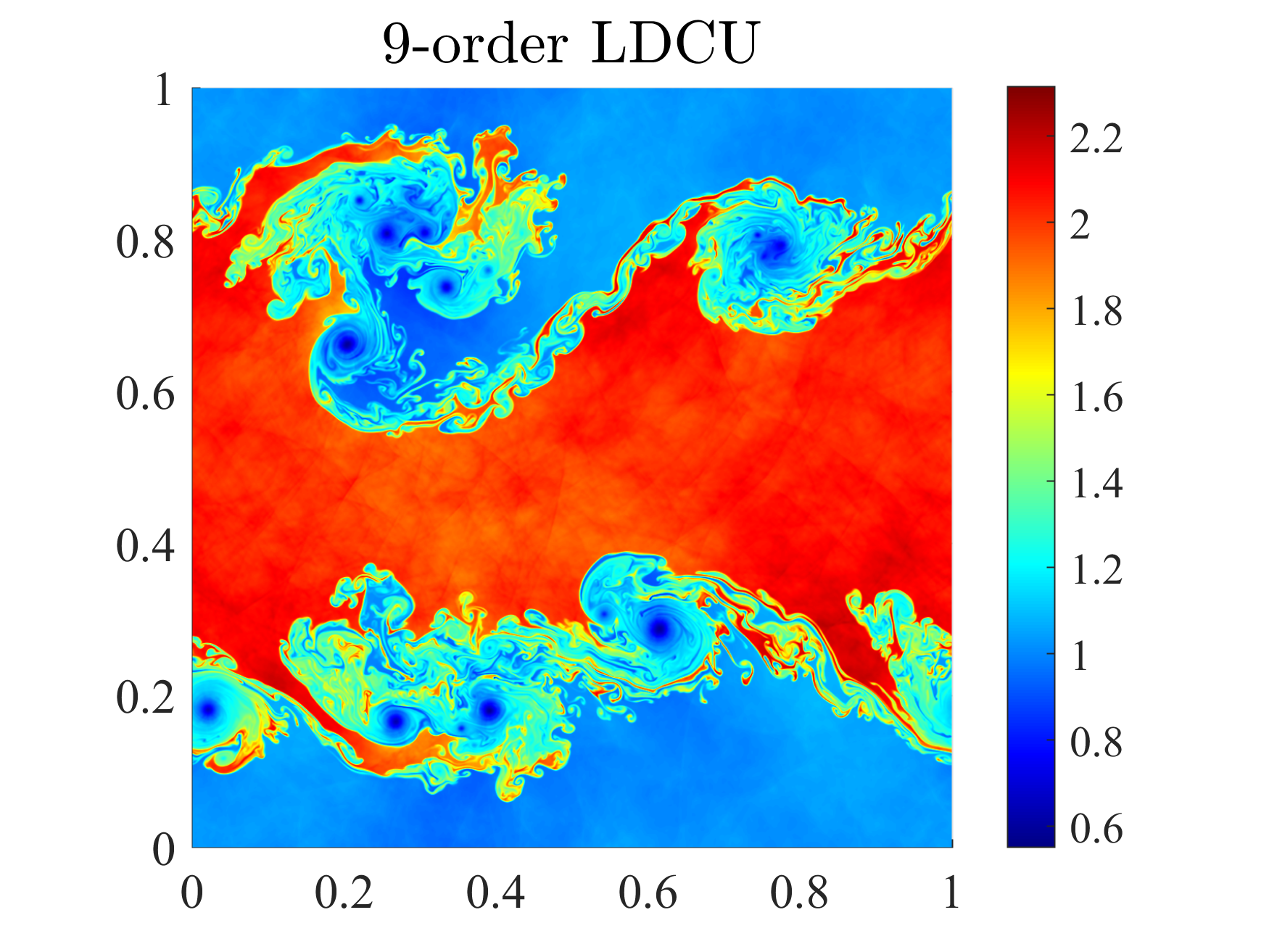}\hspace*{0.0cm}
            \includegraphics[trim=1.3cm 0.5cm 1.6cm 0.2cm, clip, width=4.2cm]{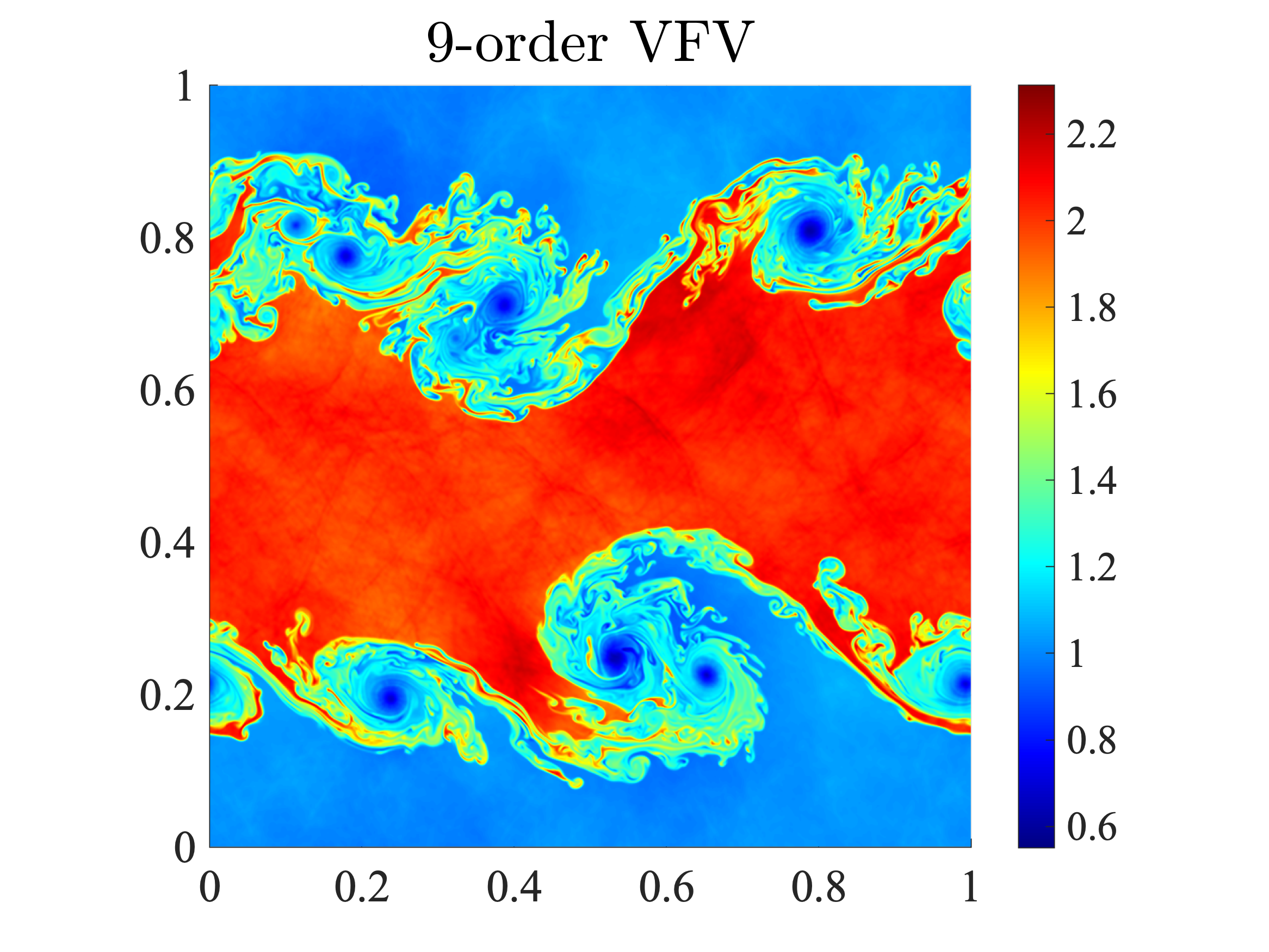}}
\caption{\sf KH Instability: Density computed by the first- (top row), third- (middle), and ninth-order (bottom row) LCDCU (left column),
LDCU (middle column), and VFV (right column) schemes.\label{fig1}}
\end{figure}

We compute the numerical solutions until the final time $T=2$ using the studied schemes and show the results, obtained by the first-, third-, and ninth-order schemes in Figure \ref{fig1}. The results indicate that the flow is unstable, with an increasing number of vortices generated when higher-order schemes are employed. Additionally, it is evident that different schemes converge to different DW solutions. This suggests that the numerical solutions will not strongly converge as the order of accuracy increases.

\begin{figure}[ht!]
\centerline{\includegraphics[trim=1.3cm 0.5cm 1.6cm 0.2cm, clip, width=4.cm]{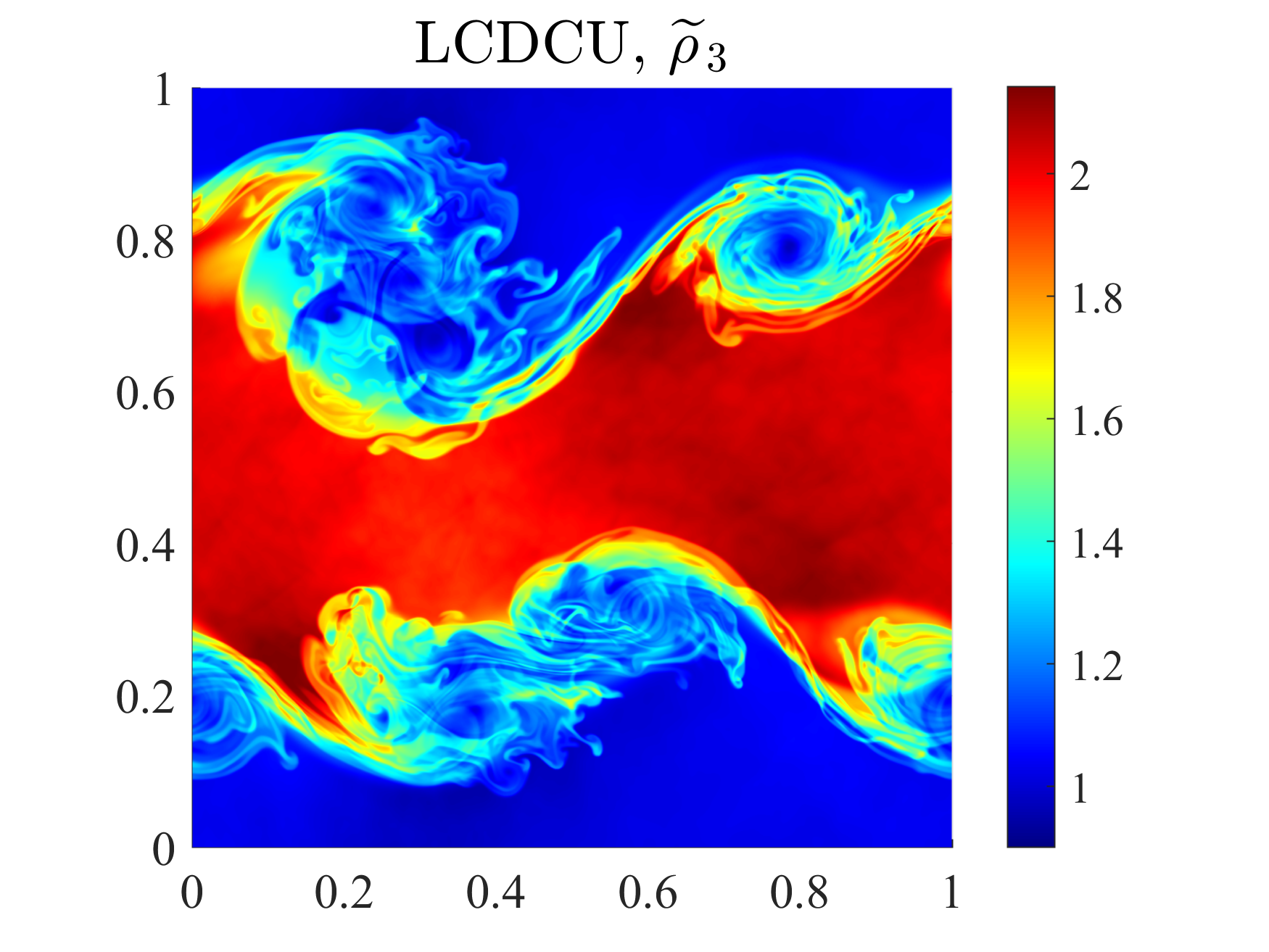}\hspace*{0.2cm}
            \includegraphics[trim=1.3cm 0.5cm 1.6cm 0.2cm, clip, width=4.cm]{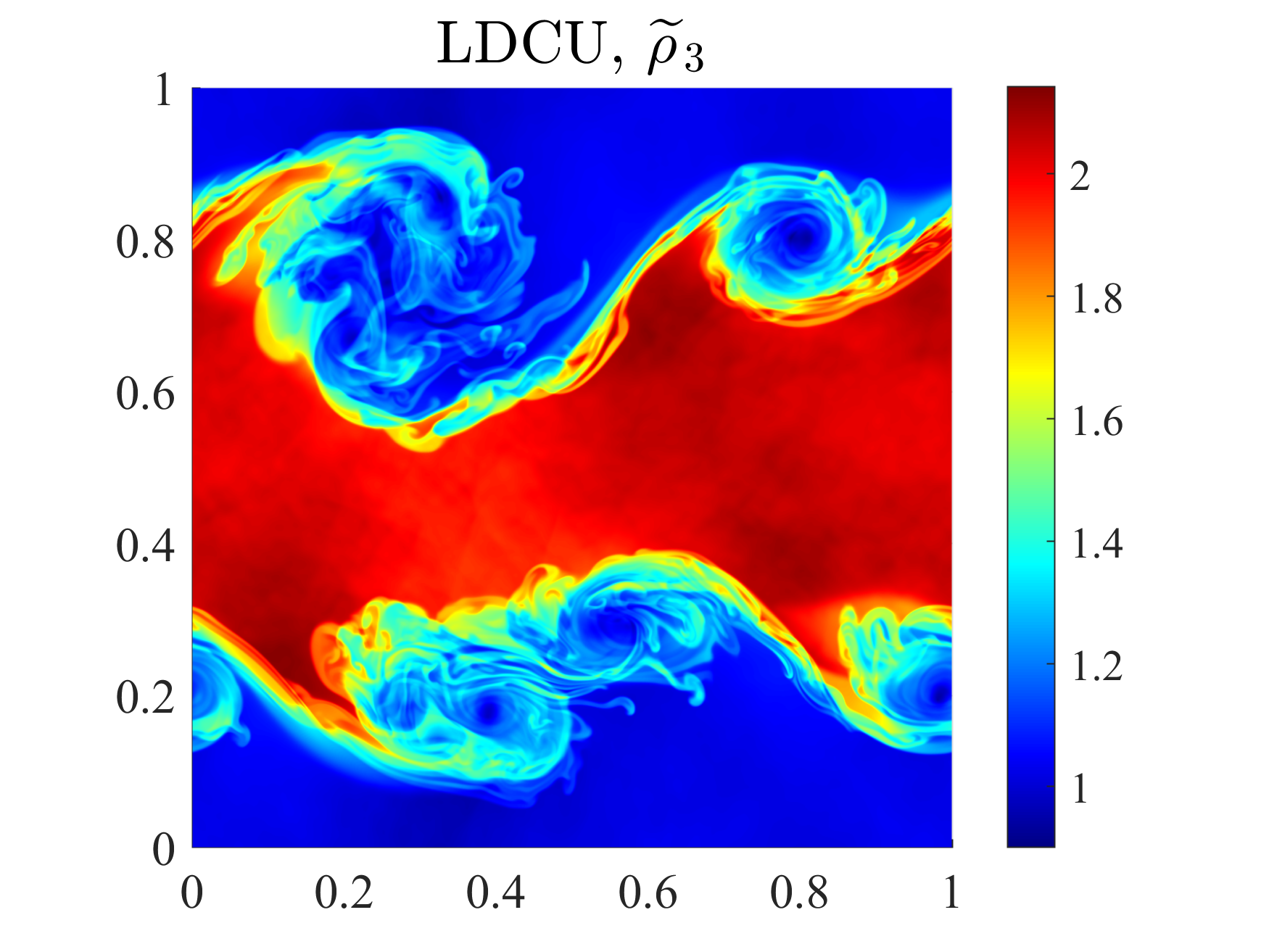}\hspace*{0.0cm}
            \includegraphics[trim=1.3cm 0.5cm 1.6cm 0.2cm, clip, width=4.2cm]{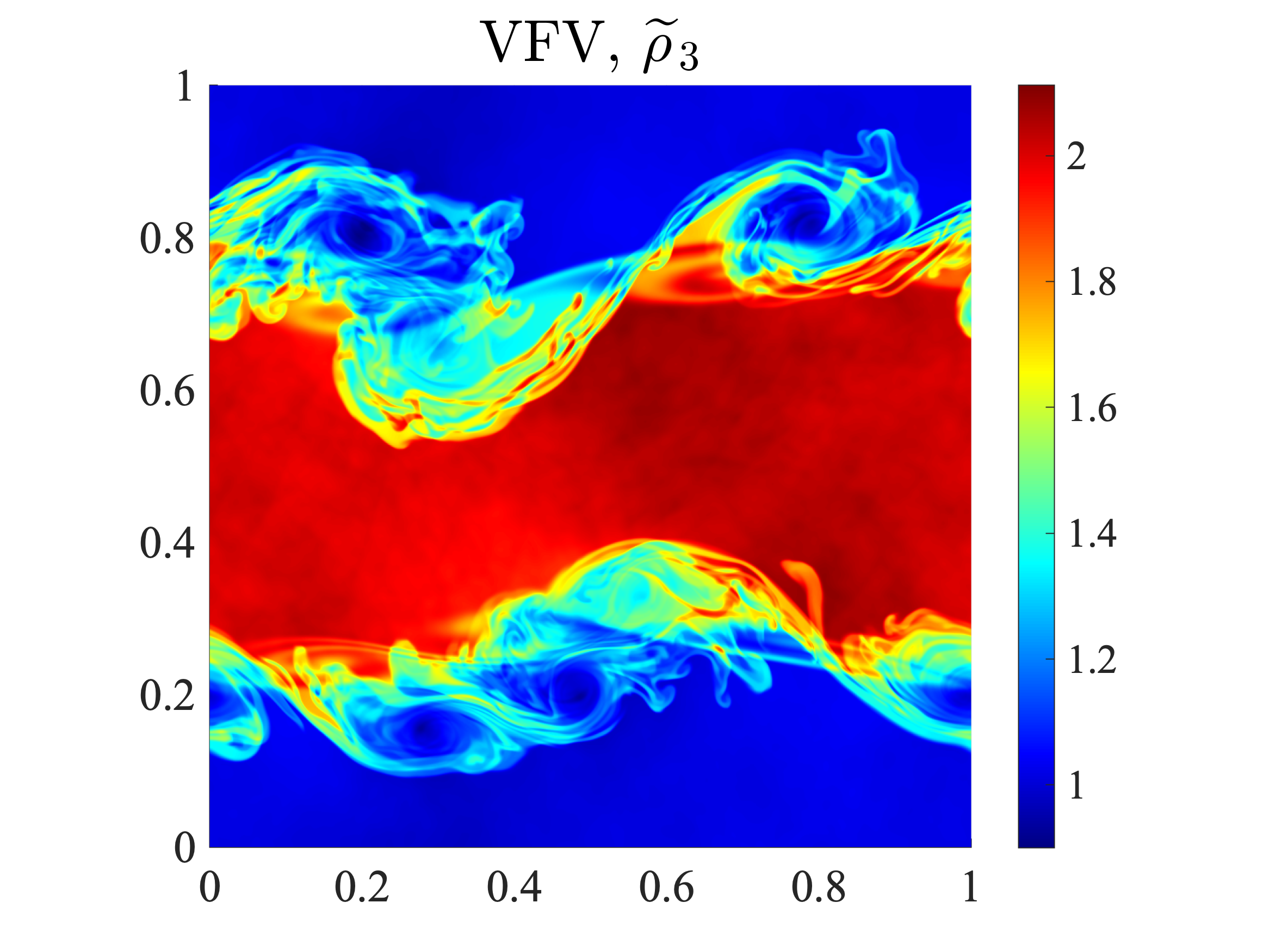}}
\vskip8pt
\centerline{\includegraphics[trim=1.3cm 0.5cm 1.6cm 0.2cm, clip, width=4.cm]{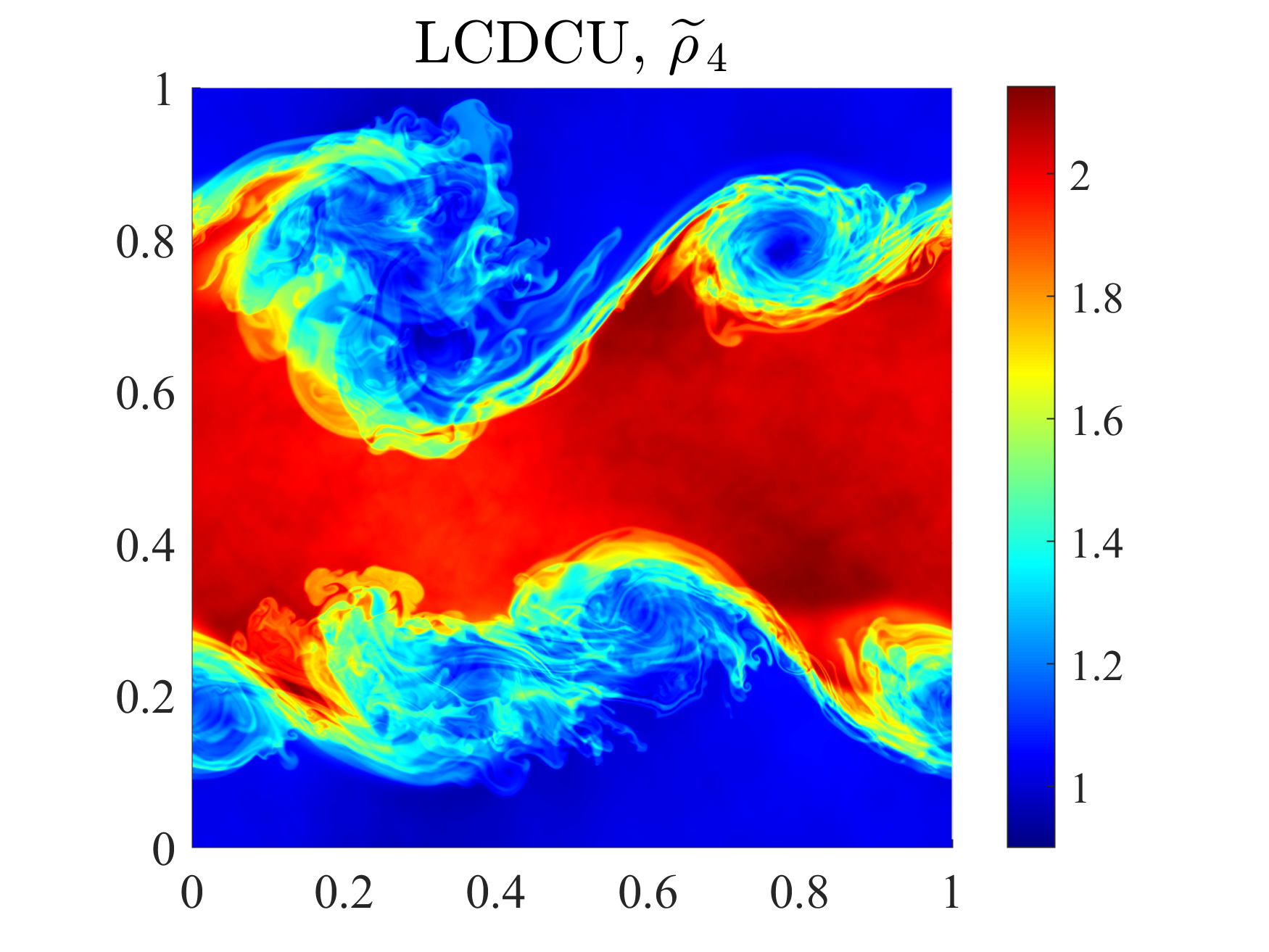}\hspace*{0.2cm}
            \includegraphics[trim=1.3cm 0.5cm 1.6cm 0.2cm, clip, width=4.cm]{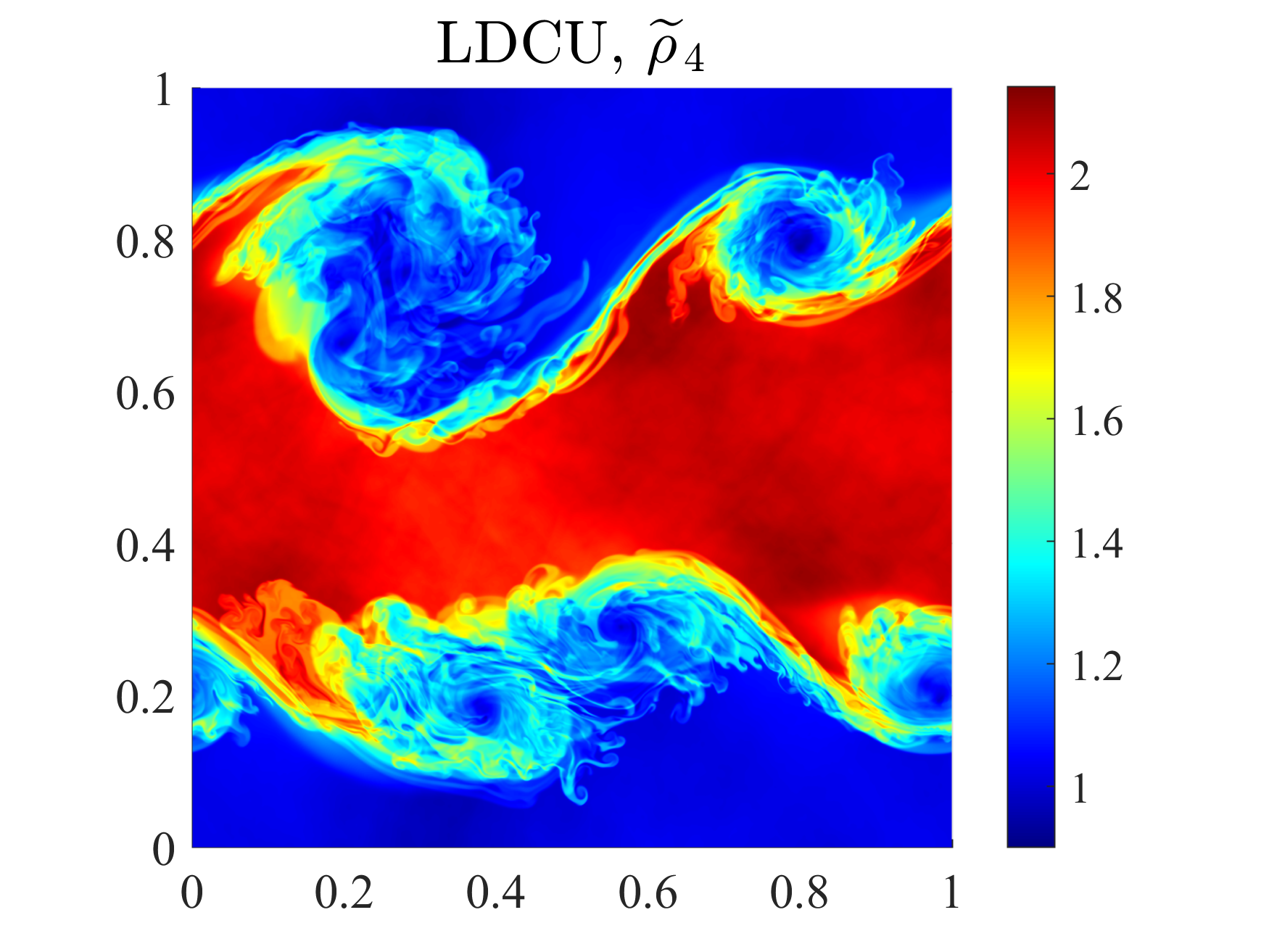}\hspace*{0.0cm}
            \includegraphics[trim=1.3cm 0.5cm 1.6cm 0.2cm, clip, width=4.2cm]{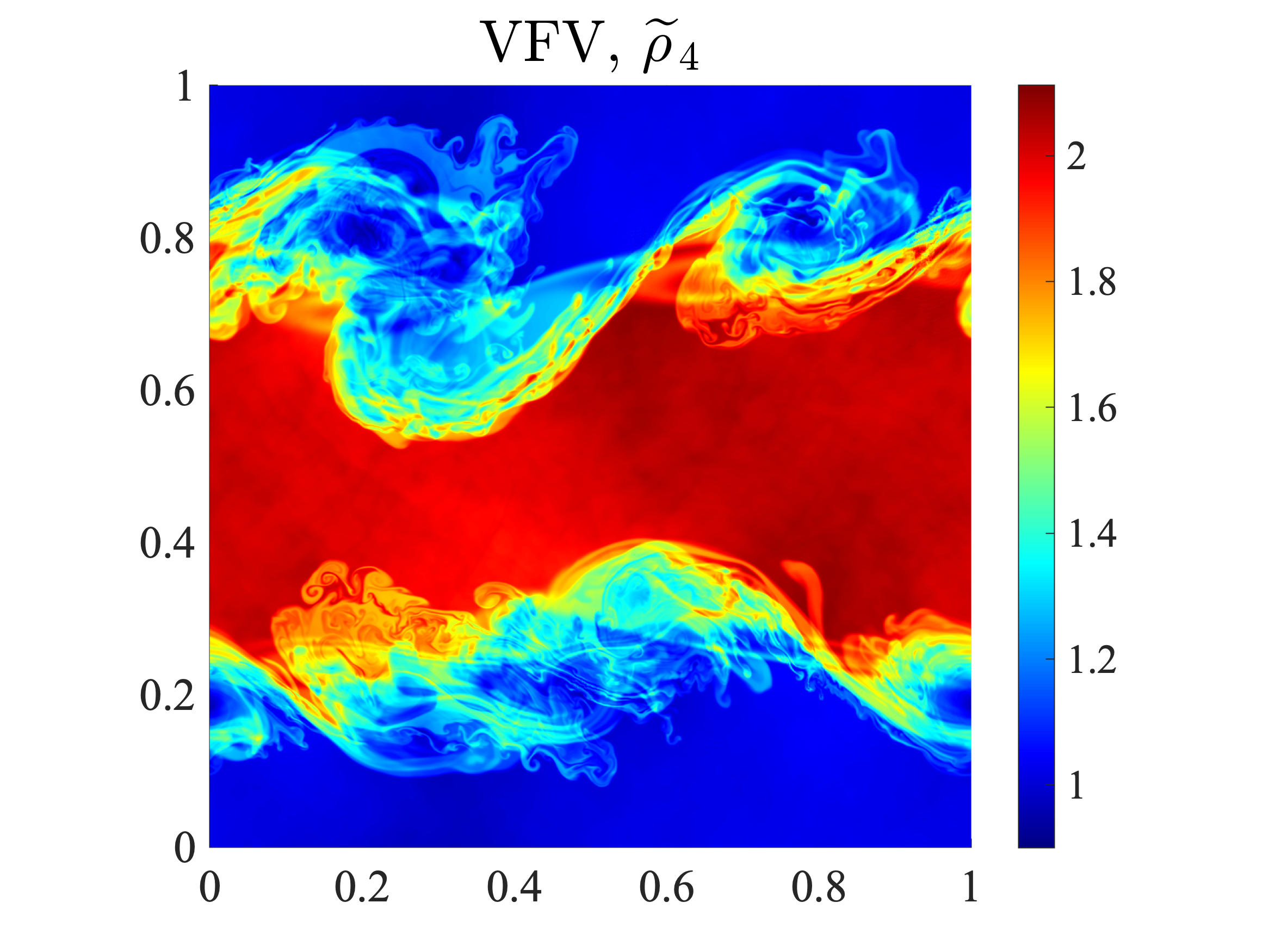}}
\vskip8pt
\centerline{\includegraphics[trim=1.3cm 0.5cm 1.6cm 0.2cm, clip, width=4.cm]{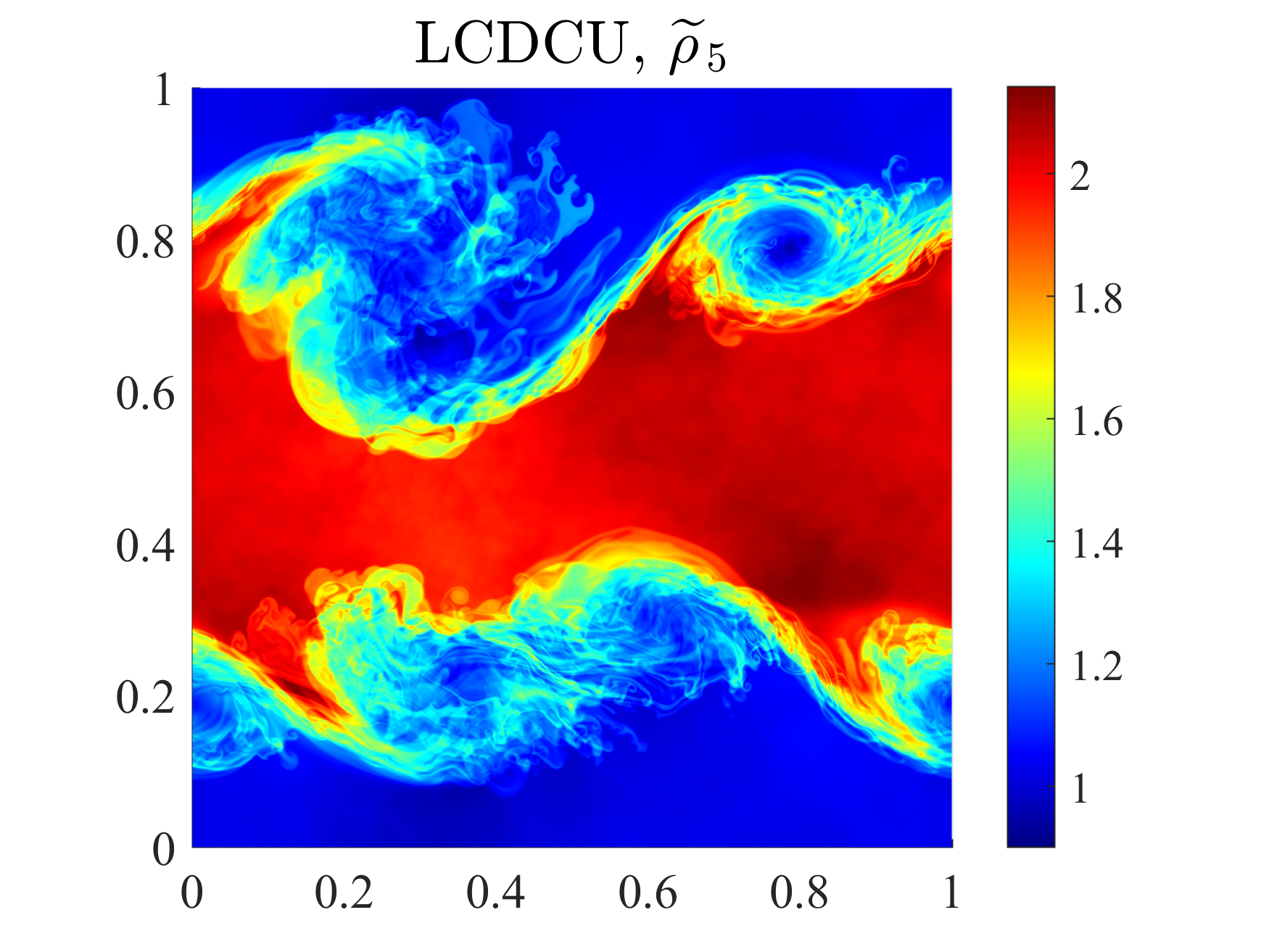}\hspace*{0.2cm}
            \includegraphics[trim=1.3cm 0.5cm 1.6cm 0.2cm, clip, width=4.cm]{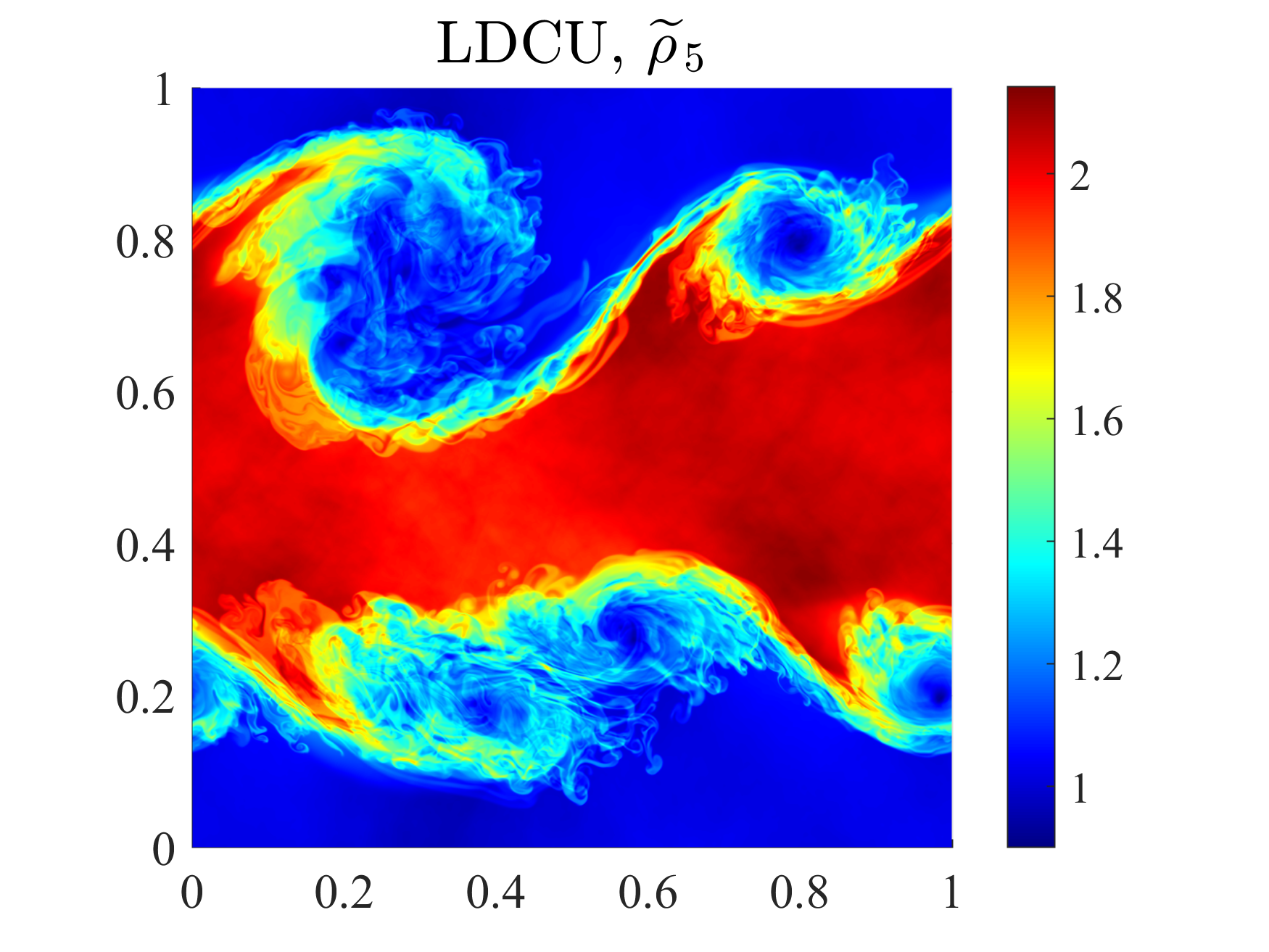}\hspace*{0.0cm}
            \includegraphics[trim=1.3cm 0.5cm 1.6cm 0.2cm, clip, width=4.2cm]{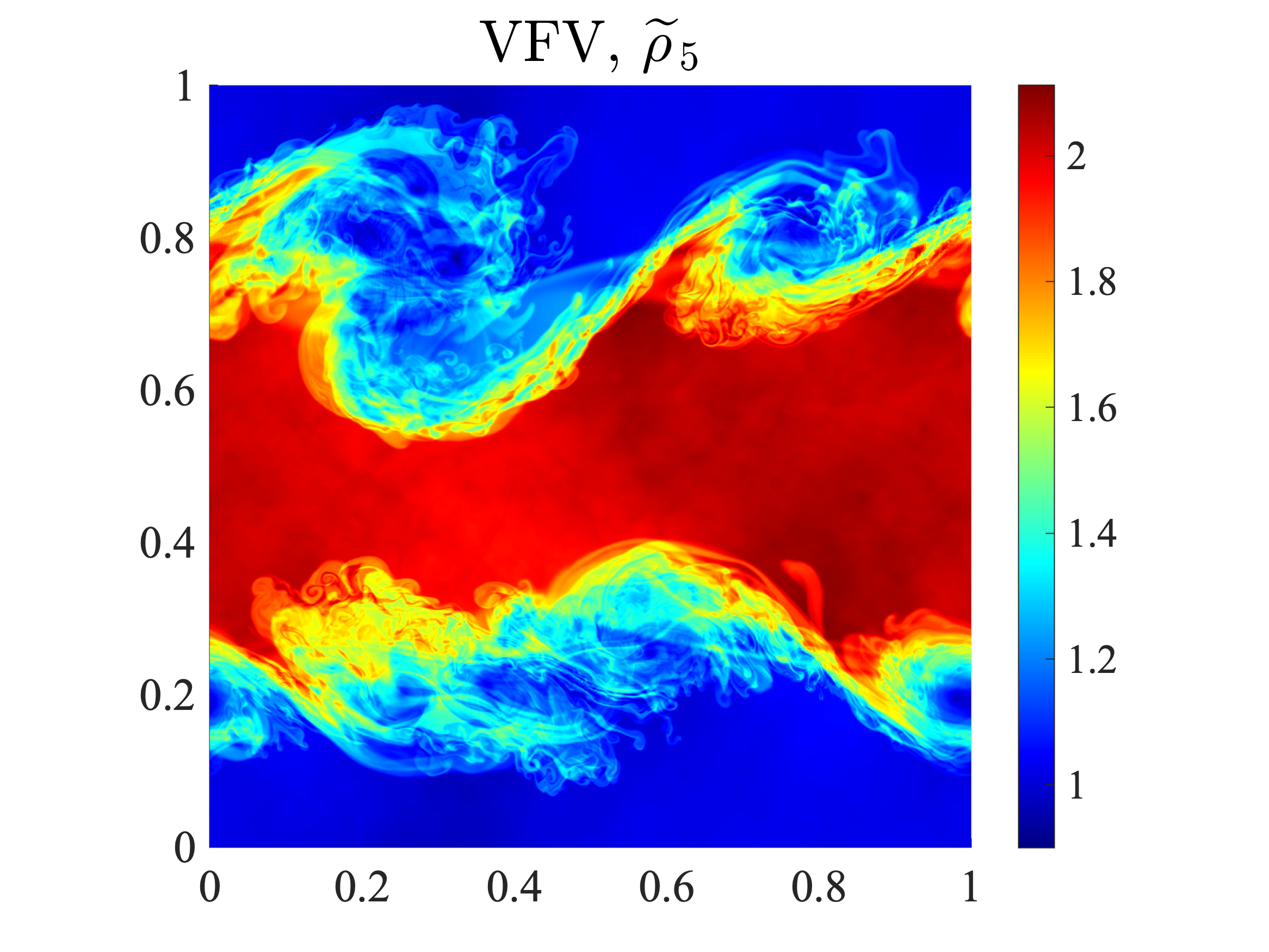}}
\vskip8pt
\centerline{\includegraphics[trim=1.3cm 0.5cm 1.6cm 0.2cm, clip, width=4.cm]{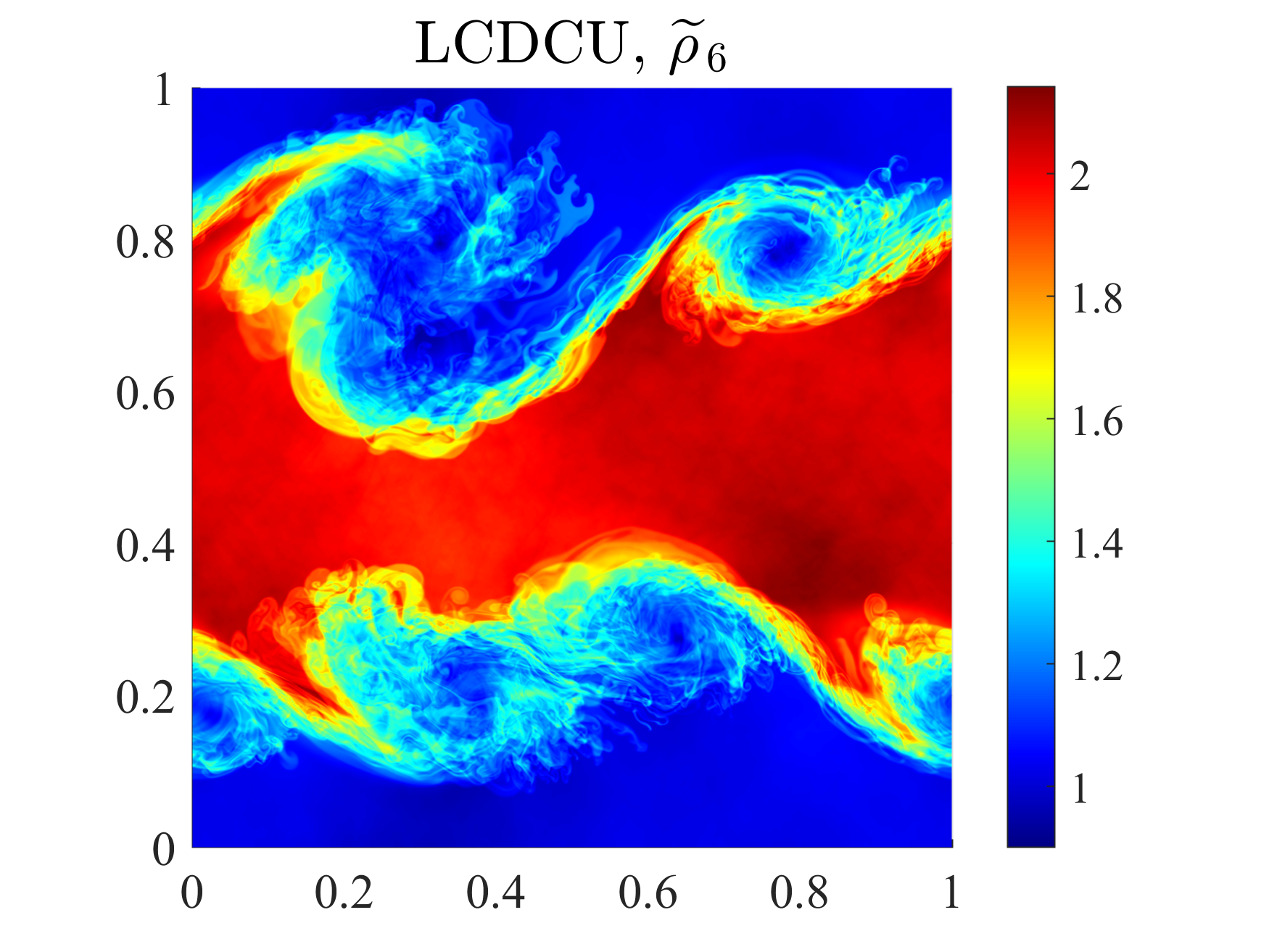}\hspace*{0.2cm}
            \includegraphics[trim=1.3cm 0.5cm 1.6cm 0.2cm, clip, width=4.cm]{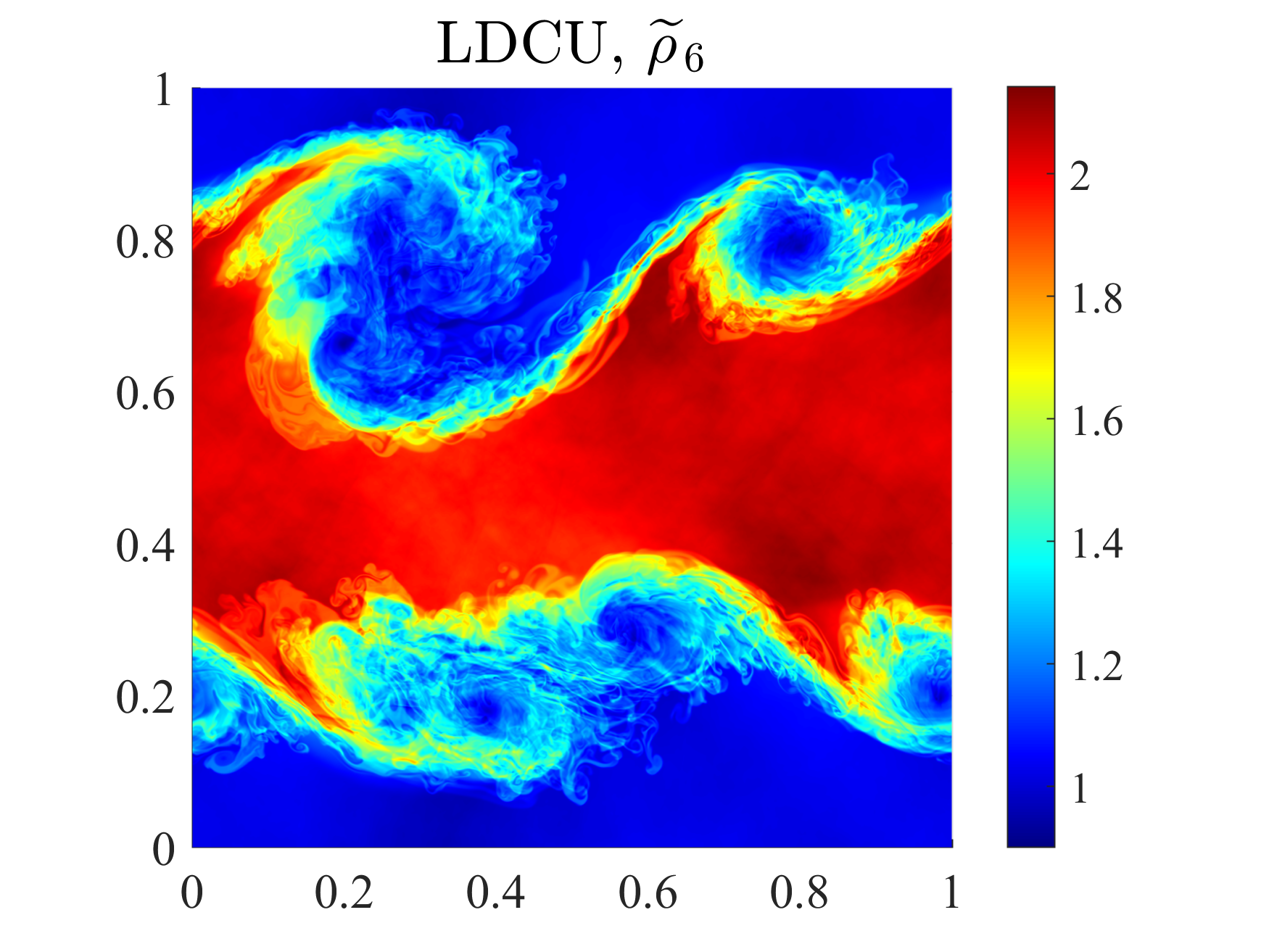}\hspace*{0.0cm}
            \includegraphics[trim=1.3cm 0.5cm 1.6cm 0.2cm, clip, width=4.2cm]{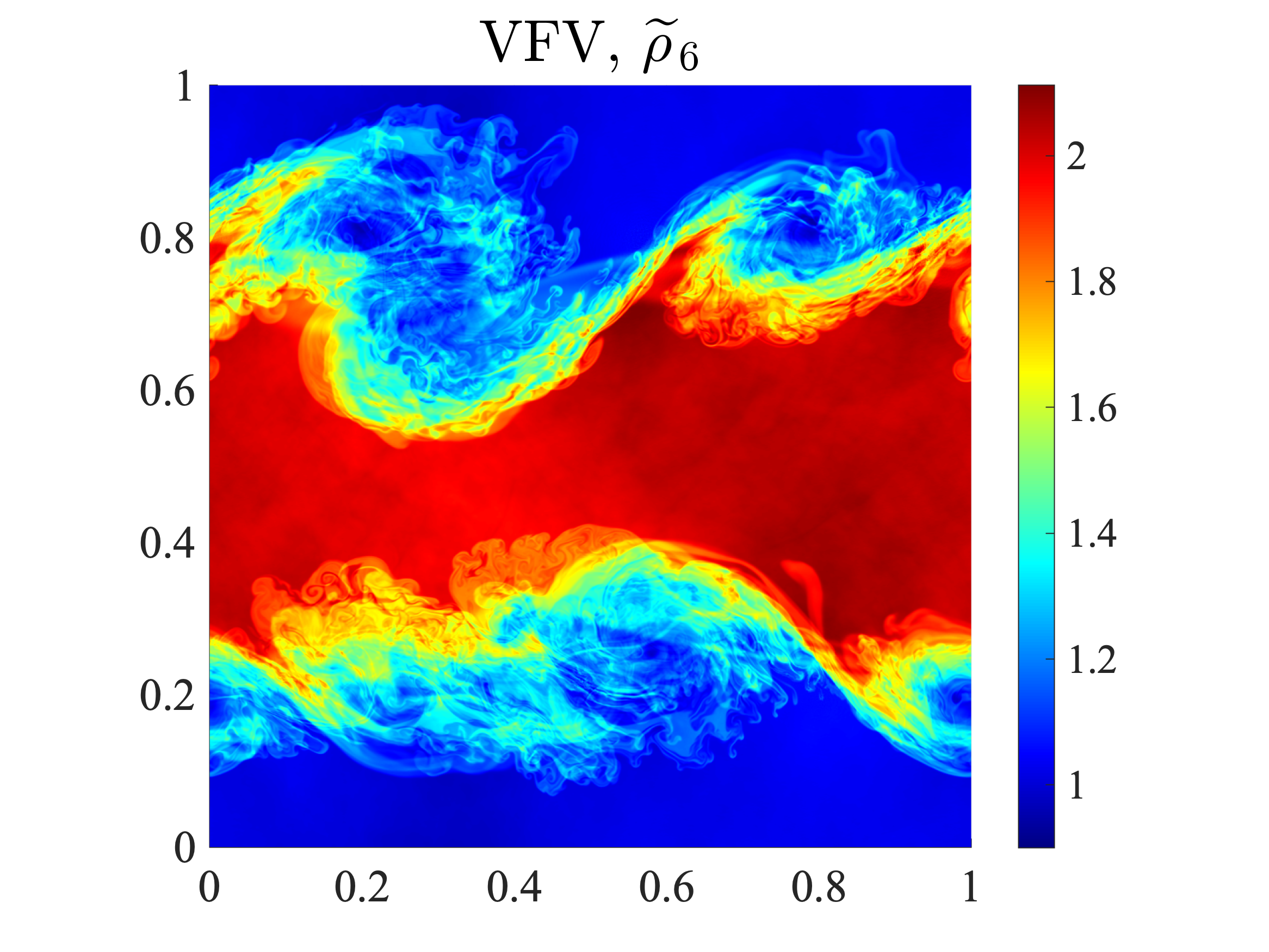}}
\caption{\sf KH Instability: $\widetilde\rho_3$ (top row), $\widetilde\rho_4$ (second row), $\widetilde\rho_5$ (third row), and $\widetilde\rho_6$ (bottom row) computed by the LCDCU (left column), LDCU (middle column), and VFV (right column) schemes.\label{fig3}}
\end{figure}

We then consider the ${\mathcal K}$-convergence of the numerical solutions over different orders. In Figure \ref{fig3}, we present $\widetilde\rho_3$, $\widetilde\rho_4$, $\widetilde\rho_5$, and $\widetilde\rho_6$, where one can see that, in contrast to Figure \ref{fig1}, the average density profiles show similarities for different $n$, suggesting that the sequence of averages is convergent. Additionally, we observe that for different methods, the limits of the average density differ significantly. This indicates that the numerical approximations generated by different methods converge to different DW solutions.

\begin{figure}[ht!]
\centerline{\includegraphics[trim=1.3cm 0.5cm 1.6cm 0.2cm, clip, width=4.cm]{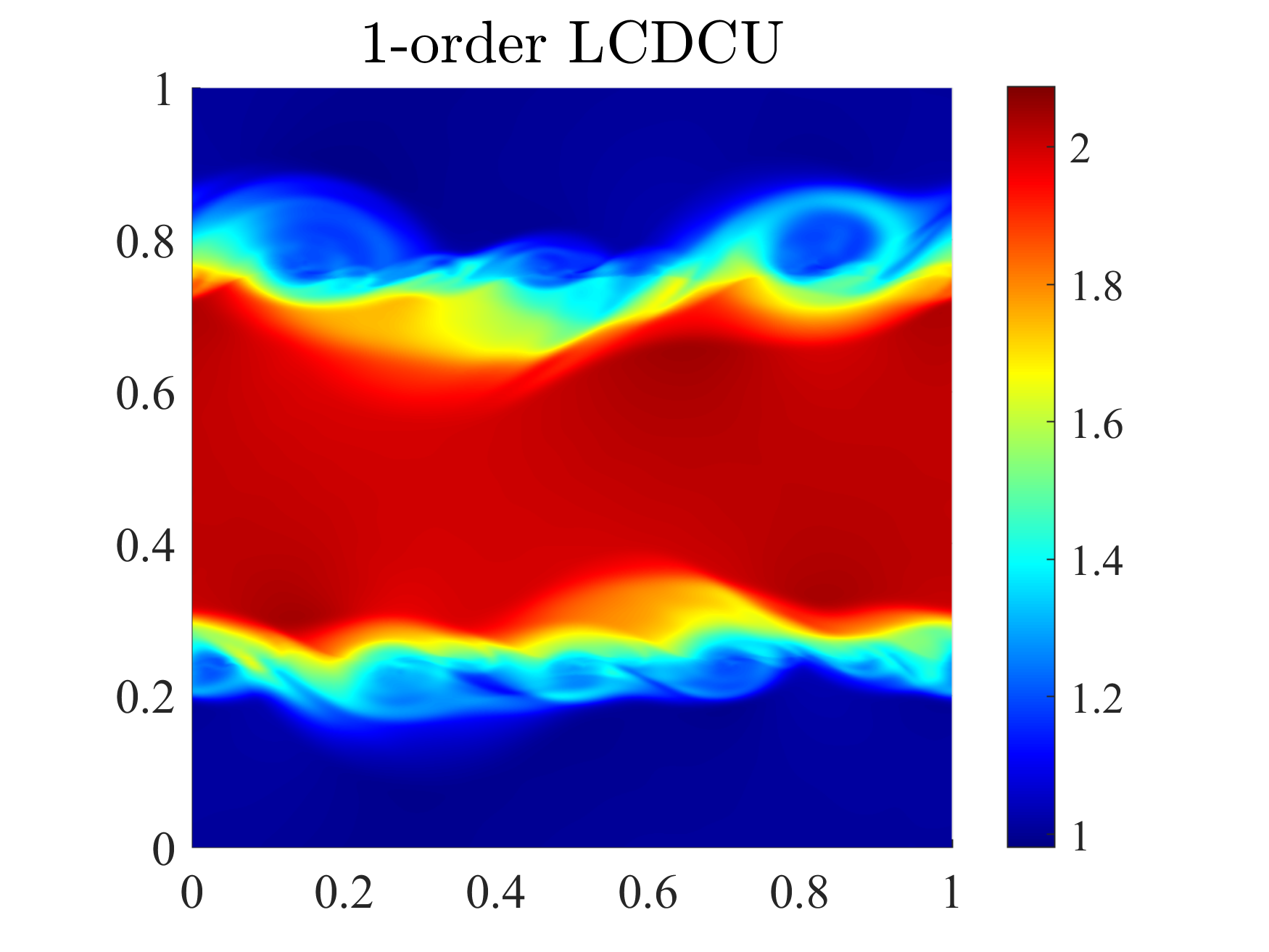}\hspace*{0.2cm}
            \includegraphics[trim=1.3cm 0.5cm 1.6cm 0.2cm, clip, width=4.cm]{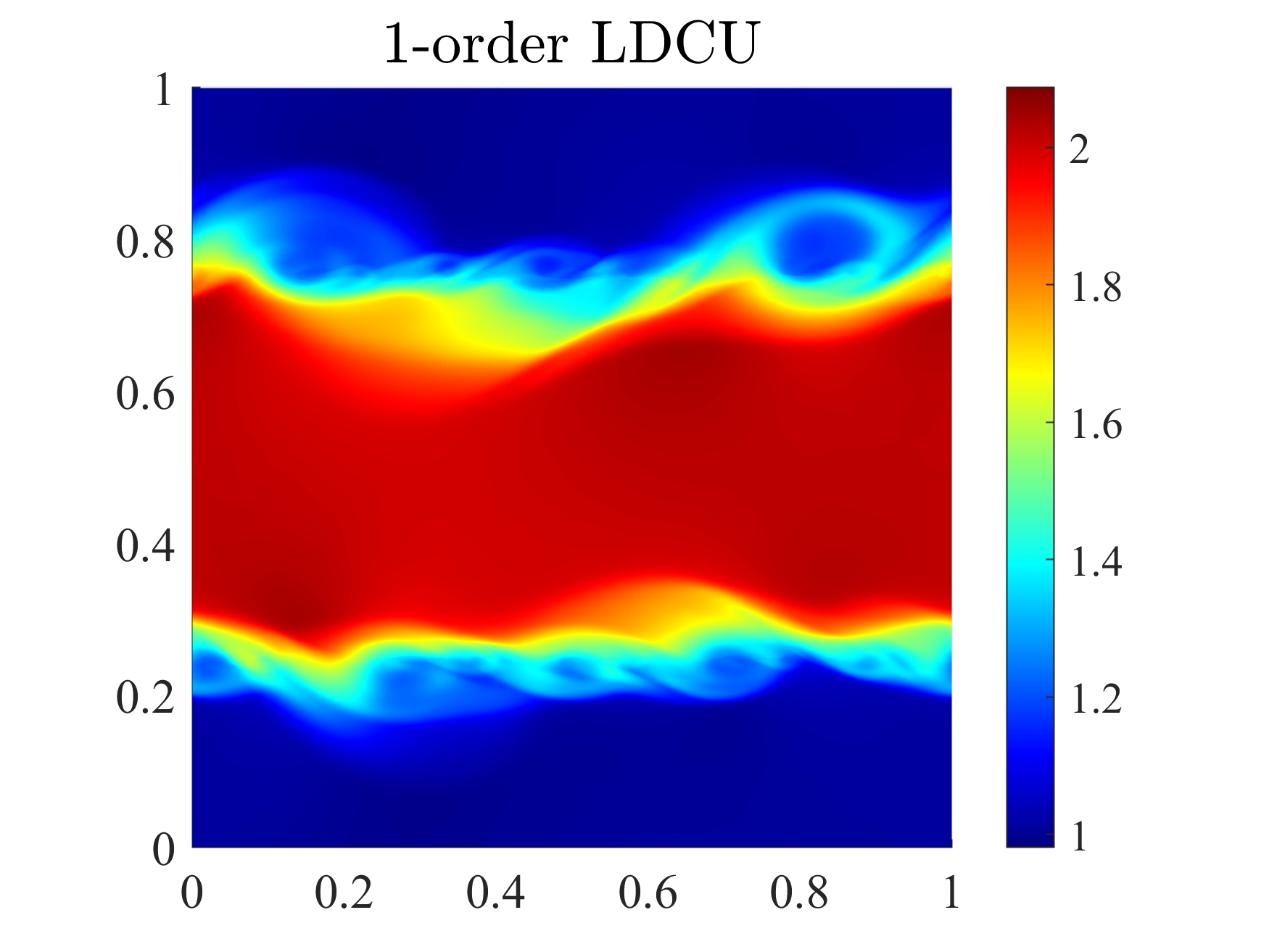}\hspace*{0.0cm}
            \includegraphics[trim=1.3cm 0.5cm 1.6cm 0.2cm, clip, width=4.2cm]{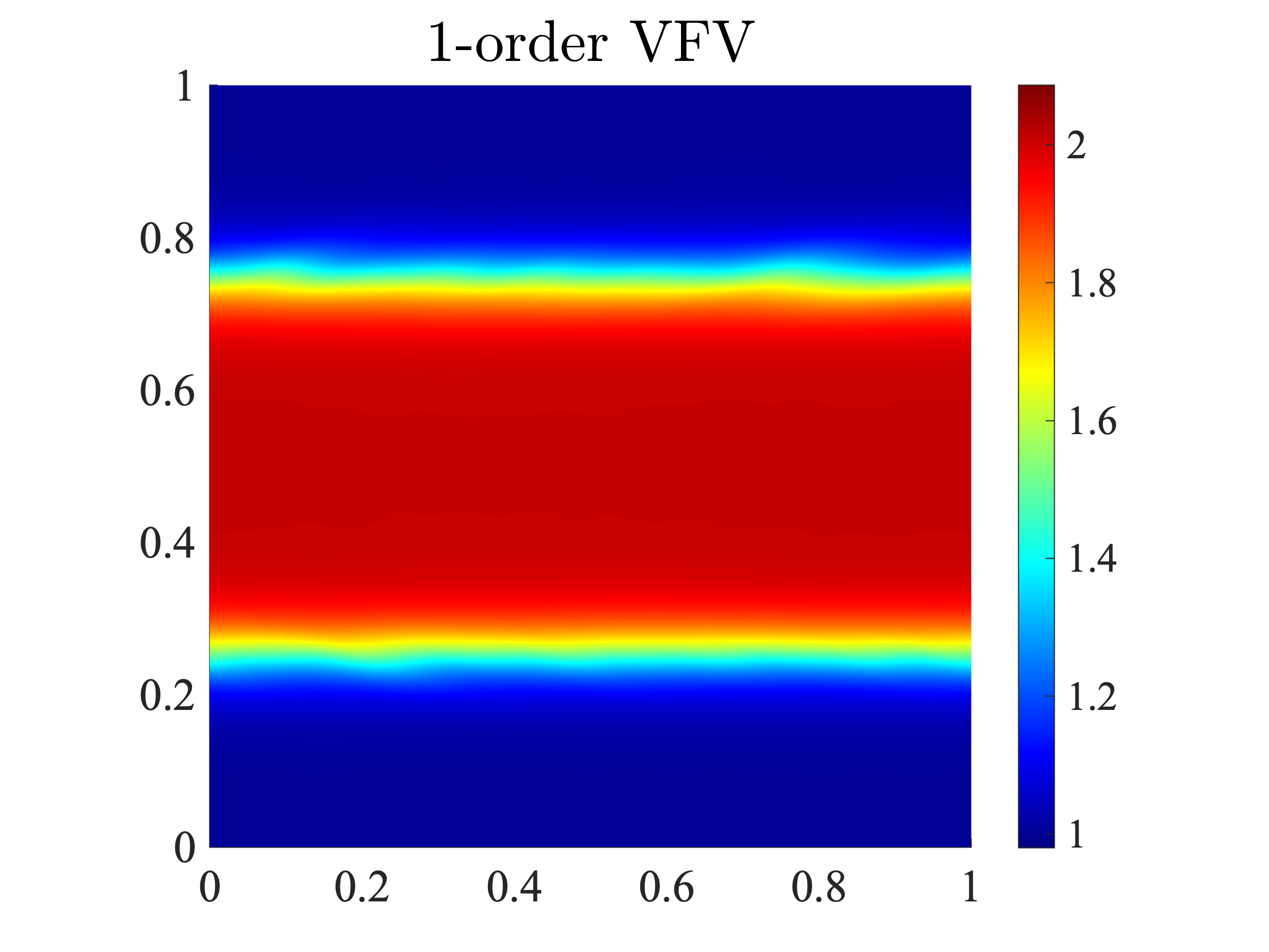}}
\vskip8pt
\centerline{\includegraphics[trim=1.3cm 0.5cm 1.6cm 0.2cm, clip, width=4.cm]{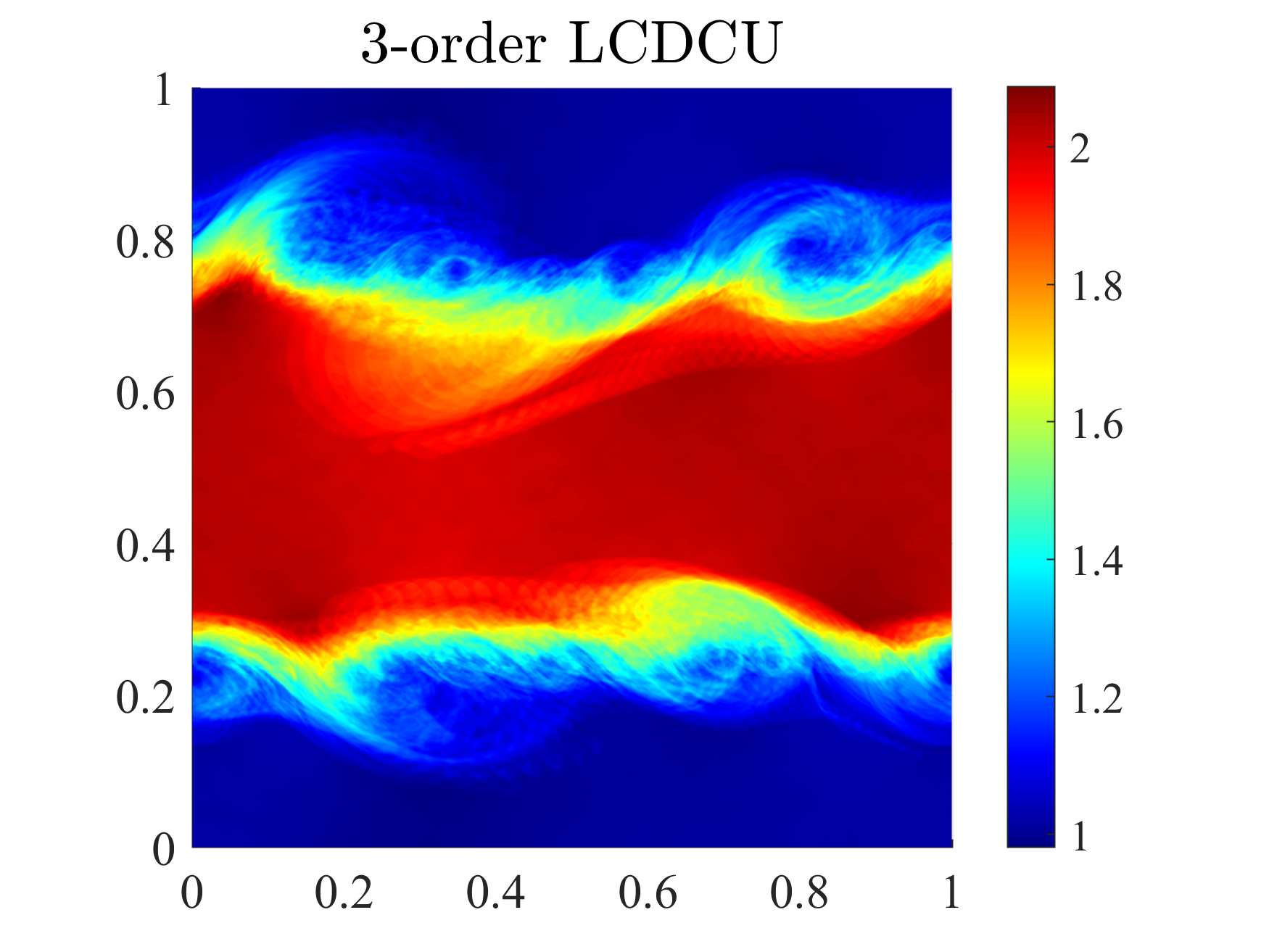}\hspace*{0.2cm}
            \includegraphics[trim=1.3cm 0.5cm 1.6cm 0.2cm, clip, width=4.cm]{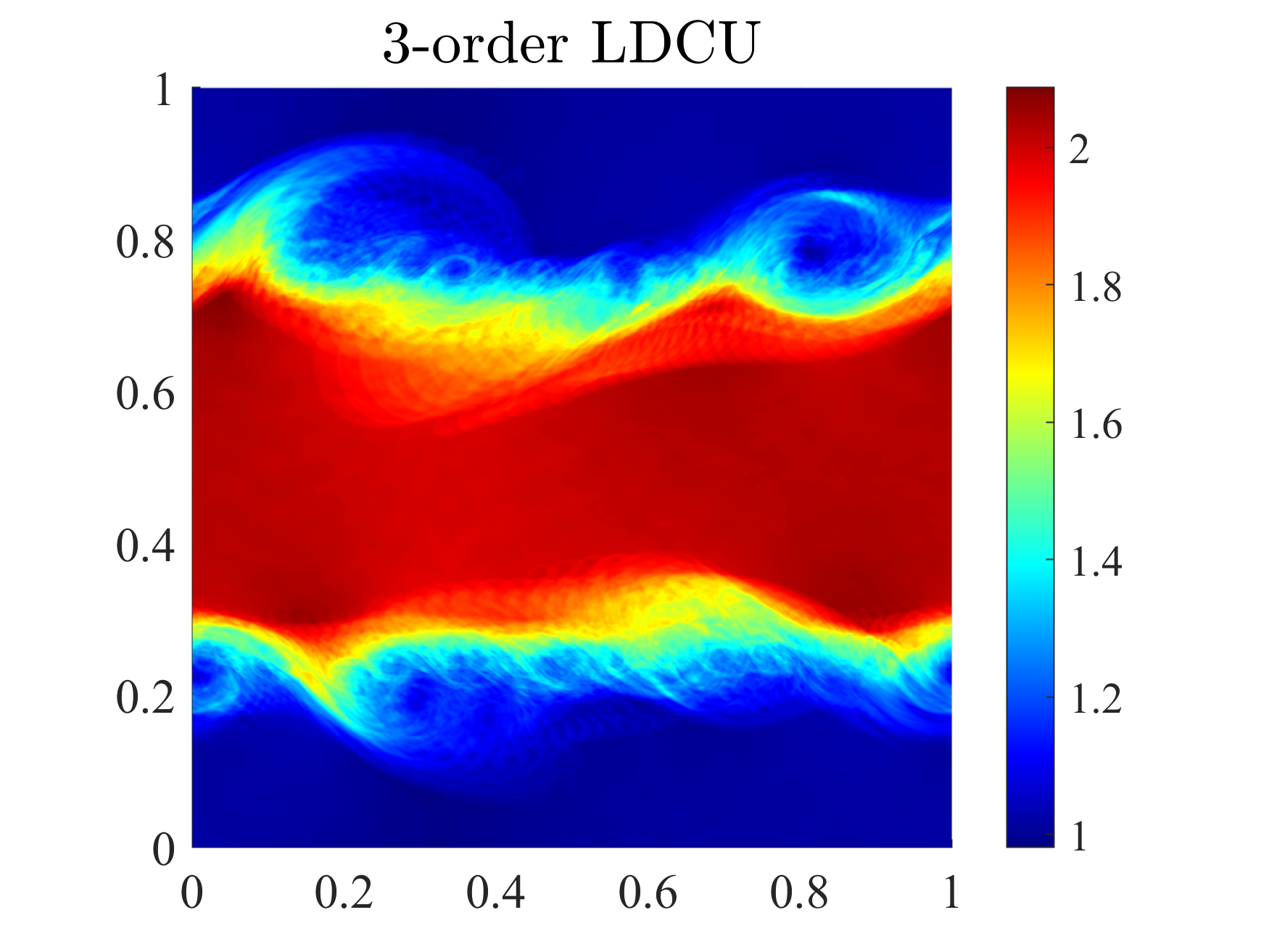}\hspace*{0.0cm}
            \includegraphics[trim=1.3cm 0.5cm 1.6cm 0.2cm, clip, width=4.2cm]{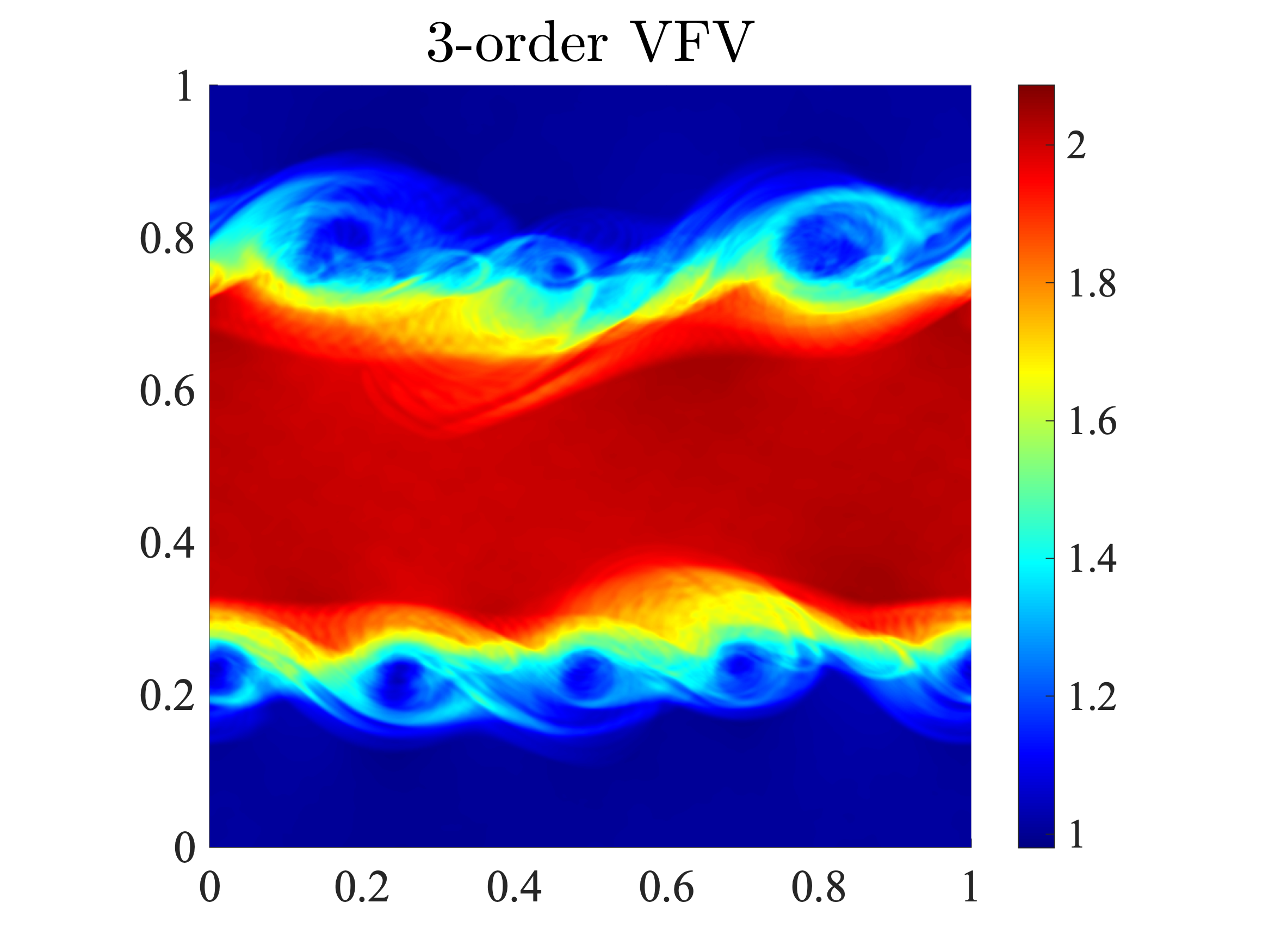}}
\vskip8pt
\centerline{\includegraphics[trim=1.3cm 0.5cm 1.6cm 0.2cm, clip, width=4.cm]{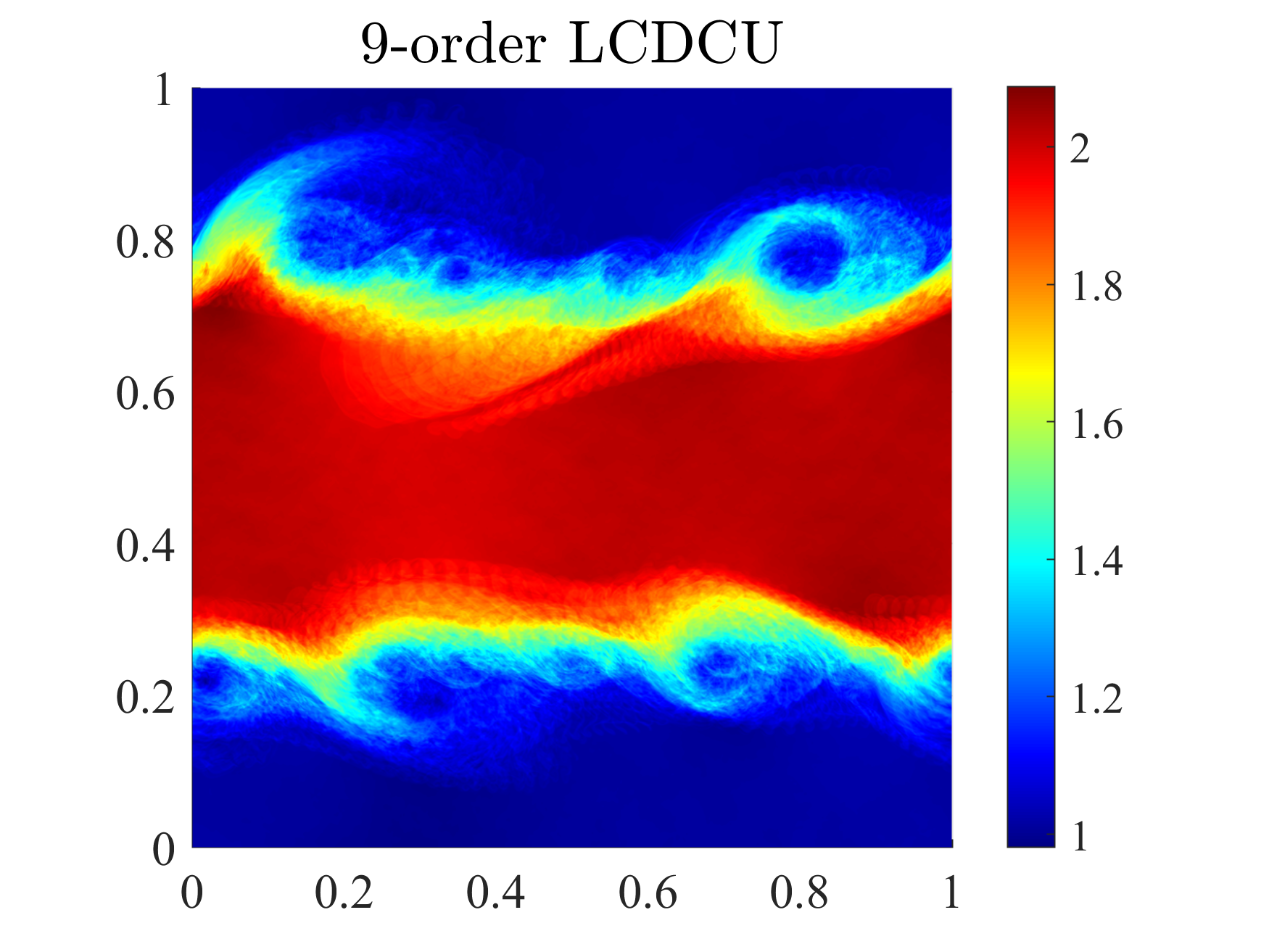}\hspace*{0.2cm}
            \includegraphics[trim=1.3cm 0.5cm 1.6cm 0.2cm, clip, width=4.cm]{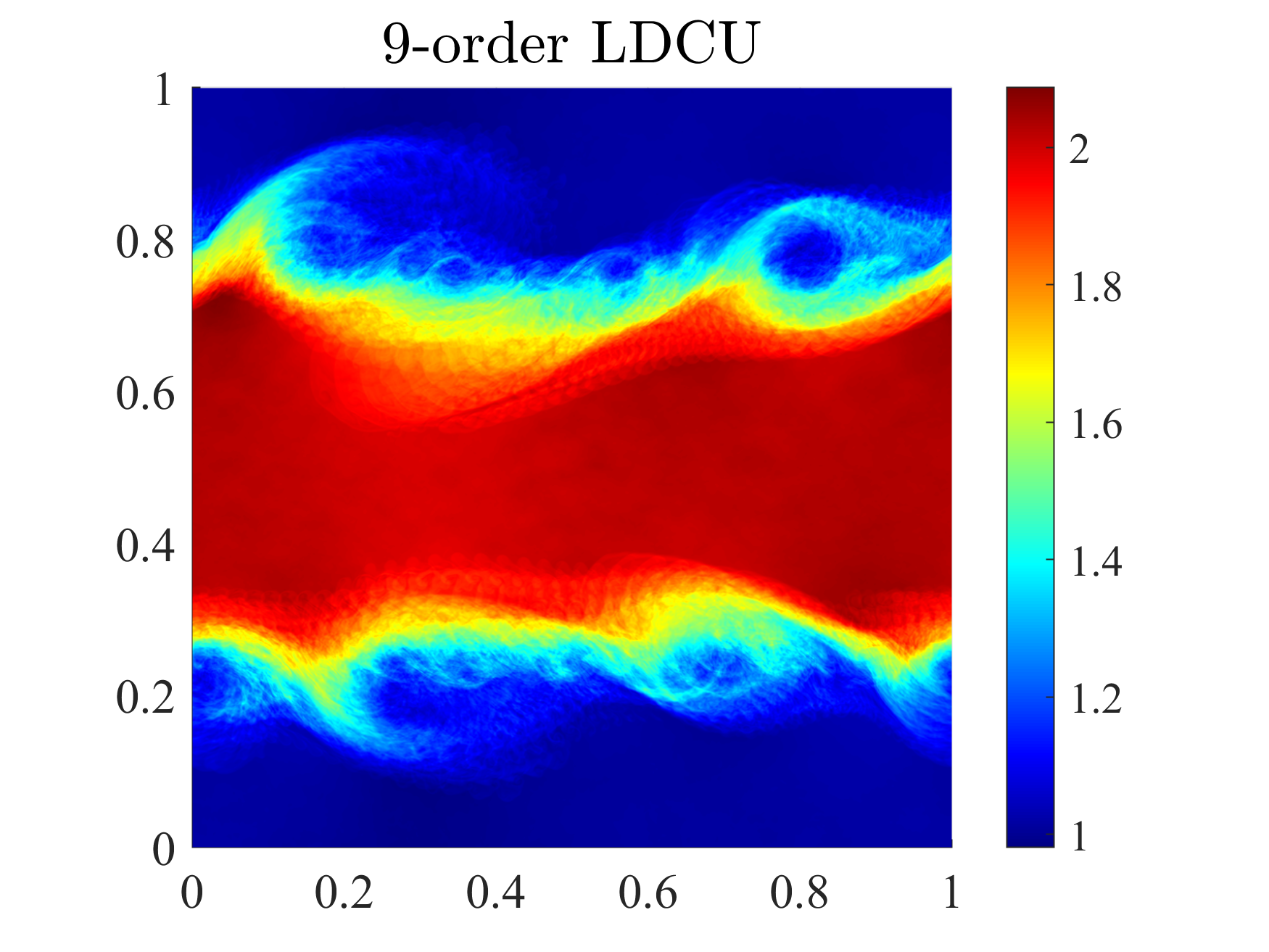}\hspace*{0.0cm}
            \includegraphics[trim=1.3cm 0.5cm 1.6cm 0.2cm, clip, width=4.2cm]{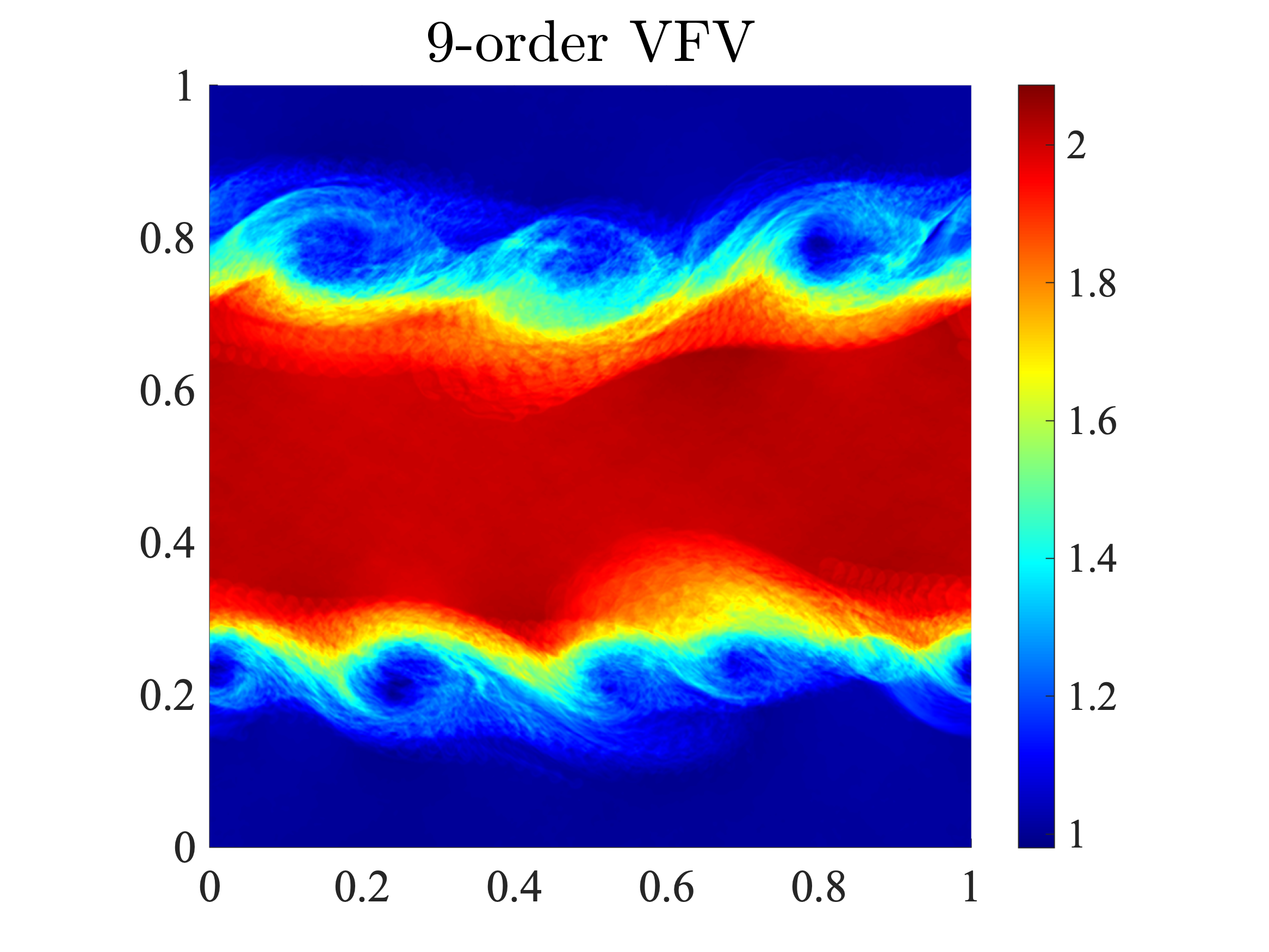}}
\caption{\sf KH Instability: $\rho^T_1$ (top row), $\rho^T_3$ (middle row), and $\rho^T_6$ (bottom row) computed by the LCDCU (left column), LDCU (middle column), and VFV (right column) schemes.\label{fig4}}
\end{figure}

\begin{figure}[ht!]
\centerline{\includegraphics[trim=1.3cm 0.5cm 1.6cm 0.2cm, clip, width=4.cm]{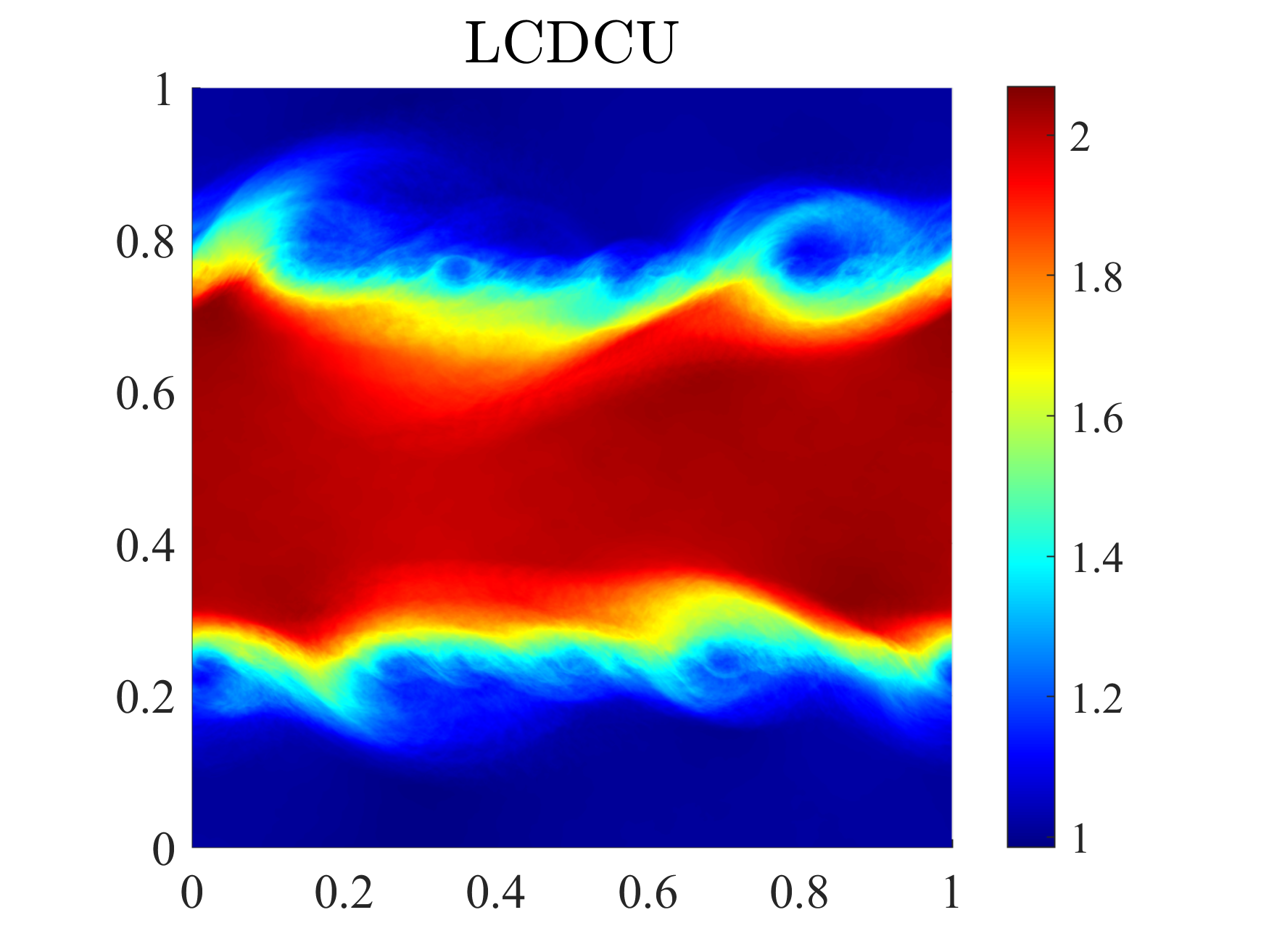}\hspace*{0.2cm}           
            \includegraphics[trim=1.3cm 0.5cm 1.6cm 0.2cm, clip, width=4.cm]{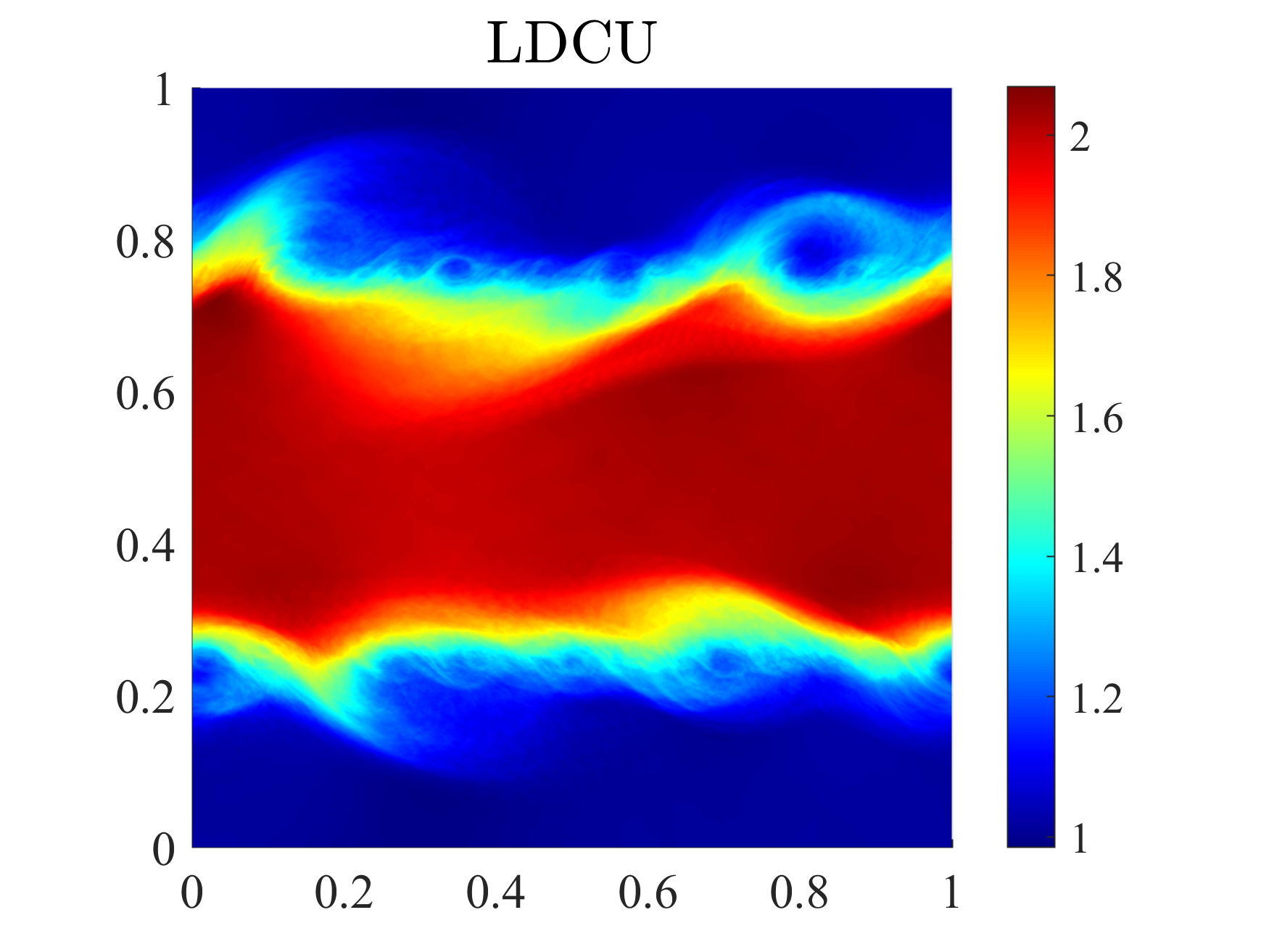}\hspace*{0.0cm}
            \includegraphics[trim=1.3cm 0.5cm 1.6cm 0.2cm, clip, width=4.2cm]{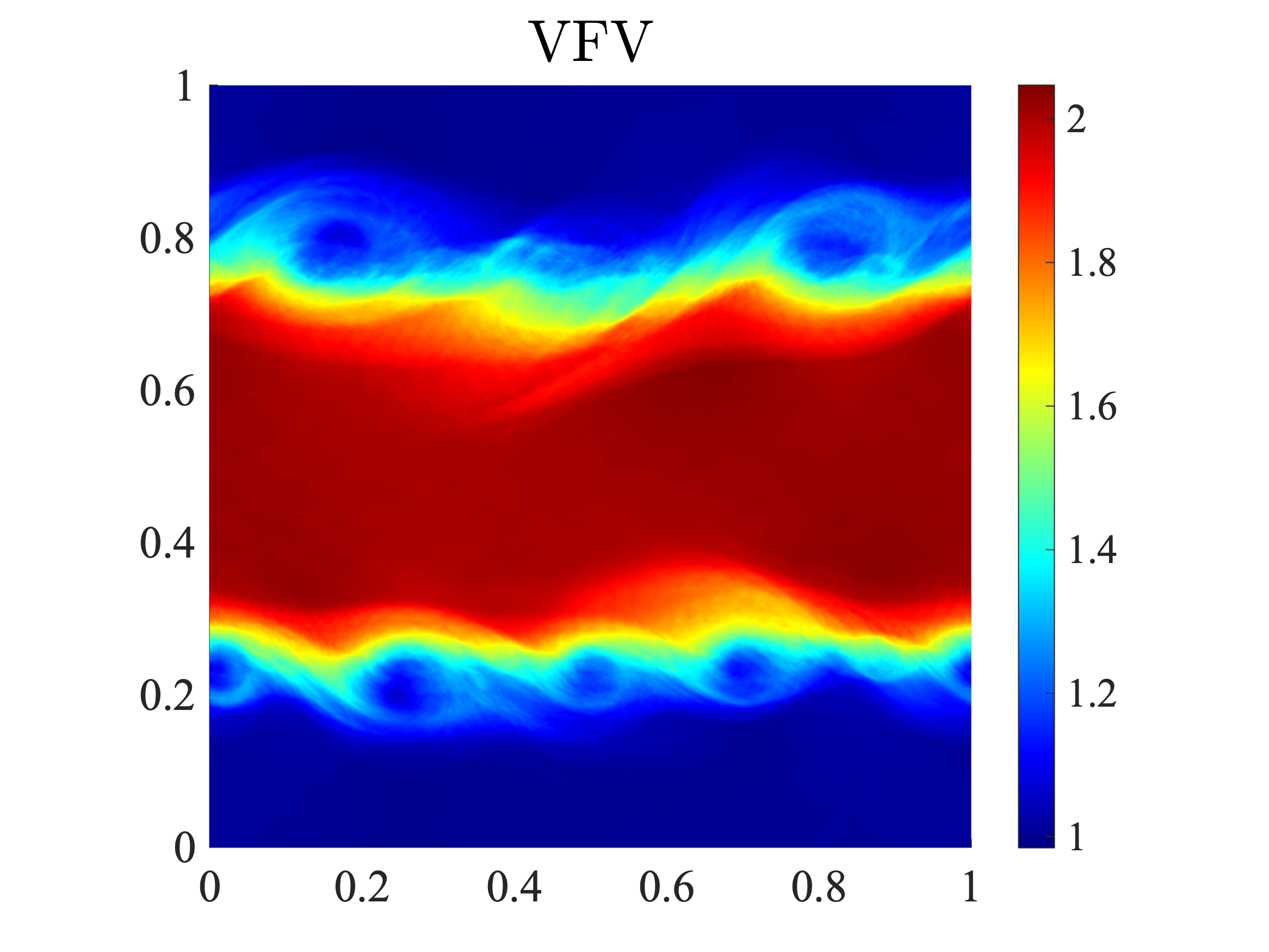}}
\caption{\sf KH Instability: $\widetilde\rho^{\,T}_6$ computed by the LCDCU (left), LDCU (middle), and VFV (right) schemes.\label{fig5}}
\end{figure}

Next, we consider the difference between the time averages of the numerical solutions. Figure \ref{fig4} presents the time averages of the density computed by the three methods. Similarly, we can observe that they are convergent as the order increases. In addition, the difference between the limiting average densities is pronounced. In Figure \ref{fig5}, we show the time averages of the mean density for the three methods. Remarkably, they look very similar. This phenomenon suggests that there is some consistency in the measure-valued solutions produced by different numerical methods.

\begin{figure}[ht!]
\centerline{\includegraphics[trim=1.2cm 0.5cm 1.4cm 0.1cm, clip, width=3.6cm]{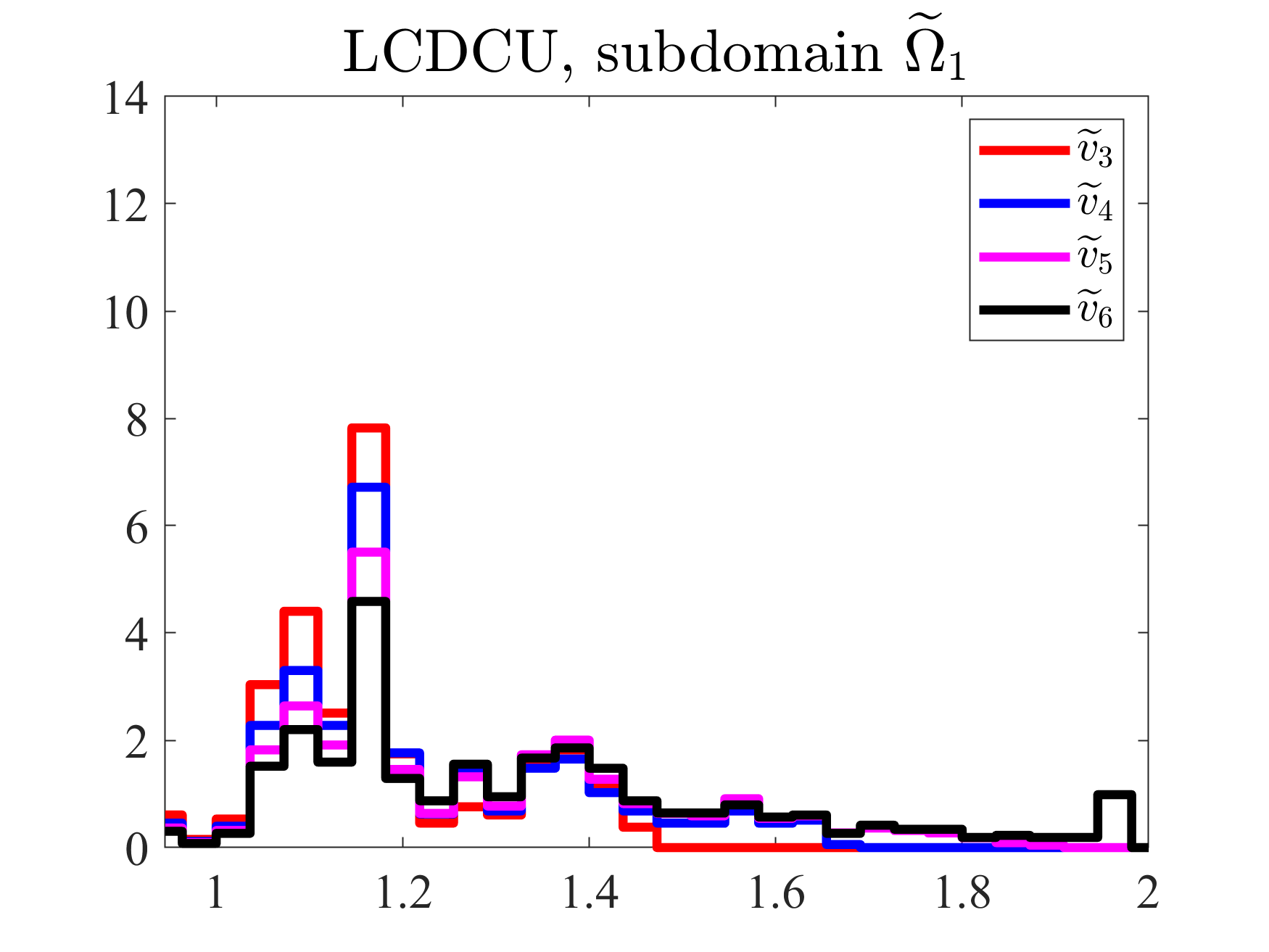}\hspace*{0.6cm}
            \includegraphics[trim=1.2cm 0.5cm 1.4cm 0.1cm, clip, width=3.6cm]{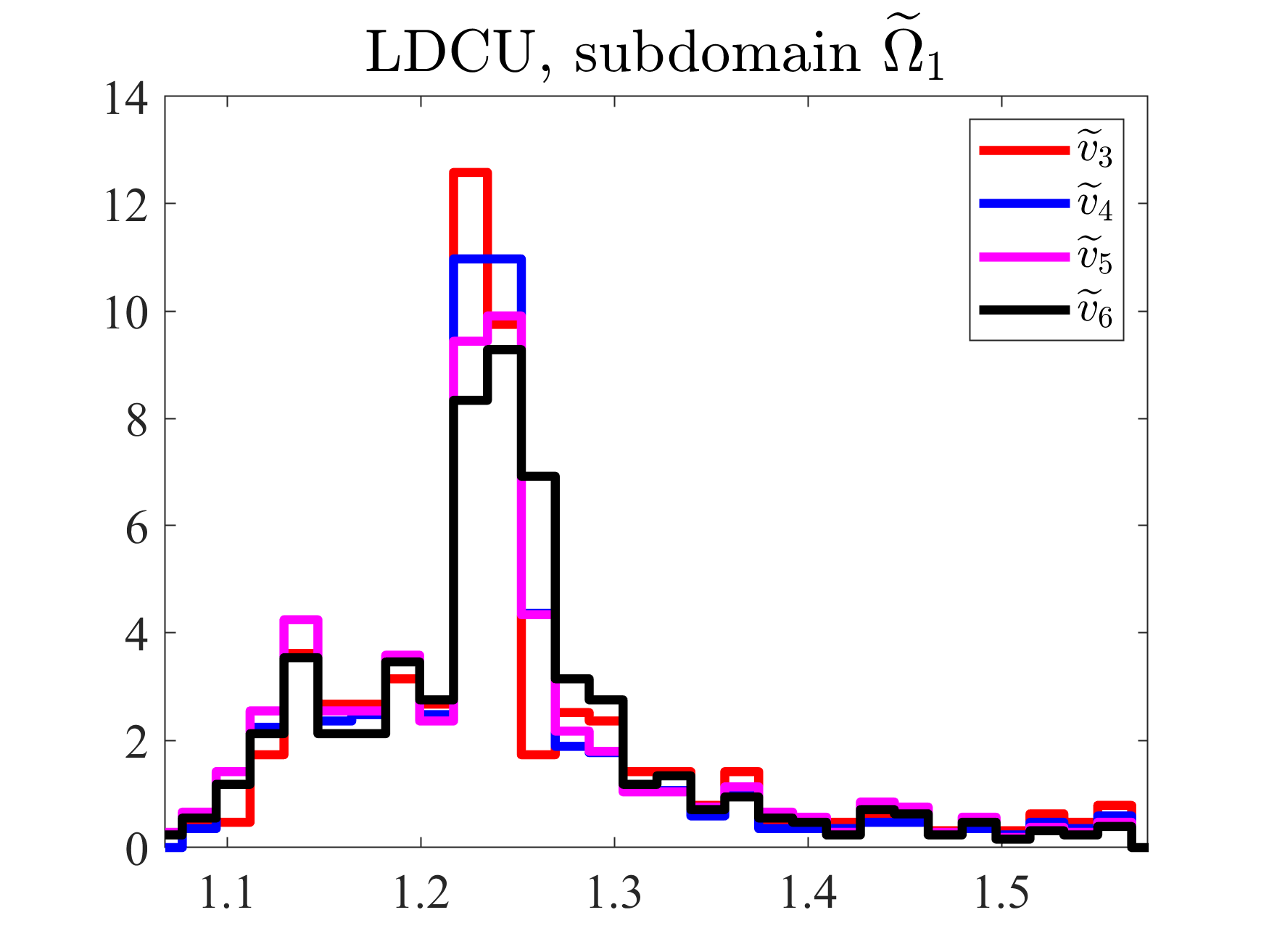}\hspace*{0.4cm}
            \includegraphics[trim=1.2cm 0.5cm 1.4cm 0.1cm, clip, width=3.8cm]{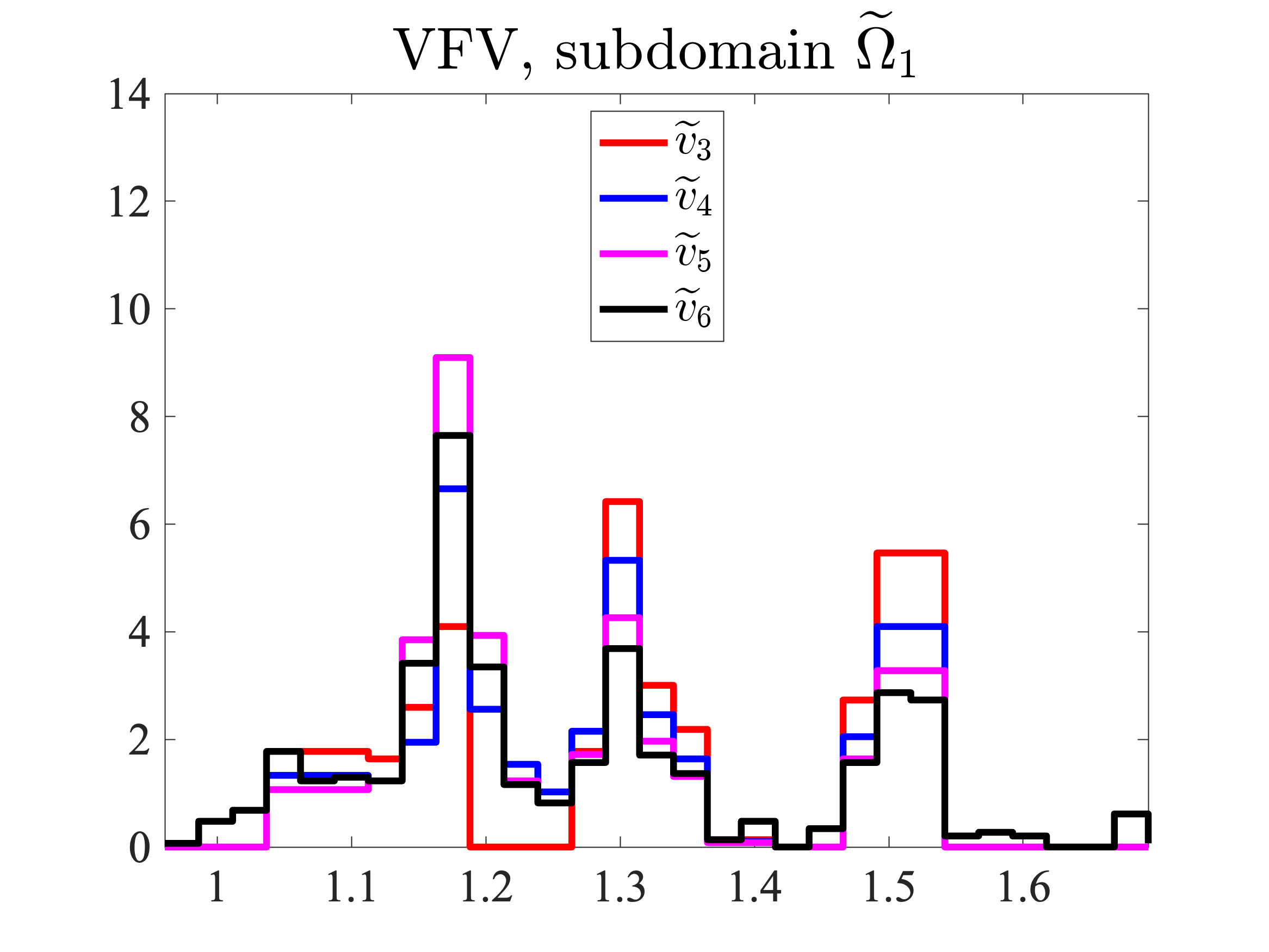}}
\vskip8pt
\centerline{\includegraphics[trim=1.2cm 0.5cm 1.4cm 0.1cm, clip, width=3.6cm]{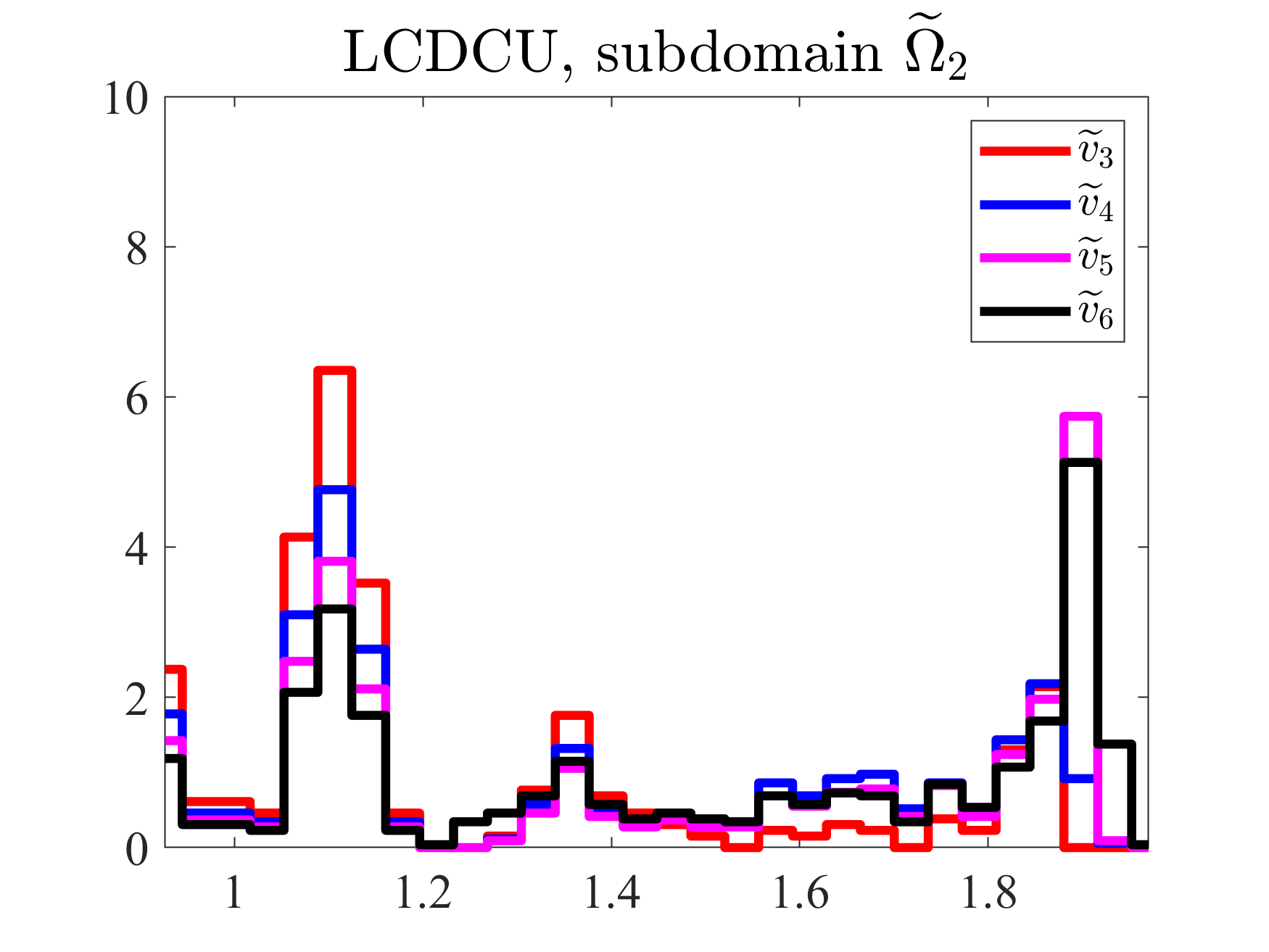}\hspace*{0.6cm}
            \includegraphics[trim=1.2cm 0.5cm 1.4cm 0.1cm, clip, width=3.6cm]{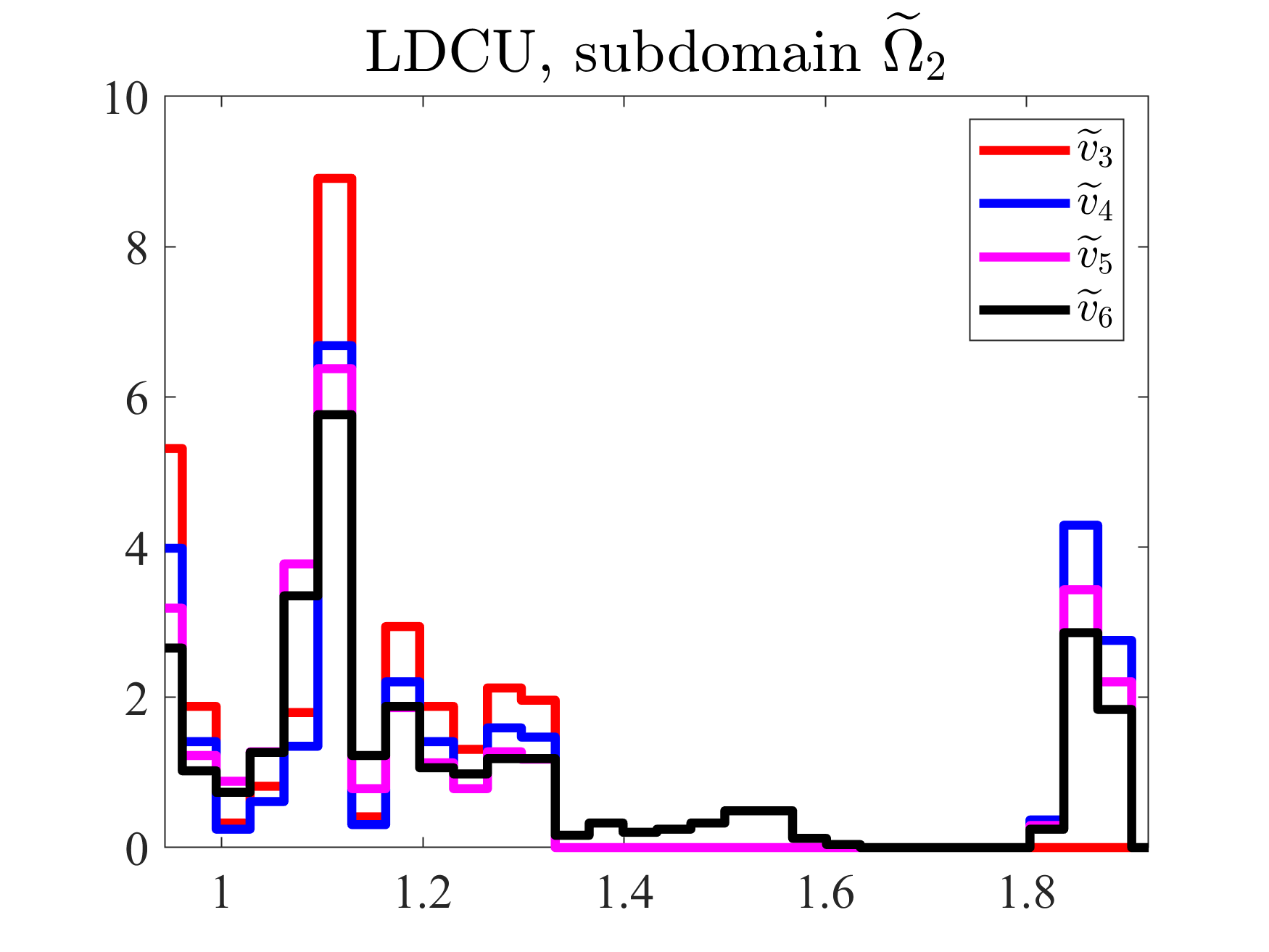}\hspace*{0.4cm}
            \includegraphics[trim=1.2cm 0.5cm 1.4cm 0.1cm, clip, width=3.8cm]{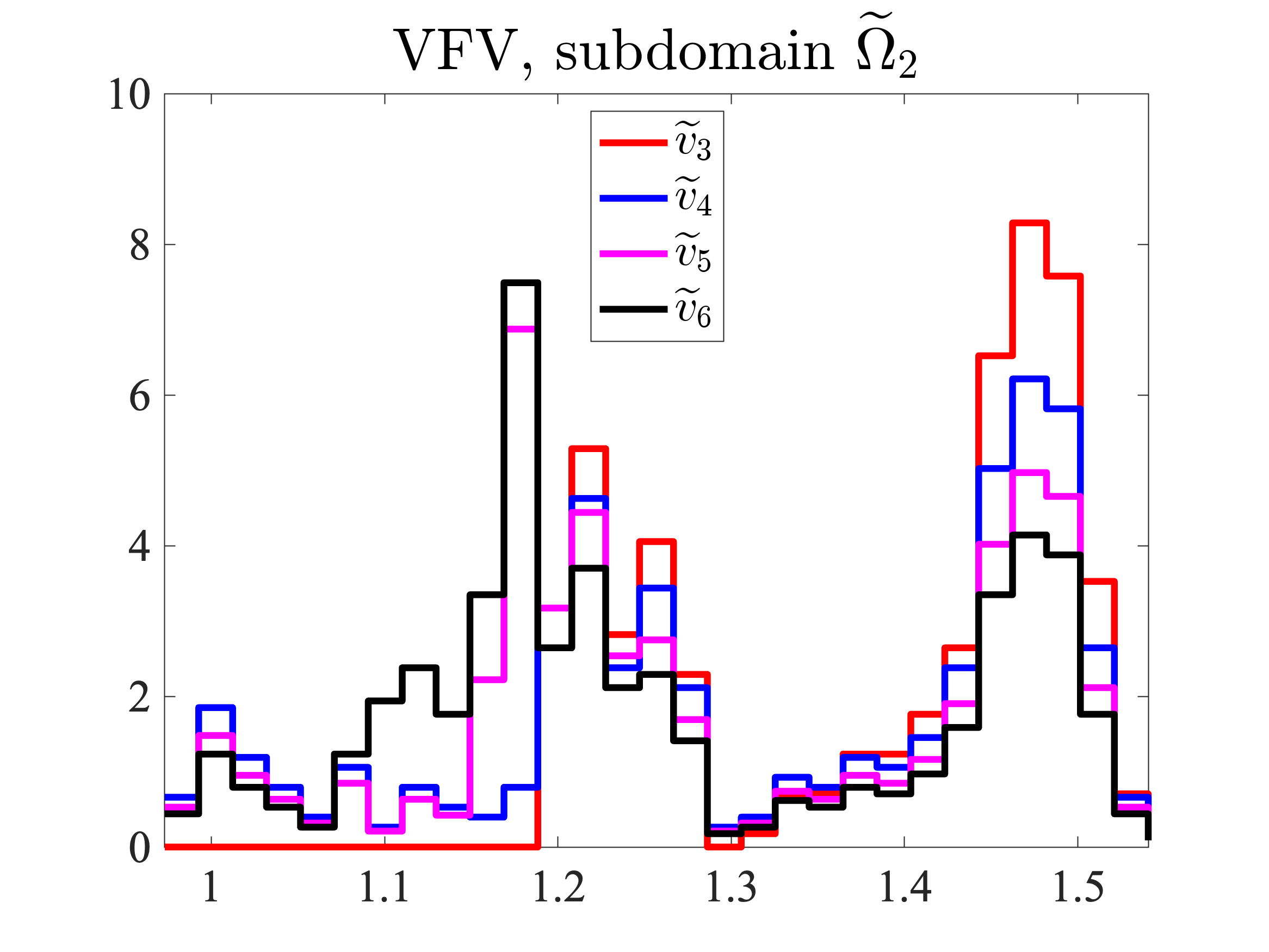}}
\caption{\sf KH Instability: $\widetilde\sigma_3$, $\widetilde\sigma_4$, $\widetilde\sigma_5$, and $\widetilde\sigma_6$ computed by the LCDCU (left column), LDCU (middle column), and VFV (right column) schemes in the subdomains $\widetilde\Omega_1$ (top row) and $\widetilde\Omega_2$ (bottom row).\label{fig6}}
\end{figure}

To illustrate the differences in numerical results, we compare the approximate PDFs of the Young measures. We select two subdomains:
\begin{equation*}
\widetilde\Omega_1=[0.50,0.51]\times[0.25,0.26]\quad\mbox{and}\quad\widetilde\Omega_2=[0.15,0.16]\times[0.75,0.76],
\end{equation*}
and present the obtained results in Figure \ref{fig6}. It is evident that these PDFs differ significantly, particularly between those
computed by the VFV method and by the other two methods. This indicates that the different methods can yield different DW solutions.

\subsection{Convergence of Averages over Orders}\label{sec43}
\begin{figure}[ht!]
\centerline{\includegraphics[trim=0.5cm 0.5cm 1.4cm 0.2cm, clip, width=3.8cm]{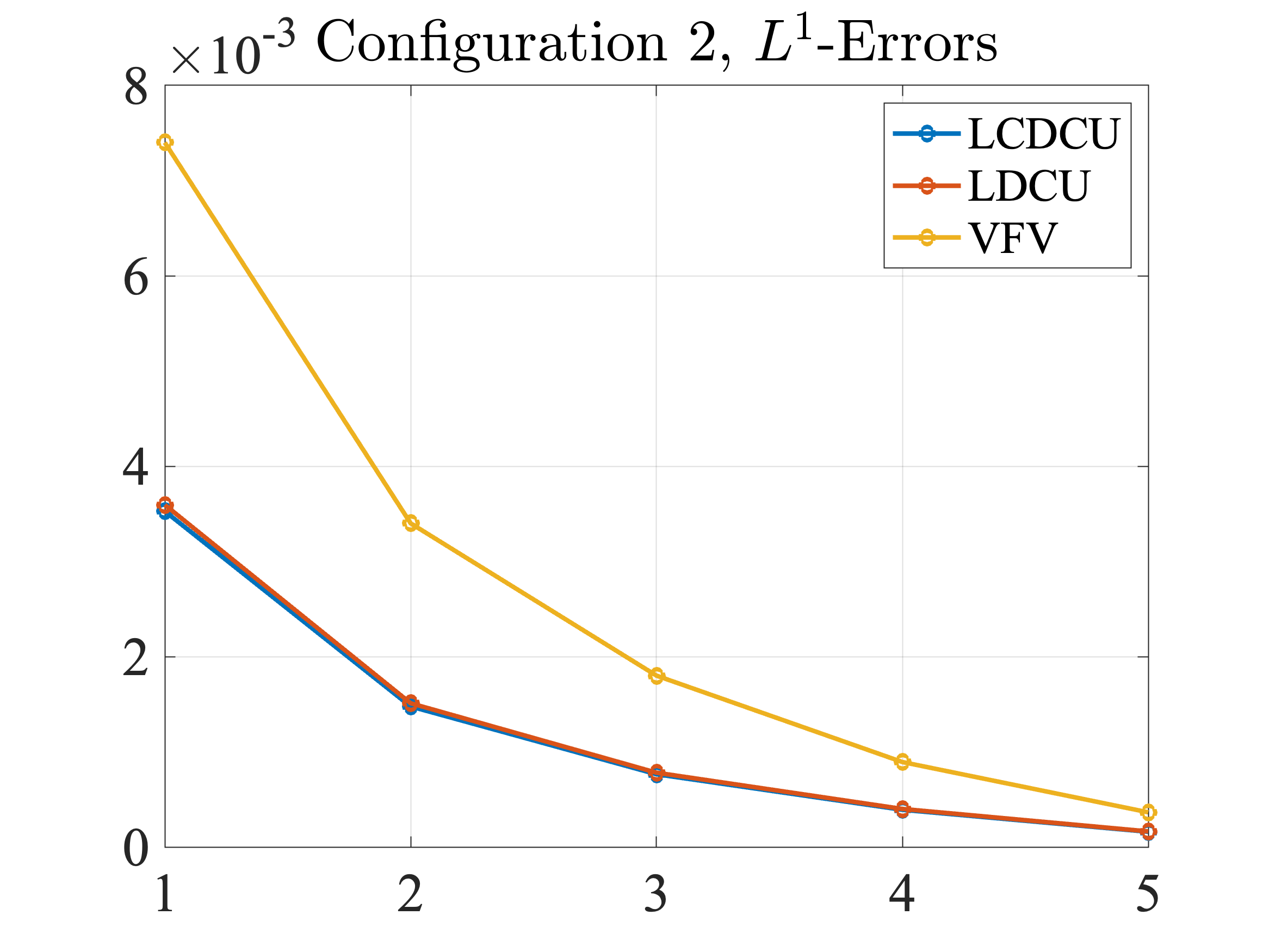}\hspace*{0.5cm}
            \includegraphics[trim=0.5cm 0.5cm 1.4cm 0.2cm, clip, width=3.8cm]{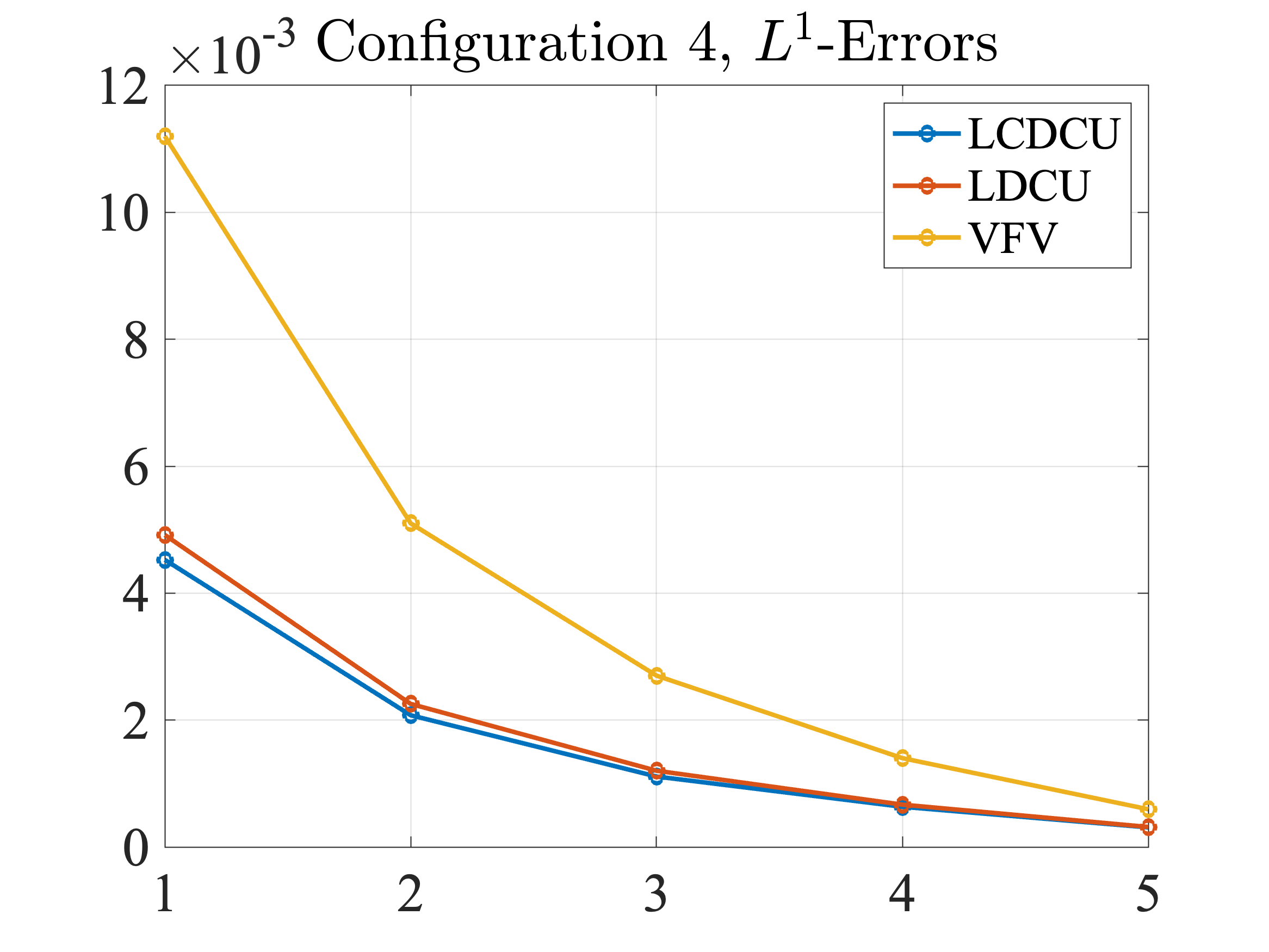}\hspace*{0.5cm}
            \includegraphics[trim=0.5cm 0.5cm 1.4cm 0.2cm, clip, width=3.8cm]{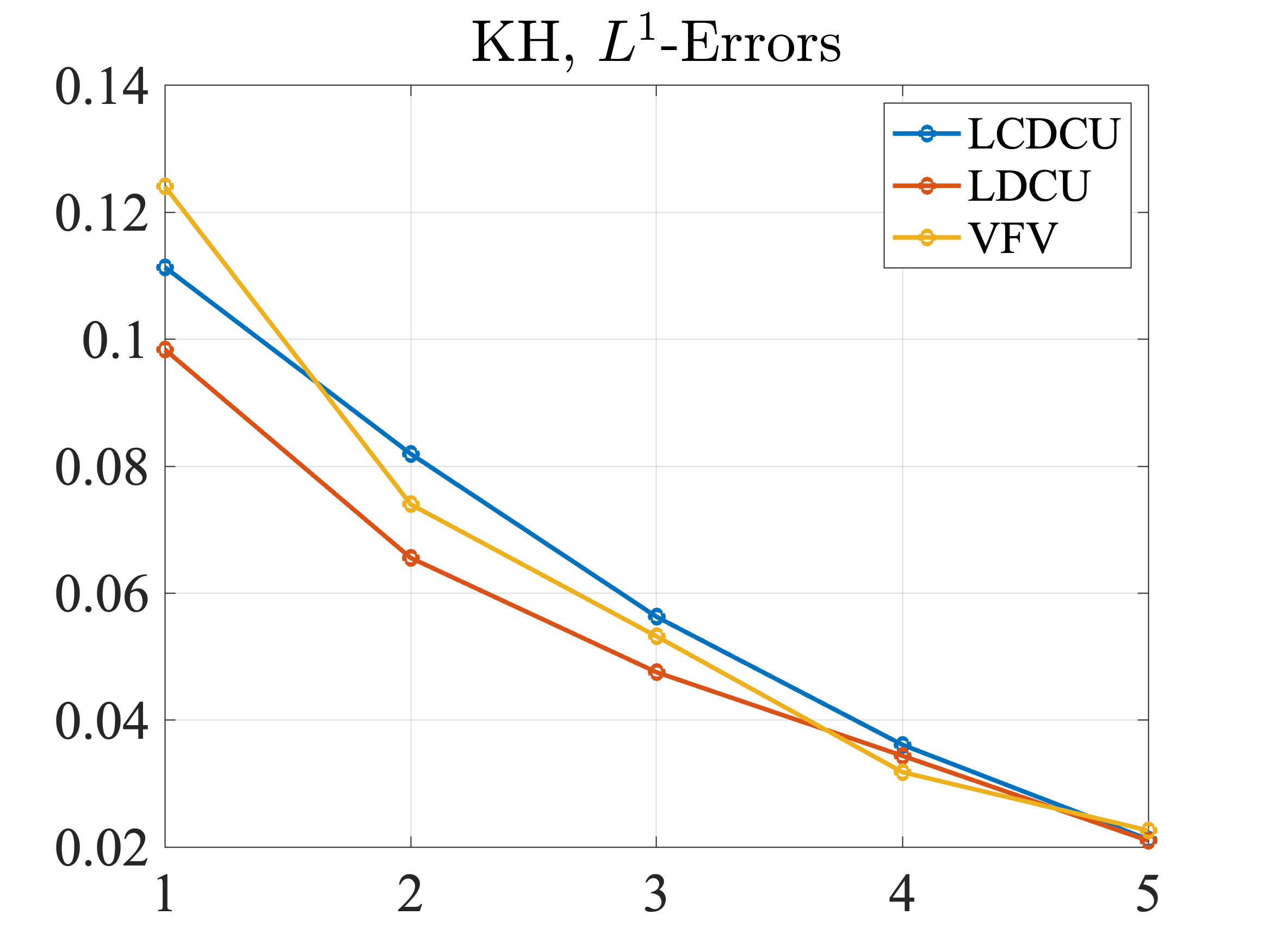}}
\caption{\sf The $L_1$-errors of the average density computed by the studied schemes for Configuration 2 (left), Configuration 4 (middle),
and KH instability (right) problems.\label{fig15}}
\end{figure}

In this section, we evaluate the convergence of the averages over orders of accuracy. We estimate the errors of these averages using the
quantities
\begin{equation*}
\left\|\widetilde\rho_n-\widetilde\rho_6\right\|_{L_1(\Omega)},\quad n=1,\dots,5,
\end{equation*}
where $\widetilde\rho_6$ serves as a reference solution. In Figure \ref{fig15}, we present the errors associated with the two Riemann
problems (Configurations 2 and 4) and KH instability problem. It is evident that the errors decrease as $n$ increases: This behavior
demonstrates the ${\mathcal K}$-convergence of the sequence of averages. 

To further analyze the behavior of errors with respect to the orders, we utilize the least squares data fitting to approximate the error
functions, which take the form of power functions of $n$; see Table \ref{table1}. One can observe a substantial difference in the decay
rates for different problems. For Configuration 2, in which the computed weak solution is continuous, the convergence is the fastest. It is
somewhat slower for the Configuration 4, in which the solution is discontinuous, but still not a truly DW solution. Finally, for the KH
instability problem, which admits genuine measure-valued solutions, the convergence is substantially slower and, as we can see from Figure
\ref{fig15}, the magnitude of the errors is much larger than those obtained for the considered Riemann problems. 
\begin{table}[ht!]
\centering
\begin{tabular}{|c|c|c|c|}\hline
&LCDCU&LDCU&VFV\\\hline
Configuration 2&$0.0043\,n^{-1.82}$&$0.0044\,n^{-1.82}$&$0.0093\,n^{-1.77}$\\\hline
Configuration 4&$0.0053\,n^{-1.60}$&$0.0059\,n^{-1.63}$&$0.0136\,n^{-1.75}$\\\hline
KH Instability &$0.1335\,n^{-0.98}$&$0.1104\,n^{-0.90}$&$0.1382\,n^{-1.04}$\\\hline
\end{tabular}\\[1.ex]
\caption{\sf The errors of the average density as functions of $n$ for the LCDCU, LDCU, and VFV schemes.\label{table1}}
\end{table}

\subsection{Selection Criteria}
As seen in the previous results (Configuration 3 for the 2-D Riemann problem and the KH instability problem), different numerical schemes
can converge to different DW solutions. Consequently, we need additional criteria to select the physically relevant one. In this section, we
analyze potential selection criteria based on the entropy production and energy defects. Notice that similar selection criteria were studied
theoretically in \cite{FL26,FJL26,FLY}. 

We denote the integral of the average total entropy over orders up to $n$ by 
\begin{equation*}
{\mathcal S}^n(t):=\int\limits_\Omega\widetilde S_n(x,y,t)\,{\rm d}x{\rm d}y.
\end{equation*}
The total energy and mean energy are defined by 
\begin{equation*}
\begin{aligned}
{\mathcal E}^n_1(t)&:=\int\limits_\Omega\widetilde E_n(x,y,t)\,{\rm d}x{\rm d}y=
\int\limits_\Omega\frac{1}{n}\sum_{\ell=1}^n\left[\frac{|\bm m_\ell|^2}{2\rho_\ell}+\rho_\ell e(\rho_\ell,S_\ell)\right]{\rm d}x{\rm d}y,\\
{\mathcal E}^n_2(t)&:=\int\limits_\Omega E\big(\widetilde\rho_n(x,y,t),\widetilde{\bm m}_n(x,y,t),\widetilde S_n(x,y,t)\big){\rm d}x{\rm d}y \\
&=\int\limits_\Omega\frac{|\widetilde{\bm m}_n|^2}{2\widetilde\rho_n}+\widetilde\rho_ne\big(\widetilde\rho_n,\widetilde S_n\big)
{\rm d}x{\rm d}y,
\end{aligned}
\end{equation*}
respectively, and the energy defect is defined by
\begin{equation*}
{\mathcal D}^n_{\rm E}(t):= \int\limits_\Omega\left|\widetilde E_n(x,y,t)-E\big(\widetilde\rho_n(x,y,t),\widetilde{\bm m}_n(x,y,t),
\widetilde S_n(x,y,t)\big)\right|{\rm d}x{\rm d}y.
\end{equation*}
In addition, we consider a so-called Bregman distance between a solution $(\rho, \bm m, S)$ at time $t$ and the constant equilibrium state given by $(\bar{\rho},\bm0, S(\bar{\rho},\bar{\theta}))$ that maximizes the entropy production. Here, 
$$
\bar{\rho}= \frac{1}{|\Omega|} \int\limits_\Omega \rho_0 \, \mathrm{d}x,
\quad \bar\theta=\frac{1}{ c_v \bar\rho |\Omega|} \int\limits_\Omega E_0 \, \mathrm{d}x,
$$
and the Bregman distance is given by
\begin{equation*}
{\mathcal J}^n(t):= \mathcal{E}^n_2(t) -\bar \theta \mathcal{S}^n(t);
\end{equation*}
see \cite{FL26} for details on its derivation.

Possible selection criteria for DW solutions are (i) the maximization of entropy production, (ii) the minimization of energy defect,
(iii) the maximization of energy defect, and (iv) minimizing their combination in $\mathcal{J}^n$.

et us briefly discuss the physical interpretation of the proposed selection criteria. First, criteria (i) and (iv) are related to the Second law of thermodynamics stating that physically reliable solutions either preserve or increase the entropy over time. Here, we select the solutions that either maximize the entropy production (i) or are closest to the equilibrium state with maximal entropy (iv). The latter expresses the Clausius formulation of the Second law of the thermodynamics: ``The energy of the world (closed system) is constant; its entropy tends to a maximum.'' 

Criteria (ii) and (iii)  are not contradictory since they help to select the solutions
with either the minimal or maximal oscillatory (turbulent) structures, respectively, expressed by the energy defects. As formulated in Definition~\ref{D2}, the energy defect controls the turbulent Reynolds stress defect.

As shown in \cite{FL26}, the class of DW solutions emanated from a given initial data is convex, closed, and weakly compact in a suitable topology. The energy $E(\rho,\bm m, S)$ is a strictly convex function on its domain of definition with $\rho > 0.$ Consequently, the energy defect is a strictly concave function and both criteria  (iii) or (iv) yield a unique selected solution in a DW class.

\begin{figure}[ht!]
\centerline{\includegraphics[trim=0.5cm 0.5cm 1.0cm 0.1cm, clip, width=3.8cm]{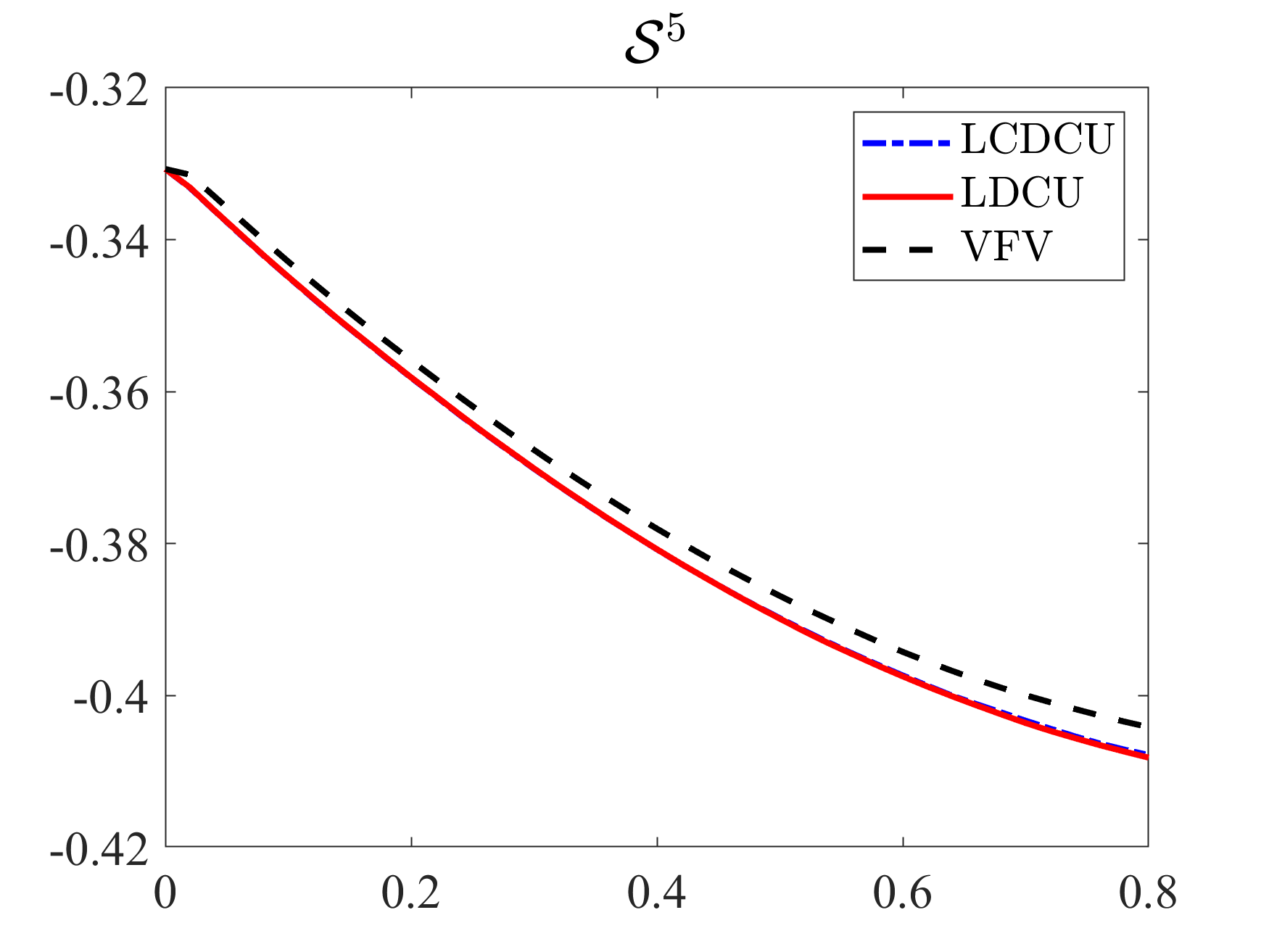}\hspace*{0.5cm}
            \includegraphics[trim=0.5cm 0.5cm 1.0cm 0.1cm, clip, width=3.8cm]{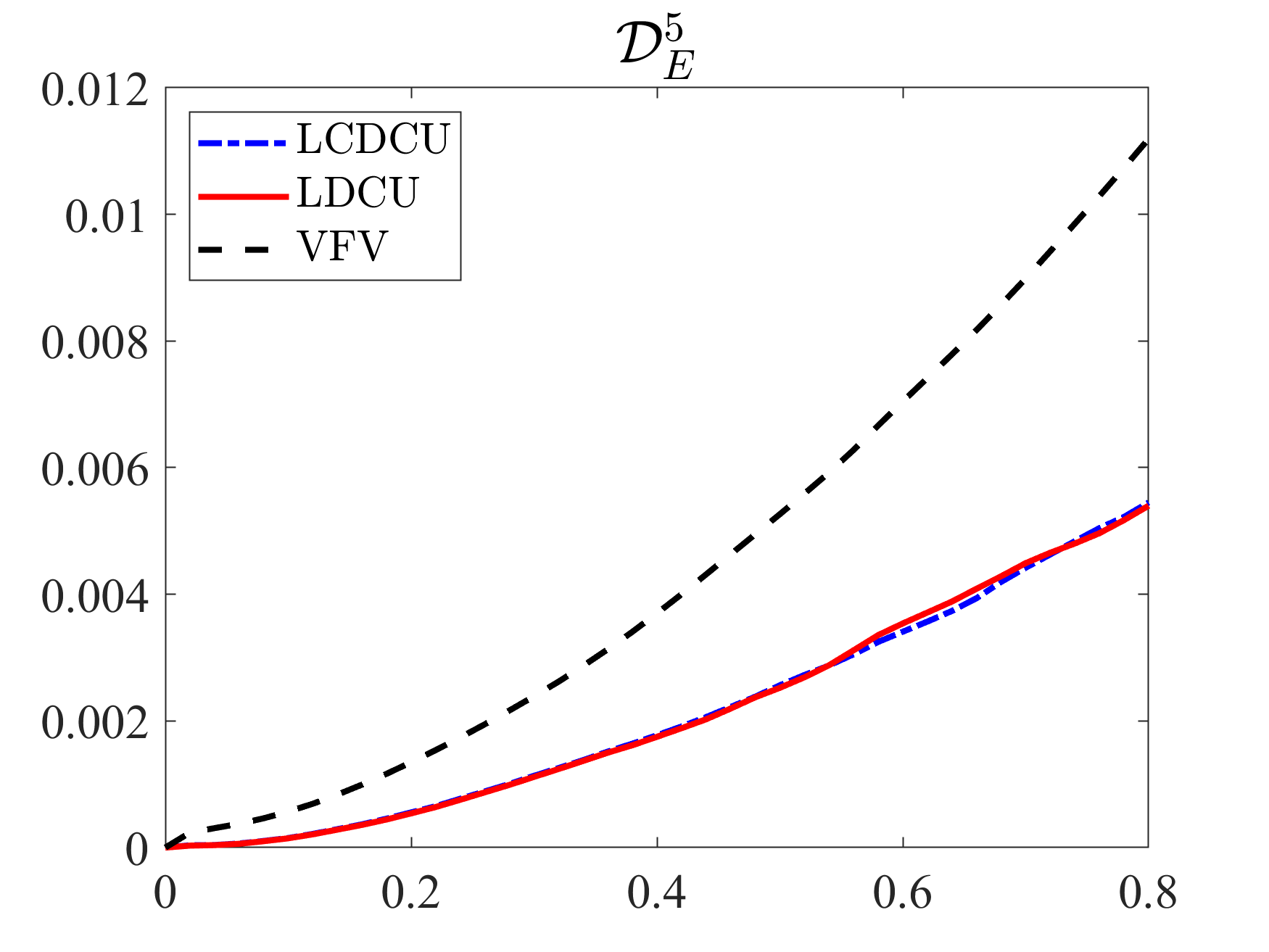}\hspace*{0.5cm}
            \includegraphics[trim=0.5cm 0.5cm 1.0cm 0.1cm, clip, width=3.8cm]{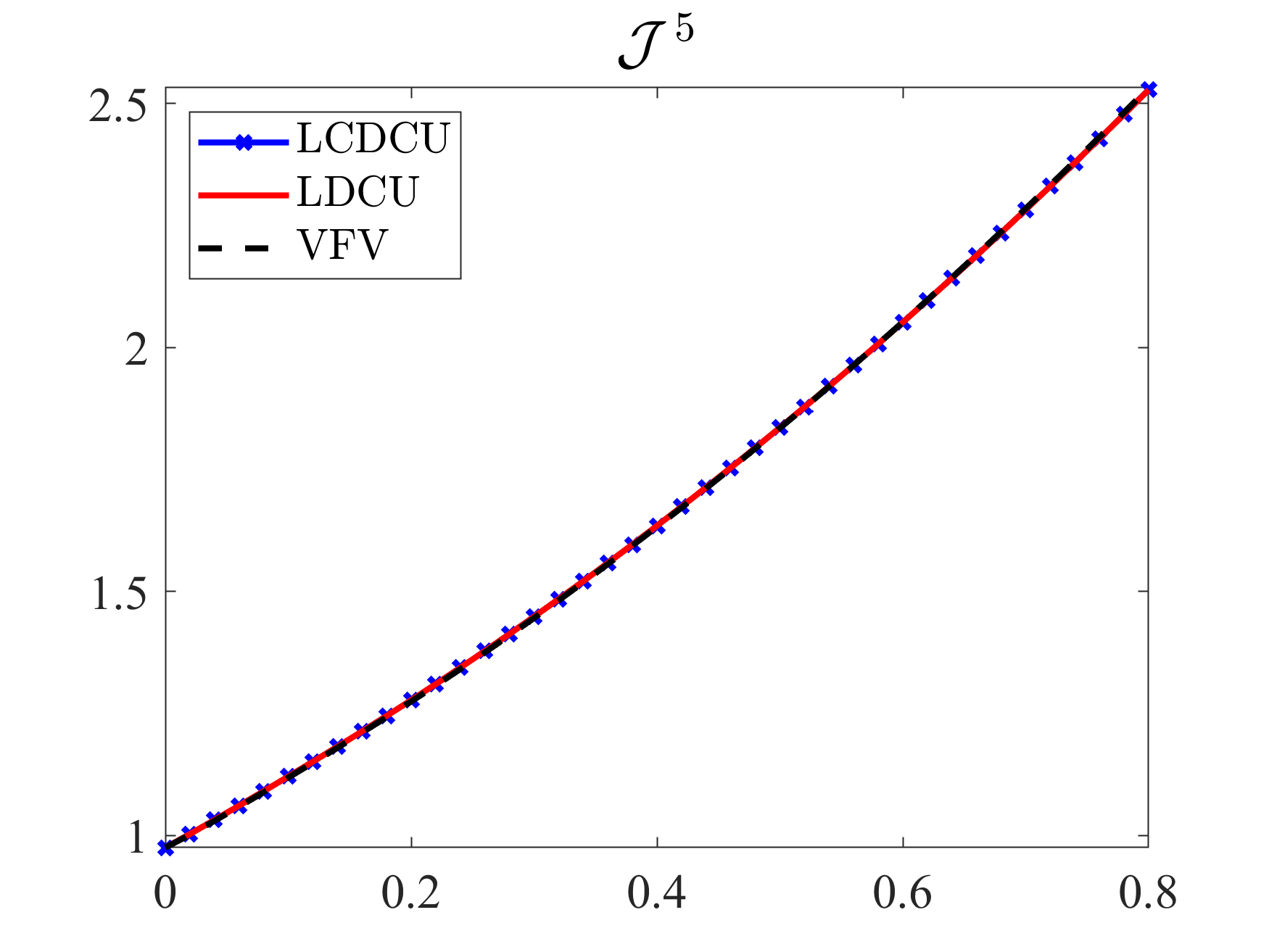}}

\caption{\sf Configuration 3: Time evolution of the entropy ${\mathcal S}^5$ (left), energy defects ${\mathcal D}^5_{\rm E}$ (middle),  and the Bregman distance ${\mathcal J}^5$ (right) computed by the LCDCU, LDCU, and VFV schemes.\label{fig18}}
\end{figure}

Based on the numerical solutions obtained by different methods, we compare their respective entropy production, energy defects, and the Bregman distance to the equilibrium state.
We begin with Configuration 3 for the 2-D Riemann problem. The time evolution of ${\mathcal S}^5$, ${\mathcal D}^5_{\rm E}$, and $\mathcal{J}^5$ are shown in Figure \ref{fig18}, and their time averages are presented in Table \ref{table4}. Here, $\bar{\theta}=0.864596$ is computed from the piecewise constant projection of the initial energy and mass.
As expected, the oscillations in specific regions of the flow, observed in \S\ref{sec41}, lead to energy defects. It can be observed that the LCDCU and LDCU
schemes produce smaller defects compared to the VFV schemes. At the same time, these two methods result in smaller entropy productions.
Consequently, criteria (i), (iii), and (iv) prioritize the VFV schemes, while criterion (ii) prioritizes the LCDCU schemes. Note that the differences between the LCDCU and LDCU schemes with respect to all selection criteria are very small.
We also emphasize that since the computed solutions develop oscillatory structures in a small part of the computational domain only, the magnitude
of the energy defects is quite small in this example.

\begin{table}[ht!]
\centering
\begin{tabular}{|c|c|c|c|}\hline
&LCDCU&LDCU&VFV\\\hline
$\frac{1}{T}\int_0^T{\mathcal S}^5(t)\,{\rm d}t$&-0.3778&-0.3778&-0.3741\\\hline
$\frac{1}{T}\int_0^T{\mathcal E}^5_1(t)\,{\rm d}t$&1.3690&1.3690&1.3541\\\hline
$\frac{1}{T}\int_0^T{\mathcal E}^5_2(t)\,{\rm d}t$&1.3669&1.3668&1.3498\\\hline
$\frac{1}{T}\int_0^T{\mathcal D}^5_{\rm E}(t)\,{\rm d}t$&0.0021&0.0022&0.0045\\\hline
$\frac{1}{T}\int_0^T{\mathcal J}^5(t)\,{\rm d}t$ &1.7179&1.7178&1.6732\\\hline
\end{tabular}\\[1.ex]
\caption{\sf Configuration 3: Time averages of the entropy,  total and mean energy,  the energy defects,  and the Bregman distance computed by the LCDCU, LDCU, and VFV schemes.\label{table4}}
\end{table}

\begin{figure}[ht!]
\centerline{\includegraphics[trim=0.5cm 0.5cm 1.0cm 0.1cm, clip, width=3.8cm]{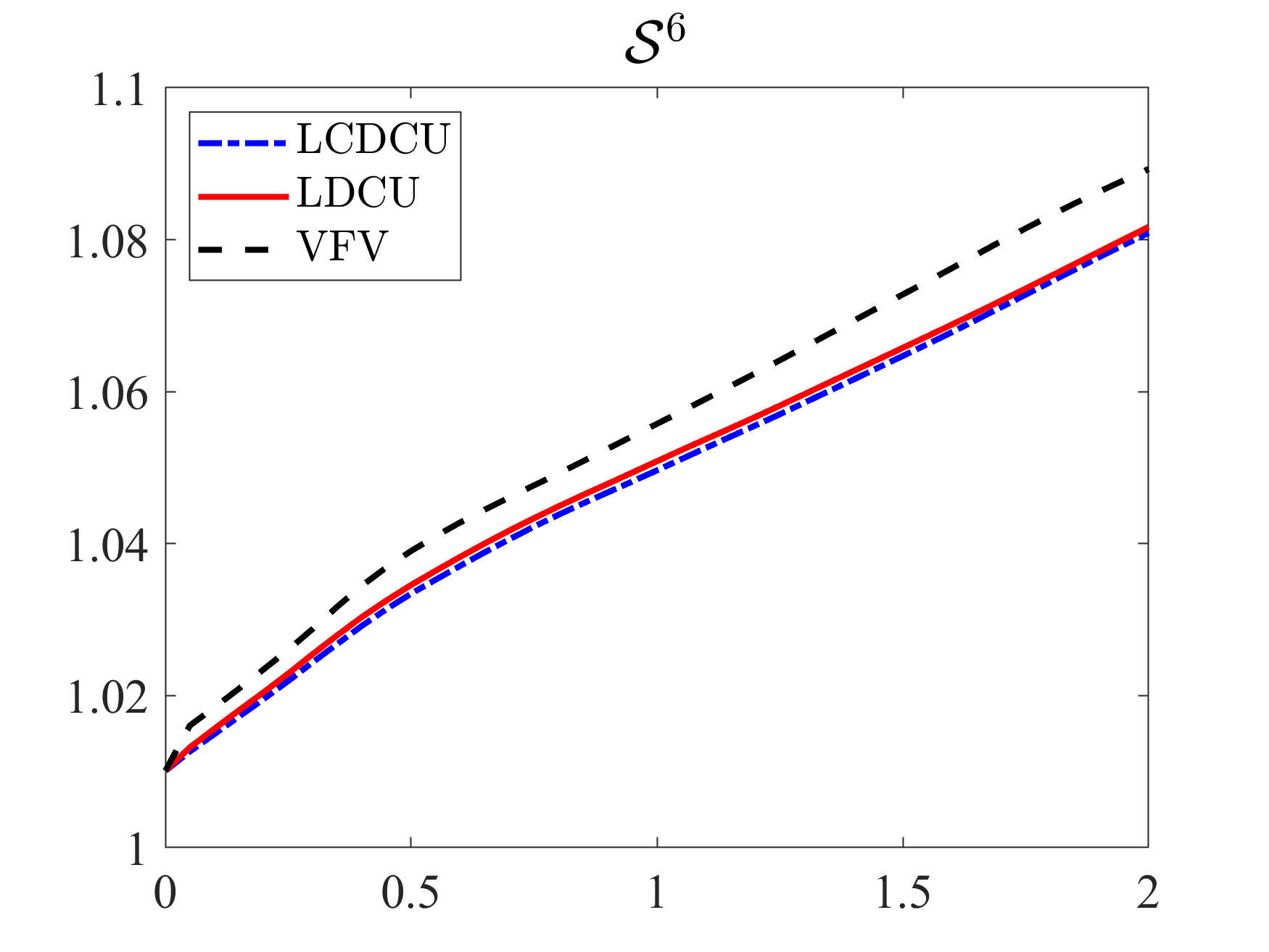}\hspace*{0.5cm}
            \includegraphics[trim=0.5cm 0.5cm 1.0cm 0.1cm, clip, width=3.8cm]{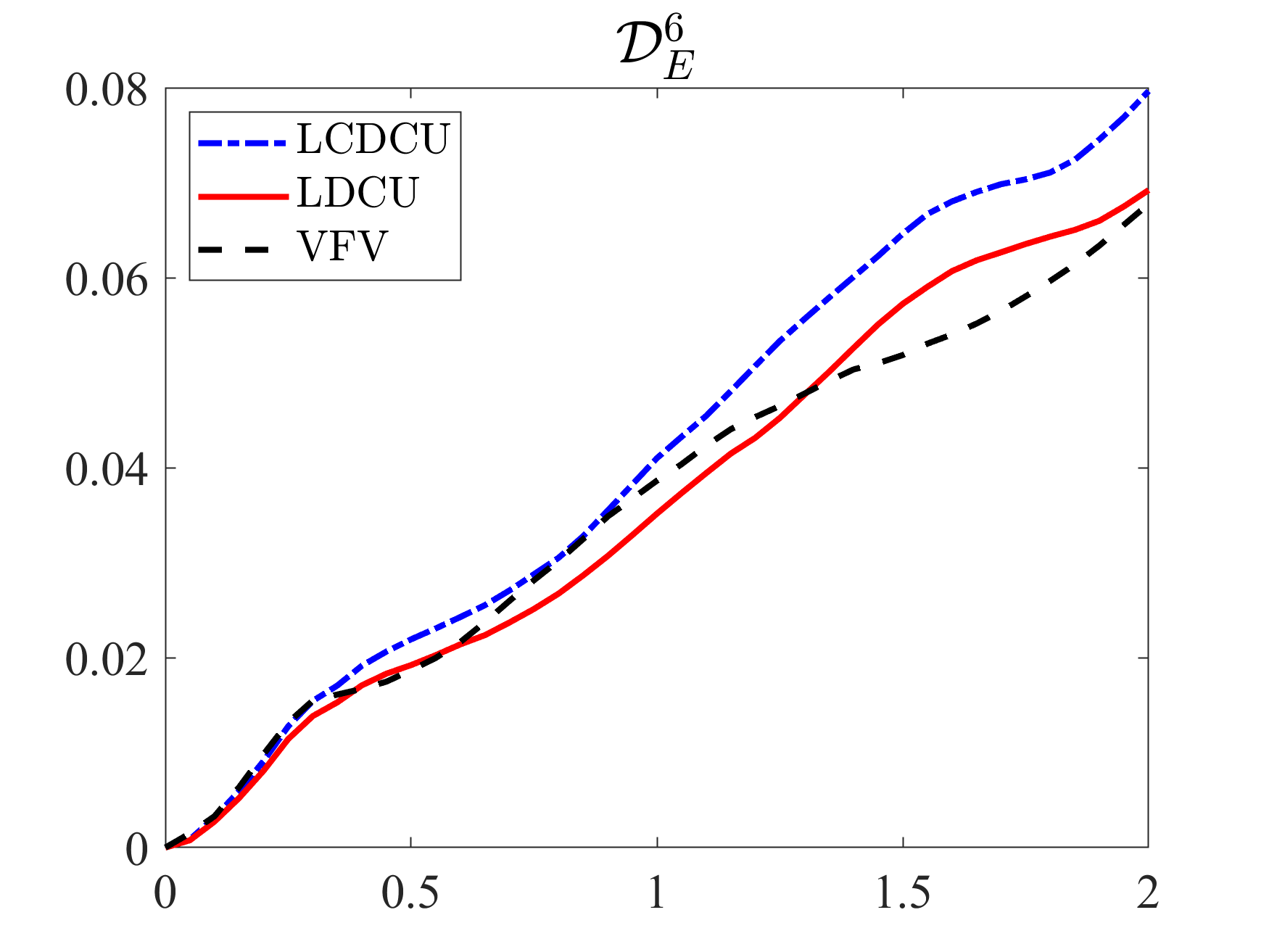}\hspace*{0.5cm}
            \includegraphics[trim=0.5cm 0.5cm 1.0cm 0.1cm, clip, width=3.8cm]{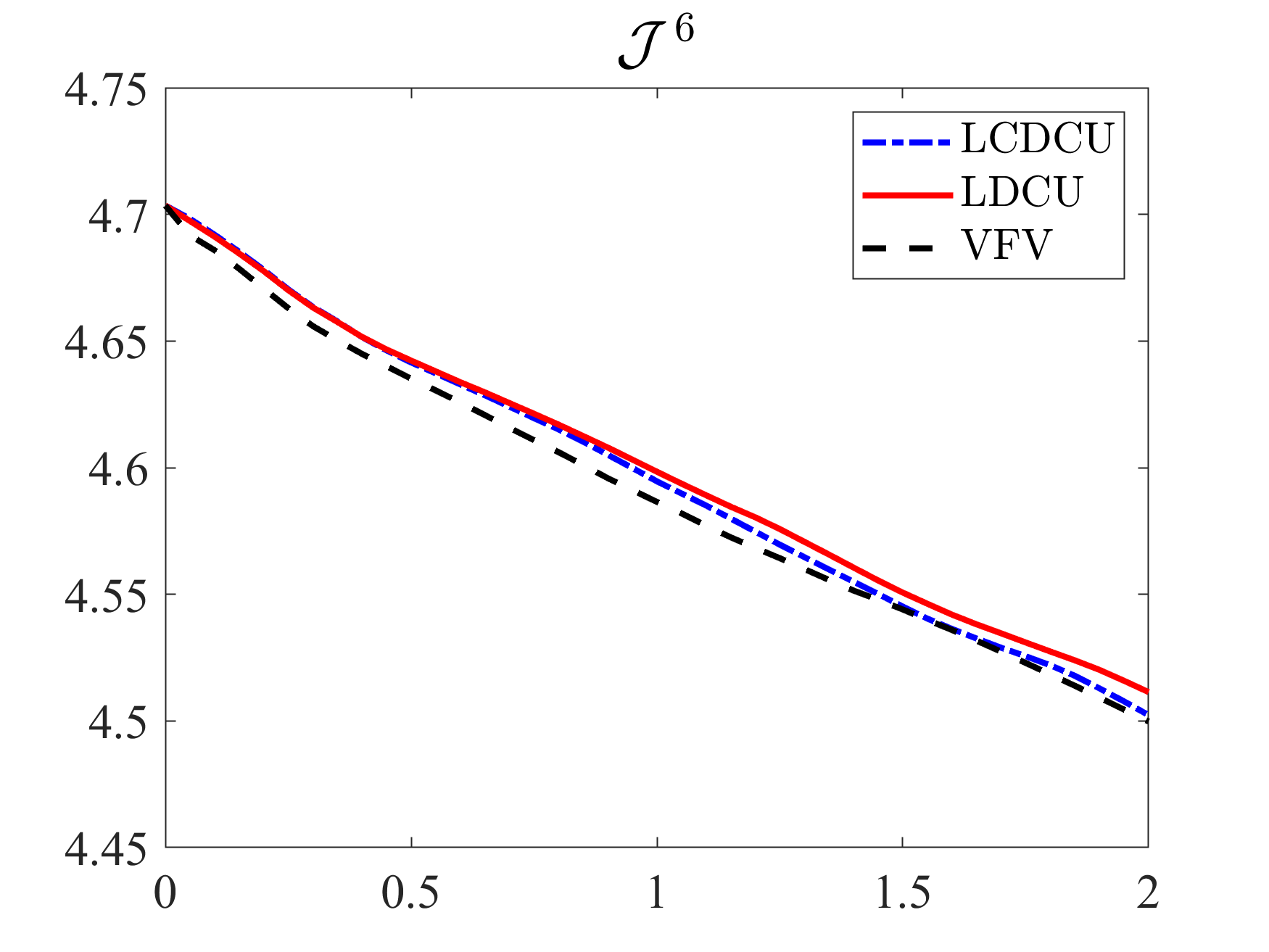}}           
\caption{\sf KH Instability: Time evolution of the entropy ${\mathcal S}^5$ (left), energy defects ${\mathcal D}^6_{\rm E}$ (middle), and the Bregman distance ${\mathcal J}^6$ (right) computed by the LCDCU, LDCU, and VFV schemes.\label{fig16}}
\end{figure}

Finally, we focus on the total entropy ${\mathcal S}^6$, the energy defect ${\mathcal D}_{\rm E}^6$, and the Bregman distance $\mathcal{J}^6$ for the KH instability problem. Here, the initial data yield $\bar{\theta}=1.716667$.
The time evolutions and time averages of these quantities are shown in Figure \ref{fig16} and Table \ref{table2}, respectively. The magnitude of
energy defects is substantially larger than in the Configuration 3 example, indicating that the limiting DW solutions are clearly not weak
solutions in the sense of distributions. This observation is consistent with the figures of the averaged numerical solutions presented in
\S\ref{sec43}. 
Based on the values reported in
Table \ref{table2}, criterion (iii) prioritizes the LCDCU schemes, while the entropy production criterion (i) and the Bregman distance to the equilibrium maximizing the entropy (iv) prioritize the VFV
schemes.

\begin{table}[ht!]
\centering
\begin{tabular}{|c|c|c|c|}\hline
&LCDCU&LDCU&VFV\\\hline
$\frac{1}{T}\int_0^T{\mathcal S}^6(t)\,{\rm d}t$&1.0493&1.0504&1.0548\\\hline
$\frac{1}{T}\int_0^T{\mathcal E}^6_1(t)\,{\rm d}t$&6.4375&6.4375&6.4375\\\hline
$\frac{1}{T}\int_0^T{\mathcal E}^6_2(t)\,{\rm d}t$&6.3952&6.4003&6.4015\\\hline
$\frac{1}{T}\int_0^T{\mathcal D}^6_E(t)\,{\rm d}t$&0.042321&0.037206&0.036857\\\hline
$\frac{1}{T}\int_0^T{\mathcal J}^6(t)\,{\rm d}t$&4.7114&4.7148&4.7058\\\hline
\end{tabular}\\[1.ex]
\caption{\sf KH Instability: Time averages of the entropy, total and mean energy, the energy defect, and the Bregman distance computed by the LCDCU, LDCU, and VFV schemes.\label{table2}}
\end{table}

\section{Conclusions}\label{sec5}
In this paper, we have presented a numerical study of dissipative weak (DW) solutions for the Euler equations of gas dynamics. We have considered different benchmark examples, 2-D Riemann problems with continuous or discontinuous solutions, and the Kelvin-Helmholtz (KH) instability problem. Three different numerical methods have been investigated, namely the local characteristic decomposition-based central-upwind scheme (LCDCU), the low-dissipation central-upwind scheme (LDCU), and the viscous finite volume method (VFV). The main result of our numerical study is that different numerical methods may converge to different DW solutions. This fact is evident in the approximation of density and other state variables as well as in the averages over different orders of accuracy, which have been confirmed to converge strongly, even when the convergence of one particular numerical realization is weak only. Furthermore, the averages over different orders of accuracy have been employed to distinguish between weak solutions in the sense of distributions and genuine DW solutions.

We have also computed the approximated error functions for each problem, which, according to our numerical results, behave similarly for different numerical methods. Moreover, we have discovered that the time averages of the average numerical solutions obtained by different methods exhibit similar behavior. In addition, we have approximated the probability density functions of Young measures. In cases where oscillations arise in numerical approximation, the limits are DW solutions and the corresponding Young measures are not Dirac $\delta$-functions and they clearly depend on the scheme. Finally, recognizing that DW solutions are generically not unique, we have explored four  potential selection criteria, which are based on the maximization of entropy production, either minimization or maximization of energy defect and the minimization of the Bregman distance to the constant equilibrium maximizing the entropy.

\begin{DA}
\noindent\textbf{Funding.} The work of S. Chu and  M. Herty was funded by the Deutsche Forschungsgemeinschaft (DFG, German Research Foundation)
-- SPP 2410 Hyperbolic Balance Laws in Fluid Mechanics: Complexity, Scales, Randomness (CoScaRa) within the Project(s) HE5386/26-1
(Numerische Verfahren f\"ur gekoppelte Mehrskalenprobleme,525842915) and (Zuf\"allige kompressible Euler Gleichungen: Numerik und ihre
Analysis, 525853336) HE5386/27-1, and the Deutsche Forschungsgemeinschaft (DFG, German Research Foundation) -- SPP 2183:
Eigenschaftsgeregelte Umformprozesse with the Project(s) HE5386/19-2,19-3 Entwicklung eines flexiblen isothermen Reckschmiedeprozesses f\"ur die eigenschaftsgeregelte Herstellung von Turbinenschaufeln aus Hochtemperaturwerkstoffen (424334423). The work of A. Kurganov was supported in part by NSFC grant W2431004. The work of  M. Luk\'{a}\v{c}ov\'{a}-Medvi{\softd}ov\'{a} was supported by the Gutenberg Research College and by the Deutsche Forschungsgemeinschaft (DFG, German Research Foundation) -- project number 233630050 -- TRR 146 and project number 525853336 -- SPP 2410 ``Hyperbolic Balance Laws: Complexity, Scales and Randomness''. She is also grateful to the Mainz Institute of Multiscale Modelling for supporting her research. C. Yu was supported by the DFG project number 525853336 -- SPP 2410 ``Hyperbolic Balance Laws: Complexity, Scales and Randomness''.\\

\noindent\textbf{Conflicts of interest.} On behalf of all authors, the corresponding author states that there is no conflict of interest.\\

\noindent\textbf{Data and software availability.} The data that support the findings of this study and FORTRAN codes developed by the authors and
used to obtain all of the presented numerical results are available from the corresponding author upon reasonable request.
\end{DA}


\end{document}